\documentclass{amsart}
\usepackage[utf8]{inputenc}
\usepackage[T1]{fontenc}
\usepackage{lmodern}
\usepackage{amsmath,amsfonts,amssymb,amscd,amsthm,amsbsy,upref}
\usepackage[dvipsnames]{xcolor}
\usepackage{mathtools,mathrsfs}
\usepackage{graphicx}
\usepackage{booktabs,array,longtable}
\usepackage{enumitem}
\usepackage{microtype}
\usepackage[colorlinks=true,linkcolor=blue!55!black,citecolor=blue!55!black,urlcolor=blue!55!black]{hyperref}
\usepackage{aliascnt}
\usepackage[nameinlink,capitalize,noabbrev]{cleveref}
\hypersetup{
  pdftitle={Magnetohydrostatic equilibria with exact finite cyclic symmetry},
  pdfsubject={Smooth unforced magnetohydrostatic equilibria on solid tori whose Euclidean symmetry group is exactly a finite cyclic group},
  pdfkeywords={magnetohydrostatics, Grad's conjecture, Nash-Moser theorem, symmetry breaking, toroidal equilibria}
}
\renewcommand{\baselinestretch}{1.2}
\numberwithin{equation}{section}
\newcommand{\R}{\mathbb R}

\newcommand{\Z}{\mathbb Z}
\newcommand{\T}{\mathbb T}
\newcommand{\D}{\mathbb D}
\newcommand{\cA}{\mathcal A}
\newcommand{\cB}{\mathcal B}
\newcommand{\cC}{\mathcal C}
\newcommand{\cD}{\mathcal D}
\newcommand{\cE}{\mathcal E}
\newcommand{\cF}{\mathcal F}
\newcommand{\cG}{\mathcal G}
\newcommand{\cH}{\mathcal H}
\newcommand{\cL}{\mathcal L}

\newcommand{\cQ}{\mathcal Q}
\newcommand{\cR}{\mathcal R}
\newcommand{\cS}{\mathcal S}
\newcommand{\cX}{\mathcal X}
\newcommand{\cY}{\mathcal Y}
\newcommand{\bX}{\mathbb X}
\newcommand{\bY}{\mathbb Y}
\newcommand{\eps}{\varepsilon}
\newcommand{\pa}{\partial}
\newcommand{\grad}{\nabla}
\newcommand{\curl}{\operatorname{curl}}
\newcommand{\diver}{\operatorname{div}}
\newcommand{\Emb}{\operatorname{Emb}}
\newcommand{\Diff}{\operatorname{Diff}}
\newcommand{\Stab}{\operatorname{Stab}}

\newcommand{\tr}{\operatorname{tr}}
\newcommand{\diag}{\operatorname{diag}}

\newcommand{\dd}{\,\mathrm d}
\newcommand{\inner}[2]{\left\langle #1,#2\right\rangle}
\newcommand{\norm}[1]{\left\lVert #1\right\rVert}
\newcommand{\abs}[1]{\left\lvert #1\right\rvert}

\newcommand{\natg}{\natural}

\theoremstyle{plain}
\newtheorem{theorem}{Theorem}[section]
\crefname{theorem}{Theorem}{Theorems}
\Crefname{theorem}{Theorem}{Theorems}
\newaliascnt{proposition}{theorem}
\newtheorem{proposition}[proposition]{Proposition}
\aliascntresetthe{proposition}
\crefname{proposition}{Proposition}{Propositions}
\Crefname{proposition}{Proposition}{Propositions}
\newaliascnt{lemma}{theorem}
\newtheorem{lemma}[lemma]{Lemma}
\aliascntresetthe{lemma}
\crefname{lemma}{Lemma}{Lemmas}
\Crefname{lemma}{Lemma}{Lemmas}
\newaliascnt{corollary}{theorem}
\newtheorem{corollary}[corollary]{Corollary}
\aliascntresetthe{corollary}
\crefname{corollary}{Corollary}{Corollaries}
\Crefname{corollary}{Corollary}{Corollaries}
\newaliascnt{claim}{theorem}

\aliascntresetthe{claim}
\crefname{claim}{Claim}{Claims}
\Crefname{claim}{Claim}{Claims}
\newaliascnt{sublemma}{theorem}

\aliascntresetthe{sublemma}
\crefname{sublemma}{Sublemma}{Sublemmas}
\Crefname{sublemma}{Sublemma}{Sublemmas}
\newaliascnt{conjecture}{theorem}

\aliascntresetthe{conjecture}
\crefname{conjecture}{Conjecture}{Conjectures}
\Crefname{conjecture}{Conjecture}{Conjectures}

\theoremstyle{definition}
\newaliascnt{definition}{theorem}

\aliascntresetthe{definition}
\crefname{definition}{Definition}{Definitions}
\Crefname{definition}{Definition}{Definitions}
\newaliascnt{assumption}{theorem}

\aliascntresetthe{assumption}
\crefname{assumption}{Assumption}{Assumptions}
\Crefname{assumption}{Assumption}{Assumptions}
\newaliascnt{remark}{theorem}
\newtheorem{remark}[remark]{Remark}
\aliascntresetthe{remark}
\crefname{remark}{Remark}{Remarks}
\Crefname{remark}{Remark}{Remarks}
\newaliascnt{warning}{theorem}

\aliascntresetthe{warning}
\crefname{warning}{Warning}{Warnings}
\Crefname{warning}{Warning}{Warnings}
\newaliascnt{question}{theorem}

\aliascntresetthe{question}
\crefname{question}{Question}{Questions}
\Crefname{question}{Question}{Questions}
\newaliascnt{objective}{theorem}

\aliascntresetthe{objective}
\crefname{objective}{Objective}{Objectives}
\Crefname{objective}{Objective}{Objectives}

\renewcommand{\labelenumi}{(\theenumi)}
\renewcommand{\le}{\leqslant}
\renewcommand{\ge}{\geqslant}

\usepackage[T1]{fontenc}
\usepackage[utf8]{inputenc}

\usepackage{amsmath,amssymb,amsthm,mathrsfs,graphicx,xcolor,pdfrender,listings}
\usepackage[scaled=0.83]{DejaVuSansMono}
\definecolor{leanblue}{RGB}{33,74,135}
\definecolor{leanop}{RGB}{24,45,75}
\definecolor{leancomment}{RGB}{112,83,36}
\newcommand{\leansymbol}[1]{{\color{leanop}\textpdfrender{
  TextRenderingMode=FillStroke,LineWidth=0.07pt}{\ensuremath{\boldsymbol{#1}}}}}
\newcommand{\leanletter}[1]{\textpdfrender{
  TextRenderingMode=FillStroke,LineWidth=0.05pt}{\ensuremath{#1}}}
\newcommand{\leanOrderType}{\leanletter{\mathbb N\mkern2mu\boldsymbol\infty\mkern2mu\omega}}
\lstdefinelanguage{LeanShowcase}{
  morekeywords={abbrev,def,inductive,structure,where,instance,theorem,by,sorry,
    let,fun,if,then,else,exact,classical,noncomputable,section,namespace,end,
    import,open,scoped,Prop,Type},
  sensitive=true,alsoletter={_},
  morecomment=[l]{--},morecomment=[s]{/-}{-/},morestring=[b]"
}
\lstdefinestyle{lean}{
  language=LeanShowcase,basicstyle=\small\ttfamily,
  keywordstyle=\color{leanblue}\bfseries,
  commentstyle=\color{leancomment}\itshape,stringstyle=\color{leancomment},
  columns=fullflexible,keepspaces=true,showstringspaces=false,
  breaklines=true,breakatwhitespace=true,breakindent=1.5em,
  aboveskip=7pt,belowskip=7pt,tabsize=2,
  literate=
    {ℝ}{{\leanletter{\mathbb R}}}1 {ℕ}{{\leanletter{\mathbb N}}}1
    {∞}{{\leansymbol{\infty}}}1 {ω}{{\leanletter{\omega}}}1
    {∀}{{\leansymbol{\forall}}}1 {∃}{{\leansymbol{\exists}}}1
    {→}{{\leansymbol{\to}}}1 {↔}{{\leansymbol{\leftrightarrow}}}1
    {∧}{{\leansymbol{\wedge}}}1 {∨}{{\leansymbol{\vee}}}1
    {¬}{{\leansymbol{\neg}}}1 {∈}{{\leansymbol{\in}}}1
    {⊆}{{\leansymbol{\subseteq}}}1 {∩}{{\leansymbol{\cap}}}1
    {×}{{\leansymbol{\times}}}1 {•}{{\leansymbol{\bullet}}}1
    {≠}{{\leansymbol{\ne}}}1 {≤}{{\leansymbol{\le}}}1
    {≥}{{\leansymbol{\ge}}}1 {‖}{{\leansymbol{\Vert}}}1
    {≃}{{\leansymbol{\simeq}}}1
    {⟨}{{\leansymbol{\langle}}}1 {⟩}{{\leansymbol{\rangle}}}1
    {∑}{{\leansymbol{\sum}}}1 {∘}{{\leansymbol{\circ}}}1
    {Ω}{{\leanletter{\Omega}}}1 {Γ}{{\leanletter{\Gamma}}}1
    {Φ}{{\leanletter{\Phi}}}1 {κ}{{\leanletter{\kappa}}}1
    {σ}{{\leanletter{\sigma}}}1 {θ}{{\leanletter{\theta}}}1
    {γ}{{\leanletter{\gamma}}}1 {λ}{{\leanletter{\lambda}}}1
    {₀}{{\leanletter{{}_0}}}1
    {ℕ∞ω}{{\leanOrderType}}3
    {×ˢ}{{\leansymbol{\times}\leanletter{{}^{\mathrm s}}}}2
    {≃ₗᵢ}{{\leansymbol{\simeq}\leanletter{{}_{\mathrm{li}}}}}3
}
\lstnewenvironment{leancode}[1][]{\lstset{style=lean,#1}}{}
\newcommand{\lean}{\lstinline[style=lean]}
\newtheorem*{leanmainresult}{Theorem (Equilibria with exact cyclic symmetry)}

\title[Counterexamples to Grad's conjecture]{Counterexamples to Grad's conjecture}

\author{Javier G\'omez-Serrano}
\address{Department of Mathematics\\ Brown University, 151 Thayer Street, 02912 Providence RI, USA.}
\email{javier\_gomez\_serrano@brown.edu}

\author[Lukas Liehr]{Lukas Liehr}
\address{Department of Mathematics, Bar-Ilan University, Ramat-Gan 5290002, Israel}
\email{lukas.liehr@biu.ac.il}

\author{Mitchell~A.\ Taylor}
\address{Mathematical Institute, University of Oxford, Andrew Wiles Building, Radcliffe Observatory Quarter, Woodstock Road, Oxford, OX2 6GG, United Kingdom.}
\email{mitchtaylor@shaw.ca}

\date{}


\begin{document}

\begin{abstract}

 For each sufficiently large $N$, we construct smooth solutions of the magnetohydrostatic equations
on embedded solid tori whose regular pressure levels are nested tori, whose magnetic
field vanishes exactly on a round magnetic axis, and whose group of
Euclidean symmetries, even when isometries reversing the sign of the field
are admitted, is exactly the cyclic group $C_N$ of rotations through
multiples of $2\pi/N$. For each such $N$, these equilibria form a smooth one-parameter family
which is locally nontrivial in the moduli space of embedded configurations and have none of the plane-reflection, axial, or helical
symmetries conjectured by Grad. A Lean 4 certification of the above results is also provided. %The starting point is the family of exact
\end{abstract}

\maketitle
\tableofcontents

\section{Introduction}\label{sec:introduction}

The objective of this article is to construct exotic smooth solutions of the magnetohydrostatic equations
\begin{equation}\label{eq:mhs}
 B\times\curl B+\grad P=0,
 \qquad \diver B=0
 \quad\text{in }\Omega,
 \qquad B\cdot n=0
 \quad\text{on }\pa\Omega ,
\end{equation}
on solid tori $\Omega\subset\R^3$, with a prescribed symmetry group.  Here, $B$ is the
magnetic field, $P$ is the pressure, $n$ is the outward unit normal, and by a \emph{solid torus} we mean a compact domain diffeomorphic to $\overline{\D^2}\times\T$. The
system \eqref{eq:mhs} describes the static equilibria of an ideal,
perfectly conducting plasma. Indeed, the first equation states that the Lorentz force balances the pressure
gradient, the second states that the field is \emph{solenoidal}, and  the third states that the boundary $\pa\Omega$ is a \emph{magnetic surface}.
Since the first equation carries no external force, we call
\eqref{eq:mhs} the \emph{unforced magnetohydrostatic (MHS) equations}. 

Taking the inner product of the first equation in \eqref{eq:mhs} with $B$ gives
$B\cdot\grad P=0$, so field lines lie on the level sets of the pressure.
In plasma confinement, one is interested in equilibria whose regular
pressure levels are nested tori filling $\Omega$ and collapsing in the
interior onto a closed curve, the \emph{magnetic axis}, where $\grad P=0$.
On each torus, the first equation imposes a solvability condition along the
field lines \cite{Hamada1962,KruskalKulsrud1958,Newcomb1959}, and this
condition obstructs the construction of nonsymmetric equilibria with
nonconstant pressure. Grad \cite{Grad1967,Grad1985} argued that, outside the
classical symmetric configurations, no smooth equilibrium with nested
toroidal pressure surfaces should belong to a family depending smoothly on
a parameter. A precise formulation of this conjecture was given by Constantin, Drivas and Ginsberg
\cite{CDG2021}, and it states that a smooth unforced MHS equilibrium on a domain of toroidal
or cylindrical topology, whose regular pressure surfaces foliate the
domain, whose field does not vanish away from the magnetic axis, and which
is not isolated, has plane-reflection, axial, or helical symmetry. Here,
\emph{non-isolation} means that in some suitable topology there are nearby
equilibria other than trivial rescalings and translations of the given
one. %Note that the first of the three symmetries reverses orientation, and the other
%two are continuous.

In this article, we construct equilibria for which none of the three alternatives holds.
More precisely, for every sufficiently large integer $N$, we produce a smooth one-parameter
family of solutions of \eqref{eq:mhs} on embedded solid tori. The
pressure levels of these solutions are nested tori, their field vanishes exactly on a round
magnetic axis, and their full group of Euclidean symmetries, even allowing
isometries that reverse the sign of $B$, is the cyclic group $C_N$ of
rotations through multiples of $2\pi/N$ about the axis of the torus. %Here, the perturbation
%parameter is $1/N$, which measures the curvature of one period cell.
%Solving for all small real values of this parameter gives the result for
%every $N\ge N_0$. %Since the domains vary, comparison with the conjecture
%depends on the topology used to define non-isolation. We make this topology explicit in
%\cref{ss:moduli} and draw the conditional consequence for Grad's
%conjecture in \cref{cor:grad}.

\subsection{Embedded configurations and their moduli topology}\label{ss:moduli}

We compare equilibria on different domains by pulling them back to a
common reference domain. Let
\[
 Q=\overline{\D^2}\times\T,
 \qquad \T=\R/(2\pi\Z),
\]
and, for an integer $k\ge3$, represent an equilibrium by the triple
\[
 (X,b,p)=(X,B\circ X,P\circ X)
 \in \Emb^k(Q,\R^3)\times C^{k-1}(Q,\R^3)\times C^k(Q),
\]
where $X$ is an embedding of $Q$ onto $\Omega=X(Q)$. The above indices are
matched so that $\curl B$ is continuously differentiable, and the product
carries its usual topology. Two triples describe the same equilibrium if
they differ by a reparametrization of the reference domain,
\begin{equation}\label{eq:reparam-action}
 (X,b,p)\longmapsto (X\circ\Phi,b\circ\Phi,p\circ\Phi),
 \qquad \Phi\in\Diff^k(Q),
\end{equation}
or by a spatial similarity, a rescaling of the magnetic amplitude, or a
shift of the pressure, all of which preserve \eqref{eq:mhs} and can be represented by
\begin{equation}\label{eq:physical-moduli-action}
 (X,b,p)\longmapsto
 \bigl(sOX\circ\Phi+c,\;\kappa Ob\circ\Phi,\;
       \kappa^2p\circ\Phi+c_0\bigr),
\end{equation}
with $s>0$, $O\in O(3)$, $c\in\R^3$, $\kappa\in\R\setminus\{0\}$, and
$c_0\in\R$. We write $\mathscr M^k_{\rm emb}(Q)$ for the quotient with
the natural quotient topology, and call it the \emph{moduli space of embedded
configurations}. For smooth configurations, we use the usual $C^\infty$ topology
and quotient by smooth reparametrizations and the same physical action.
The embedded domain is part of the configuration.

For a Euclidean isometry $g(x)=Ox+c$, we set $(g_*B)(x)=OB(g^{-1}x)$ and
define the \emph{extended stabilizer} by
\begin{equation}\label{eq:stabilizer-pm}
 \Stab_{\pm}(\Omega,B,P)
 =\{g:\ g\Omega=\Omega,\ P\circ g^{-1}=P,\ g_*B=\pm B\}.
\end{equation}
Thus, this stabilizer includes isometries that reverse the field. We will
prove that distinct members of each family we construct are inequivalent
even under the larger group \eqref{eq:physical-moduli-action}.

% To
%distinguish members of this family, we allow the larger equivalence group
%in \eqref{eq:physical-moduli-action}.

\subsection{Historical comments}\label{ss:history}

The existence of three-dimensional equilibria with nested pressure surfaces
is the central mathematical question of stellarator theory. We refer to
Freidberg \cite{Freidberg2014} and Helander \cite{Helander2014} for an account of the underlying
physics. Grad's argument \cite{Grad1967} rests on the analysis of the real characteristics
of the MHS system along field lines. Near a torus with irrational
rotational transform, the solvability condition of
\cite{Hamada1962, KruskalKulsrud1958,Newcomb1959} must hold on a dense set
of rational tori, and this appears to force the pressure to flatten. One
way around this obstruction is to admit current sheets. Bruno and Laurence
\cite{BrunoLaurence1996} proved the existence of equilibria with piecewise
constant pressure in tori that are small perturbations of an axisymmetric
one. In their setting, the pressure jumps across finitely many interfaces carrying current
sheets, and the rotational transform on each interface must be
sufficiently irrational (see also Kaiser and Salat \cite{KaiserSalat1994}).
This  picture underlies the equilibrium program of
Hudson, Dewar, and their colleagues
\cite{Hudson2012, HudsonHoleDewar2007,HudsonKraus2017,KrausHudson2017,MacKay2025}.
Notably, Enciso, Luque and Peralta-Salas \cite{EncisoLuquePeraltaSalas2025}
recently constructed piecewise smooth equilibria with arbitrarily many
current sheets in toroidal domains that need not be close to axisymmetric
and may be knotted. In \cite{BrunoLaurence1996,EncisoLuquePeraltaSalas2025}, the pressure is constant
between the interfaces and the solutions are not $C^1$.

Smooth equilibria without continuous symmetries are known to exist in a few special
situations. Lortz \cite{Lortz1970} constructed toroidal equilibria with
nonconstant pressure whose domains are symmetric with respect to a plane
but otherwise have no symmetry. In this construction, all field lines are closed and the rotational
transform is zero. The planar symmetry makes it consistent with Grad's conjecture. Salat and Kaiser
\cite{SalatKaiser1995} produced exact equilibria in an infinite cylinder
with a straight axis but no cylindrical, helical, or mirror symmetry. Their
field lines are closed curves in the planes orthogonal to the axis, and
their pressure surfaces have elliptical cross-sections whose eccentricity
and orientation vary freely along the axis. Sato \cite{Sato2019} obtained
equilibria in bounded domains without continuous Euclidean symmetry, and
Sato and Yamada \cite{SatoYamada2023} constructed nested invariant tori without such
symmetry for anisotropic magnetohydrodynamics. Weitzner obtained formal expansions of nonsymmetric MHS equilibria to all
orders in doubly periodic domains \cite{Weitzner2014,Weitzner2016} and
about a circular magnetic axis \cite{Weitzner2016}. Weitzner and Sengupta
\cite{WeitznerSengupta2021} treated nearly parallel plasma flows in a
doubly periodic domain. For related work on MHS equilibria and their symmetry properties, we also
refer to
\cite{BurbyKallinikosMacKay2020,CaryShasharina1997,GarrenBoozer1991,LandremanSengupta2018}.

On the numerical side, three-dimensional MHS equilibria have been computed
for decades, and the difficulty behind Grad's conjecture is visible in the
design of the codes. VMEC \cite{HirshmanWhitson1983} and, more recently,
DESC \cite{DudtKolemen2020} assume nested flux surfaces from the outset and
minimize the energy or the force residual within this class. Hudson and
Kraus \cite{HudsonKraus2017,KrausHudson2017} have argued that nested
surfaces with a smooth pressure profile produce infinite currents at the
rational surfaces, so that the pressure gradient must vanish on an interval
around each of them and the resulting profile is fractal. PIES
\cite{ReimanGreenside1986}, HINT \cite{HarafujiHayashiSato1989} and SIESTA
\cite{HirshmanSanchezCook2011} drop the assumption of nested surfaces and
resolve islands and stochastic regions, while SPEC \cite{Hudson2012}
implements the stepped-pressure picture described above. Loizu, Hudson,
Bhattacharjee and their collaborators
\cite{LoizuHudsonBhattacharjeeHelander2015,LoizuHudsonBhattacharjeeLazersonHelander2015}
computed the singular currents that form at rational surfaces in the
presence of a pressure gradient, and Zhou, Huang, Qin and Bhattacharjee
\cite{ZhouHuangQinBhattacharjee2016} followed the formation of current
singularities in topologically constrained plasmas. More recently,
physics-informed neural networks \cite{RaissiPerdikarisKarniadakis2019},
which minimize the residual of the equations over a mesh-free
parametrization of the unknowns, have been applied to the Grad--Shafranov
equation \cite{JangKaptanogluGaurPanLandremanDorland2024,ZhouZhu2025} and to
fully three-dimensional stellarator equilibria with nested surfaces
\cite{ThunMerloConlinPaniciBockenhoff2026}, where they reach lower force
residuals than the conventional solvers.

The formulation of Grad's conjecture stated above is due to Constantin, Drivas
and Ginsberg \cite{CDG2021}. In their paper, they also proved rigidity and flexibility theorems
for steady Euler, Boussinesq, and MHS flows. Their flexibility results
deform a given equilibrium to nearby domains, which is the point of view
we take in \cref{ss:moduli}. In \cite{CDG2021jpp}, the same authors
constructed smooth nonsymmetric equilibria with nested flux surfaces at
the price of a small forcing term, and in \cite{CDG2022} they studied the
free boundary version of the problem. Cardona, Duignan and Perrella
 showed in \cite{CardonaDuignanPerrella2025} that on a Riemannian manifold the
absence of Killing symmetries is generic among adapted metrics. On the other hand, in 
\cite{Perrella2025} it is proved that an orientation-reversing isometry of an
equilibrium with toroidally nested pressure surfaces forces all field
lines to be periodic, and in 
\cite{ClellandKlotz2020, EncisoPeraltaSalas2016,PeraltaSalasVaquero2024} it is shown that Beltrami fields with
nonconstant proportionality factor are rare. Interestingly, Enciso, Kepplinger and
Peralta-Salas \cite{EncisoKepplingerPeraltaSalas2025} found steady Euler
flows on certain closed three-manifolds that are isolated in $C^1$, suggesting that the
non-isolation hypothesis in the conjecture is not vacuous. For 
structure theorems and steady Euler flows of low regularity, we refer to
\cite{ArnoldKhesin2021,EncisoPenafielPeraltaSalas2025}.

The closest result to ours is the construction by Drivas, Elgindi and
Ginsberg \cite{DrivasElgindiGinsberg2025} of smooth steady Euler flows in a
periodic cylinder, fibered by invariant cylinders that are level sets of
the Bernoulli function. Their flows have an $m$-fold rotational symmetry
and a plane-reflection symmetry but no continuous Euclidean symmetry. The
solutions that they construct come in infinite-dimensional families, the perturbation parameter is
$1/m$, and the iteration scheme is motivated by the approach of Lortz. They ask whether all Euclidean
symmetries, including the reflection, can be broken, and they propose to
refine Grad's conjecture by requiring ergodic orbits on almost every
invariant torus.

%Our equilibria have no reflection but retain a finite
%cyclic group, so the first question remains open. We do not address the
%refined conjecture: as in \cite{Lortz1970,DrivasElgindiGinsberg2025},
%every nonconstant field line of our equilibria is closed. This removes
%the obstruction of \cite{Grad1967}.
In the present article, we construct new MHS solutions that exhibit the following properties.
\begin{enumerate}
\item \label{i} a smooth pressure with nested embedded tori in a solid
torus and a magnetic axis on which the field vanishes;
\item \label{ii} a prescribed symmetry group, excluding
reflections as well as continuous symmetries, even when isometries
reversing the field are admitted;
%\item \label{iii} a conclusion for every sufficiently large field period;
\item \label{iv} for each of the prescribed symmetry groups, a locally nontrivial smooth family of solutions obtained by a Nash-Moser argument.
\end{enumerate}

\begin{figure}[htbp]
\centering
\includegraphics[width=\textwidth]{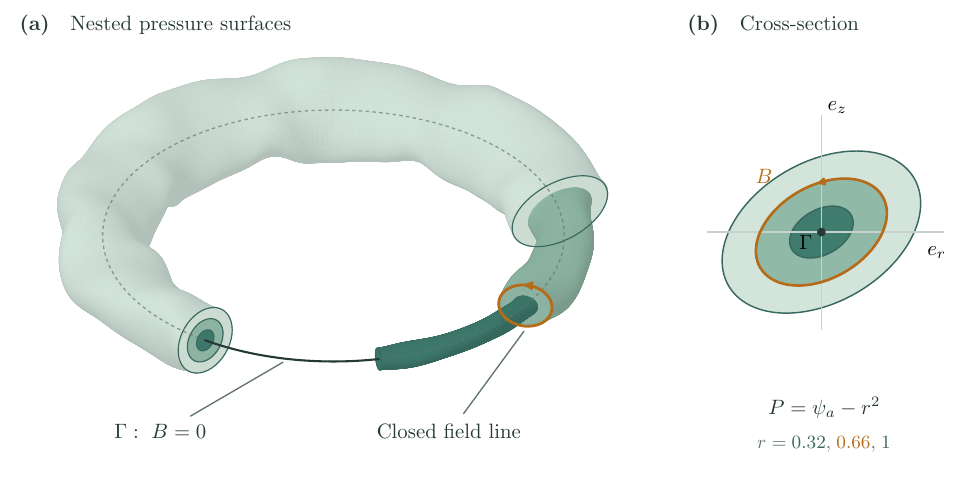}
\caption{Nested pressure surfaces for the model
$X^0(y,\phi)=R_\phi(NLe_x+\iota M_\lambda(N\phi)y)$.
(a) A cutaway showing three pressure surfaces, a closed field line,
and the magnetic axis $\Gamma$. The axis is dashed where hidden.
(b) A cross-section through the highlighted field line.
We use $N=8$, $L=0.65$, $\rho_*=0.40$, $\delta_*=0.425$,
$\alpha_0=\pi/6$, and $\lambda=0.375$.
The cross-sections have axis ratio approximately $1.53:1$.
The eccentricity is enlarged for visibility and lies outside the theorem's
small-eccentricity regime. The nonlinear correction in \eqref{eq:intro-X}
is omitted.}
\label{fig:grad}
\end{figure}

\textbf{LLM usage disclosure:} The work on this project began in July 2026. In July, the authors set up the problem and provided a detailed roadmap
to construct the solutions. After this, they used a combination of the LLMs
GPT-5.6 Sol, Claude Fable 5, and Claude Opus 5 to fill in the technical details, assist with calculations, and provide mathematical feedback
on the parts of the blueprint that were not detailed enough. The Lean code was produced by the same models with continuous guidance from the authors. During the latter stages of the Lean verification, the authors switched to the newer models GPT-6 Astra and 
Claude Fable 5.1. These models were also used for final proofreading but did not contribute to any of the mathematical content of the article. All mathematical statements and proofs have
been checked by the authors, who take full responsibility for their correctness.

We also remark that our examples share some resemblances to the recent and spectacular counterexamples \cite{cao2026counterexamples,colbrook2026computer} to Schiffer's conjecture. Our approach is closer in spirit to that of \cite{cao2026counterexamples}, but we also ran an independent computer-assisted approach. This latter attempt was abandoned once it became clear that the approach we present here would close, but it would nonetheless be interesting to systematically pursue such methods as they should have immense value.

\subsection{The main result}\label{ss:main}

Let $e_x,e_y,e_z$ be the standard basis of $\R^3$ and $R_\phi$ the
rotation through $\phi$ about the $e_z$-axis. Throughout, $N$ is the
number of periods of the field. %and $J$ an unspecified constant that will appear later as a cutoff in the number of Fourier modes needed in the linear analysis.
Note that the pressure is determined by \eqref{eq:mhs} up to an additive constant,
which we fix by prescribing its value $\psi_a$ on the magnetic axis.

\Cref{fig:grad} illustrates the geometry and terminology. The $C_N$
symmetry consists of rotations about the $e_z$-axis. The magnetic axis is
the circle $\Gamma$ where the nested tori collapse and the magnetic field
vanishes.
A \emph{cell} is one of the $N$ repeated pieces of the torus. Consecutive
cells are related by a rotation through $2\pi/N$. The magnetic axis has
radius $NL$, so its arc within each cell has length $2\pi L$. At leading
order, $\lambda$ controls how the elliptical cross-sections turn along a cell.

Our main result may be stated as follows.
% AUTHOR QUERY: The existence claim depends on the unresolved gluing and
% inverse estimates marked in Sections 7, 9--12.
\begin{theorem}[Equilibria with exact cyclic symmetry]\label{thm:main}
Fix $L>0$. There exist a compact interval $I\subset(0,\tfrac12)$ with nonempty interior and
an integer $N_0<\infty$ such that the following holds. For every
integer $N\ge N_0$ and every $\lambda\in I$, there is a smooth
embedding
\[
 X_{N,\lambda}:Q\longrightarrow\R^3
\]
and smooth fields $B_{N,\lambda},P_{N,\lambda}$ on
$\Omega_{N,\lambda}=X_{N,\lambda}(Q)$ satisfying \eqref{eq:mhs}, with the
following properties.
\begin{enumerate}
\item (Nested pressure surfaces). The unique critical set of $P_{N,\lambda}$
is the round circle
\begin{equation}\label{eq:main-axis}
 \Gamma_N=\{NL(\cos\phi,\sin\phi,0):\phi\in\T\},
\end{equation}
every regular pressure level is an embedded torus, and these levels
foliate $\Omega_{N,\lambda}\setminus\Gamma_N$.
\item (Magnetic axis). The field vanishes exactly on the axis, that is,
\begin{equation}\label{eq:main-zero-set}
 B_{N,\lambda}(x)=0\quad\Longleftrightarrow\quad x\in\Gamma_N.
\end{equation}
\item (Exact stabilizer). The extended stabilizer is the
cyclic group of rotations,
\begin{equation}\label{eq:main-stabilizer}
 \Stab_{\pm}(\Omega_{N,\lambda},B_{N,\lambda},P_{N,\lambda})
 =\{R_{2\pi j/N}:j=0,\ldots,N-1\}\cong C_N.
\end{equation}
In particular, the equilibrium has neither a continuous Euclidean symmetry
nor an orientation-reversing Euclidean symmetry.
\item (Non-isolation). For each fixed $N\ge N_0$, the representatives may
be chosen to depend smoothly on $\lambda$, and the map
\begin{equation}\label{eq:main-moduli-curve}
 \lambda\longmapsto
 [X_{N,\lambda},B_{N,\lambda}\circ X_{N,\lambda},
                    P_{N,\lambda}\circ X_{N,\lambda}]
 \in\mathscr M^k_{\rm emb}(Q)
\end{equation}
is continuous and injective for every finite $k\ge3$, as well as in the
$C^\infty$ topology. Consequently, every member with
$\lambda\in\operatorname{int}I$ is non-isolated in the moduli space of
embedded configurations.
\end{enumerate}
\end{theorem}

The rotations in \eqref{eq:main-stabilizer} exhaust the symmetry group,
including isometries sending $B$ to $-B$. This group has none of the
plane-reflection, axial, or helical symmetries posited in \cite{CDG2021}, and thus \cref{thm:main} refutes Grad's conjecture.

\begin{figure}[htbp]
\centering
\includegraphics[width=0.94\textwidth]{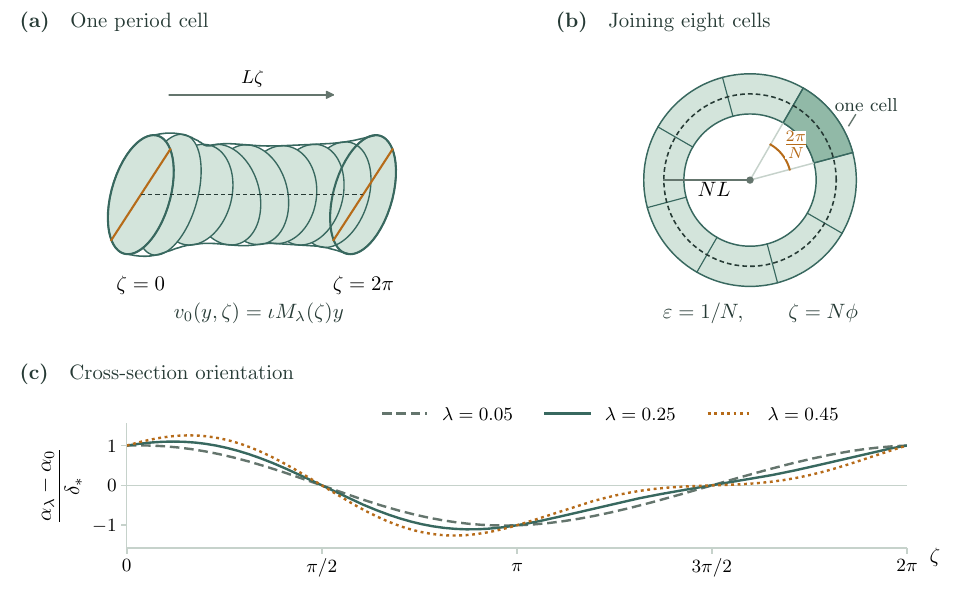}
\caption{From one period cell to a torus.
(a) A straight cell at zero curvature, with matching end cross-sections
and the parameters of \cref{fig:grad}.
(b) A schematic view of eight cells joined around a circle of radius
$NL=5.2$. Each cell spans an angle $2\pi/N$, and the radial lines mark
where the cells meet.
(c) The cross-section angle over one period, plotted as
$(\alpha_\lambda-\alpha_0)/\delta_*=\cos\zeta+\lambda\sin2\zeta$
for $\lambda=0.05$, $0.25$, and $0.45$.}
\label{fig:cell-geometry}
\end{figure}

As illustrated in \Cref{fig:grad,fig:cell-geometry}, our construction depends on the three constants
\[
 0<\rho_*<\frac14,\qquad \delta_*\ne0,\qquad
 \alpha_0\notin\frac\pi2\Z,
\]
which we will later fix. These constants allow us to make
the equilibria completely explicit to leading order in $1/N$. Writing $\zeta=N\phi$ for the coordinate along a cell, we set
\begin{align}
 \alpha_\lambda(\zeta)
   &=\alpha_0+\delta_*(\cos\zeta+\lambda\sin2\zeta),
                                                        \label{eq:intro-alpha}\\
 M_\lambda(\zeta)
   &=\mathsf R_{\alpha_\lambda(\zeta)}
     \begin{pmatrix}\sqrt{1+\rho_*}&0\\0&\sqrt{1-\rho_*}\end{pmatrix}
     \mathsf R_{-\alpha_\lambda(\zeta)},               \label{eq:intro-M}
\end{align}
where $\mathsf R_\alpha$ is the planar rotation through $\alpha$, so that
$M_\lambda(\zeta)$ is a symmetric matrix with eigenvalues
$\sqrt{1\pm\rho_*}$ whose principal axis turns by $\alpha_\lambda(\zeta)$
around the circle. We set $\iota(a,b)=ae_x+be_z$ and
$v_{0,\lambda}(y,\zeta)=\iota M_\lambda(\zeta)y$. When the axis is straightened, that
is, in the limit $N\to\infty$, the map
$(y,\zeta)\mapsto L\zeta e_y+v_{0,\lambda}(y,\zeta)$, with $\zeta\in\R$, parametrizes an exact equilibrium.
Its pressure surfaces have elliptical cross-sections whose orientation
varies along the axis, and its nonconstant field lines are the images
under this map of the circles $|y|=r>0$ at fixed $\zeta$. The parameter $\rho_*$ governs the eccentricity of the
ellipse, $\alpha_0$ represents the mean tilt of the major axis of the ellipse measured from the plane of the torus, and $\delta_*$ controls the amplitude of the twisting. We refer the reader to  the associated GitHub repository for an interactive plot that may be used to visualize the geometry. Up to
normalization, these maps describe a subfamily of solutions constructed by Salat and Kaiser
\cite{SalatKaiser1995}. \Cref{thm:main} says that these solutions can be bent into a
closed torus of major radius $NL$ while keeping the pressure surfaces
nested and the field smooth. Indeed, for every small curvature parameter $\eps$,
the proof constructs a $2\pi$-periodic map $v_{\eps,\lambda}(y,\zeta)$
close to $v_{0,\lambda}$. Setting $\eps=1/N$, the equilibrium is
\begin{align}
 X_{N,\lambda}(y,\phi)
   &=R_\phi\bigl(NLe_x+v_{\eps,\lambda}(y,N\phi)\bigr),
                                                        \label{eq:intro-X}\\
 P_{N,\lambda}\circ X_{N,\lambda}
   &=\psi_a-|y|^2,
 &B_{N,\lambda}\circ X_{N,\lambda}
   &=\pa_\theta X_{N,\lambda},                         \label{eq:intro-BP}
\end{align}
with $\pa_\theta=-y_2\pa_{y_1}+y_1\pa_{y_2}$. In the normalized chart,
we write the perturbed map as its linear Taylor term in $y$ at the axis
$y=0$ plus a remainder,
\begin{equation}\label{eq:intro-jet}
 v_{\eps,\lambda}(y,\zeta)
 =P_{M_\lambda,\tau_{\eps,\lambda}}(\zeta)y
   +\widetilde v_{\eps,\lambda}(y,\zeta),
 \qquad j_y^1\widetilde v_{\eps,\lambda}|_{y=0}=0.
\end{equation}
Here,
\[
 P_{M,\tau}=a(\tau)\iota M+e_y\otimes\tau,
 \qquad a(\tau)=\sqrt{1-|\tau|^2/2},
\]
and the unknown function $\tau=\tau_{\eps,\lambda}:\T\to\R^2$ determines the component of this
linear term along the straightened axis, namely $(\tau(\zeta)\cdot y)e_y$.
The condition $j_y^1\widetilde v_{\eps,\lambda}|_{y=0}=0$ means that
$\widetilde v_{\eps,\lambda}$ and its first derivatives in $y$ vanish at $y=0$.
We will also prove that, uniformly for $\lambda\in I$,
\begin{equation}\label{eq:intro-Oeps}
 \norm{\tau_{\eps,\lambda}}_{C^2}
 +\norm{\widetilde v_{\eps,\lambda}}_{C^2}
 \le C|\eps|.
\end{equation}
The term $\cos\zeta$ in \eqref{eq:intro-alpha} fixes the period of the
tilt at $2\pi/N$ in the variable $\phi$. The term $\lambda\sin2\zeta$,
which is odd in $\zeta$, distinguishes the equilibrium from its mirror
images and from nearby members of the family. The vanishing condition in \eqref{eq:intro-jet} makes the Hessian of the pressure
on the magnetic axis independent of $\widetilde v$. Restricted to the normal
plane of the axis, its inverse can be normalized to remove the dependence
on $\tau$. This normalized inverse Hessian determines
$\alpha_\lambda(N\phi)$ modulo $\pi$, which gives the exact symmetry group.

\subsection{Strategy of the proof}\label{ss:strategy}

The proof has three parts. First, we reduce \eqref{eq:mhs} to a problem
on one period, posed on a fixed domain with the
curvature as a parameter. Second, we invert the linearization of this
problem, with estimates that are uniform on an open set of states. Third,
we apply the Nash-Moser theorem to solve the nonlinear equations, close the
cell into a torus, and compute the symmetry group. Most of
the paper is devoted to the second part of the argument.

\textbf{Choosing coordinates.} We look for an immersion
$X=X(\psi,\theta,\zeta)$ with $B\circ X=X_\theta$ and $P\circ X=\psi$, so
that pressure surfaces and field lines are built into the coordinates.
This makes \eqref{eq:mhs} equivalent to three scalar equations for $X$
(see \cref{sec:fixed-cell}). Two of these equations are in conservation form, and
the third says that the Jacobian of $X$ is independent of $\theta$. To
isolate one period of the field, we write
\[
 X(\psi,\theta,\zeta)
 =R_{\eps\zeta}\Bigl(\frac{L}{\eps}e_x+v(\psi,\theta,\zeta)\Bigr),
\]
where $\eps$ is a curvature parameter. Although $L/\eps$ diverges as $\eps\to0$, the equations for $v$ remain
well defined at $\eps=0$, since the
curvature enters only through $D_\eps v=v_\zeta+\eps e_z\times v+Le_y$,
which replaces $v_\zeta$. The unknown $v$ lives on the fixed domain
$Q=\overline{\D^2}\times\T$, and we solve for all small $\eps$ before
setting $\eps=1/N$ and joining $N$ cells into a torus. This is why
\cref{thm:main} holds for every large $N$.

At $\eps=0$, the equations are solved exactly by every
$v=\iota M(\zeta)y$ with $M$ a smooth periodic $2\times2$ matrix
satisfying $\det M>0$ and $\tr(M^TM)=2$ (see
\cref{lem:flat-affine-seeds}). We perturb
this infinite-dimensional family rather than a single state,
and we choose the subfamily \eqref{eq:intro-M} so that the symmetry
group can be read off at the end. The degeneracies of the problem, namely
reparametrizations of the cell and motion along the axis, both lie
in the kernel of the linearization. We remove them with the normalized
chart \eqref{eq:intro-jet}, two gauge conditions, and a condition on the
outer boundary (see \cref{sec:tame-problem}).

\textbf{Cell frequency.} The variable $\zeta\in\T$ parametrizes the
period cell, and we expand every unknown in a Fourier series in it,
\[
 u(y,\zeta)=\sum_{n\in\Z}\widehat u_n(y)\,e^{in\zeta}.
\]
We call the index $n$ the \emph{cell frequency}. Once $\eps=1/N$, the
\emph{cell variable} is $\zeta=N\phi$, so $e^{in\zeta}$ is the toroidal mode
$e^{inN\phi}$. Thus, cell frequency $n$ corresponds to toroidal mode
number $nN$, and the modes with $n=0$ are the axisymmetric ones. We write
$\langle n\rangle=(1+n^2)^{1/2}$ and $\Lambda=\langle D_\zeta\rangle$.
The operator $\Lambda$ multiplies the $n$-th Fourier coefficient by
$\langle n\rangle$. Since the coefficients of the linearized operator
depend on $\zeta$ through $M$, different cell frequencies are coupled. Nevertheless, the whole linear analysis is organized
by the size of $\langle n\rangle$. A cutoff $2^J$ separates the
\emph{high frequencies} $|n|\ge2^J$, which we treat by energy estimates
and a Neumann series, from the finitely many \emph{low frequencies},
which we treat by comparison with an explicit reference.

The cell frequency also determines the scale at which we split the disk.
We write $y=re^{i\theta}$. In polar coordinates, the equations have
singular coefficients at the axis $r=0$. The transition from the behavior
near the axis to a radial evolution occurs when $r\langle n\rangle\sim1$.
Accordingly, at cell frequency $n$ we split the disk into a \emph{cap}
of radius comparable to $\langle n\rangle^{-1}$ around the axis, where
the equations are of Bessel type, and the complementary \emph{annulus},
where they are a first-order evolution in $r$ with $\langle n\rangle$ as
the large parameter. The cap shrinks as the frequency grows, and the
analysis on it is a rescaled version of a single fixed problem.

\textbf{The linearized problem.} Three difficulties arise in constructing
a uniform inverse. First, the cap and the annulus require different
methods, and we must match the solutions at their interface with
estimates uniform in the cell frequency $n$. Second, the system on the
annulus has three parts that require different estimates. At a constant
ellipse, the \emph{center part} consists of two coupled oscillators whose
speeds depend on the ellipse parameter $\rho$ and coincide at the circle
$\rho=0$. We estimate them together so that the bounds remain uniform
as the ellipse approaches a circle. The \emph{affine part} comes from
the solutions at zero curvature that are linear in the disk variable $y$,
of the form $v(y,\zeta)=\iota M(\zeta)y$. Its equations allow linear
growth in $r$.
The remaining part is elliptic and satisfies a coercive estimate.
Finally, the source must satisfy compatibility conditions, and converting
it from Cartesian variables to the components used in these estimates
requires differentiation. This conversion costs a fixed number $\ell$
of derivatives, which must be accounted for in the inverse estimate.

\begin{figure}[htbp]
\centering
\includegraphics[width=\textwidth]{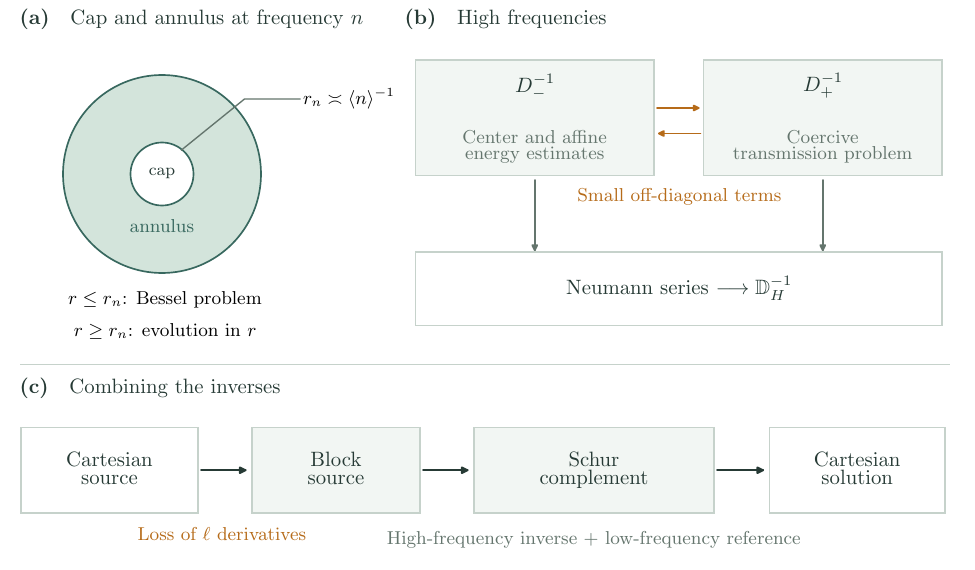}
\caption{Construction of the uniform inverse.
(a) We match the cap and annulus at $r_n\asymp\langle n\rangle^{-1}$.
(b) A Neumann series combines the diagonal inverses and the small
off-diagonal terms to give the inverse at high frequencies.
(c) A Schur complement joins this inverse to the low-frequency reference.
We convert the Cartesian source to block variables once, before applying
the inverse, to avoid incurring this loss of derivatives multiple times.}
\label{fig:inverse-scheme}
\end{figure}

\textbf{The uniform inverse.} The construction is summarized in
\cref{fig:inverse-scheme}. To control the growth of the center and affine
solutions on the annulus, we introduce an exponential weight in $r$
whose strength grows with the frequency. We work in fixed spaces for
the unknowns and the sources, equipped with the analytic weight
\[
 W_\gamma(r,n)
 =\exp\Bigl\{\sigma_0\langle n\rangle
  -\gamma\bigl(\sqrt{1+r^2\langle n\rangle^2}-1\bigr)\Bigr\}.
\]
Here, $0<\gamma<\min\{\sigma_0,1\}$. On the axis, $W_\gamma$ is the weight
$e^{\sigma_0\langle n\rangle}$ of
analyticity in $\zeta$ on a strip of width $\sigma_0$. Outside the cap,
the exponent decreases like $\gamma r\langle n\rangle$, so the analytic
width shrinks linearly in $r$ and, after applying the weight, the evolution in
$r$ acquires a damping of size $\gamma\langle n\rangle$. Near the reference state, this damping controls the growth on the annulus.
On the cap, the ratio $W_\gamma(r,n)/W_\gamma(0,n)$ is bounded above and
below by positive constants independent of $n$, so the weight does not
alter the cap estimates except by fixed factors.

We rewrite the unknowns and sources in \emph{block variables}, grouping
the center and affine components together and the elliptic components
separately. Each component has a weight chosen for its equation.
In these variables, the high-frequency part of the linearized operator
is a $2\times2$ matrix
\[
 \mathbb D_H=
 \begin{pmatrix}
  D_-&E_{-+}\\ E_{+-}&D_+
 \end{pmatrix}.
\]
The block $D_-$ contains the center and affine equations. To construct
its inverse, we match the cap and annulus solutions (see
\cref{app:minus:cap,app:minus:annulus}). On the cap, we start from the explicit inverse at
the circle, obtained from Bessel equations for the center components and
integral formulas for the affine components, with prescribed data on
the axis. The cap solution supplies the initial data for the radial
equations on the annulus. There, we estimate all four center modes
together using a positive quadratic energy, so the bounds remain uniform
even when the two oscillator speeds coincide. A separate, modified
energy controls the affine components (see \cref{app:minus:annulus}).
To combine these estimates, we use a change of variables that removes
the leading coupling between the center and affine components (see
\eqref{app:minus:eq:sylvester-transform}). The block $D_+$ contains the elliptic part. We
construct its inverse by solving on the cap and annulus and matching
the values and conormal derivatives at their interface (see
\cref{app:positive}). We correct these inverses for the full operator and estimate the coupling
between its diagonal blocks. The highest-order terms in the
off-diagonal blocks cancel, and the remaining blocks vanish at the
circular reference state and are small nearby (see
\cref{app:macro-cross}). A Neumann series then combines the diagonal
inverses to give the inverse of $\mathbb D_H$. We handle the low
frequencies by comparison with the inverse at a constant ellipse (see
\cref{app:lowref}). A Schur complement combines this with the
high-frequency inverse (see \cref{app:global}). Importantly, we only convert the
Cartesian source to block variables once, before applying the block
inverse, to avoid incurring this loss of derivatives multiple times.
The resulting two-sided inverse satisfies tame estimates on an open set
of states (see \cref{thm:uniform-current-inverse}). The constants are fixed in
\cref{sec:inverse:acyclic-choice}, with the frequency cutoff chosen
before the ellipse parameter $\rho_*>0$.

\textbf{The Nash-Moser scheme and the resulting geometry.} We apply
Hamilton's Nash-Moser theorem
\cite{Hamilton1982,Moser1966b,Moser1966a,Nash1956} with $\eps$ and $\lambda$
as parameters. The inverse estimates established above allow us to
solve the nonlinear equations for $|\eps|<\eps_H$ and $\lambda$ in a
compact interval, obtaining solutions that depend smoothly on both
parameters (see \cref{sec:hamilton}). The correction vanishes at
$\eps=0$, since the unperturbed family already solves the equations. Smooth dependence on
$\eps$ therefore gives \eqref{eq:intro-Oeps}, uniformly for $\lambda$
in the compact interval. At $\eps=1/N$, we bound the
Jacobian of $X_{N,\lambda}$ away from zero by a perturbation estimate and
obtain $|B|\ge c|y|$. Together with $B=0$ on the axis, this shows that
the field vanishes exactly there. We remark, however, that global
injectivity does not follow from the Jacobian bound, since the fast
variable $N\phi$ could produce distant self-intersections. To establish global
injectivity, we use
coordinates around the axis and allow the angular coordinate to range
over $\R$. In these coordinates, the map is close to the identity
together with its first derivatives, which rules out self-intersections.
The major radius $NL$ compensates for the factor $N$ from differentiation
of the fast variable (see \cref{subsec:geometry-injectivity}). To determine the
symmetry group, we use the normalized chart to compute
\[
 \cS_\phi
 :=\frac{(-H_\phi)^{-1}}{\tr((-H_\phi)^{-1})}
 =\frac12M_\lambda(N\phi)M_\lambda(N\phi)^T,
\]
where $H_\phi$ is the Hessian of the pressure restricted to the normal
plane of the axis, written in the radial and vertical basis $(R_\phi e_x,e_z)$. The tensor
$\cS_\phi$ is independent of $\tau$ and $\widetilde v$. An isometry
preserving the equilibrium preserves the axis. It therefore acts on the
axis by $\phi\mapsto\sigma\phi+c$ with $\sigma=\pm1$ and preserves or
flips $e_z$. Preservation of this tensor imposes a relation on the angle
$\alpha_\lambda$. Averaging this relation excludes a flip of $e_z$ because
$2\alpha_0\notin\pi\Z$. The first harmonic $\cos\zeta$ then forces
$c\in(2\pi/N)\Z$, and the odd second harmonic $\sin2\zeta$ excludes
$\sigma=-1$. This proves \eqref{eq:main-stabilizer}. The same tensor
shows that distinct values of $\lambda$ give inequivalent equilibria.

\subsection{Outline of the paper}\label{ss:outline}
The paper follows the three parts of the strategy described above.

The first part, \crefrange{sec:fixed-cell}{sec:fixed-map}, sets up the
problem on one cell. In \cref{sec:fixed-cell}, we write \eqref{eq:mhs} in
flux coordinates, pass to one rotating cell with the curvature $\eps$ as
a parameter, integrate the equations once to obtain a first-order system,
and exhibit the family of exact solutions at $\eps=0$. In
\cref{sec:tame-problem}, we introduce the normalized chart, specify the
target space of $H_\eps$ through its compatibility conditions, and
choose the scale of spaces used in our estimates. In
\cref{app:fixed-graphs}, we construct extension operators, projections,
and smoothing operators, and identify the varying function spaces with
one fixed space. In
\cref{sec:fixed-map}, we define the nonlinear map on the resulting fixed
spaces, verify that the family at zero curvature consists of exact zeros,
and compute the linearization $L_{\mathbf b}$ at a general state
$\mathbf b$.

The second part, \crefrange{sec:block}{sec:inverse}, inverts
$L_{\mathbf b}$ and is the bulk of the paper. In \cref{sec:block}, we
introduce the block variables, the transfer maps between them and the
Cartesian variables, and the augmented block operator
$\mathbb D_{\mathbf b}$. This is also where the conversion of Cartesian
sources is performed and its loss of derivatives is recorded. In
\cref{app:reference}, we construct the inverse at a fixed ellipse and
analyze the circular case. \Cref{app:minus,app:positive} construct
inverses at high frequencies for the center and affine equations and for
the elliptic boundary problem, respectively. In \cref{app:macro-cross},
we estimate the off-diagonal blocks and correct these inverses to obtain
inverses of the diagonal blocks of the full operator. In
\cref{app:lowref}, we invert the low-frequency part at the constant
ellipse, and in \cref{app:global} we assemble the inverse of
$\mathbb D_{\mathbf b}$ by a Neumann series and a Schur complement,
estimate its derivatives with respect to the state, and return to real
Cartesian variables.
In \cref{sec:inverse}, we fix the constants in their required order and
prove \cref{thm:uniform-current-inverse}, the uniform inverse theorem for
$L_{\mathbf b}$.

The third part, \cref{sec:hamilton,sec:geometry}, produces the
desired equilibria. In \cref{sec:hamilton}, we verify the hypotheses of the
Nash-Moser theorem and obtain a smooth branch of solutions in the
curvature and the parameter $\lambda$. In \cref{sec:geometry}, we set
$\eps=1/N$, close the cell into a torus, bound the Jacobian away from zero,
prove global injectivity, reconstruct the magnetic field and pressure,
verify the MHS equations, and show that the regular pressure levels are
nested tori. We then determine the symmetry group, show that distinct
values of $\lambda$ give inequivalent equilibria, and complete the proof
of \cref{thm:main}. An explanation of the Lean formalization of this theorem is
given in \Cref{app:lean}.

The article is designed to be modular. With the nonlinear setup and
parameters fixed, \cref{thm:uniform-current-inverse} supplies the inverse
and tame estimates used in the Nash-Moser argument of \cref{sec:hamilton}.
A reader interested in the linear analysis may take the block operator
of \cref{sec:block} as the starting point. \Cref{app:reference,app:lowref}
contain the explicit computations at the constant ellipse that give the
inverse formulas and estimates used later.

\subsection*{Acknowledgments}

JGS was partially supported by the NSF through Grants DMS-2247537, DMS-2434314 and DMS-2554957.
LL~is grateful to the Azrieli Foundation for the award of an Azrieli Fellowship and acknowledges the support of this research by ISF Grant No.~854/25.

\section{Coordinates and equations}
\label{sec:fixed-cell}

In this section, we derive the equations on one period cell and describe
the explicit family of solutions at zero curvature.

Throughout, \(e_x,e_y,e_z\) denote the standard basis vectors of \(\R^3\).
We write \(R_\varphi\) for rotation about \(e_z\) through the angle
\(\varphi\), and \(\mathsf R_\alpha\) for rotation in \(\R^2\) through
the angle \(\alpha\). We set
\[
 Aq=e_z\times q,\qquad q\in\R^3,
\]
so that \(\partial_\varphi R_\varphi=AR_\varphi\).
All angular variables have period \(2\pi\).

\subsection{Flux coordinates}

We first work away from the magnetic axis. Let
\(X=X(\psi,\theta,\zeta)\) be a smooth immersion, and locally prescribe
\[
 B\circ X=X_\theta,\qquad P\circ X=\psi,
 \qquad J=\det(X_\psi,X_\theta,X_\zeta).
\]
In these \emph{flux coordinates}, \(\psi\) labels pressure surfaces and
\(\theta\) parametrizes field lines. Since \(J\neq0\), the vectors
\(X_\psi,X_\theta,X_\zeta\) form a basis of
\(\R^3\), so the \emph{force imbalance}
\(\cQ=B\times\curl B+\grad P\) vanishes if and only if its inner
product with each basis vector vanishes.  Combining the identity
\[
 B\times\curl B=\grad\frac{|B|^2}{2}-(B\cdot\grad)B
\]
with the chain rule
\(\bigl((B\cdot\grad)f\bigr)\circ X=\pa_\theta(f\circ X)\),
we obtain
\begin{align}
 (\cQ\circ X)\cdot X_\theta&=0, \label{eq:fixed-Qtheta}\\
 (\cQ\circ X)\cdot X_\psi&=-F_1,\label{eq:fixed-Qpsi}\\
 (\cQ\circ X)\cdot X_\zeta&=F_2,\label{eq:fixed-Qzeta}
\end{align}
where
\begin{align}
 F_1(X)&=X_\psi\cdot X_{\theta\theta}
        -X_\theta\cdot X_{\psi\theta}-1,\label{eq:fixed-F1}\\
 F_2(X)&=X_\theta\cdot X_{\theta\zeta}
        -X_\zeta\cdot X_{\theta\theta}.\label{eq:fixed-F2}
\end{align}
Moreover,
\begin{equation}
 (\diver B)\circ X=J^{-1}\pa_\theta J,                \label{eq:fixed-div}
\end{equation}
so the magnetohydrostatic equations in these coordinates can be written as
\begin{equation}
 F_1=F_2=F_3=0,
 \qquad F_3:=\pa_\theta J.                            \label{eq:fixed-flux-system}
\end{equation}
The boundary condition \(B\cdot n=0\) holds automatically when the
physical boundary is a level set of \(\psi\), since \(X_\theta\) is tangent
to that level set. To obtain a first-order system, we write the
two force equations in conservation form,
\begin{align}
 F_1&=\pa_\theta(X_\psi\cdot X_\theta)
       -\pa_\psi|X_\theta|^2-1,\label{eq:fixed-conservation1}\\
 F_2&=\pa_\zeta|X_\theta|^2
       -\pa_\theta(X_\theta\cdot X_\zeta).
                                                        \label{eq:fixed-conservation2}
\end{align}

We now bend one period cell around a circle. We fix \(L>0\), with
\(2\pi L\) the length of the reference axis within one cell. For 
\(\eps\neq0\), we take \(v\) to be \(2\pi\)-periodic in \(\zeta\) and set
\begin{equation}
 X_{\eps,v}(\psi,\theta,\zeta)
 =R_{\eps\zeta}\left(\frac{L}{\eps}e_x
                       +v(\psi,\theta,\zeta)\right),   \label{eq:fixed-cell-ansatz}
\end{equation}
where \(\zeta\) is allowed to range over \(\R\). Increasing \(\zeta\) by
\(2\pi\) rotates the image through \(2\pi\eps\). We also set
\begin{equation}
 D_\eps v=v_\zeta+\eps Av+Le_y.                       \label{eq:fixed-Deps}
\end{equation}
Since \(AR_\varphi=R_\varphi A\) and \(Ae_x=e_y\), we obtain from
\eqref{eq:fixed-cell-ansatz} the identities
\[
 X_\psi=R_{\eps\zeta}v_\psi,\qquad
 X_\theta=R_{\eps\zeta}v_\theta,\qquad
 X_\zeta=R_{\eps\zeta}D_\eps v.
\]
Thus, \eqref{eq:fixed-flux-system} becomes
\(F_{1,\eps}=F_{2,\eps}=F_{3,\eps}=0\), where
\begin{align}
 F_{1,\eps}(v)
 &=v_\psi\cdot v_{\theta\theta}
   -v_\theta\cdot v_{\psi\theta}-1,\label{eq:fixed-cell-F1}\\
 F_{2,\eps}(v)
 &=v_\theta\cdot v_{\theta\zeta}
   -D_\eps v\cdot v_{\theta\theta},\label{eq:fixed-cell-F2}\\
 F_{3,\eps}(v)
 &=\pa_\theta\det(v_\psi,v_\theta,D_\eps v),
                                                        \label{eq:fixed-cell-F3}
\end{align}
on the cell, away from the magnetic axis. The equations remain well
defined at \(\eps=0\), although \eqref{eq:fixed-cell-ansatz} does not.
At \(\eps=0\), they are the equations for the map
\(X_{0,v}=L\zeta e_y+v\). For \(\eps=1/N\) with \(N\in\mathbb N\),
the map is periodic after \(N\) cells. Global injectivity of the
constructed solutions is proved in \cref{subsec:geometry-injectivity}.
Notably, we will solve the system for all small \(\eps\) and only restrict to
the discrete values \(\eps=1/N\) at the very end of the argument.

\subsection{A first-order system}

We use Cartesian coordinates \(y=(y_1,y_2)\) on the disk and set
\[
 y=r(\cos\theta,\sin\theta),\qquad
 \psi=\psi_a-r^2,\qquad 0\le r\le1,
\]
where \(\psi_a\) is the pressure on the axis. The unknown \(v\) is now
viewed as a map on \(Q=\overline{\D^2}\times\T\). For a scalar or vector
function \(f\), we define its \emph{angular average} by
\[
 (\Pi f)(y,\zeta)=\frac1{2\pi}\int_0^{2\pi}
 f(\mathsf R_s y,\zeta)\dd s.
\]
Thus, \(\Pi\) averages at fixed \(r\) and \(\zeta\). We set
\[
 E=|v_\theta|^2,\qquad
 g_{\psi\theta}=v_\psi\cdot v_\theta,\qquad
 g^\eps_{\theta\zeta}=v_\theta\cdot D_\eps v,\qquad
 J_\eps=\det(v_\psi,v_\theta,D_\eps v),
\]
and introduce an auxiliary scalar \(w\) to integrate the force equations
in \(\theta\). Note that its additive freedom is fixed by \(\Pi w=0\). 

For fixed
\(\eps\), we define
\begin{align}
 G_0(v,w)&=w_\theta-(E-r^2),\label{eq:fixed-G0}\\
 G_1(v,w)&=(I-\Pi)(g_{\psi\theta}-w_\psi),\label{eq:fixed-G1}\\
 G_2(v,w)&=(I-\Pi)(g^\eps_{\theta\zeta}-w_\zeta),\label{eq:fixed-G2}\\
 G_3(v,w)&=(I-\Pi)J_\eps.\label{eq:fixed-G3}
\end{align}
The equations \(G_0=G_1=G_2=G_3=0\) involve only first derivatives of
\((v,w)\). Note that only \(G_3\) contains terms of degree three. The following
lemma shows that introducing \(w\) does not change the solutions.

\begin{lemma}
\label{lem:fixed-integrated-equivalence}
Let \(\eps\in\R\) and let \(v\in C^\infty(Q;\R^3)\). The following
are equivalent.
\begin{enumerate}[label=\textup{(\roman*)}]
\item \(F_{1,\eps}(v)=F_{2,\eps}(v)=F_{3,\eps}(v)=0\) on \(0<r\le1\).
\item There is a unique \(w\in C^\infty(Q)\) with \(\Pi w=0\) such that
\(G_0(v,w)=G_1(v,w)=G_2(v,w)=G_3(v,w)=0\) on \(0<r\le1\).
\end{enumerate}
\end{lemma}

\begin{proof}
Suppose first that \(G_0=G_1=G_2=G_3=0\). The equation \(G_0=0\)
gives \(w_\theta=E-r^2\), from which we obtain
\begin{align*}
 \pa_\theta(g_{\psi\theta}-w_\psi)&=F_{1,\eps},\\
 \pa_\theta(g^\eps_{\theta\zeta}-w_\zeta)&=-F_{2,\eps}.
\end{align*}
Here, we used that \(\pa_\psi r^2=-1\) and the identity 
\[
 v_\theta\cdot\pa_\theta(D_\eps v)
 =v_\theta\cdot\bigl(v_{\theta\zeta}+\eps Av_\theta\bigr)
 =v_\theta\cdot v_{\theta\zeta},
\]
which is a consequence of the skew-symmetry of \(A\).  Hence, \(G_1=G_2=0\) gives
the two force equations, and \(G_3=0\) gives \(F_{3,\eps}=0\).

Conversely, suppose that \(F_{1,\eps}=F_{2,\eps}=F_{3,\eps}=0\).
The first conservation identity becomes
\(F_{1,\eps}=\pa_\theta g_{\psi\theta}-\pa_\psi E-1\), so averaging in
\(\theta\) gives \(\pa_\psi\Pi E=-1\).  Since \(v\) is smooth in \(y\), we
have \(v_\theta=0\) and hence \(E=0\) on the axis. Integrating from
\(\psi=\psi_a\), we obtain \(\Pi E=\psi_a-\psi=r^2\). We therefore set
\[
 w(y,\zeta)=\frac{1}{2\pi}\int_0^{2\pi}
 s\,(E-r^2)(\mathsf R_sy,\zeta)\dd s.
\]
This function is smooth on \(Q\). Since \(\Pi(E-r^2)=0\), integration
by parts gives \(w_\theta=E-r^2\), and averaging the formula gives
\(\Pi w=0\). Any two such functions differ by a function independent
of \(\theta\), so the mean zero condition gives uniqueness.
The preceding derivative identities and \(F_{3,\eps}=0\) show that
\(g_{\psi\theta}-w_\psi\), \(g^\eps_{\theta\zeta}-w_\zeta\), and
\(J_\eps\) are independent of \(\theta\). Applying \(I-\Pi\) therefore
gives \(G_1=G_2=G_3=0\).
\end{proof}

The equations \(G_1=0\) and \(G_3=0\) still contain the derivative
\(\partial_\psi\), which is singular at the axis. To obtain equations
that are smooth in \(y\), we set
\begin{equation}
 \cD=y_1\pa_{y_1}+y_2\pa_{y_2},\qquad
 \cR=-y_2\pa_{y_1}+y_1\pa_{y_2},                    \label{eq:fixed-disk-generators}
\end{equation}
so that when \(r>0\) we have
\[
 \pa_\theta=\cR,\qquad \pa_\psi=-\frac{1}{2r^2}\cD.
\]
We multiply the two singular equations by the appropriate factors of
\(r^2\) and define
\begin{align}
 H_{0,\eps}(v,w)
 &=\cR w-\bigl(|\cR v|^2-r^2\bigr),\label{eq:fixed-H0}\\
 H_{1,\eps}(v,w)
 &=(I-\Pi)\bigl(\cD w-\cD v\cdot\cR v\bigr),\label{eq:fixed-H1}\\
 H_{2,\eps}(v,w)
 &=(I-\Pi)\bigl(\cR v\cdot D_\eps v-w_\zeta\bigr),\label{eq:fixed-H2}\\
 H_{3,\eps}(v,w)
 &=(I-\Pi)\det(\cD v,\cR v,D_\eps v).               \label{eq:fixed-H3}
\end{align}
We unify the above by defining \(H_\eps=(H_{0,\eps},H_{1,\eps},H_{2,\eps},H_{3,\eps})\).
On \(r>0\), these equations satisfy
\begin{equation}
 H_{0,\eps}=G_0,\qquad H_{1,\eps}=2r^2G_1,\qquad
 H_{2,\eps}=G_2,\qquad H_{3,\eps}=-2r^2G_3.           \label{eq:fixed-HG}
\end{equation}
The components \(G_1\) and \(G_3\) need not extend smoothly to the axis,
whereas every component of \(H_\eps\) does. Thus, \(H_\eps=0\) is a
system on the whole disk with the same solutions as the original
equations on \(r>0\). When the reconstructed map is smooth and has
nonzero Jacobian in Cartesian coordinates, the physical MHS equations also hold on the
axis by continuity.

\subsection{A family of solutions at zero curvature}

At zero curvature, the equations admit an explicit family of solutions
that are linear in \(y\), including the straight-axis family of Salat and
Kaiser \cite{SalatKaiser1995}. We recall the inclusion
\(\iota:\R^2\to\R^3\),
\[
 \iota(q_1,q_2)=q_1e_x+q_2e_z,
\]
whose ordering fixes the sign of the Jacobian.  For a
smooth periodic matrix \(M:\T\to GL^+(2,\R)\), we set
\begin{equation}
 v_M(y,\zeta)=\iota M(\zeta)y.                         \label{eq:fixed-affine-seed}
\end{equation}
The following trace condition makes \(v_M\) an exact solution.\nopagebreak
\begin{lemma}
\label{lem:flat-affine-seeds}
Suppose that
\begin{equation}
 \tr(M^TM)=2                             \label{eq:fixed-trace-condition}
\end{equation}
for every \(\zeta\). Then \(F_{1,0}(v_M)=F_{2,0}(v_M)=F_{3,0}(v_M)=0\)
on \(r>0\). Let \(w_M\in C^\infty(Q)\) be the unique function with
\(\Pi w_M=0\) satisfying
\begin{equation}
 (w_M)_\theta=|(v_M)_\theta|^2-r^2.                \label{eq:fixed-wM-def}
\end{equation}
Then \(H_0(v_M,w_M)=0\) on the whole disk.
\end{lemma}

\begin{proof}
With \(e_r=(\cos\theta,\sin\theta)\), \(e_\theta=(-\sin\theta,\cos\theta)\),
\(rr_\psi=-1/2\), and \(M'=\pa_\zeta M\), we have
\begin{align*}
 (v_M)_\psi&=r_\psi\iota Me_r,&
 (v_M)_\theta&=r\iota Me_\theta,&
 (v_M)_{\theta\theta}&=-r\iota Me_r,\\
 (v_M)_\zeta&=r\iota M'e_r,&
 (v_M)_{\theta\zeta}&=r\iota M'e_\theta.
\end{align*}
Inserting these at \(\eps=0\), where \(D_0v_M=r\iota M'e_r+Le_y\), we obtain
\begin{align}
 F_{1,0}(v_M)
 &=\frac12\bigl(|Me_r|^2+|Me_\theta|^2\bigr)-1
   =\frac12\tr(M^TM)-1,\label{eq:fixed-affine-F1}\\
 F_{2,0}(v_M)
 &=r^2\bigl(Me_\theta\cdot M'e_\theta
             +M'e_r\cdot Me_r\bigr)
   =\frac{r^2}{2}\pa_\zeta\tr(M^TM),\label{eq:fixed-affine-F2}\\
 J_0(v_M)&=\frac{L}{2}\det M(\zeta),                  \label{eq:fixed-affine-J}
\end{align}
by orthonormality of \(e_r,e_\theta\) and by the identity
\(\det(\iota a,\iota b,e_y)=-\det[a\ b]\).
The first two lines vanish by \eqref{eq:fixed-trace-condition}, the last is
independent of \(\theta\), and the assertion for \(w_M\) now follows from
\cref{lem:fixed-integrated-equivalence}.
\end{proof}

The trace condition alone makes the equations vanish. To reconstruct
smooth physical fields, we also need \(M\) to be invertible. Our
assumption \(M(\zeta)\in GL^+(2,\R)\), together with compactness of
\(\T\), gives
\begin{equation}
 \inf_{\zeta\in\T}\det M(\zeta)>0.                  \label{eq:fixed-det-condition}
\end{equation}
For \(\zeta\in\R\), the map
\(X_{0,v_M}=L\zeta e_y+\iota M(\zeta)y\) is injective, since its
\(e_y\)-component determines \(\zeta\) and \(M(\zeta)\) is invertible.
Its Jacobian in Cartesian coordinates is \(-L\det M\), so it defines smooth physical
fields through the axis.

We use the following family in the construction. For \(0<\rho<1\), we set
\begin{align}
 \alpha_\lambda(\zeta)
 &=\alpha_0+\delta\bigl(\cos\zeta+\lambda\sin2\zeta\bigr),
                                                        \label{eq:fixed-alpha-lambda}\\
 M_\lambda(\zeta)
 &=\mathsf R_{\alpha_\lambda(\zeta)}
   \begin{pmatrix}\sqrt{1+\rho}&0\\0&\sqrt{1-\rho}\end{pmatrix}
   \mathsf R_{-\alpha_\lambda(\zeta)},               \label{eq:fixed-M-lambda}
\end{align}
so that
\begin{equation}
 \tr(M_\lambda^TM_\lambda)=2,\qquad
 \det M_\lambda=\sqrt{1-\rho^2}                      \label{eq:fixed-M-invariants}
\end{equation}
and both \eqref{eq:fixed-trace-condition} and \eqref{eq:fixed-det-condition}
hold identically.
We require \(\delta\neq0\), \(0<\lambda<1/2\), and
\(\alpha_0\notin(\pi/2)\Z\). The orientation \(\alpha_0\) is fixed at the beginning of the argument, and
the small parameter \(\delta\neq0\) is chosen only after the linear
analysis. As in \cref{ss:main}, the term \(\cos\zeta\) fixes the period
of the tilt at \(2\pi\) in \(\zeta\), or \(2\pi/N\) in \(\phi\), when
\(\eps=1/N\). The term \(\lambda\sin2\zeta\) excludes the reversal
\(\zeta\mapsto-\zeta\). Both are used to determine the symmetry group.

We finally compute the error caused by bending the cell. For nonzero
curvature, \(v_M\) still satisfies the two force equations. Indeed,
\(\eps Av_M=\eps r(Me_r)_1e_y\) is orthogonal to the range of \(\iota\),
which contains \((v_M)_{\theta\theta}\), so that
\(F_{1,\eps}(v_M)=F_{2,\eps}(v_M)=0\) while the Jacobian becomes
\begin{align}
 J_\eps(v_M)
 &=\frac12\det M(\zeta)
   \bigl(L+\eps r(M(\zeta)e_r)_1\bigr),
                                                        \label{eq:fixed-bent-J}\\
 F_{3,\eps}(v_M)
 &=\frac{\eps}{2}\det M(\zeta)\,
   \pa_\theta\bigl(r(M(\zeta)e_r)_1\bigr).          \label{eq:fixed-bent-residual}
\end{align}
Thus, only the divergence equation fails. Since
\(r(M(\zeta)e_r)_1=(M(\zeta)y)_1\), its residual is smooth through the
axis and is \(O(|\eps|)\) in every fixed \(C^k(Q)\) norm. This is the
error that the nonlinear argument must correct.

\section{Function spaces and constraints}\label{sec:tame-problem}
We now choose coordinates for the unknowns, specify the constraints on
these coordinates and on the right-hand sides, and define the norms used
in the analysis. The constraints remove the freedom described in
\cref{ss:strategy}, while the norms control both smoothness and
analyticity in the periodic variable.
\subsection{The normalized chart}\label{subsec:normalized-slice}

We let \(e_T=e_y\) be tangent to the limiting straight axis and let
\(\iota:\R^2\to e_T^\perp\) be the inclusion into the \((e_x,e_z)\)-plane.
The unperturbed solution at zero curvature is the pair given by
\[
 v_M(y,\zeta)=\iota M(\zeta)y,
 \qquad \tr(M^TM)=2,
\]
together with the auxiliary scalar \(w_M\) of \eqref{eq:fixed-wM-def}.
To allow the linear Taylor term of \(v\) to vary along the axis while
preserving its trace normalization, we introduce an unknown function
\(\tau:\T\to\R^2\). For \(|\tau(\zeta)|^2<2\), we define
\begin{equation}\label{eq:normalized-affine-map}
 a(\tau)=\sqrt{1-\frac{|\tau|^2}{2}},
 \qquad P_{M,\tau}=a(\tau)\iota M+e_T\otimes\tau.
\end{equation}
Since \(\iota^Te_T=0\) gives
\(P_{M,\tau}^TP_{M,\tau}=a(\tau)^2M^TM+\tau\tau^T\), we obtain
\begin{equation}\label{eq:normalized-affine-trace}
 \tr(P_{M,\tau}^TP_{M,\tau})=2.
\end{equation}
%so that the first equation satisfies its nonlinear mean compatibility
%condition.
%We require \(|\tau|^2<2\) pointwise and will choose the state
%ball in \cref{subsec:fixed-map} accordingly. 
We use the triple \((\tau,\widetilde v,s)\) as coordinates for the
unknown pair \((v,w)\), writing
\begin{equation}\label{eq:normalized-state}
 v=P_{M,\tau}(\zeta)y+\widetilde v(y,\zeta),
 \qquad w=w_M+s,
 \qquad j_y^1\widetilde v|_{y=0}=0.
\end{equation}
Here, $j_y^1u|_{y=0}$ denotes the value and first derivatives of $u$ in $y$
on the axis, which we call its \emph{first jet}. These are the coefficients
of its Taylor expansion through first order. Thus, $P_{M,\tau}(\zeta)y$
is the linear term in the Taylor expansion of $v$ at $y=0$, and
$\widetilde v$ and its first derivatives in $y$ vanish there.
At $\tau=0$, the derivative of $P_{M,\tau}$ in $\tau$ is
$e_T\otimes\dot\tau$, so a variation of $\tau$ changes this linear term
only along the axis direction $e_T$. Its component in the normal plane
is unchanged to first order in $\tau$.

With \(\Pi\) denoting the angular average in \(y\) and
\(\cR=-y_2\pa_{y_1}+y_1\pa_{y_2}\) as in \eqref{eq:fixed-disk-generators},
we fix the remaining freedom to change coordinates by imposing the two
\emph{gauge conditions}
\begin{equation}\label{eq:two-gauges}
 \Pi(\widetilde v\cdot\cR v_M)=0,
 \qquad
 \Pi(\widetilde v\cdot D_0v_M)=0,
\end{equation}
where we use \(D_0\) in place of \(D_\eps\) so that the constraints do not
depend on the curvature. We also impose the following condition on the
outer boundary:
\begin{equation}\label{eq:outer-slice}
 \Pi_{|m|\ge3}
 \left[e_r\cdot M(\zeta)^{-1}\iota^T
             \widetilde v(1,\theta,\zeta)\right]=0.
\end{equation}
Here, \(P_m\) projects onto the angular Fourier mode \(e^{im\theta}\),
and \(\Pi_{|m|\ge3}=\sum_{|m|\ge3}P_m\). Thus, the boundary condition
removes the modes \(|m|\ge3\) of the displayed radial component.
Finally, we require \(\Pi s=0\), so that \(w=w_M+s\) retains the
mean zero normalization of \cref{lem:fixed-integrated-equivalence}.

For fixed \(M\), we let \(\cX_M^\infty\) be the linear space of smooth
triples \(x=(\tau,\widetilde v,s)\) satisfying these constraints. With
\begin{equation}\label{eq:physical-slice-row}
 B_Mu:=\Pi_{|m|\ge3}
 \left[e_r\cdot M(\zeta)^{-1}\iota^Tu|_{r=1}\right],
\end{equation}
we put
\begin{equation}\label{eq:complete-domain-constraint}
 \mathscr C_M(\tau,\widetilde v,s)
 =\left(
 j_y^1\widetilde v|_{y=0},\ \Pi s,\
 \Pi(\widetilde v\cdot\cR v_M),\
 \Pi(\widetilde v\cdot D_0v_M),\ B_M\widetilde v
 \right),
\end{equation}
so that \(\cX_M^\infty=\ker\mathscr C_M\), while the tangential coefficient \(\tau\)
is unrestricted in this linear space. The nonlinear chart
\eqref{eq:normalized-state} is defined on the open subset where
\(|\tau(\zeta)|^2<2\) for every \(\zeta\).
Thus, the conditions on the axis, the two gauges, the angular mean of
\(s\), and the outer boundary are all included in the domain.% Their values
%will serve as auxiliary coordinates when we represent the restricted linear
%operator.

\subsection{The target space of $H_\eps$}\label{subsec:compatible-range}

The components of \(H_\eps(v,w)\) satisfy identities on the axis and
have prescribed angular means. We build these conditions into the target
space. We call a quadruple \(f=(f_0,f_1,f_2,f_3)\) a \emph{residual}
and introduce its \emph{quotient variables}, initially for \(y\ne0\), by
\begin{equation}\label{eq:quotient-variables}
 g_+=\frac{f_1+if_0}{\bar z},
 \qquad
 g_-=\frac{f_1-if_0}{z},
 \qquad
 g_3=\frac{f_3}{z\bar z},
 \qquad z=y_1+iy_2.
\end{equation}
We call \(f\) \emph{compatible} if these quotients are smooth in Cartesian
coordinates and
\begin{equation}\label{eq:range-means}
 \Pi f_1=\Pi f_2=\Pi f_3=0,
 \qquad
 \frac14\Delta_y\Pi f_0(0,\zeta)=0.
\end{equation}
We write \(\cY_{\rm comp}^\infty\) for the space of compatible
residuals.  We take \(Z=(g_+,g_-,g_3,h)\), with \(h=f_2\), as the range
coordinate and reconstruct
\begin{equation}\label{eq:range-reconstruction}
 J_YZ=\left(
 \frac{\bar zg_+-zg_-}{2i},
 \frac{\bar zg_++zg_-}{2},
 h,z\bar zg_3\right).
\end{equation}
Thus, passing from \(Z\) to \(f\) only multiplies by smooth factors.
Passing from \(f\) to \(Z\), however, requires the divisibility conditions
in \eqref{eq:quotient-variables}. We write
\begin{equation}\label{eq:compatible-extraction}
 \mathscr Q_Y(f_0,f_1,f_2,f_3)=(g_+,g_-,g_3,f_2)
\end{equation}
for the inverse of \(J_Y\) on such residuals. We measure regularity in
these quotient coordinates. In particular, the norm of \(f\) is the
product norm of \(\mathscr Q_Yf\), rather than the product norm of its
four original components. In \cref{subsec:range-projection}, we construct
a bounded projection that imposes \eqref{eq:range-means} in these
coordinates.

We now check that $H_\eps$ takes values in this target space when \((v,w)\) is
written in the chart \eqref{eq:normalized-state}. The components \(H_{1,\eps}\), \(H_{2,\eps}\), and \(H_{3,\eps}\)
carry the factor \(I-\Pi\), so their angular means vanish. For
the zeroth entry, we have \(\Pi f_0=r^2-\Pi|\cR v|^2\). Since \(\widetilde v\) vanishes
to second order on the axis, so does \(\cR\widetilde v\), and only the
linear Taylor term \(P_{M,\tau}y\) contributes to the coefficient of \(r^2\) in
\(\Pi|\cR v|^2\). That coefficient is
\(\tfrac12\tr(P_{M,\tau}^TP_{M,\tau})=1\) by
\eqref{eq:normalized-affine-trace}, so \(\Pi f_0=O(|y|^3)\). Its quadratic Taylor coefficient therefore
vanishes, as required by the last condition in \eqref{eq:range-means}.

It remains to see that the quotients \eqref{eq:quotient-variables} of
\(H_\eps(v,w)\) are smooth across the axis. With
\[
 \mathfrak m_1(v,w):=\Pi(\cD w-\cD v\cdot\cR v),
\]
and with the curvature subscript suppressed, we have
\begin{align}
 H_1+iH_0&=(\cD+i\cR)w
   -(\cD+i\cR)v\cdot\cR v+ir^2
   -\mathfrak m_1(v,w),\label{eq:literal-quotient-plus}\\
 H_1-iH_0&=(\cD-i\cR)w
   -(\cD-i\cR)v\cdot\cR v-ir^2
   -\mathfrak m_1(v,w),\label{eq:literal-quotient-minus}\\
 \det(\cD v,\cR v,D_\eps v)
 &=-2iz\bar z\det(\pa_zv,\pa_{\bar z}v,D_\eps v).
                                                        \label{eq:literal-quotient-J}
\end{align}
Here, \(\pa_z=\tfrac12(\pa_{y_1}-i\pa_{y_2})\) and
\(\pa_{\bar z}=\tfrac12(\pa_{y_1}+i\pa_{y_2})\), so
\(\cD+i\cR=2\bar z\pa_{\bar z}\) and \(\cD-i\cR=2z\pa_z\).
Since \(\mathfrak m_1\) is an angular average of a smooth function, it
is smooth in \((r^2,\zeta)\). It vanishes on the axis and is therefore
divisible by \(r^2=z\bar z\). The first two identities then give the
factors \(\bar z\) and \(z\) in \(H_1+iH_0\) and \(H_1-iH_0\).
The last identity gives a factor \(r^2\) in the determinant. Applying
\(I-\Pi\) preserves this factor because \(r^2\) is independent of
\(\theta\). Thus, these identities define the quotients \eqref{eq:quotient-variables} of \(H_\eps(v,w)\) on the axis by smooth
expressions, so \(H_\eps(v,w)\) is compatible whenever \((v,w)\) is written in the
chart \eqref{eq:normalized-state} and is smooth in \(y\).

\subsection{A fixed scale of spaces}\label{subsec:analytic-scale}

We use a weight that controls analyticity in \(\zeta\) and supplies
radial damping in the linear estimates. Let
\(\Lambda=\langle D_\zeta\rangle\), whose symbol is
\(\langle n\rangle=(1+n^2)^{1/2}\), and fix \(\sigma_0>0\) and
\(0<\gamma<\min\{\sigma_0,1\}\). For \(n\in\Z\), we set
\begin{equation}\label{eq:phase-weight}
 \Phi_\gamma(r,\langle n\rangle)
 =\sigma_0\langle n\rangle
 -\gamma\bigl(\sqrt{1+r^2\langle n\rangle^2}-1\bigr),
 \qquad W_\gamma=e^{\Phi_\gamma(r,\Lambda)}.
\end{equation}
We require \(\cR v_M\), \(D_0v_M\), and \(M^{-1}\) to be analytic in
\(\zeta\) on a strip wider than \(\sigma_0\). Note that our unperturbed solutions satisfy
this condition.
We define the norms in the interior by applying the relevant
Sobolev norms to \(W_\gamma u\). Thus, we set
\begin{equation}\label{eq:analytic-grade}
 \norm{u}_{\cA_\gamma^s}=\norm{W_\gamma u}_{H^s(Q)}.
\end{equation}
Here, \(W_\gamma\) acts on the \(n\)-th Fourier coefficient by
multiplication by \(W_\gamma(r,n)=e^{\Phi_\gamma(r,\langle n\rangle)}\).
We weight functions on the axis, including \(\tau\), by
\(e^{\sigma_0\Lambda}\), and outer boundary traces by
\(e^{\Phi_\gamma(1,\Lambda)}\), which is equivalent to
\(e^{(\sigma_0-\gamma)\Lambda}\). For the target space, we apply these
norms to the quotient coordinates \(Z\). The three properties of the
weight used below are
\begin{align}
 e^{(\sigma_0-\gamma r)\langle n\rangle}
 &\le W_\gamma(r,n)
 \le e^\gamma e^{(\sigma_0-\gamma r)\langle n\rangle},
                                                        \label{eq:phase-linear-equivalence}\\
 W_\gamma(r,n_1+n_2)
 &\le W_\gamma(r,n_1)W_\gamma(r,n_2),             \label{eq:phase-submultiplicative}\\
 \Phi_\gamma(0,\langle n\rangle)&=\sigma_0\langle n\rangle,
 \qquad \grad_y\Phi_\gamma(0,\langle n\rangle)=0.    \label{eq:phase-axis}
\end{align}
% AUTHOR QUERY: Prove the mixed Cartesian/cell weighted tame product bound
% below. Fourier submultiplicativity alone does not control derivatives of
% the y-dependent weight; specify the base index and any fixed loss.
We also use the following product estimate: For a fixed base index
\(s_b>3/2\) and \(s\ge s_b\), we have
\begin{equation}\label{eq:analytic-product}
 \norm{uv}_{\cA_\gamma^s}
 \le C_s\bigl(
 \norm{u}_{\cA_\gamma^s}\norm{v}_{\cA_\gamma^{s_b}}
 +\norm{u}_{\cA_\gamma^{s_b}}\norm{v}_{\cA_\gamma^s}
 \bigr).
\end{equation}
Note that the weight is smooth in \(|y|^2\), so
multiplication by it preserves vanishing Cartesian jets on the axis.
Its negative logarithmic radial derivative,
\begin{equation}\label{eq:phase-damping}
 d_\gamma(r,\langle n\rangle)
 =-\pa_r\Phi_\gamma
 =\frac{\gamma r\langle n\rangle^2}
        {\sqrt{1+r^2\langle n\rangle^2}},
\end{equation}
is the damping introduced when we apply the weight to an outward radial
equation. For \(r\langle n\rangle\ge1\), it is comparable to
\(\gamma\langle n\rangle\). We use this damping to control growth on the
annulus without changing the chosen spaces during the iteration.

\section{Relations between different spaces}\label{app:fixed-graphs}
We now construct projections that impose the constraints of
\cref{sec:tame-problem}, identify the domains for different \(M\), and
build smoothing operators that preserve the constraints. These
constructions allow us to apply the Nash-Moser theorem on fixed spaces.

Throughout this section, we fix a compact family of periodic matrix functions
\begin{equation}\label{eq:A-seed-family}
 M\in C^\infty(\T;GL(2,\R)),
 \qquad \tr(M^TM)=2,
 \qquad \norm{M}+\norm{M^{-1}}\le C_M,
\end{equation}
and every estimate below is locally uniform in this family. For the
weighted estimates, we also require the matrix functions to satisfy the analyticity
bounds on the strip of \cref{subsec:analytic-scale}. We first work in the
standard Sobolev spaces on $Q=\overline{\D^2}\times\T$, and we then
conjugate by the weight \eqref{eq:phase-weight}. We also fix the
projections
\begin{equation}\label{eq:A-exact-cell-split}
 H_J:=\mathbf 1_{\{|n|\ge 2^J\}}(D_\zeta),
 \qquad
 P_J:=I-H_J=\mathbf 1_{\{|n|<2^J\}}(D_\zeta),
\end{equation}
onto the high and low cell frequencies, where the integer \(J\) sets
the cutoff at \(2^J\). Superscripts indicate the fixed coordinate space
on which a cutoff acts. We write \(P_J^X,H_J^X\) for the domain and
\(P_J^Y,H_J^Y\) for the range.

For the unconstrained unknowns, we use the product space
\[
 \mathbb E^s=H^s(\T;\R^2)\times H^s(Q;\R^3)\times H^s(Q)
\]
for \((\tau,\widetilde v,s)\). In the range, we set
\[
 \mathbb B^s=J_Y\bigl((H^s(Q;\mathbb C))^4\bigr),
 \qquad
 \norm{f}_{\mathbb B^s}=\norm{\mathscr Q_Yf}_{(H^s)^4}.
\]
Thus, \(\cX_M^s\) and \(\cY_{\rm comp}^s\) are the subspaces of
\(\mathbb E^s\) and \(\mathbb B^s\) satisfying the respective
constraints. We use the same notation for their complexifications and
for the weighted spaces, specifying which norms are in use. The
superscript \(\infty\) denotes the intersection over all regularity
indices. Throughout the projection constructions, we take \(s>2\), so
that the first derivatives have well-defined traces on the axis.

We recall the terminology of Hamilton \cite{Hamilton1982} that we use
below. A \emph{graded Fr\'echet space} is a Fr\'echet space whose topology
is defined by an increasing sequence of norms $\norm{\cdot}_s$, such as
$\bigcap_s\cX_M^s$ with the norms of \cref{subsec:analytic-scale}.
We prove that our spaces are \emph{tame} in Hamilton's sense by realizing
them as ranges of bounded projections on standard spaces of smooth
functions. We also construct \emph{smoothing operators} \(S_\nu\),
\(\nu\ge1\), satisfying, whenever \(b\ge a\) are above the base index,
\begin{align}
 \norm{S_\nu q}_b&\le C_{a,b}\nu^{(b-a)_+}\norm q_a,
                                                        \label{eq:A-smoothing-high}\\
 \norm{(I-S_\nu)q}_a&\le C_{a,b}\nu^{a-b}\norm q_b.  \label{eq:A-smoothing-low}
\end{align}
Thus, $S_\nu q$ is as regular as we wish, at the price of a power of
$\nu$, and it approximates $q$ with an error that is small in a weaker
norm when $q$ is regular. A map $F$ between graded spaces is called
\emph{tame} if, on a neighborhood of each point of its domain, it satisfies
an estimate of the form
\[
 \norm{F(x)}_s\le C_s\bigl(1+\norm x_{s+r}\bigr)
\]
for all \(s\) above a fixed base index, with a loss of derivatives
\(r\) independent of \(s\). A smooth map is called \emph{smooth tame}
if every derivative is tame, viewed as a map of the base point and the
direction variables together. The
main result of this section may be stated as follows.

\begin{proposition}
\label{prop:fixed-graphs}
There is an index $s_*$ such that the following hold for every $s\ge s_*$,
with constants locally uniform in the family \eqref{eq:A-seed-family}.
\begin{enumerate}[label=\textup{(\roman*)}]
\item There are bounded projections
\[
 P_{X,M}:\mathbb E^s\longrightarrow\cX_M^s,
 \qquad
 P_Y:\mathbb B^s\longrightarrow\cY_{\rm comp}^s,
\]
and the derivatives of \(P_{X,M}\) in \(M\) of every fixed order are
bounded at the same regularity level. The constants may depend on the
regularity index and the number of parameter derivatives.
\item For every $M$ in the family and a fixed reference matrix function $M_*$ in it,
there is a linear isomorphism
\[
 T_M:\cX_{M_*}^\infty\xrightarrow{\ \sim\ }\cX_M^\infty,
\]
which is bounded with bounded inverse at every regularity level, depends
analytically on $M$ with every derivative tame, and commutes with complex
conjugation.
\item The constrained scales $\cX_M^s$ and $\cY_{\rm comp}^s$ admit
smoothing operators \(S_{\nu,M}^X\) and \(S_\nu^Y\) satisfying
\eqref{eq:A-smoothing-high}--\eqref{eq:A-smoothing-low}. Their
intersections over all regularity indices are tame Fr\'echet spaces.
\end{enumerate}
The same conclusions hold on the weighted analytic scales of
\cref{subsec:analytic-scale}. In particular, neither the projections nor
the smoothing operators lose derivatives.
\end{proposition}

We prove \cref{prop:fixed-graphs} in the next five subsections. We begin
with the extension operators from which the projections are built.

\subsection{Extension operators on the axis and on the boundary}

The constraints \eqref{eq:complete-domain-constraint} prescribe the first
jet on the axis and a trace on the outer boundary. To project onto them, we
need operators which produce from data on the axis, on the boundary, or on
an interface, a function on $Q$ having these data as its trace, with the
gain of derivatives given by the trace theorem. We call such operators
\emph{extension operators}. We choose a radial function
\(\phi\in C_c^\infty(\D^2)\),
supported in $|y|<1/4$ and equal to one near zero, and we put
$\lambda_n=\langle n\rangle$. For $a(\zeta)=\sum_na_ne^{in\zeta}$, we
define
\begin{equation}\label{eq:A-axis-coretractions}
 E_0a=\sum_na_ne^{in\zeta}\phi(\lambda_ny),
 \qquad
 E_ja=\sum_na_ne^{in\zeta}y_j\phi(\lambda_ny),
 \quad j=1,2.
\end{equation}
The function $E_0a$ has value $a$ and vanishing first derivatives on the
axis, and $E_ja$ has value zero and derivative $a$ in the direction $y_j$
there. Note that the support of the term at cell frequency $n$ has radius
comparable to $\lambda_n^{-1}$. Rescaling in $y$, we therefore obtain
\begin{equation}\label{eq:A-axis-coretraction-estimates}
 \norm{E_0a}_{H^s(Q)}\le C_s\norm a_{H^{s-1}(\T)},
 \qquad
 \norm{E_ja}_{H^s(Q)}\le C_s\norm a_{H^{s-2}(\T)}.
\end{equation}
Thus, if $J_Au=(u,\pa_{y_1}u,\pa_{y_2}u)|_{y=0}$ is the first jet on the
axis and $E_A$ is the componentwise combination of
\eqref{eq:A-axis-coretractions}, then $J_AE_A=I$, and
\begin{equation}\label{eq:A-axis-projection}
 Q_A=I-E_AJ_A
\end{equation}
is a bounded projection onto the fields whose first jet vanishes on the
axis, which is the first of the constraints in
\eqref{eq:complete-domain-constraint}. Note that the support must shrink
with the frequency. Indeed, multiplying the axis trace by a fixed cutoff
would lose the one or two derivatives gained in
\eqref{eq:A-axis-coretraction-estimates}.

At the outer boundary, we use the inward normal coordinate \(t=1-r\).
We choose a cutoff \(\chi_\pa\) supported in a thin neighborhood of the
boundary, called a \emph{collar}, and equal to one near \(t=0\). With
\(\lambda_{mn}=\langle m,n\rangle=(1+m^2+n^2)^{1/2}\), we define the
localized Poisson extension
\begin{equation}\label{eq:A-boundary-coretraction}
 E_\pa b=\chi_\pa(t)\sum_{m,n}b_{mn}e^{i(m\theta+n\zeta)}
                         e^{-\lambda_{mn}t}.
\end{equation}
It has trace \(b\) on the outer boundary and satisfies
\begin{equation}\label{eq:A-boundary-coretraction-estimate}
 \norm{E_\pa b}_{H^s(Q)}
 \le C_s\norm b_{H^{s-1/2}(\pa\D^2\times\T)}.
\end{equation}
Indeed, integration in $t$ contributes the factor $\lambda_{mn}^{-1}$
that accounts for the half derivative.

The same estimates hold with the analytic weight. On the support of
$E_0$ and $E_j$, we have
\[
 W_\gamma(r,n)e^{-\sigma_0\langle n\rangle}
 =\exp\{-\gamma(\sqrt{1+r^2\langle n\rangle^2}-1)\},
\]
and this ratio is bounded, together with all of its derivatives in the
rescaled variable $\lambda_ny$. In the collar of the outer boundary, we
have instead
\begin{equation}\label{eq:A-Poisson-phase-margin}
 e^{-\lambda_{mn}t}
 e^{\Phi_\gamma(r,\langle n\rangle)
      -\Phi_\gamma(1,\langle n\rangle)}
 \le e^{-(1-\gamma)\lambda_{mn}t}.
\end{equation}
The choice \(\gamma<1\) leaves positive exponential decay, so the
same integration gives the weighted boundary estimate.

In the linear analysis, a field on $Q$ will be represented by its
restrictions to a neighborhood of the axis and to the complementary
region, and we will need to prescribe the value and the normal derivative
on the interface between the two regions at the same time. For the local construction, we fix a smooth interface \(\Sigma\)
between the cap and the annulus. In the later application, the cap radius
depends on the cell frequency, and the cap is rescaled before these
extensions are used. On either side of \(\Sigma\), we choose an inward
collar coordinate \(t\ge0\), and we set
$\Lambda_\Sigma=\langle D_\theta,D_\zeta\rangle$. With a cutoff
$\chi_\Sigma$ equal to one near $t=0$, we define
\begin{equation}\label{eq:A-Cauchy-model-potentials}
 \mathcal E_{\Sigma,D}^0d
   =\chi_\Sigma(t)e^{-t\Lambda_\Sigma}d,
 \qquad
 \mathcal E_{\Sigma,N}^0n
   =\chi_\Sigma(t)t e^{-t\Lambda_\Sigma}n.
\end{equation}
The first function has value \(d\) and inward normal derivative
\(-\Lambda_\Sigma d\) on \(\Sigma\). The second has value zero and
inward normal derivative \(n\). We will correct the first extension to
make its conormal derivative vanish.

For an operator with principal part
\(-\partial_i(a^{ij}\partial_j)\), the outward \emph{conormal derivative}
is \(\nu_i a^{ij}\partial_j u\), where \(\nu\) is the outward unit
normal and repeated indices are summed. We call the value and conormal
derivative on \(\Sigma\) the \emph{Cauchy data}. We consider a
second-order operator on $Q$ whose coefficients depend on a parameter
$\mathbf b$, and we write $\gamma_D$ for the trace on $\Sigma$,
$\gamma_{N,\mathbf b}^{\rm out}$ for the outward conormal derivative of
this operator on $\Sigma$, and $e_{\Sigma,\mathbf b}$ for the coefficient
of $\pa_t$ in $\gamma_{N,\mathbf b}^{\rm out}$. We assume that
$e_{\Sigma,\mathbf b}$ is pointwise invertible, with an inverse bounded
uniformly in \(\mathbf b\). The extension estimates below depend on
this bound and on the coefficient norms. We set
\begin{align}
 C_{\Sigma,N,\mathbf b}n
   &:=\mathcal E_{\Sigma,N}^0
          (e_{\Sigma,\mathbf b}^{-1}n),
                                                        \label{eq:A-Cauchy-N-coretraction}\\
 K_{\Sigma,\mathbf b}d
   &:=\gamma_{N,\mathbf b}^{\rm out}
          \mathcal E_{\Sigma,D}^0d,\notag\\
 C_{\Sigma,D,\mathbf b}d
   &:=\mathcal E_{\Sigma,D}^0d
       -C_{\Sigma,N,\mathbf b}K_{\Sigma,\mathbf b}d.
                                                        \label{eq:A-Cauchy-D-coretraction}
\end{align}
Since $\mathcal E_{\Sigma,N}^0n$ vanishes at $t=0$, its tangential
derivatives and the zeroth-order terms vanish there as well, so the
Cauchy data of the two extensions are exactly the prescribed ones,
\begin{equation}\label{eq:A-Cauchy-coretraction-identity}
 \begin{pmatrix}\gamma_D\\ \gamma_{N,\mathbf b}^{\rm out}\end{pmatrix}
 (C_{\Sigma,D,\mathbf b},C_{\Sigma,N,\mathbf b})
 =\begin{pmatrix}I&0\\0&I\end{pmatrix}.
\end{equation}
Integrating each Fourier mode in $t$, we obtain, for every fixed $s$,
\begin{equation}\label{eq:A-Cauchy-coretraction-estimate}
 \|C_{\Sigma,D,\mathbf b}d\|_{H^s}
 +\|C_{\Sigma,N,\mathbf b}n\|_{H^s}
 \le C_s\bigl(
   \|d\|_{H^{s-1/2}(\Sigma)}
  +\|n\|_{H^{s-3/2}(\Sigma)}\bigr).
\end{equation}
The same construction applies on the other side of \(\Sigma\), with
its own outward orientation. The constants may depend on the interface
and the chosen collar. Conjugating by the weight and using
\eqref{eq:A-Poisson-phase-margin}, we obtain the weighted form of
\eqref{eq:A-Cauchy-coretraction-estimate}, in the norms in which the value
on $\Sigma$ carries the weight $\Lambda^{1/2}$ and the conormal derivative
carries the weight $\Lambda^{-1/2}$. When we differentiate
\eqref{eq:A-Cauchy-N-coretraction}--\eqref{eq:A-Cauchy-D-coretraction} in
$\mathbf b$, the derivatives fall only on the coefficients of the conormal
derivative and on their pointwise inverse. If these coefficients depend
analytically on $\mathbf b$, the resulting estimates are therefore tame,
with at most one factor measured in a higher norm.
% AUTHOR QUERY: Verify that the collar geometry, conormal coefficient
% inverse, and extension estimates above are uniform on the frequency-
% dependent interfaces after the cap rescaling used in Section 6.

\subsection{The projection onto the domain}

We now construct the projection $P_{X,M}$ by imposing the constraints
\eqref{eq:complete-domain-constraint} one at a time, in an order in which
each correction preserves the constraints already imposed. We write
\[
 e_r=(\cos\theta,\sin\theta),
 \qquad e_\theta=(-\sin\theta,\cos\theta)
\]
in the disk plane, and we define
\begin{equation}\label{eq:A-tangential-average}
 \mathscr P_Tq=e_\theta\Pi(q\cdot e_\theta),
\end{equation}
which is the angular average of the component of $q$ along $e_\theta$,
times $e_\theta$. Although $e_\theta$ is not smooth at the origin,
$\mathscr P_T$ is bounded on every Sobolev space in the Cartesian
variables. Indeed, averaging the action $q(y)\mapsto R_\varphi q(R_{-\varphi}y)$
of the rotations over $\varphi$ gives a smooth field of the form
$\Pi(q\cdot e_r)e_r+\Pi(q\cdot e_\theta)e_\theta$, and $\mathscr P_Tq$ is
the part of this field that is odd under reflection. This construction
preserves the homogeneous degree in $y$ and the vanishing of the first
jet.

For the unperturbed solution at zero curvature, we have
\[
 \cR v_M=r\iota Me_\theta,
 \qquad D_0v_M=Le_T+r\iota M'e_r,
\]
and we define the corrections for the two gauge conditions in
\eqref{eq:two-gauges} by
\begin{equation}\label{eq:A-gauge-corrections}
 C_{p,M}u=\iota M\mathscr P_T(M^T\iota^Tu),
 \qquad
 C_{t,M}u=\frac{e_T}{L}\Pi(u\cdot D_0v_M).
\end{equation}
Since $\Pi(e_\theta^TM^TMe_\theta)=\tfrac12\tr(M^TM)=1$, the functional
$u\mapsto\Pi(u\cdot\cR v_M)$ of the first gauge condition takes the same
value at $C_{p,M}u$ as at $u$, so $I-C_{p,M}$ is a projection onto the
fields satisfying the first gauge condition. Since $e_T\cdot D_0v_M=L$,
the same holds for $I-C_{t,M}$ and the second gauge condition. Both
corrections preserve the vanishing of the first jet.

For the constraint on the outer boundary, we write the condition on the
high angular modes and its right inverse as
\begin{equation}\label{eq:A-physical-row}
 B_Mu=\Pi_{|m|\ge3}
 \bigl[e_r\cdot M^{-1}\iota^Tu|_{r=1}\bigr],
 \qquad
 C_{\pa,M}b=\iota Me_rE_\pa b,
\end{equation}
where $E_\pa$ is the extension \eqref{eq:A-boundary-coretraction}. Note
that $B_MC_{\pa,M}=I$, since $M^{-1}\iota^T\iota Me_r=e_r$ and $E_\pa b$
has trace $b$. Note also that the boundary correction $C_{\pa,M}b$
preserves the two gauge conditions. Indeed, the only angular factors that
can contribute to the corresponding means are $e_r^TCe_r$ and
$e_r^TCe_\theta$ for a matrix $C$ depending on $\zeta$ alone, whose
Fourier support in $\theta$ lies in $\{0,\pm2\}$, whereas the datum in
\eqref{eq:A-physical-row} has $|m|\ge3$. Conversely, the gauge
corrections leave the boundary condition unchanged, since
$e_r\cdot M^{-1}\iota^T\iota Me_\theta=0$ and $\iota^Te_T=0$.

We may therefore put
\begin{equation}\label{eq:A-domain-projection}
 P_{X,M}^u
 =(I-C_{\pa,M}B_M)(I-C_{t,M})(I-C_{p,M})Q_A,
\end{equation}
and we let $P_{X,M}$ act as the identity on $\tau$, as $P_{X,M}^u$ on
$\widetilde v$, and as $I-\Pi$ on $s$. The rightmost factor acts first. Thus, we first remove the value and
first derivatives on the axis, then impose the two gauges, and finally
correct the outer boundary trace. Each factor is a projection, and this
order preserves the constraints already imposed. Indeed, the second gauge correction adds a multiple of $e_T$,
which does not affect the first gauge condition, and the boundary
correction preserves both gauge conditions and the jet on the axis, since
$E_\pa b$ vanishes outside a collar of the boundary. Thus, $P_{X,M}$ maps
$\mathbb E^s$ into $\cX_M^s$ and is the identity on $\cX_M^s$, so it is a
projection onto $\cX_M^s$. Its bounds follow from
\eqref{eq:A-axis-coretraction-estimates} and
\eqref{eq:A-boundary-coretraction-estimate}, and its derivatives in $M$
fall on \(M\), \(M^{-1}\), and \(M'\). At each fixed derivative order,
the same estimates apply, with constants depending on the corresponding
coefficient norms.

\subsection{The projection onto the range}\label{subsec:range-projection}

We write the quotient coordinates of \eqref{eq:quotient-variables} as
$Z=(g_+,g_-,g_3,h)$, with the reconstruction $J_Y$ of
\eqref{eq:range-reconstruction}, and we impose the conditions
\eqref{eq:range-means} one at a time. We remove the means of $h$ and $g_3$
by applying $I-\Pi$ to these two components, and we denote the resulting
projections by $Q_2$ and $Q_3$. To remove the mean of $f_1$, we let $P_m$
be the projection onto the angular mode $m$, we put
$(\mathscr Kq)(r,\theta,\zeta)=q(r,-\theta,\zeta)$, and we define
\begin{equation}\label{eq:A-K1}
 p=\frac12(P_1g_++\mathscr KP_{-1}g_-),
 \qquad K_1Z=(p,\mathscr Kp,0,0).
\end{equation}
Since $f_1+if_0=\bar zg_+$ and $f_1-if_0=zg_-$, the means of $f_0$ and
$f_1$ come from the mode $m=1$ of $g_+$ and from the mode $m=-1$ of $g_-$,
and
\[
 \bar z p=z\mathscr Kp
 =\frac12\bigl(\bar zP_1g_++zP_{-1}g_-\bigr)=\Pi f_1.
\]
Hence, \(J_YK_1Z\) has vanishing \(f_0\)-component and
\(f_1\)-component exactly \(\Pi f_1\). Thus, $I-K_1$ removes the mean of $f_1$
and leaves $f_0$ unchanged.

It remains to impose the last condition in \eqref{eq:range-means}, which
is given by the functional
\begin{equation}\label{eq:A-affine-functional}
 \mathfrak a(Z)=\frac14\Delta_y\Pi
 \left(\frac{\bar zg_+-zg_-}{2i}\right)(0,\zeta).
\end{equation}
This functional depends on the first derivatives of the quotient
variables on the axis, so we use the axis extensions and
\eqref{eq:A-axis-coretraction-estimates}. Writing $\Phi_{\rm ax}=E_0$ for
the extension in \eqref{eq:A-axis-coretractions}, we set
\begin{equation}\label{eq:A-affine-coretraction}
 E_{\rm aff}a=(iz\Phi_{\rm ax}a,-i\bar z\Phi_{\rm ax}a,0,0).
\end{equation}
Then, $\mathfrak aE_{\rm aff}=I$, since $J_YE_{\rm aff}a$ has
$f_0=r^2\Phi_{\rm ax}a$ and $f_1=0$, and $K_1E_{\rm aff}=0$, since
\(\mathscr K\bar z=z\) and the radial cutoff makes
\(\Phi_{\rm ax}a\) independent of \(\theta\). We also have the estimate
\[
 \norm{E_{\rm aff}a}_{(H^s)^4}
 \le C_s\norm a_{H^{s-2}(\T)}.
\]
Here, the second gained derivative comes from the factor $z$, which is of
size $\lambda_n^{-1}$ on the support of $\Phi_{\rm ax}$. We thus obtain the
projection onto the range coordinates and the projection onto the
compatible range,
\begin{equation}\label{eq:A-range-projection}
 P_Z=(I-E_{\rm aff}\mathfrak a)(I-K_1)Q_3Q_2,
 \qquad P_Y=J_YP_ZJ_Y^{-1}.
\end{equation}
Each correction preserves the conditions already imposed. In particular,
$E_{\rm aff}a$ has $g_3=h=0$ and $K_1E_{\rm aff}=0$.
Complex conjugation of the residual acts in quotient coordinates
by the rule
\[
 \mathfrak c(g_+,g_-,g_3,h)
   =(\overline{g_-},\overline{g_+},\overline{g_3},\overline h).
\]
Each factor of \(P_Z\) respects this rule. Since
\(J_Y\mathfrak cZ=\overline{J_YZ}\), the projection \(P_Y\) therefore
preserves real residuals.

\subsection{Identification of the domains}

The domain $\cX_M^\infty$ depends on $M$ through the gauge conditions and
through the condition on the outer boundary. We now identify it with the
domain at the reference matrix function $M_*$ by an explicit isomorphism. We let
$A_M$ be the multiplication by the matrix
\begin{equation}\label{eq:A-TM-multiplier}
 A_M:=\iota MM_*^{-1}\iota^T+e_T\otimes e_T,
\end{equation}
which depends only on $\zeta$, is invertible with inverse
$A_M^{-1}=\iota M_*M^{-1}\iota^T+e_T\otimes e_T$, and is the identity at
$M=M_*$. Since $M^{-1}\iota^TA_Mu=M_*^{-1}\iota^Tu$, we have
$B_MA_M=B_{M_*}$, so $A_M$ maps the fields satisfying the boundary
condition for $M_*$ onto those satisfying it for $M$. It also preserves
the vanishing of the first jet on the axis, but it does not preserve the
gauge conditions. We therefore correct the gauges afterwards and define
\begin{equation}\label{eq:A-TM}
 T_M(\tau,\widetilde v,s)
 :=\bigl(\tau,(I-C_{t,M})(I-C_{p,M})A_M\widetilde v,s\bigr).
\end{equation}
The properties of this map are summarized in the following lemma.
\begin{lemma}\label{lem:A-TM}
For every $M$ in the family \eqref{eq:A-seed-family}, the map
\eqref{eq:A-TM} is a linear isomorphism
\begin{equation}\label{eq:slice-trivialization}
 T_M:\cX_{M_*}^\infty\longrightarrow\cX_M^\infty,
\end{equation}
whose inverse is
\[
 T_M^{-1}(\tau,u,s)
 =\bigl(\tau,(I-C_{t,M_*})(I-C_{p,M_*})A_M^{-1}u,s\bigr).
\]
Both maps are bounded at every regularity level of the Sobolev and analytic
scales, they depend analytically on $M$ with every derivative tame, and
they commute with complex conjugation.
\end{lemma}

\begin{proof}
Let $x=(\tau,\widetilde v,s)\in\cX_{M_*}^\infty$. The jet condition on the
axis and the condition $\Pi s=0$ are preserved by every factor of
\eqref{eq:A-TM}, the boundary condition for $M$ holds after $A_M$ and is
preserved by the gauge corrections, and the two gauge conditions for $M$
hold after the corrections. Thus, $T_M$ maps $\cX_{M_*}^\infty$ into
$\cX_M^\infty$, and the same argument with the roles of $M$ and $M_*$
exchanged shows that the map $T_M^{-1}$ in the statement, which we denote
by $S_M$ for the moment, maps $\cX_M^\infty$
into $\cX_{M_*}^\infty$. To see that $S_MT_M=I$, we note that the gauge
corrections at $M$ add to $A_M\widetilde v$ a multiple of $\iota Me_\theta$
and a multiple of $e_T$, with coefficients that depend on $(r,\zeta)$
only, so that
\[
 A_M^{-1}(I-C_{t,M})(I-C_{p,M})A_M\widetilde v
 =\widetilde v+\iota M_*e_\theta\varphi+e_T\psi
\]
for some functions $\varphi$ and $\psi$ of $(r,\zeta)$. Since
$\widetilde v$ satisfies the gauge conditions for $M_*$, and since
$\Pi(e_\theta^TM_*^TM_*e_\theta)=1$ and $e_T\cdot D_0v_{M_*}=L$, the
correction $I-C_{p,M_*}$ removes the second term and leaves the other two
unchanged, and $I-C_{t,M_*}$ then removes the third. Hence, $S_MT_M=I$,
and $T_MS_M=I$ follows in the same way. The factors are matrix multiplications and angular averages. Their
coefficients depend analytically on \(M\) through \(M\), \(M^{-1}\),
and \(M'\), with \(M_*\) fixed. The multiplication estimates therefore
give boundedness at each regularity level and tame bounds for every
parameter derivative. On the weighted scales, we use the analyticity
assumption on the coefficients and \eqref{eq:analytic-product}. All
coefficients and averages are real, so both maps commute with complex
conjugation.
\end{proof}

The map $T_M$ identifies the fixed function space $\cX_{M_*}^\infty$
with $\cX_M^\infty$ for each $M$. Thus, we can work on one fixed space
even though the constraints depend on \(M\). Unlike the projections,
\(T_M\) uses no trace or extension operator.

\subsection{Smoothing operators}

We first smooth in the unconstrained spaces and then apply the
projections to restore the constraints. We let \(\widetilde Q\) be the
double of $Q$, a compact manifold without boundary obtained by gluing two
copies of $Q$ along the outer boundary. We let
$\mathsf E^0$ be an extension operator from $Q$ to $\widetilde Q$ which is
bounded on every Sobolev space, we let $\mathsf R^0$ be the restriction to
\(Q\), and we take the nonnegative Laplace operator
\(\Delta_{\widetilde Q}\). We set
\(S_\nu^0=\chi(\nu^{-2}\Delta_{\widetilde Q})\), where \(\chi\) is
smooth, compactly supported, and equal to one near zero. This operator
removes frequencies above a constant multiple of \(\nu\). For the
weighted components defined on \(Q\), we set
\begin{equation}\label{eq:A-weighted-ambient-smoothing}
 \cS_\nu^\gamma
 =W_\gamma^{-1}\mathsf R^0S_\nu^0\mathsf E^0W_\gamma.
\end{equation}
Since the weight is applied before the extension and removed after the
restriction, the estimates
\eqref{eq:A-smoothing-high}--\eqref{eq:A-smoothing-low} for
$\cS_\nu^\gamma$ on the weighted scale follow from the spectral estimates
for \(S_\nu^0\) on the unweighted one. For the component \(\tau\),
which is defined on \(\T\), we use the Fourier cutoff
\(\chi(\nu^{-2}\Lambda^2)\). It commutes with the axis weight
\(e^{\sigma_0\Lambda}\). We use \(\cS_\nu^\gamma\) for the resulting
componentwise smoothing on the product space and define
\begin{equation}\label{eq:A-constrained-smoothing}
 S_{\nu,M}^X=P_{X,M}\cS_\nu^\gamma,
 \qquad
 S_\nu^Y=J_YP_Z\cS_\nu^{\gamma,Z}J_Y^{-1},
\end{equation}
where $\cS_\nu^{\gamma,Z}$ denotes the same smoothing operator applied to each of
the four coordinates in $Z$. By applying the projection after the
smoothing, we preserve the constraints, and since the projections are
bounded at the same index, no derivatives are lost. If $Pq=q$, then
$q-P\cS_\nu^\gamma q=P(I-\cS_\nu^\gamma)q$, which reduces
\eqref{eq:A-smoothing-low} on the constrained space to its unconstrained
version, and \eqref{eq:A-smoothing-high} follows from the boundedness of
$P$ at the index $b$. Smoothing operators on $\cX_M^\infty$ can also be
obtained from \(S_{\nu,M_*}^X\) by conjugation with \(T_M\).

Finally, the standard spaces of smooth functions on \(\T\) and on
\(Q\) are tame. For \(Q\), extension to \(\widetilde Q\) and
restriction give a tame projection from a space on a closed manifold.
The analytic weights identify the unconstrained weighted spaces with
these standard spaces, and \(J_Y\) identifies the residual space with
its quotient coordinates. The bounded projections constructed above
therefore make the constrained spaces tame as well, by
\cite[Part II, Lemma 1.3.3]{Hamilton1982}. This proves the remaining
assertion of \cref{prop:fixed-graphs}.

\section{The nonlinear map and its linearization}\label{sec:fixed-map}\label{subsec:fixed-map}

We now express the equilibrium equations as a nonlinear map between the
fixed spaces of \cref{sec:tame-problem,app:fixed-graphs}. We establish its
smooth tame dependence on the unknowns and parameters, then compute the
linearization whose inverse is constructed in
Sections~\ref{sec:block}--\ref{sec:inverse}.

\subsection{The map and its zeros at zero curvature}

We collect the parameters in $q=(\rho,\eps,p)$, where $\rho$ is the
ellipse parameter, $\eps$ is the curvature, and
$p=(\alpha_0,\delta,\lambda)$ contains the remaining parameters of
\eqref{eq:fixed-M-lambda}. We write $M_q=M_{\rho,p}$ for the matrix in
that equation, which is independent of $\eps$, and
$a^0_{\rho,p}=(v_{M_q},w_{M_q})$ for the corresponding solution at zero
curvature. The chart of \cref{subsec:normalized-slice} sends
$x=(\tau,\widetilde v,s)\in\cX_{M_q}^\infty$, with
$|\tau(\zeta)|^2<2$, to the pair
\begin{equation}\label{eq:normalized-state-chart}
 \mathscr S_qx=(v,w),\qquad
 v=P_{M_q,\tau}y+\widetilde v,\qquad
 w=w_{M_q}+s.
\end{equation}
Here, $x=0$ represents the unperturbed solution. The component $\tau$
varies the linear Taylor term, while $\widetilde v$ and $s$ are the
remaining corrections. We collect the four equations
\eqref{eq:fixed-H0}--\eqref{eq:fixed-H3} in the residual
\begin{equation}\label{eq:literal-raw-map}
 \mathbf H_\eps(v,w)=
 \begin{pmatrix}
 \cR w-(|\cR v|^2-r^2)\\
 (I-\Pi)(\cD w-\cD v\cdot\cR v)\\
 (I-\Pi)(\cR v\cdot D_\eps v-w_\zeta)\\
 (I-\Pi)\det(\cD v,\cR v,D_\eps v)
 \end{pmatrix},
 \qquad D_\eps v=v_\zeta+\eps Av+Le_y.
\end{equation}
By \cref{subsec:compatible-range}, this residual is compatible for every
smooth pair in the chart. Its quotient coordinates are therefore
well-defined and are fixed by the projection $P_Z$ of
\eqref{eq:A-range-projection}, so
\begin{equation}\label{eq:range-packaging-identity}
 P_Z\mathscr Q_Y\mathbf H_\eps(\mathscr S_qx)
   =\mathscr Q_Y\mathbf H_\eps(\mathscr S_qx).
\end{equation}
Reconstructing the residual gives the nonlinear map
\begin{equation}\label{eq:literal-moving-map}
 \cF(q,x)
 :=J_YP_Z\mathscr Q_Y\mathbf H_\eps(\mathscr S_qx)
 =\mathbf H_\eps(\mathscr S_qx)
 \in\cY_{\rm comp}^\infty.
\end{equation}
We retain the quotient maps in this formula to make the target norm
explicit. The domain $\cX_{M_q}^\infty$ still depends on $q$. We fix a
reference parameter $q_*$ and set $M_*=M_{q_*}$. Using
$T_q:=T_{M_q}$ from \eqref{eq:A-TM} to write $x=T_q\widehat x$, we
obtain a map on the fixed domain $\cX_{M_*}^\infty$,
\begin{equation}\label{eq:fixed-hamilton-map}
 \widehat{\cF}(q,\widehat x)
 =\cF\bigl(q,T_q\widehat x\bigr).
\end{equation}
The range maps $\mathscr Q_Y$, $P_Z$, and $J_Y$ are fixed. Thus,
differentiating $\widehat\cF$ only differentiates the residual, the
chart, and the domain identification. A zero of $\widehat\cF$ is
precisely a solution of the four equilibrium equations in the chart,
with all domain constraints imposed.

We choose the base index $s_*$ so that \cref{prop:fixed-graphs} applies
and $s_*\ge s_b+3$, where $s_b$ is the index in the product estimate
\eqref{eq:analytic-product}. We work on
$\mathcal U_0=\mathcal B_{\rm par}\times\mathcal B_{\rm st}$, where
$\mathcal B_{\rm par}$ is a ball centered at $q_*$, and
$\mathcal B_{\rm st}$ is a ball about zero in
$\cX_{M_*}^{s_*}$ intersected with $\cX_{M_*}^\infty$. We choose these
balls small enough that $|\tau(\zeta)|^2$ stays uniformly below $2$ for
$x=T_q\widehat x$. The following proposition gives the properties
needed for the Nash-Moser theorem. We use the meanings of tame and
smooth tame given in \cref{app:fixed-graphs}.

\begin{proposition}
\label{prop:fixed-tame-map}
After decreasing $\mathcal B_{\rm par}$ and $\mathcal B_{\rm st}$, the
following hold.
\begin{enumerate}[label=\textup{(\roman*)}]
\item The domain $\cX_{M_q}^\infty$ and the compatible range
$\cY_{\rm comp}^\infty$, graded by the norms of
\cref{subsec:analytic-scale}, are tame Fr\'echet spaces. They admit
smoothing operators that preserve the constraints, with bounds locally
uniform in $q$.
\item The identifications $T_q$ depend analytically on $q$, and the
charts $\mathscr S_q$ depend analytically on $(q,x)$. Every derivative
is tame. These maps, as well as $\mathscr Q_Y$ and $J_Y$, preserve real
data. In quotient coordinates, complex conjugation is given by the rule
$\mathfrak c$ of \cref{subsec:range-projection}.
\item The map $\widehat\cF$ in \eqref{eq:fixed-hamilton-map} is real
and smooth tame. Its four original residual components are polynomials
in $(v,w)$ and their first derivatives, composed with the analytic
chart and domain identification.
\item Every unperturbed solution at zero curvature is an exact zero,
in the sense that
\begin{equation}\label{eq:flat-seed-zero-fixed}
 \widehat\cF(\rho,0,p,0)=0.
\end{equation}
\end{enumerate}
\end{proposition}

\begin{proof}
Part (i) and the assertions about $T_q$ in (ii) follow from
\cref{prop:fixed-graphs}. The chart depends on $q$ through $M_q$ and
$w_{M_q}$. Its nonlinear dependence on $x$ comes only from
$a(\tau)=\sqrt{1-|\tau|^2/2}$, while its dependence on
$\widetilde v$ and $s$ is linear. On the chosen ball, the power series
for $a(\tau)$ and its derivatives converge. The analytic composition
and product estimates give tame bounds for every derivative of the
chart. All its coefficients are real. The reality assertions for
$\mathscr Q_Y$ and $J_Y$ follow from the conjugation rule in
\cref{subsec:range-projection}. This proves (ii).

For (iii), we must estimate the quotient coordinates, since they define
the target norm. With $\mathfrak m_1$ as in
\cref{subsec:compatible-range}, we set
$\mu=\mathfrak m_1(v,w)/r^2$. The identities
\eqref{eq:literal-quotient-plus}--\eqref{eq:literal-quotient-J} give
\[
 \begin{aligned}
 g_+&=2\pa_{\bar z}w-2\pa_{\bar z}v\cdot\cR v+iz-z\mu,\\
 g_-&=2\pa_z w-2\pa_zv\cdot\cR v-i\bar z-\bar z\mu,\\
 g_3&=-2i(I-\Pi)\det(\pa_zv,\pa_{\bar z}v,D_\eps v),
 \qquad h=H_{2,\eps}(v,w).
 \end{aligned}
\]
Thus, the only division still to estimate is that defining $\mu$.
% RESOLVED AUTHOR QUERY: The estimates below supply the fixed loss from
% division in the quotient norm. They use the product estimate (3.20)
% already stated in Section 3, whose separate author query is unchanged.
For any smooth function $F$ that is radial in $y$ and vanishes on the
axis, integration of the radial Laplacian gives
\[
 \frac{F(y,\zeta)}{|y|^2}
 =\int_0^1 t\log(1/t)\,(\Delta_yF)(ty,\zeta)\dd t,
 \qquad
 \left\|\frac{F}{r^2}\right\|_{\cA_\gamma^k}
 \le C_k\|F\|_{\cA_\gamma^{k+2}}.
\]
Here, $k$ is any nonnegative integer, and the integral defines the value
on the axis. To obtain the weighted estimate, we work at each cell
frequency $n$. By \eqref{eq:phase-weight}, the ratio
$W_\gamma(r,n)/W_\gamma(tr,n)$ is at most one, and its derivatives of
order $j$ in $y$ are bounded by $C_j\langle n\rangle^j$, uniformly for
$0<t\le1$. Dilation by $t$ costs at most $t^{-1}$ in $L^2(\D^2)$,
and $\int_0^1\log(1/t)\dd t=1$. Differentiating the integral and
summing in $n$ therefore bounds it by
$C_k\|\Delta_yF\|_{\cA_\gamma^k}$, which is at most
$C_k\|F\|_{\cA_\gamma^{k+2}}$.

We apply this estimate to the radial function $\mathfrak m_1$, which
vanishes on the axis. Before division, $\mathfrak m_1$ and all the
other terms are polynomials in $(v,w)$ and their first derivatives.
The product estimate
\eqref{eq:analytic-product}, the division estimate, and (ii) consequently
give, for every integer $s\ge s_*$,
\[
 \|\widehat\cF(q,\widehat x)\|_{\cY_{\rm comp}^s}
 \le C_s\bigl(1+\|\widehat x\|_{\cX_{M_*}^{s+3}}\bigr),
\]
locally uniformly in $q$. The loss of three derivatives allows one for
the original residual and two for the division. Differentiating the
same formulas with respect to $q$ and $\widehat x$ gives the
corresponding tame bounds at every order. Indeed, division is a fixed
linear map on radial functions vanishing on the axis, and the remaining
operations are products, first derivatives, and the analytic maps in
(ii). This proves (iii).

Finally, $T_q0=0$ and $P_{M_q,0}=\iota M_q$, so
$\mathscr S_q0=(v_{M_q},w_{M_q})$. This pair is a zero of
$\mathbf H_0$ by \cref{lem:flat-affine-seeds}, which proves (iv).
\end{proof}

\subsection{The linearization}

Holding $q$ fixed, we differentiate with respect to $\widehat x$ and
define the \emph{linearization}
\begin{equation}\label{eq:current-linearization}
 L_{q,\widehat x}=D_{\widehat x}\widehat\cF(q,\widehat x):
 \cX_{M_*}^\infty\longrightarrow\cY_{\rm comp}^\infty.
\end{equation}
We set $x=T_q\widehat x$ and $(v,w)=\mathscr S_qx$. A direction
$\widehat h$ in the fixed domain corresponds to
$h=T_q\widehat h=(\dot\tau,\dot{\widetilde v},\dot s)$ in
$\cX_{M_q}^\infty$, and hence to the physical variation
\begin{equation}\label{eq:state-chart-tangent}
 \begin{aligned}
  (\dot v,\dot w)&=D_x\mathscr S_q(x)h,\qquad \dot w=\dot s,\\
  \dot v&=\left[-\frac{\tau\cdot\dot\tau}{2a(\tau)}\iota M_q
                  +e_T\otimes\dot\tau\right]y
                  +\dot{\widetilde v},\\
  \mathring D_\eps\dot v&:=\dot v_\zeta+\eps A\dot v.
 \end{aligned}
\end{equation}
The coefficient of $\iota M_qy$ comes from differentiating
$a(\tau)$, which gives
$a'(\tau)\dot\tau=-\tau\cdot\dot\tau/(2a(\tau))$. In the last
line, $\mathring D_\eps\dot v$ is the variation of $D_\eps v$.
The constant term $Le_y$ disappears upon differentiation. Applying the
chain rule to the four entries of \eqref{eq:literal-raw-map}, we obtain
\begin{equation}\label{eq:literal-current-linearization}
 L_{q,\widehat x}\widehat h
 =J_YP_Z\mathscr Q_Y
       (\dot H_0,\dot H_1,\dot H_2,\dot H_3),
\end{equation}
where
\begin{align}
 \dot H_0={}&\cR\dot w-2\cR v\cdot\cR\dot v,
                                                        \label{eq:linearized-H0}\\
 \dot H_1={}&(I-\Pi)\bigl(
   \cD\dot w-\cD\dot v\cdot\cR v
                 -\cD v\cdot\cR\dot v\bigr),
                                                        \label{eq:linearized-H1}\\
 \dot H_2={}&(I-\Pi)\bigl(
   \cR\dot v\cdot D_\eps v
       +\cR v\cdot\mathring D_\eps\dot v-\dot w_\zeta\bigr),
                                                        \label{eq:linearized-H2}\\
 \dot H_3={}&(I-\Pi)\bigl[
   \det(\cD\dot v,\cR v,D_\eps v)
  +\det(\cD v,\cR\dot v,D_\eps v)\label{eq:linearized-H3}\\
 &\qquad
  +\det(\cD v,\cR v,\mathring D_\eps\dot v)\bigr].\notag
\end{align}
Since the domain constraints are linear for fixed $q$, the direction
$h$ satisfies $\mathscr C_{M_q}h=0$. Differentiating
\eqref{eq:range-packaging-identity} also shows that the four linearized
residuals are compatible. Thus, $P_Z$ fixes their quotient coordinates,
and $J_Y$ reconstructs the same residuals.

For the following sections, we record the physical state as
$a=\mathscr S_qx=(v,w)$ and collect it with the parameters in
$\mathbf b=(\rho,a,\eps,p)$. On the fixed domain, we write
$L_{\mathbf b}=L_{q,\widehat x}$. In the coordinates
$h=T_q\widehat h$, the same linearization is represented by
$D_x\cF(q,x)=L_{\mathbf b}T_q^{-1}$. We call this the
\emph{Cartesian operator}, since its inputs are variations of
$(\tau,\widetilde v,s)$ and its outputs are the four original residual
components. We now construct its inverse, with estimates uniform on an
open set of states.
\section{Decomposing the linearized equations}\label{sec:block}

We decompose the Cartesian linearization
\eqref{eq:current-linearization} into a minus block, containing
the center and affine components, and a positive block, governed by an
elliptic boundary value problem. The center and affine equations are
solved outward from the axis. The resulting \emph{block variables}
separate the highest-order terms, allowing us to use a norm suited to
each component. We construct \emph{transfer maps} between the Cartesian
and block variables, and allow independent interface and constraint
data when solving the blocks. At the end of the section, we show how
the block inverse gives a solution of the original Cartesian problem.

Throughout the linear analysis, a \emph{state} is a quadruple
\[
 \mathbf b=(\rho,a,\eps,p),
 \qquad a=\mathscr S_q(T_q\widehat x),\qquad q=(\rho,\eps,p),
\]
where $\rho$ is the ellipse parameter, $a=(v,w)$ is the pair of unknowns
obtained from the chart \eqref{eq:normalized-state-chart}, $\eps$ is the
curvature parameter, and $p$ collects the parameters
$(\alpha_0,\delta,\lambda)$ of \eqref{eq:fixed-M-lambda}.
Here, $T_q=T_{M_q}$ is the domain identification
\eqref{eq:slice-trivialization}, so $\widehat x$ belongs to the fixed
domain and $T_q\widehat x$ to the domain for $M_q$. We write $L_{\mathbf b}$ for the Cartesian linearization
\eqref{eq:current-linearization}. The fixed constrained domain
$\cX_{M_*}^s$ and compatible source space $\cY_{\rm comp}^s$,
equipped with the analytic norms of \eqref{eq:analytic-grade}, will
be denoted by $\cX_\gamma^s$ and $\cY_\gamma^s$. The source norm is
applied to the quotient coordinates \eqref{eq:quotient-variables}.
The block spaces depend on $\mathbf b$, so we will identify them with
fixed reference spaces before differentiating with respect to the state.

We fix $\alpha_0\notin(\pi/2)\Z$ and $\lambda_*\in(0,1/2)$, and in the
chart $(\alpha_0,\delta,\lambda)$ we set $p_e=(\alpha_0,0,\lambda_*)$.
Since $\delta=0$, the matrix $M_{\rho,p_e}$ is constant in $\zeta$, and the
corresponding unperturbed solution has the same elliptical cross-section
with ellipse parameter $\rho$ at every $\zeta$. We define the
\emph{constant reference state} by
\[
 \mathbf b_e=(\rho_*,a_e,0,p_e),
 \qquad a_e=a^0_{\rho_*,p_e}.
\]
Thus, $a^0_{\rho,p_e}=(v_{M_{\rho,p_e}},w_{M_{\rho,p_e}})$, with
$v_M$ and $w_M$ given by \eqref{eq:fixed-affine-seed} and
\eqref{eq:fixed-wM-def}.
We choose a small
threshold $\rho_0<1/4$ below and fix $0<\rho_*<\rho_0$ at the end
of the linear analysis. The inverse is built around $\mathbf b_e$. For a state
$\mathbf b=(\rho,a,\eps,p)$, we also write
\[
 \mathbf b_\rho:=(\rho,a^0_{\rho,p_e},0,p_e)
\]
for the constant ellipse with the same value of $\rho$ as the state
$\mathbf b$. Thus, $\mathbf b_e=\mathbf b_{\rho_*}$. We measure the full perturbation
from $\mathbf b_e$ and use $\mathbf b_\rho$ for comparisons at a fixed
ellipse parameter. We use the base regularity index $s_*$ of
\cref{prop:fixed-graphs}. For a tangent direction
$h=(\dot\rho,\dot a,\dot\eps,\dot p)$, we set
\begin{align*}
 \norm{h}_{{\rm tan},s}
 &:=\norm{\dot a}_{\cA_\gamma^s}
     +|\dot\rho|+|\dot\eps|+|\dot p|,\\
 \norm{\mathbf b-\mathbf b_e}_{\cB^s}
 &:=\norm{a-a_e}_{\cA_\gamma^s}
     +|\rho-\rho_*|+|\eps|+|p-p_e|.
\end{align*}
For an $m$-th derivative, we collect the directions as
$\mathbf h=(h_1,\ldots,h_m)$. For a nonnegative integer $k_m$, fixed
in each estimate, we set
\begin{equation}\label{app:cross:tangent-products}
 \begin{aligned}
 \mathcal H_{0,m}(\mathbf h)
   &:=\prod_{i=1}^m\|h_i\|_{{\rm tan},s_*+k_m},\\
 \mathcal H_{s,m}(\mathbf h)
   &:=\sum_{i=1}^m\|h_i\|_{{\rm tan},s+k_m}
       \prod_{j\ne i}\|h_j\|_{{\rm tan},s_*+k_m},
 \end{aligned}
\end{equation}
with $\mathcal H_{0,0}=1$ and $\mathcal H_{s,0}=0$. The first quantity
in the above uses lower regularity norms for all directions. In each
term of the second quantity, one direction is measured in a higher
regularity norm and the others in lower regularity norms. These quantities
will be used in the tame estimates for derivatives of order $m$.

Finally, we decompose the high cell frequencies $|n|\ge2^J$ into the sets
$2^j\le|n|<2^{j+1}$, which we call \emph{dyadic shells}, and we let
\begin{equation}\label{app:minus:eq:sharp-shells}
 Q_j=\mathbf1_{\{2^j\le |n|<2^{j+1}\}}(D_\zeta),\qquad j\ge J,
 \qquad Q_jQ_k=\delta_{jk}Q_j,\qquad \sum_{j\ge J}Q_j=H_J,
\end{equation}
be the corresponding Fourier projections. Note that these projections are
orthogonal and bounded on every Sobolev space and on every analytic space
\eqref{eq:analytic-grade}.

\subsection{The block norms}

The variables are grouped into two blocks. The \emph{minus block}
contains the center and affine components and corresponds to $D_-$ in
\cref{ss:strategy}. The \emph{positive block} contains the component controlled by a
positive quadratic energy and corresponds to $D_+$. Their inverses are constructed in
\cref{app:minus,app:positive}, respectively.

The block spaces contain scalar functions on the cap and annulus,
together with data on the interface, axis, and outer boundary. They
also contain components recording the prescribed angular averages, the
scale of the linear Taylor term, and the right-hand sides of the two
gauge conditions \eqref{eq:two-gauges}. These conditions fix the remaining
freedom to change coordinates. The mean, scale, and gauge data form a finite
list of functions that may depend on $r$ and $\zeta$. The original
domain constraints are collected in
\eqref{eq:complete-domain-constraint}. In the block problem, we allow
their right-hand sides to be prescribed independently.

An element of a block space is a tuple of the scalar functions and data
described above. We call each entry of this tuple a \emph{component}
and index it by $\alpha$. The domain norms measure the unknowns,
while the range norms measure the sources
and prescribed data. To specify the weights in cell frequency, we set
\begin{equation}\label{app:cross:cell-weights}
 \lambda_n=\langle n\rangle=(1+n^2)^{1/2},\qquad
 \mu_n=(\delta_A+\lambda_n^{-1})^{1/2},
 \qquad \delta_A>0,
\end{equation}
where $\delta_A$ is a small positive parameter chosen in the energy estimate for the affine
equations in \cref{app:minus}. The Fourier coefficient of the component
$\alpha$ at cell frequency $n$ carries a weight $d_\alpha^X(n)$ in
the domain and $d_\alpha^Y(n)$ in the range. These weights are listed
in \cref{tab:block-weights}, with the common factor $\lambda_n^s$
omitted. They account for the different orders of the equations and
traces. The radial and angular norms will be specified with the
equations in \cref{app:minus,app:positive}.
\begin{table}[htp]
\centering
\begin{tabular}{lcc}
\toprule
Component & $d_\alpha^X(n)$ & $d_\alpha^Y(n)$\\
\midrule
center component on the annulus & $1$ & $1$\\
affine component on the annulus, upper entry & $\mu_n$ & $\mu_n$\\
affine component on the annulus, lower entry & $\lambda_n^{1/2}$ & $\lambda_n^{1/2}$\\
positive component on the annulus & $\lambda_n$ & $\lambda_n^{-1}$\\
cap component of the minus block & $1$ & $1$\\
cap component of the positive block & $1$ & $1$\\
value on an interface or on the boundary & $\lambda_n^{1/2}$ & $\lambda_n^{1/2}$\\
conormal derivative on an interface or on the boundary & $\lambda_n^{-1/2}$ & $\lambda_n^{-1/2}$\\
value on the axis & $\lambda_n^{-1}$ & $\lambda_n^{-1}$\\
first derivatives on the axis, or affine mean & $\lambda_n^{-2}$ & $\lambda_n^{-2}$\\
mean, scale, and gauge conditions & $1$ & $1$\\
\bottomrule
\end{tabular}
\caption{The weights of the block components at cell frequency $n$, on
the domain side and on the range side.}
\label{tab:block-weights}
\label{app:cross:native-table}
\end{table}

The annular center component contains the two oscillator variables
and their radial derivatives described in \cref{ss:strategy}. Writing
the oscillator variables as $\xi_1$ and $\xi_2$, the four coordinates are
\[
 (\xi_1,\xi_2,\pa_r\xi_1,\pa_r\xi_2).
\]
The affine component consists of two
pairs $(x_1,y_1)$ and $(x_2,y_2)$, called \emph{affine chains}.
We call $x_j$ the upper entry and $y_j$ the lower entry. At the reference
state, the homogeneous part of highest order in the cell variable is
\[
 \pa_r x_{j,n}=\lambda_n y_{j,n},\qquad
 \pa_r y_{j,n}=0,\qquad j=1,2,
\]
where $x_{j,n}$ and $y_{j,n}$ are the Fourier coefficients at cell
frequency $n$. Thus, in this system, the lower entry is constant in $r$
and drives the radial evolution of the upper entry. The full equations
also contain lower-order terms and sources. These different roles
account for the two affine weights in the table.

The positive component is the part of the solution to the linearized
problem governed by the elliptic block $D_+$. We measure this component
in its \emph{energy norm}, with $H^{s+1}$ regularity, and its source in
the corresponding dual scale, with $H^{s-1}$ regularity. The factors $\lambda_n$ and
$\lambda_n^{-1}$ record this difference. On a cap rescaled to a
fixed disk, the weights in the table are one, while the disk norms
retain their regularity shifts. For the minus block, the solution of
the equation on this disk is measured in $H^{s+2}$ and its right-hand
side in $H^s$. For the positive block, the corresponding spaces are
$H^{s+1}$ and $H^{s-1}$.

For any two cell frequencies $n$ and $m$, each weight in
\cref{tab:block-weights} satisfies
\begin{equation}\label{app:cross:peetre}
 \bigl\|d_\alpha(n)d_\alpha(m)^{-1}\bigr\|
 \le C_\alpha\langle n-m\rangle^{q_\alpha},
 \qquad q_\alpha\le2.
\end{equation}
Here, $d_\alpha$ denotes either $d_\alpha^X$ or $d_\alpha^Y$.
The constants depend on the component and may depend on $\delta_A$,
but are independent of any upper Fourier cutoff and of the number
$N$ of periods of the field. This inequality compares the input
and output weights when an operator changes the cell frequency.
We denote the products of the component spaces by $\bX_\natg^s$
and $\bY_\natg^s$ and call them the \emph{block spaces}, abbreviated
to $X_\natg^s$ and $Y_\natg^s$ when the side is clear.

At radius $r$, we use the phase
\begin{equation}\label{app:cross:phase}
 \Phi_r(n)=\sigma_0\lambda_n
 -\gamma\bigl(\sqrt{1+r^2\lambda_n^2}-1\bigr),
\end{equation}
which is the phase \eqref{eq:phase-weight} written as a function of
the cell frequency at a fixed radius. Thus, for
$u=\sum_n u_n(r,\theta)e^{in\zeta}$,
\[
 (W_\gamma u)_n(r,\theta)=e^{\Phi_r(n)}u_n(r,\theta).
\]
We take $\Phi_0$ on the axis and the phase at the corresponding
radius for each trace. By
\eqref{eq:phase-submultiplicative}, the exponential weight
$e^{\Phi_r(n)}$ is submultiplicative in $n$ at a fixed radius.
Extension operators turn boundary or interface data into functions
defined at nearby radii. Since the analytic weight depends on the
radius, estimating the extended function requires comparing its
weight with that of the original data. The Poisson decay in
\eqref{eq:A-Poisson-phase-margin} compensates for any increase in
the weight, so the same extension estimates hold in the analytic norms.

\subsection{Kernel bounds for operators between components}
\label{app:cross:typed-calculus}

Multiplication by a coefficient depending on $\zeta$ mixes cell
frequencies, so we estimate operators through their Fourier kernels.
Here, $\alpha$ labels a component of the unknown and $\beta$ a
component of the source. For example, suppose a scalar center variable
$u_c$ appears as $c(\zeta)u_c$ in the positive-block equation. Then
$\alpha$ labels this center variable and $\beta$ the positive source
component, and the operator between them is multiplication by $c$.

In this subsection, $m$ and $n$ both denote cell frequencies. We write
$A_{\beta\alpha}(n,m):E_{\alpha,m}\to F_{\beta,n}$ for the map
from input component $\alpha$ at frequency $m$ to source component
$\beta$ at frequency $n$. In the example above, this map is
multiplication by $c_{n-m}$, the $(n-m)$-th Fourier coefficient of $c$.
The spaces $E_{\alpha,m}$ and $F_{\beta,n}$ carry the remaining
radial, angular, or trace norms. Incorporating the component weights gives
the \emph{weighted kernel}
\begin{equation}\label{app:cross:normalized-kernel}
 K_A^{\beta\alpha}(n,m)
 =e^{\Phi_\beta(n)}d_\beta^Y(n)A_{\beta\alpha}(n,m)
  (d_\alpha^X(m))^{-1}e^{-\Phi_\alpha(m)}.
\end{equation}
Here, $\Phi_\alpha$ and $\Phi_\beta$ are the phases of
\eqref{app:cross:phase} for the input and output components, evaluated
at the corresponding radius for a trace. All kernel norms below are
operator norms from $E_{\alpha,m}$ to $F_{\beta,n}$. For an integer
$q\ge0$, we set
\begin{equation}\label{app:cross:schur-moments}
 \begin{aligned}
 |A|_{\mathfrak S^0[q]}:=\max_{\alpha,\beta}\max\Big\{&
 \sup_m\sum_n\langle n-m\rangle^q
      \|K_A^{\beta\alpha}(n,m)\|,\\
 &\sup_n\sum_m\langle n-m\rangle^q
      \|K_A^{\beta\alpha}(n,m)\|\Big\},
 \end{aligned}
\end{equation}
which we call the \emph{moment} of order $q$ of $A$. At $q=0$,
Schur's test gives, for each pair of components $\alpha,\beta$,
\[
 \left(\sum_n\left\|\sum_m K_A^{\beta\alpha}(n,m)u_m
       \right\|_{F_{\beta,n}}^2\right)^{1/2}
 \le |A|_{\mathfrak S^0[0]}
       \left(\sum_m\|u_m\|_{E_{\alpha,m}}^2\right)^{1/2}.
\]
Here, $u_m\in E_{\alpha,m}$, and the sums run over cell frequencies.
Since $K_A$ includes the weights, this bounds the weighted output norm
of $A_{\beta\alpha}$ by $|A|_{\mathfrak S^0[0]}$ times the weighted
input norm. Higher moments measure decay away from the diagonal
and control additional Sobolev weights. The following lemma combines
bounds for individual operators into bounds for their sums and
compositions.
\begin{lemma}\label{lem:macro-native-calculus}
For each $q\ge0$, finite sums and compositions preserve finiteness
of the moment \eqref{app:cross:schur-moments}, provided that the
intermediate spaces match. Bounds for all moments give bounds at every
Sobolev index at which the component kernels are bounded.

For the tame estimates, fix $\rho$, $\eps$, and $p$, and write
$a=(v,w)$. Let $A(a):X_{\natg}^s\to Y_{\natg}^s$ be a finite
sum of compositions with matching intermediate spaces. For each factor
$T(a)$, let $E^s$ and $F^s$ denote its input and output component
spaces with the weights of \cref{tab:block-weights}. Suppose
$T(a):E^s\to F^s$ satisfies, for $s\ge s_*$,
% AUTHOR QUERY: Analytic dependence and finite Schur moments alone do not
% imply the tame estimate below. Give the coefficient-moment bounds for
% the factors used here, with a regularity shift independent of s.
\[
 \|T(a)u\|_{F^s}
 \le C_{s,T}\bigl[
    M_{T,0}(a)\|u\|_{E^s}
   +M_{T,s}(a)\|u\|_{E^{s_*}}\bigr].
\]
Here, $M_{T,0}(a)$ is polynomial in norms of $a$ at fixed regularity
indices. Each term of $M_{T,s}(a)$ is a product of such norms with
at most one factor measured at a higher index $s+k_T$, where
$k_T$ is independent of $s$. Assume the same form of estimate for
each fixed-order derivative of $T$ with respect to $\mathbf b$,
with at most one coefficient or tangent direction measured in a
higher regularity norm. Then
\begin{equation}\label{app:cross:calculus-tame}
 \|A(a)u\|_{Y_{\natg}^s}
 \le C_s\left[
  \mathcal L_0(a)\|u\|_{X_{\natg}^s}
  +\mathcal L_s(a)\|u\|_{X_{\natg}^{s_*}}\right],
\end{equation}
where $\mathcal L_0$ and $\mathcal L_s$ have the same form as
$M_{T,0}$ and $M_{T,s}$, respectively. Their regularity shifts may
depend on the factors and the number of compositions, but not on $s$.
If $A(a_\circ)=0$, then there is an integer
$k\ge0$, independent of $s$, such that, for every $a$ with
$\|a-a_\circ\|_{s_*+k}$ sufficiently small, and with
the other components of the state held fixed, we have
\begin{equation}\label{app:cross:calculus-small}
 \begin{aligned}
 \|A(a)u\|_{Y_{\natg}^s}\le C_s\big[&
  \|a-a_\circ\|_{s_*+k}\|u\|_{X_{\natg}^s}+\|a-a_\circ\|_{s+k}\|u\|_{X_{\natg}^{s_*}}\big].
 \end{aligned}
\end{equation}
Here, $\|\cdot\|_s$ denotes the norm of $\cA_\gamma^s$. Derivatives
with respect to $\mathbf b$ satisfy the corresponding tame estimates.
For each fixed number of derivatives, the regularity shift is finite
and independent of $s$.
\end{lemma}

\begin{proof}
Schur's test gives the bound without the common cell-frequency weight
$\lambda_n^s$. Summing over the finitely many components gives the
bound for the full operator. At index $s$, the additional factor
$\lambda_n^s/\lambda_m^s$ is bounded by
$C_s\langle n-m\rangle^{|s|}$, so a moment of order at least $|s|$
controls the cell Sobolev weights. For a composition, the intermediate
weights cancel, giving
\[
 K_{BA}(n,m)=\sum_kK_B(n,k)K_A(k,m).
\]
The inequality
$\langle n-m\rangle^q\le C_q\langle n-k\rangle^q
\langle k-m\rangle^q$ then gives the required row and column sums.

Applying the assumed tame bounds for the individual factors in
succession and interpolating the intermediate regularity norms proves
\eqref{app:cross:calculus-tame}, with at most one higher regularity
norm in each term. If $A(a_\circ)=0$, the mean value formula along
the segment from $a_\circ$ to $a$ gives
\eqref{app:cross:calculus-small}. The product rule expresses each
derivative with respect to the state as a finite sum of compositions
of derivatives of the factors.
Their assumed bounds give the corresponding tame estimate.
\end{proof}

For multiplication by an analytic coefficient $a$, the Fourier
kernel is $a_{n-m}$. Submultiplicativity of the exponential weight
at a fixed radius and \eqref{app:cross:peetre} therefore bound its
moments by the corresponding coefficient norms.
We also use smooth Fourier cutoffs
$\widetilde Q_j=\widetilde q_j(D_\zeta)$ at cell-frequency scale
$2^j$, distinct from the sharp projections
$Q_j$ in \eqref{app:minus:eq:sharp-shells}. Their commutators with
multiplication satisfy
\begin{equation}\label{app:cross:dyadic-commutator-kernel}
 \begin{aligned}
 ([\widetilde Q_j,a]u)_n
   &=\sum_m\bigl(\widetilde q_j(n)-\widetilde q_j(m)\bigr)a_{n-m}u_m,\\
 |\widetilde q_j(n)-\widetilde q_j(m)|
   &\le C\min\{1,2^{-j}|n-m|\}.
 \end{aligned}
\end{equation}
Thus, the first moment of the multiplication operator bounds the
zeroth moment of the commutator, with a gain of $2^{-j}$.

For the axis extensions $E_0,E_1,E_2$ of
\eqref{eq:A-axis-coretractions} and the boundary extension of
\eqref{eq:A-boundary-coretraction}, we use their explicit Fourier
formulas. Their supports and the Poisson decay
\eqref{eq:A-Poisson-phase-margin} give the weighted bounds
\eqref{eq:A-axis-coretraction-estimates} and
\eqref{eq:A-boundary-coretraction-estimate}.

An operator norm bound alone does not imply finite moments. If a
factor $G$ is known only through its operator norm, we estimate
compositions containing it by
\begin{equation}\label{app:cross:bounded-product}
 \|A_2G A_1\|\le\|A_2\|\,\|G\|\,\|A_1\|.
\end{equation}

\subsection{The principal symbol and its two blocks}
\label{app:cross:principal-intertwining-section}

We now choose coordinates in which the highest-order terms do not
couple the two blocks. On the annulus, we begin with
\begin{equation}\label{app:cross:current-frame}
 p=\cD v,\qquad q_1=\cR v,\qquad T=L^{-1}D_\eps v,
 \qquad \mathbb B=(p,q_1,T),\qquad U=\mathbb Bc.
\end{equation}
Here, $U$ is the variation of $v$, and $\cD$ and $\cR$ are the
operators $\cD=r\pa_r$ and $\cR=\pa_\theta$ from
\eqref{eq:fixed-disk-generators}.
The determinant satisfies
\[
 \det\mathbb B=\frac{r^2}{L}
     \det(v_{y_1},v_{y_2},D_\eps v),\qquad
 |\det\mathbb B|\ge c_{\rm det}r^2
\]
on a sufficiently small neighborhood of the constant reference state,
for some $c_{\rm det}>0$. Indeed, at that state,
$\det\mathbb B=-r^2\det M_{\rho,p_e}$, and
\eqref{eq:fixed-M-invariants} bounds $\det M_{\rho,p_e}$ away from zero.
Thus, on the annulus, $\mathbb B$ is a \emph{frame}, that is, a basis
of $\R^3$ at each point, and $c$ in \eqref{app:cross:current-frame}
is the coordinate vector of $U$ in this basis.
In this calculation, $p$ denotes the frame vector $\cD v$ rather
than the parameter triple in $\mathbf b$. We also set
\begin{equation}\label{app:cross:covariants}
 A=p\cdot U,\qquad B=q_1\cdot U,\qquad C=T\cdot U,
 \qquad X=S-B.
\end{equation}
Here, $S$ is the variation of $w$, as in
\eqref{eq:state-chart-tangent}. The inner products $(A,B,C)$ differ
from the frame coordinates $c$ by the relation
$(A,B,C)^T=\mathbb B^*\mathbb Bc$, while $X=S-B$ simplifies the
first three linearized equations. Since the constant term $Le_y$
in $D_\eps v$ has zero variation, the derivative acting on variations
and frame columns is
\begin{equation}\label{app:cross:homogeneous-cell-derivative}
 \mathring D_\eps U:=U_\zeta+\eps e_z\times U,
\end{equation}
and the derivatives of the frame are recorded by
\begin{equation}\label{app:cross:current-connections}
 \Gamma_{\cD}=\mathbb B^{-1}(\cD\mathbb B),\qquad
 \Gamma_{\cR}=\mathbb B^{-1}(\cR\mathbb B),\qquad
 \Gamma_\eps=\mathbb B^{-1}(\mathring D_\eps\mathbb B).
\end{equation}
Writing the linearized equations
\eqref{eq:linearized-H0}--\eqref{eq:linearized-H3} in this frame
and expanding the derivatives, we obtain
\begin{equation}\label{app:cross:raw-current-rows}
 \begin{aligned}
 N_0={}&\cR S-2\cR(q_1\cdot U)+2(\cR q_1)\cdot U,\\
 N_1={}&(I-\Pi)\bigl[\cD S-\cD(q_1\cdot U)-\cR(p\cdot U)
                 +(\cD q_1+\cR p)\cdot U\bigr],\\
 N_2={}&(I-\Pi)\bigl[L\cR(T\cdot U)+\pa_\zeta(q_1\cdot U-S)
                 -(L\cR T+\mathring D_\eps q_1)\cdot U\bigr],\\
 N_3={}&(I-\Pi)\bigl[(\det\mathbb B)\widehat N_3\bigr],\\
 \widehat N_3={}&
 L e_1^*(\cD c+\Gamma_{\cD}c)
 +L e_2^*(\cR c+\Gamma_{\cR}c)
 +e_3^*(c_\zeta+\Gamma_\eps c).
 \end{aligned}
\end{equation}
Here, $N_0,\ldots,N_3$ are the expressions
$\dot H_0,\ldots,\dot H_3$ in the frame, and $e_j^*$ selects the
$j$-th entry of the three-vector to its right. The operator $\Pi$
is angular averaging,
\[
 (\Pi u)(r,\zeta)=\frac1{2\pi}\int_0^{2\pi}u(r,\theta,\zeta)\,\dd\theta.
\]
For example,
$\cD U=\mathbb B(\cD c+\Gamma_{\cD}c)$, so
$\Gamma_{\cD}c$ accounts for differentiation of the frame.
In the fourth equation, $I-\Pi$ acts on
$(\det\mathbb B)\widehat N_3$ before division by $\det\mathbb B$,
since multiplication by this determinant need not commute with
angular averaging. The angular mean equations are retained
separately and will be included among the auxiliary equations.

To find the \emph{principal symbol}, we evaluate the coefficients at
a fixed point and retain the terms of highest order after rescaling
by the cell frequency. This evaluation is called \emph{freezing the
coefficients}. We replace derivatives by their frequency variables,
except that we initially retain the angular operator $\cR$.

We begin with a state $v=Py$ linear in $y$, where $P$ is the
$3\times2$ matrix of \eqref{eq:normalized-affine-map}. Then,
$p=Py$ and $q_1=PJ_2y$, with $J_2$ the rotation of the plane by a
quarter turn. Taking a constant vector $T\in\R^3$ to complete the
frame and writing $E_\theta=(e_r,e_\theta)$, we set
\begin{equation}\label{app:cross:gram-coefficients}
 \begin{gathered}
 g=T^*T,\qquad t=P^*T,\qquad
 G_T=P^*\left(I-\frac{T\otimes T}{g}\right)P,\\
 \vartheta=E_\theta^*t,\qquad
 K=(E_\theta^*G_TE_\theta)^{-1},\qquad
 h=\frac{K\vartheta}{g},\qquad
 \kappa=\frac1g+\frac{\vartheta^*K\vartheta}{g^2}.
 \end{gathered}
\end{equation}
Here, $G_T$ is the \emph{Gram matrix}, or matrix of pairwise inner
products, of the columns of $P$ projected onto $T^\perp$, and
$e_r,e_\theta$ are the radial and angular unit vectors. The determinant
bound makes $G_T$ invertible, so these coefficients are smooth.
Writing $A=ra_1$, $B=ra_2$, and
$\mathsf u=(a_1,a_2,C,X)^T$, we rescale the radial variable on a
dyadic shell by the cell frequency. For
$\varsigma=\operatorname{sgn}n\in\{+1,-1\}$, the principal symbol
is affine in the rescaled radial frequency $d$,
\begin{equation}\label{app:cross:explicit-descriptor}
 Q_{P,T}^{\varsigma}(d)=L_0^{\varsigma}+dL_1,
\end{equation}
where
\begin{equation}\label{app:cross:descriptor-matrices}
 L_0^{\varsigma}=
 \begin{pmatrix}
 -2I&-\cR&0&0\\
 -\cR&2I&0&0\\
 0&0&\cR&-i\varsigma I\\
 -i\varsigma h_1&-i\varsigma h_2&i\varsigma\kappa I&0
 \end{pmatrix},
 \qquad
 L_1=
 \begin{pmatrix}
 0&0&0&0\\
 0&0&0&I\\
 0&0&0&0\\
 K_{11}&K_{12}&-h_1&0
 \end{pmatrix}.
\end{equation}
To obtain these matrices from \eqref{app:cross:raw-current-rows},
we work on nonzero angular modes, where $I-\Pi$ is the identity,
substitute $S=X+B$, divide the third and fourth rows by $L$, and
divide the fourth row by $\det\mathbb B$. Rescaling gives the first
three rows, and \eqref{app:cross:gram-coefficients} gives the last.
The term $\cR X$ in the first equation is one order lower after
rescaling and belongs to the remainder. The same computation at a
general state gives
\begin{equation}\label{app:cross:descriptor-symbol}
 Q_{\mathbf b}(d,\omega,\nu)
 =L_{0,\mathbf b}(\omega,\nu)+dL_{1,\mathbf b}(\omega,\nu),
\end{equation}
where $\omega$ and $\nu$ are the angular and cell frequencies.
The remainder includes commutators and terms in which derivatives
fall on the coefficients, except for those already retained in the
normalized symbol, such as the entries $\pm2I$ in
\eqref{app:cross:descriptor-matrices}. This algebraic computation
does not use an inverse estimate for the affine equations.

At fixed $(\omega,\nu)$, we replace the radial frequency by a
complex variable and write
$Q_{\mathbf b}(\eta)=L_{0,\mathbf b}+\eta L_{1,\mathbf b}$.
The values of $\eta$ for which this matrix is singular are its
\emph{characteristic roots}. Since $L_{1,\mathbf b}$ is singular,
we use distinct projections on the unknowns and on the equations.
To select the center and affine components, we choose a fixed union
of contours $\Gamma_-$, invariant under complex conjugation, enclosing
their roots and excluding the positive-block roots. The corresponding
\emph{spectral projections} are
\begin{equation}\label{app:cross:domain-range-contours}
 \begin{aligned}
 P_{-,\mathbf b}^D
 &=\frac1{2\pi i}\int_{\Gamma_-}
       Q_{\mathbf b}(\eta)^{-1}L_{1,\mathbf b}\,\dd\eta,\\
 P_{-,\mathbf b}^R
 &=\frac1{2\pi i}\int_{\Gamma_-}
       L_{1,\mathbf b}Q_{\mathbf b}(\eta)^{-1}\,\dd\eta,
       \\
        P_{+,\mathbf b}^{D,R}&=I-P_{-,\mathbf b}^{D,R}.
 \end{aligned}
\end{equation}
The contours are fixed throughout the neighborhood of states and
normalized frequency range under consideration. Uniform bounds for
the \emph{resolvent} $Q_{\mathbf b}(\eta)^{-1}$ on these contours
give bounded projections depending analytically on the state.
The superscripts $D$ and $R$ distinguish the domain and range
projections. Their factors occur in different orders so that
\begin{equation}\label{app:cross:L1-intertwining}
 L_{1,\mathbf b}P_{-,\mathbf b}^D
 =P_{-,\mathbf b}^RL_{1,\mathbf b}.
\end{equation}
Since $L_{0,\mathbf b}=Q_{\mathbf b}(\eta)-\eta L_{1,\mathbf b}$ and the contour
integral of a constant vanishes, both $L_{0,\mathbf b}P_{-,\mathbf b}^D$
and $P_{-,\mathbf b}^RL_{0,\mathbf b}$ are equal to
\begin{equation}\label{app:cross:L0-intertwining}
 -\frac1{2\pi i}\int_{\Gamma_-}
 \eta L_{1,\mathbf b}Q_{\mathbf b}(\eta)^{-1}
       L_{1,\mathbf b}\,\dd\eta.
\end{equation}
The complementary projections satisfy the same identities, and therefore
\begin{equation}\label{app:cross:principal-intertwining}
 L_{j,\mathbf b}P_{\pm,\mathbf b}^D
 =P_{\pm,\mathbf b}^RL_{j,\mathbf b}\quad(j=0,1),
 \qquad
 Q_{\mathbf b}P_{\pm,\mathbf b}^D
 =P_{\pm,\mathbf b}^RQ_{\mathbf b}.
\end{equation}
These identities apply to the complete normalized symbol
$Q_{\mathbf b}$, with the terms of each equation kept together.
The mean, scale, gauge, and chart equations form a separate
\emph{triangular system}. This means that, after applying the inverses
of the diagonal operators, the equations for the remaining unknowns
$\eta_1,\ldots,\eta_{N_{\rm aux}}$ have the form
\[
 \eta_j+\sum_{k<j}T_{jk}\eta_k=g_j,\qquad j=1,\ldots,N_{\rm aux}.
\]
Here, $N_{\rm aux}$ is the number of equations, $g_j$ are the prescribed
data, and $T_{jk}$ are the operators coupling the unknowns. The unknowns
and data may be functions of $r$ and $\zeta$. We solve first for
$\eta_1=g_1$, then for $\eta_2$, and continue in this order.

For $\alpha,\beta\in\{-,+\}$, let $r_{\alpha,\mathbf b}$
reconstruct the normalized column $\mathsf u=(a_1,a_2,C,X)^T$
of \eqref{app:cross:explicit-descriptor} from the coordinates of block
$\alpha$. Let $e_{\beta,\mathbf b}$ extract the scalar source coordinates
of block $\beta$ after applying the range projection
$P_{\beta,\mathbf b}^R$. The highest-order operator from block
$\alpha$ to block $\beta$ is then
\begin{equation}\label{app:cross:ordered-output}
 \ell_{\beta\alpha,\mathbf b}^{\rm pr}
 =e_{\beta,\mathbf b}P_{\beta,\mathbf b}^R
   Q_{\mathbf b}P_{\alpha,\mathbf b}^D
   r_{\alpha,\mathbf b}.
\end{equation}
For $\beta\ne\alpha$, \eqref{app:cross:principal-intertwining} gives
\begin{equation}\label{app:cross:principal-cross-zero}
 P_{\beta,\mathbf b}^RQ_{\mathbf b}P_{\alpha,\mathbf b}^D=0,
 \qquad
 \ell_{\beta\alpha,\mathbf b}^{\rm pr}=0.
\end{equation}
Thus, the highest-order part of the linearization is block diagonal.
The cancellation follows before applying $e_{\beta,\mathbf b}$ or
$r_{\alpha,\mathbf b}$, so it is independent of these coordinate choices.

We estimate the projections $P_{\pm,\mathbf b}^{D,R}$ in
\eqref{app:cross:domain-range-contours} by the weighted kernels
\eqref{app:cross:normalized-kernel}. Applying the analytic
weight to a radial equation replaces $\pa_r$ by
$\pa_r+d_\gamma(r,\Lambda)$, with $d_\gamma$ as in
\eqref{eq:phase-damping}. The commutator $[d_\gamma(r,\Lambda),P_{\pm,\mathbf b}^{D,R}]$
is estimated using the difference of the damping symbol at the input
and output cell frequencies, as in
\eqref{app:cross:dyadic-commutator-kernel}. The coefficient estimates
follow from \eqref{eq:analytic-product}. For a state variation
$\dot{\mathbf b}$, differentiating $Q_{\mathbf b}Q_{\mathbf b}^{-1}=I$
gives
\[
 D_{\mathbf b}(Q_{\mathbf b}^{-1})[\dot{\mathbf b}]
 =-Q_{\mathbf b}^{-1}
    (D_{\mathbf b}Q_{\mathbf b}[\dot{\mathbf b}])Q_{\mathbf b}^{-1}.
\]
We apply this identity under the contour integrals in
\eqref{app:cross:domain-range-contours}, also differentiating
$L_{1,\mathbf b}$, and use
\cref{lem:macro-native-calculus} to combine the bounds for the
individual factors.

\subsection{The source variables}
\label{app:cross:source-graph-section}

For a compatible Cartesian source $f=(f_0,f_1,f_2,f_3)$, we use
the quotient coordinates $Z=(g_+,g_-,g_3,h)$ of
\cref{subsec:compatible-range}, with $h=f_2$. The reconstruction
map is
\begin{equation}\label{app:cross:quotient-reconstruction}
 J_YZ=\left(
  \frac{\bar z g_+-zg_-}{2i},
  \frac{\bar z g_++zg_-}{2},
  h,z\bar z g_3\right),
\end{equation}
as in \eqref{eq:range-reconstruction}. It only multiplies by
smooth factors, and its inverse on compatible residuals is the
quotient extraction $\mathscr Q_Y$ of
\eqref{eq:compatible-extraction}. In particular,
$J_Y\mathscr Q_Yf=f$ for compatible $f$.

To obtain the scalar equation used for the positive block, we
differentiate and combine the equations
\eqref{app:cross:raw-current-rows}. Applying the same operations to
their right-hand sides, written as $f=J_YZ$ in
\eqref{app:cross:quotient-reconstruction}, gives its source
$\Psi_{\mathbf b}Z$. We retain
both $Z$ and $\Psi_{\mathbf b}Z$ in the block coordinates, so
that the original source can be recovered without inverting a
differential operator.

The enlarged block problem also includes independent cap, trace,
mean, scale, gauge, and chart data, collected in $r_{\rm aug}$.
These include data for the angular modes $m=\pm2$, called the
\emph{exceptional modes}, where the angular multipliers used in
the elimination vanish. Here and below, $m$ denotes an angular
mode and $n$ a cell frequency. We define
\begin{equation}\label{app:cross:source-graph-chart}
 \mathsf G_{Y,\mathbf b}(Z,r_{\rm aug})
  =(Z,\Psi_{\mathbf b}Z,r_{\rm aug}),
\end{equation}
and its inverse on its image is
\begin{equation}\label{app:cross:source-graph-inverse}
 (Z,\Psi_{\mathbf b}Z,r_{\rm aug})
       \longmapsto (Z,r_{\rm aug}).
\end{equation}
The image of $\mathsf G_{Y,\mathbf b}$ is the
\emph{graph} of $Z\mapsto\Psi_{\mathbf b}Z$,
together with the independent data $r_{\rm aug}$. Thus,
$\Psi_{\mathbf b}Z$ is determined by $Z$, and the inverse simply
discards this derived entry. In the triangular elimination leading to the scalar equations,
one Cartesian source component is replaced by $\Psi_{\mathbf b}Z$,
with a nonzero coefficient on that component.
% AUTHOR QUERY: Display the triangular row operation and its inverse, and
% identify this graph with the forcing ranges used by the diagonal inverses.
% The forgetful inverse on the graph alone does not establish the equivalence.

The derivative orders in $\Psi_{\mathbf b}Z$, before applying
the angular multipliers, are listed in \cref{tab:source-orders}.

\begin{table}[htp]
\centering
\begin{tabular}{lcccc}
\toprule
 & $g_+$ & $g_-$ & $g_3$ & $h$\\
\midrule
derivatives $(D_y,D_\zeta)$ at zero curvature & $(1,1)$ & $(1,1)$ & $(0,1)$ & $(2,0)$\\
total order at a general state & $\le2$ & $\le2$ & $\le1$ & $\le2$\\
\bottomrule
\end{tabular}
\caption{Derivative orders in the scalar source $\Psi_{\mathbf b}Z$,
before the angular multipliers are applied.}
\label{tab:source-orders}
\label{app:cross:source-order-matrix}
\end{table}

At the unperturbed solution at zero curvature, the terms are, up
to fixed nonzero factors, $\pa_\zeta\pa_zg_+$,
$\pa_\zeta\pa_{\bar z}g_-$, $\pa_\zeta g_3$, and $\Delta_yh$.
This gives the first row of the table. The total orders in the
second row follow from \eqref{app:cross:raw-current-rows}.
Each entry is a composition of at most two first-order operators,
and differentiation of a coefficient leaves one fewer derivative
on the source.

Since $\cR=\pa_\theta$, the angular multipliers act by
\[
 (\cR-2i)e^{im\theta}=i(m-2)e^{im\theta},\qquad
 (\cR+2i)e^{im\theta}=i(m+2)e^{im\theta}.
\]
They can therefore be inverted on the odd
modes and on the even modes with $|m|\ge4$. At $m=2$ and $m=-2$,
respectively, one multiplier vanishes. We retain the corresponding
equation and its source as part of the affine component, instead of
dividing by that multiplier. The detailed reconstruction appears in
\cref{app:reference}.

We retain the equations for the angular mean and zero cell
frequency, together with the remaining mean, scale, gauge, chart,
and compatibility conditions, and solve them successively in the
triangular system. Their right-hand sides, together with the
coordinates in \eqref{app:cross:source-graph-chart}, give the
source variables used by the block inverses.

\subsection{The order of the reduced scalar equation}

We eliminate $A$, $B$, and $X$ from
\eqref{app:cross:raw-current-rows} to obtain the scalar equation for
$C=T\cdot U$ introduced in \eqref{app:cross:covariants}. Its positive
spectral part supplies the equation for the positive block. Although the calculation appears to produce
third derivatives, an angular factor cancels and leaves an
equation of order two. We first verify this at zero curvature,
using
\[
 D=r\pa_r,\qquad R=\pa_\theta,
 \qquad E=D^2+R^2.
\]
On a nonzero angular mode, $R$ is invertible. Keeping the terms of
\eqref{app:cross:raw-current-rows} in which a derivative falls on an
unknown gives
\begin{equation}\label{app:cross:canonical-rows}
 \begin{aligned}
 n_0&=RX-RB,& n_1&=DX-RA,\\
 n_2&=LRC-\pa_\zeta X,& n_3&=LDA+LRB+r^2\pa_\zeta C.
 \end{aligned}
\end{equation}
We write $n=(n_0,n_1,n_2,n_3)$ for these four expressions. Up to a
nonzero normalization, the highest-order part of the source
conversion is the combination
\begin{equation}\label{app:cross:canonical-source-row}
 \Phi^{\rm pr}n=(LR)^{-1}\left\{
  En_2+\pa_\zeta\left(Dn_1+Rn_0+\frac RL n_3\right)\right\}.
\end{equation}
Since $D$, $R$, and $\pa_\zeta$ commute, substitution gives
\begin{equation}\label{app:cross:canonical-syzygy}
 Dn_1+Rn_0+\frac RL n_3
 =EX+\frac{r^2}{L}R\pa_\zeta C,
 \qquad
 \Phi^{\rm pr}n=\left(E+\frac{r^2}{L^2}\pa_\zeta^2\right)C.
\end{equation}
Thus, the $A$, $B$, and $X$ terms cancel, and the third-order
numerator in \eqref{app:cross:canonical-source-row} contains the angular
factor $R$ that appears in the denominator. Canceling this factor gives the
second-order equation for $C$.

For a general frame, we first absorb the nonzero factors of $L$
into the third and fourth equations. Thus, $\nu$ below is the
frequency corresponding to $L^{-1}\pa_\zeta$. We set
\begin{equation}\label{app:cross:gram}
 G_{\mathbf b}=\mathbb B_{\mathbf b}^*\mathbb B_{\mathbf b},
 \qquad \xi=(d,\omega,\nu),\qquad
 (\alpha,\beta,\chi)=\xi^*G_{\mathbf b}^{-1},
\end{equation}
so that the principal symbol at the state $\mathbf b$ has the rows
\[
 n_0=\omega(X-B),\quad n_1=dX-\omega A,\quad
 n_2=\omega C-\nu X,\quad
 n_3=\alpha A+\beta B+\chi C,
\]
where $d$ and $\omega$ are the radial and angular frequencies.
For $\omega\ne0$, we set
\begin{equation}\label{app:cross:current-source-row}
 s_{G_{\mathbf b}}(\xi)
 =\left(\frac{\nu\beta}{\omega},
         \frac{\nu\alpha}{\omega},
         \frac{\alpha d+\beta\omega}{\omega},\nu\right).
\end{equation}
We let $\ell_{1,G_{\mathbf b}}(\xi)$ be the matrix of
$n_0,\ldots,n_3$ in the unknowns $(A,B,C,X)$, and let $e_C^*$ select
the $C$ entry. Multiplying the four equations by the entries of
$s_{G_{\mathbf b}}(\xi)$ and adding, we obtain
\begin{equation}\label{app:cross:current-syzygy}
 s_{G_{\mathbf b}}(\xi)\ell_{1,G_{\mathbf b}}(\xi)
   =(\xi^*G_{\mathbf b}^{-1}\xi)e_C^*.
\end{equation}
Indeed, the coefficients of $A$, $B$, and $X$ cancel, while the
coefficient of $C$ is
$\alpha d+\beta\omega+\chi\nu=\xi^*G_{\mathbf b}^{-1}\xi$.
Multiplying by $\omega$ gives a polynomial identity. The angular mean
is kept in the triangular system, and the exceptional modes are treated
by the affine equations. On the compatible Cartesian spaces on the cap, that is,
on sources whose quotients \eqref{eq:quotient-variables} are smooth,
\eqref{app:cross:current-syzygy} extends across the axis.
% AUTHOR QUERY: Prove the cap extension and principal cap cancellation in
% Cartesian variables. The annular Gram inverse is singular at the axis;
% the frame map of eq:A-frame-map is defined to be the identity on the cap.
When the coefficients vary, a derivative falling on the frame, on the
inverse of $G_{\mathbf b}$, or on a projection, leaves one fewer derivative
on the unknown. Thus, the reduced scalar operator is still of order two,
and after canceling its highest-order off-diagonal part, the remaining
terms have order at most one.

\subsection{From Cartesian to block variables}
\label{app:fixed-graphs:adapters}

We now combine the chart derivative, the annular frame map, and
the spectral projections to express a Cartesian variation in
block variables. We give the inverse of each step so that the
variation can be reconstructed. The restrictions to the cap and
annulus, and the interface data, are treated in the next subsection.

We write $\mathbf b=(\rho,a,\eps,p)$, with
$a=(v,w)=\mathscr S_q(x)$, $x=T_q\widehat x$, and
$q=(\rho,\eps,p)$. For a direction $\widehat h$ in the fixed
domain $\cX_\gamma^s$, the chart derivative gives the physical
variation $(U,S)$,
\begin{equation}\label{eq:A-slice-tangent-factor}
 \Theta_{\rm sl,\mathbf b}:=D_x\mathscr S_q(x)T_q,
 \qquad
 \Theta_{\rm sl,\mathbf b}\widehat h=(U,S),
\end{equation}
as in \eqref{eq:state-chart-tangent}. Its inverse first extracts the
tangential coefficient $\dot\tau$ of the linear Taylor term and then sets
\begin{equation}\label{eq:A-slice-tangent-inverse}
 \dot{\widetilde v}
 =U-\left[-\frac{\tau\cdot\dot\tau}{2a(\tau)}\iota M_q
                    +e_T\otimes\dot\tau\right]y,
 \qquad
 \widehat h=T_q^{-1}(\dot\tau,\dot{\widetilde v},S).
\end{equation}
Here, $a(\tau)=\sqrt{1-|\tau|^2/2}$ is the scalar factor in
\eqref{eq:normalized-affine-map}, whereas $a=(v,w)$ denotes the state.
The tangential first derivatives on the axis determine
$\dot\tau_j=e_T\cdot\pa_{y_j}U(0,\zeta)$ for $j=1,2$.
Subtracting their linear Taylor contribution leaves
$\dot{\widetilde v}$ with vanishing value and first derivatives on the axis.

On the annulus, we express this variation using the frame
\eqref{app:cross:current-frame} and its Gram matrix,
\begin{equation}\label{eq:A-frame-factor}
 p=\cD v,\qquad q_1=\cR v,\qquad T=L^{-1}D_\eps v,
 \qquad \mathbb B=(p,q_1,T),\qquad G_{\mathbf b}=\mathbb B^*\mathbb B,
\end{equation}
and we define the pointwise map and its inverse by
\begin{align}
 \Theta_{\rm fr,\mathbf b}(U,S)
   &=(A,B,C,X)
     =(p\cdot U,q_1\cdot U,T\cdot U,S-q_1\cdot U),
                                                        \label{eq:A-frame-map}\\
 \Theta_{\rm fr,\mathbf b}^{-1}(A,B,C,X)
   &=(\mathbb Bc,X+B),\qquad
     c=G_{\mathbf b}^{-1}(A,B,C)^T.                 \label{eq:A-frame-inverse}
\end{align}
On the cap, we keep the Cartesian components $(U,S)$, so the
frame map acts as the identity there. On the annulus, the bound
$|\det\mathbb B|\ge c_{\rm det}r^2$ following
\eqref{app:cross:current-frame} makes
\eqref{eq:A-frame-map}--\eqref{eq:A-frame-inverse} mutually inverse
analytic maps of order zero.

Finally, the projections \eqref{app:cross:domain-range-contours}
split the unknowns and sources into the two blocks,
\begin{equation}\label{eq:A-contour-factors}
 \begin{aligned}
  \Theta_{\rm sp,\mathbf b}^Xz
    &=(P_{-,\mathbf b}^Dz,P_{+,\mathbf b}^Dz),&
  (\Theta_{\rm sp,\mathbf b}^X)^{-1}(z_-,z_+)&=z_-+z_+,\\
  \Theta_{\rm sp,\mathbf b}^Yf
    &=(P_{-,\mathbf b}^Rf,P_{+,\mathbf b}^Rf),&
  (\Theta_{\rm sp,\mathbf b}^Y)^{-1}(f_-,f_+)&=f_-+f_+.
 \end{aligned}
\end{equation}
The inverses are defined on the images of the corresponding
projections. We then write the minus component in center and affine coordinates
and the positive component in scalar coordinates, retaining the
exceptional angular modes $m=\pm2$ separately. The
source graph map \eqref{app:cross:source-graph-chart} is applied
after this splitting. For the mean, scale, gauge, chart, and zero
cell frequency equations, we apply the inverse of each diagonal
operator, giving the triangular form described after
\eqref{app:cross:principal-intertwining}.

\subsection{Cutting and gluing at the interfaces}

When solving the cap and annulus equations separately, we allow
their solutions to have different values and conormal derivatives
on the common interface. These values and derivatives are the
\emph{Cauchy data}. Their difference, with the signs specified below,
is the \emph{mismatch}, which we keep as an independent unknown. Subtracting an extension of the
mismatch from the annular field gives a matching pair. Retaining
the mismatch as a separate coordinate makes this change invertible.

For the value, conormal, and outer-boundary coordinates, we use
$\Lambda=\langle D_\zeta\rangle$, so that
$(\Lambda^t u)_n=\lambda_n^t u_n$ with $\lambda_n$ as in
\eqref{app:cross:cell-weights}. The change of trace variables is
\begin{align}
 \mathcal T_\Lambda(z_D,z_N,z_B)
  &=(\Lambda^{1/2}z_D,\Lambda^{-1/2}z_N,
                         \Lambda^{-1/2}z_B),\notag\\
 \mathcal T_\Lambda^{-1}(d,n,b)
  &=(\Lambda^{-1/2}d,\Lambda^{1/2}n,\Lambda^{1/2}b).
                                                        \label{eq:A-balanced-trace-factor}
\end{align}
Here, $z_D$ and $z_N$ are the interface value and conormal data,
and $z_B$ is the datum for the first-order outer boundary condition.
The factors $\Lambda^{1/2}$ and $\Lambda^{-1/2}$ account for the
one-derivative difference in \eqref{eq:A-Cauchy-coretraction-estimate}.
We call the resulting norms the \emph{balanced trace norms}.

At the interface $\Sigma$ between a cap and an annulus, we write
$u_c$ and $u_a$
for the two restrictions, and we define their \emph{value and outward conormal
mismatches} by
\begin{equation}\label{eq:A-Cauchy-mismatches}
 d_{\mathbf b}(u_c,u_a)
   =\gamma_D^au_a-\gamma_D^cu_c,
 \qquad
 n_{\mathbf b}(u_c,u_a)
   =\gamma_{N,\mathbf b}^{a,{\rm out}}u_a
      +\gamma_{N,\mathbf b}^{c,{\rm out}}u_c.
\end{equation}
Here, $\gamma_D$ is the trace,
$\gamma_{N,\mathbf b}^{\rm out}$ is the outward conormal derivative
at $\mathbf b$, and the superscripts $a$ and $c$ indicate the
annular and cap sides. The value mismatch is independent of the
state and will be denoted by $d$. To compare the conormal data
with those at the constant ellipse $\mathbf b_\rho$, we also define
the \emph{reference conormal mismatch}
\begin{equation}\label{eq:A-reference-conormal-mismatch}
 h_\rho(u_c,u_a)
 :=\gamma_{N,\mathbf b_\rho}^{a,{\rm out}}u_a
       +\gamma_{N,\mathbf b_\rho}^{c,{\rm out}}u_c.
\end{equation}
Using the annular outward normal $\nu$ as a common direction,
we set $q_\Sigma=\pa_\nu u_a-\pa_\nu u_c$. The opposite outward
normals explain the sum in \eqref{eq:A-Cauchy-mismatches}.
Since the coefficients agree on $\Sigma$, both conormal mismatches
can be expressed in terms of $q_\Sigma$ and $d$,
\begin{equation}\label{eq:A-current-reference-conormal-decomposition}
 n_{\mathbf b}=e_{\Sigma,\mathbf b}q_\Sigma
       +T_{\Sigma,\mathbf b}d,
 \qquad
 h_\rho=e_{\Sigma,\rho}q_\Sigma+T_{\Sigma,\rho}d,
\end{equation}
where $e_{\Sigma,\mathbf b}$ and $e_{\Sigma,\rho}$ are the
invertible coefficients of the common normal derivative, and
$T_{\Sigma,\mathbf b}$ and $T_{\Sigma,\rho}$ are tangential
operators of order at most one. Eliminating $q_\Sigma$ gives
\begin{align}
 \binom d{n_{\mathbf b}}
 &=\mathcal T_{\Sigma,\mathbf b,\rho}^{CN}
       \binom d{h_\rho},
 &\mathcal T_{\Sigma,\mathbf b,\rho}^{CN}
 &=\begin{pmatrix}I&0\\K_{\Sigma,\mathbf b,\rho}
                         &E_{\Sigma,\mathbf b,\rho}\end{pmatrix},
                                                        \label{eq:A-current-reference-row-map}\\
 E_{\Sigma,\mathbf b,\rho}
 &:=e_{\Sigma,\mathbf b}e_{\Sigma,\rho}^{-1},
 &K_{\Sigma,\mathbf b,\rho}
 &:=T_{\Sigma,\mathbf b}
       -E_{\Sigma,\mathbf b,\rho}T_{\Sigma,\rho},\notag\\
 (\mathcal T_{\Sigma,\mathbf b,\rho}^{CN})^{-1}
 &=\begin{pmatrix}I&0\\
   -E_{\Sigma,\mathbf b,\rho}^{-1}K_{\Sigma,\mathbf b,\rho}
       &E_{\Sigma,\mathbf b,\rho}^{-1}\end{pmatrix}.
                                                        \label{eq:A-current-reference-row-inverse}
\end{align}
The map and its inverse are uniformly bounded on the balanced
trace spaces. Indeed, the one-derivative difference between the
value and conormal norms absorbs the order of
$K_{\Sigma,\mathbf b,\rho}$. In particular,
\begin{equation}\label{eq:A-identical-glue-kernels}
 d=0,\ n_{\mathbf b}=0
 \quad\Longleftrightarrow\quad d=0,\ h_\rho=0,
\end{equation}
so using the conormal derivative at $\mathbf b$ or at
$\mathbf b_\rho$ gives the same matching condition. For the minus block,
the cap uses Cartesian coordinates and the annulus uses the frame and
spectral coordinates of \eqref{eq:A-frame-map} and
\eqref{eq:A-contour-factors}. We denote the induced change of Cauchy
coordinates by $T_{-,\mathbf b}$. If $\Gamma_{c,-,\mathbf b}u$ and
$\Gamma_{a,-,\mathbf b}u$ represent the Cauchy data of the same field
$u$ in these two coordinate systems, then
\[
 T_{-,\mathbf b}\Gamma_{c,-,\mathbf b}u
   =\Gamma_{a,-,\mathbf b}u.
\]
This map has order zero in the balanced trace norms. We apply it to the
cap data before subtracting them from the annular data. For the positive
block, we use the value and conormal mismatches directly.

We now subtract their annular extensions,
\eqref{eq:A-Cauchy-N-coretraction}--\eqref{eq:A-Cauchy-D-coretraction},
by setting
\begin{equation}\label{eq:A-Cauchy-zero-remainder}
 u_a^0=u_a-C_{\Sigma,D,\mathbf b}^a
                 d_{\mathbf b}(u_c,u_a)
             -C_{\Sigma,N,\mathbf b}^a
                 n_{\mathbf b}(u_c,u_a).
\end{equation}
% AUTHOR QUERY: Value and conormal matching do not imply global H^s
% regularity for s >= 5/2: a second normal derivative may jump. Specify
% higher-jet matching or a restricted graph domain that supplies it, and
% the smooth-core completions on which the all-order reconstruction,
% eq:A-ambient-kernel and finite-index descent identities hold.
By \eqref{eq:A-Cauchy-coretraction-identity}, the pair $(u_c,u_a^0)$
has matching values and opposite outward conormal derivatives.
It therefore agrees to first order across $\Sigma$. Keeping
the original mismatch with this pair gives the mutually inverse maps
\begin{align}
 \mathfrak B_{\Sigma,\mathbf b}(u_c,u_a)
   &=\bigl(u_c,u_a^0,
          d_{\mathbf b}(u_c,u_a),n_{\mathbf b}(u_c,u_a)\bigr),
                                                        \label{eq:A-Cauchy-break-map}\\
 \mathfrak B_{\Sigma,\mathbf b}^{-1}(u_c,u_a^0,d,n)
   &=\bigl(u_c,u_a^0+C_{\Sigma,D,\mathbf b}^ad
                         +C_{\Sigma,N,\mathbf b}^an\bigr).
                                                        \label{eq:A-Cauchy-glue-map}
\end{align}
On the axis, $J_Au=(u,\pa_{y_1}u,\pa_{y_2}u)|_{y=0}$ records
the value and first derivatives, and $E_A$ is their extension from
\eqref{eq:A-axis-coretractions}. The identity $J_AE_A=I$ gives
\[
 u\longmapsto\bigl((I-E_AJ_A)u,J_Au\bigr),\qquad
 (u_0,z)\longmapsto u_0+E_Az,\qquad J_Au_0=0.
\]
For the outer boundary condition, the corresponding maps are
\begin{equation}\label{eq:A-physical-break-glue}
 u\longmapsto\bigl((I-C_{\pa,M}B_M)u,B_Mu\bigr),
 \qquad
 (u_0,b)\longmapsto u_0+C_{\pa,M}b,\qquad B_Mu_0=0,
\end{equation}
where $B_M$ is the boundary operator
\eqref{eq:physical-slice-row} and $C_{\pa,M}$ its extension
\eqref{eq:A-physical-row}. The identity $B_MC_{\pa,M}=I$ ensures
that the remainder satisfies the boundary condition and that
the two maps are inverse.

\smallskip\noindent\textit{The shell decomposition.}
On each dyadic shell, we use the extensions
\eqref{eq:A-Cauchy-N-coretraction}--\eqref{eq:A-Cauchy-D-coretraction} at
the constant ellipse $\mathbf b_\rho$ with the same ellipse parameter as the
state, and we write
\begin{equation}\label{eq:A-reference-Cauchy-coretractions}
 C_{\Sigma,N,\rho}:=C_{\Sigma,N,\mathbf b_\rho},
 \qquad C_{\Sigma,D,\rho}:=C_{\Sigma,D,\mathbf b_\rho}.
\end{equation}
The reference coefficients are independent of $\zeta$, so these
extensions commute with $Q_j$. At a general state, coefficients
may mix cell frequencies. We therefore reconstruct the full field
before converting the reference conormal data to those at the state.

For the shell $Q_j$ of \eqref{app:minus:eq:sharp-shells}, we use
\[
 r_j=2TL2^{-j},\qquad
 0\le r\le r_j\ \text{on the cap},\qquad
 t_j=2^jr/L\in[0,2T].
\]
Here, $T>0$ is the fixed matching scale, and $J$ is large enough that
$r_J<1$. We restrict $Q_ju$ to the cap and its complementary annulus
$r_j\le r\le1$, rescale the cap by $t_j=2^jr/L$, and apply
\eqref{eq:A-Cauchy-break-map} at $\mathbf b_\rho$. Taking
the square sum over $j\ge J$ gives the spaces
\begin{equation}\label{eq:A-infinite-broken-bundle}
 \mathscr X_{{\rm br},H}^s
   =\bigoplus_{j\geq J}^{\ell_s^2}
       \bigl(X_{c,j}^s\oplus X_{a,j}^s\bigr),
 \qquad
 \mathscr Z_{{\rm aux},H}^s
   =\bigoplus_{j\geq J}^{\ell_s^2}Z_j^s .
\end{equation}
Here, $X_{c,j}^s$ and $X_{a,j}^s$ are the cap and annulus spaces
on shell $j$, while $Z_j^s$ contains the interface pair
$(d,h_\rho)$. For the minus block, this pair is replaced by the
equivalent Cauchy datum for the outward evolution on the annulus,
as described in \cref{app:minus}.
The notation $\ell_s^2$ means that, for a sequence $(u_j)$ of shell
components, its squared norm is $\sum_{j\ge J}\|u_j\|_s^2$.
The component norm $\|\cdot\|_s$ already includes the cell-frequency
weights of \cref{tab:block-weights} and the common factor
$\lambda_n^s$.

On each rescaled cap at the state $\mathbf b$, let
$\mathcal L_{c,-}$ and $\mathcal L_{c,+}$ denote the differential operators for the minus and positive equations
in the rescaled coordinates. A \emph{graph norm} measures both a
solution $u$ and its source $\mathcal L_{c,\pm}u$. The two norms have
the form
\[
 \begin{aligned}
 \|u\|_{c,-,s}
   &=\|u\|_{H^{s+2}}+\|\mathcal L_{c,-}u\|_{H^s},\\
 \|u\|_{c,+,s}
   &=\|u\|_{H^{s+1}}+\|\mathcal L_{c,+}u\|_{H^{s-1}}.
 \end{aligned}
\]
Here, the Sobolev norms are taken on the fixed rescaled cap, with the
analytic weight \eqref{app:cross:phase}.
The values and conormal derivatives of the positive component
on the interface have regularities $H^{s+1/2}$ and $H^{s-1/2}$,
respectively, on both sides.
We return each cap to the original coordinates and join it to its
own complementary annulus before summing over the shells. Variable
coefficients then act on the reconstructed field. This avoids
comparing norms on disks rescaled by different shell frequencies.

For a field $u$ and source $F$, we write
$\mathscr A_{a,j}^{X}u$ and $\mathscr A_{a,j}^{Y}F$ for the
restrictions of $Q_ju$ and $Q_jF$ to the annulus of shell $j$.
Together with the rescaled cap restrictions, these give the
domain and source maps for the shell decomposition,
\begin{equation}\label{eq:A-sharp-shell-break-insertion}
 \begin{aligned}
 \mathfrak B_{{\rm sh},J,\rho}^{X}
  &:=\bigoplus_{j\geq J}^{\ell_s^2}
      \mathfrak B_{\Sigma_j,\mathbf b_\rho}
       \bigl(\mathscr A_{X,j},\mathscr A_{a,j}^{X}\bigr),\\
 \mathfrak B_{{\rm sh},J,\rho}^{Y,\rm raw}
  &:=\bigoplus_{j\geq J}^{\ell_s^2}
       \bigl(\mathscr A_{Y,j},\mathscr A_{a,j}^{Y},0_{Z_j}\bigr).
 \end{aligned}
\end{equation}
Here, $\mathscr A_{X,j}u=\mathscr R_j^XQ_ju$ and
$\mathscr A_{Y,j}F=\mathscr R_j^YQ_jF$. The maps
$\mathscr R_j^X$ and $\mathscr R_j^Y$ restrict an unknown and a source
to $r\le r_j$ and pull them back to the fixed cap by
$y=L2^{-j}\widetilde y$, with $|\widetilde y|\le2T$.
The zero entry $0_{Z_j}$ records that an original Cartesian interior
source prescribes no independent interface data.
The reference extensions preserve each shell, and we denote
the reconstruction maps by
\begin{align}
 \mathfrak R_{{\rm sh},J,\rho}^{X}
   &:=(\mathfrak B_{{\rm sh},J,\rho}^{X})^{-1},\notag\\
 \mathfrak R_{{\rm sh},J,\rho}^{Y,\rm raw}
   &:=(\mathfrak B_{{\rm sh},J,\rho}^{Y,\rm raw})^{-1}
       \quad\hbox{on the image of the Cartesian interior sources.}
                                                        \label{eq:A-sharp-shell-reconstruction}
\end{align}
On the domain, we apply
$\mathfrak B_{\Sigma_j,\mathbf b_\rho}^{-1}$ and join each cap
to its annulus before summing over the shells. On the range,
we join the interior entries and discard the zero auxiliary entry.
The local inverse identities and $\sum_{j\ge J}Q_j=H_J$ show
that these maps invert the restrictions of Cartesian fields.

\smallskip\noindent\textit{The remaining constraints.}
We add the low-frequency components $|n|<2^J$ and the chart
coordinates to \eqref{eq:A-infinite-broken-bundle}. Allowing
arbitrary entries gives the product
$\cX_{\gamma,{\rm amb},\mathbf b}^s$, before matching and the
other constraints are imposed. We let $\mathscr C_{X,\mathbf b}$
record the reference mismatches
$(d,h_\rho)$ on each shell, the low-frequency interface data,
and the constraints in \eqref{eq:complete-domain-constraint}.
Its target $\mathscr Z_{{\rm aux},\mathbf b}^s$ contains all
these constraint values. At fixed $\rho$, the shell extensions
are fixed, and the remaining state dependence comes from the
frame and chart maps. Since the conormal derivatives at
$\mathbf b_\rho$ and $\mathbf b$ impose the same matching condition by
\eqref{eq:A-identical-glue-kernels}, we identify the Cartesian
domain with
\begin{equation}\label{eq:A-ambient-kernel}
 \cX_\gamma^s\simeq
 \ker\bigl(\mathscr C_{X,\mathbf b}:
       \cX_{\gamma,{\rm amb},\mathbf b}^s
       \longrightarrow\mathscr Z_{{\rm aux},\mathbf b}^s\bigr).
\end{equation}
Subtracting the extensions in
\eqref{eq:A-Cauchy-zero-remainder},
\eqref{eq:A-axis-coretractions}, and
\eqref{eq:A-physical-break-glue}, and reconstructing the field,
gives a bounded map
$\mathscr R_{X,\mathbf b}^{\rm br}:\cX_{\gamma,{\rm amb},\mathbf b}^s
\to\cX_\gamma^s$. It removes the mismatch and constraint data recorded
by $\mathscr C_{X,\mathbf b}$ before reconstructing a constrained
Cartesian field. Conversely,
$\mathscr I_{X,\mathbf b}^{\rm br}$ restricts a constrained
field to the caps and annuli with zero constraint data.
The matching identities and $\sum_{j\ge J}Q_j=H_J$ give
\begin{align}
 \mathscr R_{X,\mathbf b}^{\rm br}
      \mathscr I_{X,\mathbf b}^{\rm br}&=I,
 &\mathscr C_{X,\mathbf b}
      \mathscr I_{X,\mathbf b}^{\rm br}&=0,             \label{eq:A-break-left}\\
 \mathscr I_{X,\mathbf b}^{\rm br}
      \mathscr R_{X,\mathbf b}^{\rm br}U&=U
 &&\text{if }\mathscr C_{X,\mathbf b}U=0.               \label{eq:A-break-zero-right}
\end{align}
Including the constraint data, we obtain the isomorphism
\begin{equation}\label{eq:A-stabilized-break-chart}
 \mathscr B_{X,\mathbf b}U
 :=\bigl(\mathscr R_{X,\mathbf b}^{\rm br}U,
          \mathscr C_{X,\mathbf b}U\bigr):
 \mathscr X_{{\rm br},\mathbf b}^s
 \xrightarrow{\ \sim\ }
 \mathscr X_{{\rm sp},\mathbf b}^s
       \oplus\mathscr Z_{{\rm aux},\mathbf b}^s,
\end{equation}
whose inverse is
\begin{equation}\label{eq:A-stabilized-break-inverse}
 \mathscr B_{X,\mathbf b}^{-1}(u,z)
 =\mathscr I_{X,\mathbf b}^{\rm br}u
       +\mathscr E_{X,\mathbf b}^{\rm br}z.
\end{equation}
Here, $\mathscr X_{{\rm br},\mathbf b}^s$ denotes the
unconstrained product $\cX_{\gamma,{\rm amb},\mathbf b}^s$,
and $\mathscr X_{{\rm sp},\mathbf b}^s$ its constrained subspace,
identified with $\cX_\gamma^s$. The map
$\mathscr E_{X,\mathbf b}^{\rm br}$ extends the constraint data
using the reference shell extensions and the axis and boundary
extensions of \cref{app:fixed-graphs}. We normalize it by
$\mathscr R_{X,\mathbf b}^{\rm br}\mathscr E_{X,\mathbf b}^{\rm br}=0$
and $\mathscr C_{X,\mathbf b}\mathscr E_{X,\mathbf b}^{\rm br}=I$,
which gives
\begin{equation}\label{eq:A-cap-trace-inverses}
 \mathscr B_{X,\mathbf b}^{-1}\mathscr B_{X,\mathbf b}=I,
 \qquad
 \mathscr B_{X,\mathbf b}\mathscr B_{X,\mathbf b}^{-1}=I.
\end{equation}

\subsection{The augmented block operator}

We now assemble the transfer maps and define the operator
$\mathbb D_{\mathbf b}$ to be inverted in the following sections.
Besides the interior equations, it records the interface, boundary,
axis, and remaining constraints. We call these additional equations
the \emph{auxiliary equations} and the full operator the
\emph{augmented block operator}. Their right-hand sides vanish
for the original Cartesian problem.

\smallskip\noindent\textit{Fixed block spaces.}
Let $\Theta_{{\rm bulk},\mathbf b}^{Y,\rm raw}$ collect the interior
source restrictions and cap rescalings in
\eqref{eq:A-sharp-shell-break-insertion}, together with the
low-frequency source coordinates. It does not include the auxiliary
data. Before applying the trace weights and identifying
the spaces with fixed reference spaces, we denote the domain
product by $\mathscr X_{\rm ret}^s$ and the range product by
$\mathscr Y_{\rm aug}^s$. The former keeps the cap and annulus
fields separate, while the latter contains their sources and
the auxiliary data. Their low- and high-frequency parts are
\begin{equation}\label{eq:A-low-high-retained-products}
 \begin{aligned}
  \mathscr X_{{\rm ret},L}^s&=P_J^X\mathscr X_{\rm ret}^s,&
  \mathscr X_{{\rm ret},H}^s&=H_J^X\mathscr X_{\rm ret}^s,\\
  \mathscr Y_{{\rm aug},L}^s&=P_J^Y\mathscr Y_{\rm aug}^s,&
  \mathscr Y_{{\rm aug},H}^s&=H_J^Y\mathscr Y_{\rm aug}^s,
 \end{aligned}
\end{equation}
using the projections of \eqref{eq:A-exact-cell-split}. The
two diagonal blocks act on the high-frequency products.
The projections are applied in these coordinates without
commuting them through the state coefficients.

Axis values and first derivatives carry, respectively, the factors
$\Lambda^{-1}$ and $\Lambda^{-2}$ of \cref{tab:block-weights}. We write
$\Theta_{\rm tr}^{X,Y}$ for the direct sum of
\eqref{eq:A-balanced-trace-factor}, these axis factors, and the identity
on the interior components and the remaining constraint coordinates.

On the range, the extraction of the quotients is the bounded isomorphism
\begin{equation}\label{eq:A-quotient-factor}
 \Theta_Q=\mathscr Q_Y:\cY_{\rm comp}^s\longrightarrow
       \operatorname{Ran}P_Z,\qquad \Theta_Q^{-1}=J_Y,
\end{equation}
from \eqref{eq:compatible-extraction} and
\eqref{eq:range-reconstruction}. Here, $P_Z$ is the projection
\eqref{eq:A-range-projection} imposing compatibility in quotient
coordinates. We apply $\Theta_Q$ separately to the two
outputs of \eqref{eq:A-contour-factors}, followed by the source graph
map of \eqref{app:cross:source-graph-chart},
\begin{equation}\label{eq:A-source-graph-factor}
 \Theta_{\rm src,\mathbf b}(Z,r_{\rm aug})
   =(Z,\Psi_{\mathbf b}Z,r_{\rm aug}),\qquad
 \Theta_{\rm src,\mathbf b}^{-1}
   (Z,\Psi_{\mathbf b}Z,r_{\rm aug})=(Z,r_{\rm aug}).
\end{equation}
Only this step differentiates the prescribed source. Its inverse
discards the derived entry $\Psi_{\mathbf b}Z$. We record the
loss of derivatives by
\begin{equation}\label{eq:A-source-order}
 \ell_Y:=\max_{\alpha}
   \operatorname{ord}_{\cY_\gamma\to\bY_\natg}
       \Psi_{\alpha,\mathbf b}<\infty,
\end{equation}
where $\Psi_{\alpha,\mathbf b}$ gives the scalar source of component
$\alpha$, and $\operatorname{ord}$ denotes the loss between the
indicated scales. Thus, $\ell_Y$ is the number of additional derivatives
required of the Cartesian source to control every converted component.
In the component norms, this means
\[
 \|\Psi_{\alpha,\mathbf b}\mathscr Q_Y f\|_{Y_\alpha^s}
 \le C_s(\mathbf b)\|f\|_{\cY_\gamma^{s+\ell_Y}}.
\]
Here, $Y_\alpha^s$ is the source space for component $\alpha$, including
its weight from \cref{tab:block-weights}. The orders in
\cref{tab:source-orders} give a finite bound $\ell_Y$,
independent of $s$ and uniform on a sufficiently small ball of
states. We also write $\ell:=\ell_Y$. The conversion is applied
once, before the block inverse, to avoid incurring this loss
of derivatives multiple times. In particular, every term in
the Neumann series acts on sources already in block coordinates.

To express the state-dependent spaces in fixed coordinates,
we use, for each domain, range, cap, or trace projection
$P_{\mathbf b}$ with reference projection $P_*$, the map
\begin{equation}\label{eq:A-complete-Kato-factor}
 U_{P,\mathbf b}=P_{\mathbf b}P_*+
       (I-P_{\mathbf b})(I-P_*).
\end{equation}
It equals the identity at $P_{\mathbf b}=P_*$ and is invertible
when the projections are sufficiently close. Since
$P_{\mathbf b}U_{P,\mathbf b}=U_{P,\mathbf b}P_*$, its inverse
maps the range at $\mathbf b$ to the fixed reference range.
We denote the direct sums of these inverses by
$\mathcal K_{X,\mathbf b}$ and $\mathcal K_{Y,\mathbf b}$,
using the weighted $\ell^2$ sums of
\eqref{eq:A-infinite-broken-bundle} on the high-frequency
components. All subsequent direct sums over shells are understood
in this sense. The fixed contours and the uniform cap rescaling
bounds give uniform bounds for these maps.

Let $\bY_{{\rm raw},\natg}^s$ be the block source space
obtained from the compatible Cartesian sources, including the
derived entry $\Psi_{\mathbf b}Z$. The subscript $\mathrm{raw}$
indicates that the independent auxiliary data have not yet been
added. Inclusion into the full space $\bY_\natg^s$ with zero
auxiliary entries and extraction of the interior source give
the maps
\begin{equation}\label{eq:A-zero-range-splitting}
 \iota_{Y,\mathbf b}:\bY_{{\rm raw},\natg}^s
       \longrightarrow\bY_\natg^s,
 \qquad
 p_{Y,\mathbf b}:\bY_\natg^s
       \longrightarrow\bY_{{\rm raw},\natg}^s,
 \qquad p_{Y,\mathbf b}\iota_{Y,\mathbf b}=I.
\end{equation}
In coordinates $(F,z)$, where $F$ is an interior source and $z$ is
the auxiliary data, these maps act by
\[
 \iota_{Y,\mathbf b}F=(F,0),\qquad
 p_{Y,\mathbf b}(F,z)=F,\qquad
 \operatorname{aux}_{Y,\mathbf b}(F,z)=z.
\]
Thus, the reverse composition sends $(F,z)$ to $(F,0)$.

\smallskip\noindent\textit{The transfer maps.}
Composing the chart and frame maps, the spectral projections,
the restrictions, the trace weights, and the identifications
with the reference spaces gives the domain transfer
\begin{equation}\label{eq:A-JX-factorization}
 \mathsf J_{X,\mathbf b}
 :=\mathcal K_{X,\mathbf b}\Theta_{\rm tr}^X
       \mathscr I_{X,\mathbf b}^{\rm br}
       \Theta_{\rm sp,\mathbf b}^X
       \Theta_{\rm fr,\mathbf b}\Theta_{\rm sl,\mathbf b}:
 \cX_\gamma^s\longrightarrow\bX_\natg^s.
\end{equation}
Reading the composition from right to left, we apply the chart derivative
\eqref{eq:A-slice-tangent-factor}, the frame map
\eqref{eq:A-frame-map}, the projections \eqref{eq:A-contour-factors},
and the restriction map characterized by \eqref{eq:A-break-left}.
The trace weights and the inverse identifications
\eqref{eq:A-complete-Kato-factor} then put these components in the
fixed block spaces. Each factor acts as the identity on the components
it does not affect. To return to the constrained Cartesian domain, we
subtract the extensions of the constraint values and reverse
the coordinate changes,
\begin{equation}\label{eq:A-JX-inverse-factorization}
 \mathsf R_{X,\mathbf b}
 :=\Theta_{\rm sl,\mathbf b}^{-1}
       \Theta_{\rm fr,\mathbf b}^{-1}
       (\Theta_{\rm sp,\mathbf b}^X)^{-1}
       \mathscr R_{X,\mathbf b}^{\rm br}
       (\Theta_{\rm tr}^X)^{-1}\mathcal K_{X,\mathbf b}^{-1}:
 \bX_\natg^s\longrightarrow\cX_\gamma^s.
\end{equation}
Writing the constraint map $\mathscr C_{X,\mathbf b}$ in
the same fixed coordinates gives
\begin{equation}\label{eq:A-native-constraint-map}
 \mathsf C_{X,\mathbf b}:\bX_\natg^s
       \longrightarrow\mathscr Z_{{\rm aux},\natg}^s.
\end{equation}
From \eqref{eq:A-break-left}--\eqref{eq:A-stabilized-break-inverse}, we
then obtain the identities
\begin{equation}\label{eq:A-domain-retract-identities}
 \mathsf R_{X,\mathbf b}\mathsf J_{X,\mathbf b}=I,
 \qquad \mathsf C_{X,\mathbf b}\mathsf J_{X,\mathbf b}=0,
 \qquad
 \mathsf J_{X,\mathbf b}\mathsf R_{X,\mathbf b}U=U
       \quad\text{if }\mathsf C_{X,\mathbf b}U=0.
\end{equation}
We write the extension $\mathscr E_{X,\mathbf b}^{\rm br}$ of
\eqref{eq:A-stabilized-break-inverse} in the same fixed coordinates
and denote it by
$\mathsf E_{X,\mathbf b}:\mathscr Z_{{\rm aux},\natg}^s
\to\bX_\natg^s$. The normalization of this extension gives
\begin{equation}\label{eq:A-native-domain-decomposition}
 \mathsf R_{X,\mathbf b}\mathsf E_{X,\mathbf b}=0,
 \qquad
 \mathsf C_{X,\mathbf b}\mathsf E_{X,\mathbf b}=I,
 \qquad
 U=\mathsf J_{X,\mathbf b}\mathsf R_{X,\mathbf b}U
       +\mathsf E_{X,\mathbf b}\mathsf C_{X,\mathbf b}U.
\end{equation}
Thus, a block field is determined by its constrained Cartesian
part and its constraint data. Equivalently,
$(\mathsf R_{X,\mathbf b},\mathsf C_{X,\mathbf b})$ identifies
$\bX_\natg^s$ with
$\cX_\gamma^s\oplus\mathscr Z_{{\rm aux},\natg}^s$,
while $\mathsf J_{X,\mathbf b}$ identifies the Cartesian domain
with $\ker\mathsf C_{X,\mathbf b}$.

For the source transfer, we split the compatible source into
blocks, extract its quotient coordinates, and restrict to the
cap and annulus. The map $\iota_{\rm aux}^{\rm pre}$ acts at this
stage by $F\mapsto(F,0)$, before the trace weights and fixed-space
identifications are applied. Applying the source graph map,
trace weights, and identifications with the reference spaces then gives
\begin{equation}\label{eq:A-JY-factorization}
 \iota_{Y,\mathbf b}\mathsf J_{Y,\mathbf b}^{\rm raw}
 :=\mathcal K_{Y,\mathbf b}\Theta_{\rm tr}^Y
       \Theta_{\rm src,\mathbf b}\iota_{\rm aux}^{\rm pre}
       \Theta_{{\rm bulk},\mathbf b}^{Y,\rm raw}\Theta_Q
       \Theta_{\rm sp,\mathbf b}^Y:
 \cY_\gamma^{s+\ell_Y}\longrightarrow\bY_\natg^s.
\end{equation}
Reversing these steps, we discard the derived entry
$\Psi_{\mathbf b}Z$, recombine the interior Fourier components,
and recover the Cartesian source from its quotient coordinates.
This gives a left inverse
\begin{equation}\label{eq:A-JY-inverse-factorization}
 \mathsf R_{Y,\mathbf b}^{\rm raw}:
 \bY_{{\rm raw},\natg}^s\longrightarrow\cY_\gamma^s.
\end{equation}
On smooth Cartesian source coordinates, these maps are mutually
inverse. At finite Sobolev regularity, the forward conversion loses
$\ell_Y$ derivatives, so the two compositions are the inclusions
\begin{align}
 \mathsf R_{Y,\mathbf b}^{\rm raw}
       \mathsf J_{Y,\mathbf b}^{\rm raw}
  &=\iota^{\cY}_{s+\ell_Y,s},                         \label{eq:A-JY-left-shifted-identity}\\
 \mathsf J_{Y,\mathbf b}^{\rm raw}
       \mathsf R_{Y,\mathbf b}^{\rm raw}
  &=\iota^{\bY_{\rm raw}}_{s+\ell_Y,s}
       \quad\text{on the Cartesian source image},      \label{eq:A-JY-right-shifted-identity}
\end{align}
where the right-hand sides include the same source from index
$s+\ell_Y$ into index $s$. In particular,
$\iota_{Y,\mathbf b}\mathsf J_{Y,\mathbf b}^{\rm raw}$
is injective on the compatible scale.

\smallskip\noindent\textit{The augmented operator.}
The domain coordinates record the reference conormal mismatch,
whereas the interface equations use the conormal derivative at
$\mathbf b$. We therefore apply
$(\mathcal T_{\mathbf b,\rho}^{CN})^{-1}$ from
\eqref{eq:A-current-reference-row-inverse} to the interface equations.
Here, $\mathcal T_{\mathbf b,\rho}^{CN}$ denotes the collection of
the maps $\mathcal T_{\Sigma,\mathbf b,\rho}^{CN}$ over all interfaces.
This conversion acts after
$\mathfrak R_{{\rm sh},J,\rho}^{X}$ reconstructs the full field,
so its variable coefficients may couple different shells.
It and its inverse are analytic and bounded without loss on the
balanced trace product, and equal the identity at
$\mathbf b=\mathbf b_\rho$.

We let $\mathbb N_{\mathbf b}^{\rm cur}:\bX_\natg^s\to\bY_\natg^s$ be the
operator obtained by applying the interior equations at the state
$\mathbf b$ separately on the cap and on the annulus and recording the
interface pair $(d,n_{\mathbf b})$ at the state. We let
$\mathscr R_{CN,\mathbf b,\rho}$ act as the identity on the other data
and as $(\mathcal T_{\mathbf b,\rho}^{CN})^{-1}$ on this pair, and we set
\begin{equation}\label{eq:A-current-row-normalization}
 \mathbb N_{\mathbf b}:=
       \mathscr R_{CN,\mathbf b,\rho}\mathbb N_{\mathbf b}^{\rm cur},
\end{equation}
so that its auxiliary pair is the reference pair
\begin{equation}\label{eq:A-normalized-current-auxiliary-row}
 (\mathcal T_{\mathbf b,\rho}^{CN})^{-1}
       \binom d{n_{\mathbf b}}=\binom d{h_\rho}.
\end{equation}
The other auxiliary entries record the mismatch between the cap
and annular Cauchy data for the minus block, together with the
outer boundary, axis, mean, scale, gauge, and chart conditions. The separate block inverses are first
constructed with the conormal data at $\mathbf b$, and the conversion is
then included in their perturbation estimates. In particular,
an inverse $G_{\mathbf b}^{\rm cur}$ of
$\mathbb N_{\mathbf b}^{\rm cur}$ gives the inverse
\begin{equation}\label{eq:A-current-row-normalized-inverse}
 G_{\mathbf b}^{\mathbb N}
 =G_{\mathbf b}^{\rm cur}\mathscr R_{CN,\mathbf b,\rho}^{-1}
\end{equation}
of $\mathbb N_{\mathbf b}$, with the converse obtained by
composing with $\mathscr R_{CN,\mathbf b,\rho}$.

There is one further adjustment. The extension of nonzero
constraint data contributes to the interior equations, and we
subtract this contribution so that the interior output depends
only on the constrained field. Let
$\mathbb P_{{\rm bulk},\mathbf b}U$ be the output
$\mathbb N_{\mathbf b}U$ with its auxiliary entries set to zero, namely
\[
 \mathbb P_{{\rm bulk},\mathbf b}
   =\iota_{Y,\mathbf b}p_{Y,\mathbf b}\mathbb N_{\mathbf b}.
\]
We also let
$\jmath_{Z,\mathbf b}:\mathscr Z_{{\rm aux},\natg}^s\to\bY_\natg^s$
insert auxiliary data with zero interior source. Thus,
\begin{equation}\label{eq:A-auxiliary-inclusion-identities}
 \operatorname {aux}_{Y,\mathbf b}\jmath_{Z,\mathbf b}=I,
 \qquad p_{Y,\mathbf b}\jmath_{Z,\mathbf b}=0.
\end{equation}
Without using any inverse of $L_{\mathbf b}$, we define
\begin{equation}\label{eq:A-coretraction-bulk-response}
 \mathbb B_{\mathbf b}:=
       \mathbb P_{{\rm bulk},\mathbf b}\mathsf E_{X,\mathbf b}:
 \mathscr Z_{{\rm aux},\natg}^s\longrightarrow\ker
       \operatorname {aux}_{Y,\mathbf b},
\end{equation}
so that $\mathbb B_{\mathbf b}z$ is the interior source produced
by the extension. Applying $\mathbb N_{\mathbf b}$ to the decomposition
\eqref{eq:A-native-domain-decomposition} gives
\begin{equation}\label{eq:A-natural-broken-formula}
 \mathbb N_{\mathbf b}U
 =\iota_{Y,\mathbf b}\mathsf J_{Y,\mathbf b}^{\rm raw}
       L_{\mathbf b}\mathsf R_{X,\mathbf b}U
  +\mathbb B_{\mathbf b}\mathsf C_{X,\mathbf b}U
  +\jmath_{Z,\mathbf b}\mathsf C_{X,\mathbf b}U.
\end{equation}
We remove the second term by a change of range variables that
leaves the auxiliary data unchanged. For a range element
$F+\jmath_{Z,\mathbf b}z$, with
$F\in\ker\operatorname{aux}_{Y,\mathbf b}$, we set
\begin{align}
 \mathscr S_{Y,\mathbf b}(F+\jmath_{Z,\mathbf b}z)
   &=F-\mathbb B_{\mathbf b}z+\jmath_{Z,\mathbf b}z,\notag\\
 \mathscr S_{Y,\mathbf b}^{-1}(F+\jmath_{Z,\mathbf b}z)
   &=F+\mathbb B_{\mathbf b}z+\jmath_{Z,\mathbf b}z.   \label{eq:A-stabilizing-range-shear}
\end{align}
Both maps are analytic and bounded without loss of derivatives.
We can now define the augmented block operator by
\begin{align}
 \mathbb D_{\mathbf b}
 &:=\mathscr S_{Y,\mathbf b}\mathbb N_{\mathbf b}:
       \bX_\natg^s\longrightarrow\bY_\natg^s,
                                                        \label{eq:A-stabilized-native-operator}\\
 \mathbb D_{\mathbf b}U
 &=\iota_{Y,\mathbf b}\mathsf J_{Y,\mathbf b}^{\rm raw}
       L_{\mathbf b}\mathsf R_{X,\mathbf b}U
       +\jmath_{Z,\mathbf b}\mathsf C_{X,\mathbf b}U.
                                                        \label{eq:A-stabilized-native-formula}
\end{align}
The auxiliary output records the constraints, while on the
constrained domain the interior output agrees with the Cartesian
equations. For $U\in\bX_\natg^\infty$ and
$u\in\cX_\gamma^\infty$, the preceding formula gives
\begin{align}
 \operatorname {aux}_{Y,\mathbf b}\mathbb D_{\mathbf b}U
   &=\mathsf C_{X,\mathbf b}U,                         \label{eq:A-auxiliary-row-identity}\\
 \mathbb D_{\mathbf b}\mathsf J_{X,\mathbf b}u
   &=\iota_{Y,\mathbf b}\mathsf J_{Y,\mathbf b}^{\rm raw}
       L_{\mathbf b}u,\qquad u\in\cX_\gamma^\infty.  \label{eq:A-stabilized-intertwining}
\end{align}
By boundedness, these identities extend to the closures of the
smooth elements in the indicated norms. If
$\mathbb D_{\mathbf b}U=\iota_{Y,\mathbf b}
\mathsf J_{Y,\mathbf b}^{\rm raw}f$, the first identity gives
$\mathsf C_{X,\mathbf b}U=0$, so
$U=\mathsf J_{X,\mathbf b}u$ with
$u=\mathsf R_{X,\mathbf b}U$. The second identity and
injectivity of the source transfer then give $L_{\mathbf b}u=f$.
Thus, solving the augmented problem with zero auxiliary data
recovers a solution of the Cartesian problem.

The domains, targets, and roles of the transfer maps are listed in
\cref{tab:transfer-maps}.
\begin{table}[!htbp]
\centering
\small
\setlength{\tabcolsep}{3pt}
\begin{tabular}{@{}lcc>{\raggedright\arraybackslash}p{0.30\textwidth}@{}}
\toprule
Map & Domain & Target & Role\\
\midrule
$\mathsf J_{X,\mathbf b}$ & $\cX_\gamma^s$ & $\bX_\natg^s$ & Expresses a Cartesian unknown in block variables.\\
$\mathsf R_{X,\mathbf b}$ & $\bX_\natg^s$ & $\cX_\gamma^s$ & Returns its constrained Cartesian component.\\
$\mathsf C_{X,\mathbf b}$ & $\bX_\natg^s$ & $\mathscr Z_{{\rm aux},\natg}^s$ & Records the constraint data.\\
$\mathsf E_{X,\mathbf b}$ & $\mathscr Z_{{\rm aux},\natg}^s$ & $\bX_\natg^s$ & Extends the constraint data.\\
$\mathsf J_{Y,\mathbf b}^{\rm raw}$ & $\cY_\gamma^{s+\ell_Y}$ & $\bY_{{\rm raw},\natg}^s$ & Converts the Cartesian source.\\
$\mathsf R_{Y,\mathbf b}^{\rm raw}$ & $\bY_{{\rm raw},\natg}^s$ & $\cY_\gamma^s$ & Reconstructs the Cartesian source.\\
$\iota_{Y,\mathbf b}$ & $\bY_{{\rm raw},\natg}^s$ & $\bY_\natg^s$ & Adds zero auxiliary data.\\
$p_{Y,\mathbf b}$ & $\bY_\natg^s$ & $\bY_{{\rm raw},\natg}^s$ & Extracts the interior source coordinates.\\
$\operatorname{aux}_{Y,\mathbf b}$ & $\bY_\natg^s$ & $\mathscr Z_{{\rm aux},\natg}^s$ & Extracts the auxiliary data.\\
$\jmath_{Z,\mathbf b}$ & $\mathscr Z_{{\rm aux},\natg}^s$ & $\bY_\natg^s$ & Inserts auxiliary data with zero interior source.\\
$\mathscr R_{CN,\mathbf b,\rho}$ & $\bY_\natg^s$ & $\bY_\natg^s$ & Converts Cauchy data to reference coordinates.\\
$\mathscr S_{Y,\mathbf b}$ & $\bY_\natg^s$ & $\bY_\natg^s$ & Subtracts the extension's interior contribution.\\
\bottomrule
\end{tabular}
\caption{The transfer maps. The last two maps are invertible on the
indicated spaces.}
\label{tab:transfer-maps}
\label{eq:A-adapter-type-ledger}
\end{table}

\smallskip\noindent\textit{The separate high-frequency block operators.}
We construct inverses for the minus and positive blocks separately.
For $\alpha\in\{-,+\}$, let $\iota_\alpha^{X,Y}$ include block
$\alpha$ in the full domain or range, and let $\pi_\alpha^{X,Y}$
extract that block. We first take the diagonal blocks of the two
factors separately,
$N_{\alpha\alpha}:=\pi_\alpha^Y\mathbb N_{\mathbf b}\iota_\alpha^X$
and
$\mathscr S_{\alpha\alpha}:=\pi_\alpha^Y
\mathscr S_{Y,\mathbf b}\iota_\alpha^Y$.
If $G_{\alpha\alpha}^{\mathbb N}$ is the inverse of
$N_{\alpha\alpha}$, we obtain mutually inverse operators on the
spaces of block $\alpha$ by setting
\begin{equation}\label{eq:A-sheared-component-inverse}
 \widehat D_{\alpha,\mathbf b}
 :=\mathscr S_{\alpha\alpha}N_{\alpha\alpha},\qquad
 \widehat G_{\alpha,\mathbf b}
 :=G_{\alpha\alpha}^{\mathbb N}\mathscr S_{\alpha\alpha}^{-1}.
\end{equation}
The diagonal block of the full product also contains a term
passing through the complementary range,
\[
 \pi_\alpha^Y\mathbb D_{\mathbf b}\iota_\alpha^X
 =\widehat D_{\alpha,\mathbf b}
 +\pi_\alpha^Y\mathscr S_{Y,\mathbf b}
   (I-\iota_\alpha^Y\pi_\alpha^Y)
   \mathbb N_{\mathbf b}\iota_\alpha^X.
\]
The additional term passes through the other high-frequency block
or the low-frequency part. Thus, the separately constructed
$\widehat D_{\alpha,\mathbf b}$ need not equal the diagonal block
of $\mathbb D_{\mathbf b}$. We keep the auxiliary equations in
this comparison and estimate the additional term in
\cref{app:macro-cross}. On the constrained Cartesian domain, we have
\begin{equation}\label{eq:A-projected-intertwining}
 p_{Y,\mathbf b}\mathbb D_{\mathbf b}
       \mathsf J_{X,\mathbf b}
 =\mathsf J_{Y,\mathbf b}^{\rm raw}L_{\mathbf b}.
\end{equation}

The mapping properties of the transfer maps and the estimates for their
derivatives with respect to $\mathbf b$ are summarized in the following proposition.
\begin{proposition}[Transfer maps between Cartesian and block variables]
\label{prop:all-order-adapters}
On a sufficiently small ball of states in the analytic norm, the
source loss $\ell_Y$ in \eqref{eq:A-source-order} is finite and
independent of the Sobolev index. The transfer maps satisfy
\begin{align}
 \mathsf J_{X,\mathbf b},\mathsf R_{X,\mathbf b},
 \mathsf C_{X,\mathbf b},\mathsf E_{X,\mathbf b},
 \mathscr R_{CN,\mathbf b,\rho}^{\pm1},
 \mathscr S_{Y,\mathbf b}^{\pm1}
   &\ \text{are bounded between the spaces listed in \cref{tab:transfer-maps}},
                                                        \label{eq:A-JX-type}\\
 \mathsf J_{Y,\mathbf b}^{\rm raw}&:
       \cY_\gamma^{s+\ell_Y}\longrightarrow
       \bY_{{\rm raw},\natg}^s,\qquad
 \mathsf R_{Y,\mathbf b}^{\rm raw}:
       \bY_{{\rm raw},\natg}^s\longrightarrow\cY_\gamma^s.
                                                               \label{eq:A-JY-type}
\end{align}
The maps in \eqref{eq:A-zero-range-splitting} lose no derivatives.
All the maps above preserve real data, with the conjugation rule of
\cref{subsec:range-projection} in quotient coordinates. They depend
analytically on $\mathbf b$ with the domain and target regularities
listed in \cref{tab:transfer-maps} and \eqref{eq:A-JY-type}. Every derivative
of fixed order with respect to $\mathbf b$ satisfies a tame estimate
with at most one factor measured in a higher regularity norm and
with regularity shifts independent of $s$. Once the inverse
$G_{\mathbf b}$ of $\mathbb D_{\mathbf b}$ has been constructed,
the corresponding Cartesian inverse is
\begin{equation}\label{eq:A-adapter-factorization}
 \mathsf R_{X,\mathbf b}G_{\mathbf b}
       \iota_{Y,\mathbf b}\mathsf J_{Y,\mathbf b}^{\rm raw}.
\end{equation}
The source conversion contributes the loss $\ell_Y$ once,
including after differentiation with respect to $\mathbf b$,
in addition to any loss in the estimates for $G_{\mathbf b}$.
\end{proposition}

\begin{proof}
\textit{Inverse and constraint identities.}
The extension identity \eqref{eq:A-Cauchy-coretraction-identity},
together with $J_AE_A=I$ from \eqref{eq:A-axis-coretractions} and
$B_MC_{\pa,M}=I$ from \eqref{eq:A-physical-row}, shows that
subtracting the extensions imposes the constraints and adding
them back reconstructs the field. On the shells, this gives
\eqref{eq:A-break-left}--\eqref{eq:A-cap-trace-inverses}.
Composing with the chart and frame maps, the spectral projections,
and the identifications with the reference spaces gives \eqref{eq:A-domain-retract-identities} and
\eqref{eq:A-native-domain-decomposition}. The reference extensions
commute with $Q_j$, and each cap is joined to its annulus before
summation. Thus, the coefficients at the state act on the reconstructed
field in the original coordinates.

Discarding the derived entry $\Psi_{\mathbf b}Z$ and reconstructing
the Cartesian source from its quotient coordinates gives
\eqref{eq:A-JY-left-shifted-identity}--\eqref{eq:A-JY-right-shifted-identity}.
Adding zero auxiliary entries and then extracting the interior source gives
\eqref{eq:A-zero-range-splitting}. Substitution in
\eqref{eq:linearized-H0}--\eqref{eq:linearized-H3} yields the
interior part of \eqref{eq:A-natural-broken-formula}, while
\eqref{eq:A-current-reference-row-inverse} converts the interface
data to the reference mismatch without changing the matching
condition. Applying $\mathscr S_{Y,\mathbf b}$ cancels the
interior contribution of the extension and gives
\eqref{eq:A-stabilized-native-formula}.
The identities \eqref{eq:A-auxiliary-inclusion-identities} and
\eqref{eq:A-domain-retract-identities} then imply
\eqref{eq:A-auxiliary-row-identity} and
\eqref{eq:A-stabilized-intertwining}.

\smallskip\noindent\textit{Bounds for the transfer maps.}
We use the weights of \cref{tab:block-weights}. The chart,
frame, and quotient maps are bounded at the same regularity
index on their respective spaces. The spectral projections
are controlled by the uniform resolvent bounds on the fixed
contours, and \eqref{eq:A-complete-Kato-factor} is invertible
with bounded inverse because the state projections are close
to their reference projections.

The trace and extension bounds follow from
\eqref{eq:A-axis-coretraction-estimates},
\eqref{eq:A-boundary-coretraction-estimate}, and
\eqref{eq:A-Cauchy-coretraction-estimate}. Orthogonality of the
shells gives the bounds for their infinite sums. For the
conormal conversion, the one-derivative difference between
the value and conormal norms absorbs
$K_{\Sigma,\mathbf b,\rho}$, so
\eqref{eq:A-current-reference-row-map}--\eqref{eq:A-current-reference-row-inverse}
are bounded without loss, as are their state derivatives.
The map \eqref{eq:A-coretraction-bulk-response}, which gives the
interior source produced by an extension, is bounded on the same
scales. Hence, subtracting this source as in
\eqref{eq:A-stabilizing-range-shear} is also bounded. Only
$\Theta_{\rm src,\mathbf b}$ differentiates the prescribed
source, with the loss $\ell_Y$ of \eqref{eq:A-source-order}.

For derivatives with respect to $\mathbf b$, differentiating
$G_{\mathbf b}G_{\mathbf b}^{-1}=I$, with $G_{\mathbf b}$ as in
\eqref{eq:A-frame-factor}, gives
\[
 D_{\mathbf b}(G_{\mathbf b}^{-1})[h]
 =-G_{\mathbf b}^{-1}
    (D_{\mathbf b}G_{\mathbf b}[h])G_{\mathbf b}^{-1}.
\]
We combine this identity with the product estimate
\eqref{eq:analytic-product}, differentiation of the contour formulas
\eqref{app:cross:domain-range-contours}, and the Neumann series for
the inverse of \eqref{eq:A-complete-Kato-factor}.
Differentiation changes the coefficients of the source map
without increasing its differential order. Since each term
contains exactly one source conversion factor, the loss
$\ell_Y$ still occurs only once. Finally, each factor respects
complex conjugation in its coordinates, so the transfer maps
preserve real data.
\end{proof}

We close this section with the decomposition of $\mathbb D_{\mathbf b}$
into low and high cell frequencies, which organizes the rest of the linear
analysis. In the notation of \eqref{eq:A-exact-cell-split}, we set
\[
 P=P_{J,\natg}:=\mathbf1_{\{|n|<2^J\}}(D_\zeta),
 \qquad
 H=H_{J,\natg}:=\mathbf1_{\{|n|\geq2^J\}}(D_\zeta),
\]
with superscripts $X$ and $Y$ distinguishing the domain and
range. Although $P$ selects finitely many cell frequencies,
its range remains infinite-dimensional because the radial
and angular variables are unrestricted. With the low-frequency
variables first, the decomposition is
\begin{equation}\label{app:global:block-decomposition}
 \mathbb D_{\mathbf b}
 =\begin{pmatrix}A_{\mathbf b}&B_{\mathbf b}\\
                  C_{\mathbf b}&D_{H,\mathbf b}
   \end{pmatrix},
 \quad
 \begin{aligned}
 A_{\mathbf b}&=P^Y\mathbb D_{\mathbf b}P^X,&
 B_{\mathbf b}&=P^Y\mathbb D_{\mathbf b}H^X,\\
 C_{\mathbf b}&=H^Y\mathbb D_{\mathbf b}P^X,&
 D_{H,\mathbf b}&=H^Y\mathbb D_{\mathbf b}H^X.
 \end{aligned}
\end{equation}

\section{The inverse at a fixed ellipse}
\label{app:reference}

We now invert the Cartesian linearization
\eqref{eq:current-linearization} at a constant reference ellipse,
with a finite loss of derivatives of the source. Since its coefficients
are independent of $\zeta$, we can solve separately at each cell
frequency. We also show that the estimates for the even and affine
components remain uniform as the ellipse approaches a circle.

\subsection{The reference ellipse and the inverse estimate}

We fix a constant orientation $\alpha_0$ and set
\begin{equation}\label{appref:ellipse}
 M_e=\mathsf R_{\alpha_0}
 \begin{pmatrix}\sqrt{1+\rho_*}&0\\0&\sqrt{1-\rho_*}\end{pmatrix}
 \mathsf R_{-\alpha_0},
 \qquad 0<\rho_*<\rho_0<\frac14,
\end{equation}
which is the matrix \eqref{eq:fixed-M-lambda} with $\delta=0$ and
$\rho=\rho_*$. Here, $\tr(M_e^TM_e)=2$, as in
\eqref{eq:normalized-affine-trace}, so $M_e$ belongs to the family
\eqref{eq:A-seed-family}, and the corresponding unperturbed solution maps
every cross-section of the cell to the same ellipse with parameter
$\rho_*$. We let
\[
 L_e:\cX_{M_e}^\infty\longrightarrow\cY_{\rm comp}^\infty
\]
denote the linearization \eqref{eq:current-linearization} in the
normalized chart. Its domain satisfies the axis, mean, scale, gauge,
and outer boundary conditions of \cref{subsec:normalized-slice},
and its range is the compatible space $\cY_{\rm comp}^\infty$ of
\cref{subsec:compatible-range}. We impose these conditions before
separating the angular modes $e^{im\theta}$ by the \emph{parity}
of $m$, so that the domain and range have matching odd and even
components. In the estimate below, we use the analytic scale of
\cref{subsec:analytic-scale} and suppress the weight subscript.

\begin{proposition}[The reference inverse at a constant ellipse]
\label{appref:raw-constant-ellipse}
If $\rho_0$ is sufficiently small, then for every fixed
$0<\rho_*<\rho_0$, the map $L_e$ is an invertible linear map. Its inverse
$V_e$ satisfies, for every $s\geq s_*$,
\begin{equation}\label{appref:raw-estimate}
 \norm{V_ef}_{\cX^s}
 \leq C_{e,s}\norm{f}_{\cY^{s+\ell_e}}
\end{equation}
for one finite $\ell_e$, which is independent of the cell frequency and
of every angular or cell Fourier cutoff. A constant change of orientation
conjugates this inverse by a unitary map.
\end{proposition}

The proof occupies the next seven subsections. In scalar coordinates,
the problem splits into odd and even angular modes. For the odd
block, we solve on the cap and annulus, match at the interface, and
establish the analytic estimates. The even block also contains the
affine modes $m=\pm2$, which require a separate calculation.
We then combine the two inverses and return to Cartesian variables.

\subsection{Scalar variables and even and odd angular modes}

A constant rotation is unitary on every Cartesian scale, so we may take
$\alpha_0=0$, and we write $\rho=\rho_*$. For the tangent vector $U$ and
the scalar $S$, we introduce the scalar coordinates
\[
 A=y\cdot U_\perp,\qquad B=\cR y\cdot U_\perp,
 \qquad C=U\cdot e_T.
\]
Here, $S$ is the variation of $w$, $U_\perp=M_e^T\iota^TU$,
and $\cR y=(-y_2,y_1)$. Thus, $A=(\cD v_{M_e})\cdot U$
and $B=(\cR v_{M_e})\cdot U$ are the combinations appearing in
\eqref{eq:linearized-H0}--\eqref{eq:linearized-H3}, while $C$ is
the tangential component of $U$. We retain the prescribed angular averages and the scale coordinate
associated with \eqref{eq:normalized-affine-trace}. We set
$Z=A+iB$, $W=A-iB$, and $z=y_1+iy_2$. With $C^\infty$ denoting smooth
functions of $y$, the divisibility conditions
\[
 Z\in\bar zC^\infty,\qquad W\in zC^\infty,
 \qquad f_1+if_0\in\bar zC^\infty,\qquad
 f_1-if_0\in zC^\infty,\qquad f_3\in z\bar zC^\infty
\]
are part of the definitions of the domain and of the compatible range in
\cref{sec:tame-problem}, so we may use them before taking any quotient.
The maps to and from the scalar coordinates are bounded, with a finite
shift in the Sobolev index, which we call $\ell_{\rm alg}$ in the proof
of \cref{appref:raw-constant-ellipse}.

At a nonzero cell frequency $n$, we set
\[
 t=\frac{|n|r}{L},\qquad D=t\pa_t,\qquad
 \cR=\pa_\theta,\qquad \sigma=\operatorname{sgn}n.
\]
We also rescale the tangential component, writing
\(C=(L/|n|)U\cdot e_T\) for the rest of the nonzero-frequency
calculation. We keep the first two source entries unchanged, divide the
third by \(|n|\), and divide the fourth by \(-L\det M_e\). We continue to
write these normalized sources as \(f_0,\ldots,f_3\). With these
conventions, the coefficients in \eqref{appref:raw-scalar-rows} depend
on \(n\) only through \(\sigma\). Thus,
$t=1$ corresponds to $r=L/|n|$, which is the radius of the cap at cell
frequency $n$, and the outer boundary $r=1$ corresponds to $t=|n|/L$.
Below, we take Cartesian derivatives in the rescaled coordinates. If
\[
 \binom{\mathfrak a}{\mathfrak b}
 =H_\rho(\theta)\binom AB,\qquad
 H_\rho(\theta)=\frac1{1-\rho^2}
 \begin{pmatrix}
  1-\rho\cos2\theta&\rho\sin2\theta\\
  \rho\sin2\theta&1+\rho\cos2\theta
 \end{pmatrix},
\]
then the four equations for $(A,B,C,S)$ are
\begin{equation}\label{appref:raw-scalar-rows}
 \begin{aligned}
 f_0&=\cR S-2A-2\cR B,\\
 f_1&=(I-\Pi)(DS-DB-\cR A+2B),\\
 f_2&=(I-\Pi)(\cR C+i\sigma B-i\sigma S),\\
 f_3&=(I-\Pi)(D\mathfrak a+\cR\mathfrak b+i\sigma t^2C).
 \end{aligned}
\end{equation}
Here, $\Pi$ takes the angular average, so $I-\Pi$ removes the mode
$m=0$. The ellipse enters only through $H_\rho$, which is the identity
at $\rho=0$. Thus, with
\begin{equation}\label{appref:spin-variables}
 X=S-B,\qquad g_\pm=f_1\pm if_0,
\end{equation}
the first three equations, in which the ellipse does not occur, become on the
nonzero angular modes
\begin{equation}\label{appref:spin-equations}
 \begin{aligned}
 X&=i\sigma(f_2-\cR C),\\
 (\cR+2i)Z&=(D+i\cR)X-g_+,\\
 (\cR-2i)W&=(D-i\cR)X-g_-.
 \end{aligned}
\end{equation}
The nonconstant entries of $H_\rho$ are combinations of
$e^{\pm2i\theta}$, so they shift the angular mode by $\pm2$
and preserve its parity. The divisibility, axis, and boundary
conditions also preserve parity. Hence, the closed scalar
operator $\mathscr L_e$, including the mean and scale equations,
splits as
\begin{equation}\label{appref:parity-split}
 \mathscr L_e
 =\mathscr L_e^{\rm od}\oplus
   \mathscr L_e^{\rm ev,aug}.
\end{equation}
We include the angular means and the scale equation in the even summand.
They couple to the other even components only through the modes $m=\pm2$,
and the resulting system is triangular.

\subsection{The odd block on the cap}

On odd angular modes, $\cR\pm2i$ are invertible because their
multipliers are $i(m\pm2)$. We solve
\eqref{appref:spin-equations} for $Z$ and $W$ and substitute
them into the last equation of \eqref{appref:raw-scalar-rows}.
Then, $u=|\cR|C$, where $|\cR|$ has multiplier $|m|$, satisfies
\begin{equation}\label{appref:odd-cartesian-operator}
 \mathscr A_\rho u=\phi(f),\qquad
 \mathscr A_\rho
 =P\Delta_y-2\rho\mathsf H
   (\pa_z^2+\pa_{\bar z}^2)\mathsf H-(1-\rho^2)I.
\end{equation}
Here, $\phi(f)$ collects the terms depending on
$f=(f_0,f_1,f_2,f_3)$ after substituting
\[
 Z=(\cR+2i)^{-1}\bigl((D+i\cR)X-g_+\bigr),\qquad
 W=(\cR-2i)^{-1}\bigl((D-i\cR)X-g_-\bigr),
\]
with $X=i\sigma(f_2-\cR C)$ from \eqref{appref:spin-equations}.
The angular operators in \eqref{appref:odd-cartesian-operator} are
\[
 P=p(-i\pa_\theta),\qquad p_m=\frac{m^2}{m^2-4},
 \qquad \mathsf H=\operatorname{sgn}(-i\pa_\theta).
\]
The nonconstant coefficients of $H_\rho$ couple the equation for mode
$m$ to modes $m\pm2$. Here and below, $u_m$ denotes the
coefficient of $e^{im\theta}$ in the angular Fourier series of $u$, and
likewise $f_{i,m}$ denotes that coefficient of the source entry $f_i$. We
obtain \eqref{appref:odd-cartesian-operator} from the identities
\[
 \begin{aligned}
 t^{-2}\bigl(D^2+(2-2m)D+m(m-2)\bigr)u_{m-2}e^{im\theta}
  &=4\pa_{\bar z}^2(u_{m-2}e^{i(m-2)\theta}),\\
 t^{-2}\bigl(D^2+(2m+2)D+m(m+2)\bigr)u_{m+2}e^{im\theta}
  &=4\pa_z^2(u_{m+2}e^{i(m+2)\theta}).
 \end{aligned}
\]
The perturbation from the circle therefore involves the Cartesian
operators $\pa_z^2$ and $\pa_{\bar z}^2$. The operators at the circle and ellipse act on the same smooth
Cartesian functions with the same axis and boundary constraints.
% AUTHOR QUERY: Specify the completed cap domain and source space for
% this perturbation argument. Prescribing a first Cartesian coefficient
% on the axis is not a bounded condition on a general H^2 function in
% two dimensions. In particular, square-integrable center forcing alone
% does not bound the integrals in odd-center-volterra near t=0.
The divisibility of the compatible sources controls the quotients
by $z$, $\bar z$, and $z\bar z$ in this substitution. Since
at most three derivatives fall on $f$, we obtain, independently
of all cutoffs,
\begin{equation}\label{appref:odd-source-comparison}
 \norm{\phi(f)}_{H^s_y}\leq C_s
 \norm{f}_{\cY^{s+3}}.
\end{equation}
Here, $H^s_y$ denotes the Sobolev space of index $s$ in the rescaled
disk variable, and the norm on the right is that of the analytic scale
of \cref{subsec:analytic-scale}, restricted to the cell frequency $n$ as
in \eqref{appref:carrier-mode-norms} below.

On the cap, which is the region $0<t<1$ in the rescaled variable, we let
\[
 \mathscr L_{\rho,\rm ax}u
  =(\mathscr A_\rho u,\mathfrak j_cu,\gamma_hu),
 \qquad
 \gamma_hu=\Pi_{|m|\geq3}u|_{t=1}.
\]
Here, $\mathfrak j_cu$ records the coefficients $\alpha_{\pm1}$ in
$u_{\pm1}(t)=\alpha_{\pm1}t+O(t^3)$, while $\gamma_hu$ records
the values of the higher odd modes on the cap boundary. The modes
$m=\pm1$ are the \emph{center modes}, for which
$p_{\pm1}=-1/3$, while $|m|\ge3$ are the \emph{positive modes},
for which $p_m\ge1$. At the circle, the center equation is
$(\Delta_y+3)c=-3\phi_c$, with $\phi_c$ the center component
of $\phi(f)$. On each mode, this is the Bessel equation of order
one and frequency $\sqrt3$. We let
\[
 u_0(t)=\frac2{\sqrt3}J_1(\sqrt3t),\qquad
 v_0(t)=Y_1(\sqrt3t),\qquad
 \mathcal W=u_0v_0'-u_0'v_0
\]
be its regular solution, its singular solution, and their Wronskian, so
that $u_0(t)=t+O(t^3)$. The solution of the center equation with
prescribed leading coefficient $\alpha$ is then given by the variation of
constants formula
\begin{equation}\label{appref:odd-center-volterra}
 c(t)=\alpha u_0(t)
 -u_0(t)\int_0^t\frac{v_0(s)F(s)}{\mathcal W(s)}\dd s
 +v_0(t)\int_0^t\frac{u_0(s)F(s)}{\mathcal W(s)}\dd s.
\end{equation}
Here, $F=-3\phi_c$. On $|m|\ge3$, we solve the Dirichlet
problem for $p_m\Delta_y-I$. The ellipticity is uniform
in $m$ because $1\le p_m\le9/5$. We call its solution operator
the \emph{Dirichlet Green map}. We denote their direct sum by $G_{0,\rm ax}$. Thus,
$G_{0,\rm ax}(\phi,\alpha,a)$ solves the cap equation with interior
source $\phi$, center coefficients $\alpha=(\alpha_1,\alpha_{-1})$,
and higher-mode boundary value $a$. For
$\rho\ne0$, the cap inverse is
\begin{equation}\label{appref:odd-core-green}
 G_{\rho,\rm ax}
 =\bigl[I+G_{0,\rm ax}
   (\mathscr L_{\rho,\rm ax}-\mathscr L_{0,\rm ax})\bigr]^{-1}
   G_{0,\rm ax}.
\end{equation}
Indeed, we have the estimate
\[
 \norm{(\mathscr L_{\rho,\rm ax}-\mathscr L_{0,\rm ax})u}_{H^s}
 \leq C_{\rm dif}(|\rho|+\rho^2)\norm{u}_{H^{s+2}},
\]
uniformly under restriction to $|m|\le m_0$, which we call an
\emph{angular cutoff}. Decreasing $\rho_0$ makes the bracket
in \eqref{appref:odd-core-green} invertible by a Neumann series
at the base index. Cartesian difference quotients give the
higher regularity estimates with the same $\rho_0$. The
prescribed axis conditions remain part of the domain, and the
Cartesian compatibility and divisibility conditions remain
part of the source space.

\subsection{The odd block on the annulus}

On the annulus $t\geq1$, we write \eqref{appref:odd-cartesian-operator} as
an evolution in the radial variable, using $D=t\pa_t$ as the derivative,
\begin{equation}\label{appref:odd-mass-equation}
 M_\rho D^2u+K_\rho Du+N_\rho u-(1-\rho^2)t^2u=F,
\end{equation}
on the scale of sequences
\[
 \mathfrak h^s=\left\{(u_m)_{m\in2\Z+1}:
   \sum_m\langle m\rangle^{2s}|u_m|^2<\infty\right\}.
\]
Here, $u$ is the sequence of odd angular modes and
$F=t^2\phi(f)$. The matrices $M_\rho$ and $N_\rho$ are
self-adjoint, and $K_\rho$ is skew-adjoint. Each row couples
only the adjacent odd modes $m$ and $m\pm2$. At the circle,
\[
 (M_0)_{mm}=p_m,\qquad K_0=0,\qquad
 (N_0)_{mm}=-p_m m^2.
\]
We call $N_\rho u$ the \emph{angular term}, since at the circle it is the
multiplier $-p_mm^2$, and we call $(1-\rho^2)t^2u$ the \emph{mass term}.
Expanding the two off-diagonal entries of $H_\rho$, we obtain, for each
fixed $s$,
\[
 \begin{aligned}
 \norm{M_\rho-M_0}_{\mathfrak h^s\to\mathfrak h^s}
  &\leq C_s|\rho|,\\
 \norm{K_\rho}_{\mathfrak h^{s+1}\to\mathfrak h^s}
  &\leq C_s|\rho|,\\
 \norm{N_\rho-N_0}_{\mathfrak h^{s+2}\to\mathfrak h^s}
  &\leq C_s|\rho|.
 \end{aligned}
\]
Since $p_{\pm1}=-1/3$ and $p_m\geq1$ for the other odd $m$, the first
estimate and the gap in the spectrum of $M_0$ show that $M_\rho$ has
exactly two negative directions for this choice of $\rho_0$. Commuting
the contour integrals that define the spectral projections with the
weights $\langle m\rangle^s$, we obtain bounded analytic projections $P_\pm(\rho)$
onto the positive and the negative spectral subspaces of every
$\mathfrak h^s$. We choose an analytic family of orthogonal maps
$Q_\rho$ which transports the projections
at $\rho$ to those at the circle, so that
\[
 Q_\rho P_\pm(\rho)Q_\rho^*=P_\pm(0).
\]
In the \emph{fixed coordinates} given by $Q_\rho$, the spectral
subspaces and their interface trace spaces are independent of
$\rho$. We set
\[
 T_\rho=|M_\rho|^{1/2},\qquad
 J_\rho=\operatorname{sgn}M_\rho=P_+-P_-,
 \qquad \widehat u=T_\rho u.
\]
Multiplying \eqref{appref:odd-mass-equation} by $|M_\rho|^{-1/2}$, we
obtain the normalized equation
\begin{equation}\label{appref:odd-sign-normal-form}
 J_\rho D^2\widehat u+B_\rho D\widehat u+C_\rho^\sharp \widehat u
 -(1-\rho^2)t^2A_\rho^\sharp \widehat u=\widetilde f.
\end{equation}
Here,
\[
 \begin{aligned}
 B_\rho&=|M_\rho|^{-1/2}K_\rho|M_\rho|^{-1/2},&
 C_\rho^\sharp&=|M_\rho|^{-1/2}N_\rho|M_\rho|^{-1/2},\\
 A_\rho^\sharp&=|M_\rho|^{-1},&
 \widetilde f&=|M_\rho|^{-1/2}F.
 \end{aligned}
\]
The two-dimensional negative spectral subspace gives the
center component. In \eqref{appref:odd-sign-normal-form}, the coupling
between the positive and the negative subspaces contains no second-order
term. The symmetry under the reflection $\mathscr K$, which sends the mode
$m$ to $-m$ and was introduced before \eqref{eq:A-K1}, and the
skew-adjointness of $B_\rho$ also give
\begin{equation}\label{appref:odd-mass-identities}
 P_-J_\rho P_+=0,\qquad P_-B_\rho P_-=0.
\end{equation}

\smallskip\noindent\textit{The positive component.}
On the positive subspace, coercivity gives a boundary value
problem on $1\le t\le T$. We prescribe the value at $t=1$
and the outer boundary datum at $t=T$, where $T=|n|/L$ at
cell frequency $n$. We let $G_{+,T}$ be the solution operator
for this problem, including its angular and mass terms.
We denote its \emph{positive quadratic form} by
$\mathfrak a_{++,\rho,T}$. With the same normalization as
the interior equations, the outer boundary condition at
the circle is
\begin{equation}\label{appref:odd-circle-physical-row}
 (\mathcal B_{0,T}u_+)_m
 =2\sqrt{p_m}(D+2)u_{+,m}(T),\qquad |m|\geq3.
\end{equation}
For small $\rho$, the boundary operator at the ellipse is an analytic
perturbation of \eqref{appref:odd-circle-physical-row} which couples each
mode only to the adjacent ones. In particular, its conormal coefficient on
the positive modes remains invertible, while its compression to the
center modes is $O(\rho)$. We therefore separate the positive part from
the center contribution and write, in the fixed coordinates,
\begin{equation}\label{appref:odd-physical-row-split}
 \mathcal B_{\rho,T}(u_c+u_+)
 =\mathcal B_{++,\rho,T}u_++\mathcal B_{+-,\rho,T}u_c.
\end{equation}
Here, $u_c=P_-(0)Q_\rho \widehat u$ and $u_+=P_+(0)Q_\rho \widehat u$ are the center and the
positive components of the unknown in the fixed coordinates. The balanced space for the boundary datum is
\begin{equation}\label{appref:odd-balanced-trace}
 \mathfrak T_T^s=\mathfrak h_+^{s+1/2},\qquad
 \norm{b}_{\mathfrak T_T^s}
 :=T^{-1/2}\norm{b}_{\mathfrak h^{s+1/2}}.
\end{equation}
Here, $\mathfrak h_\pm^s=P_\pm(0)\mathfrak h^s$ are the positive and the
negative subspaces of $\mathfrak h^s$ in the fixed coordinates,
and the center contribution satisfies
\begin{equation}\label{appref:odd-boundary-one-way}
 \norm{\mathcal B_{+-,\rho,T}u_c}_{\mathfrak T_T^s}
 \leq C_s|\rho|\norm{u_c}_{\cC_T^s}.
\end{equation}
Here, writing $\xi=t^{1/2}u_c$, we measure the center solution by
\[
 \norm{u_c}_{\cC_T^s}
 =\sup_{1\leq t\leq T}
   \bigl(\norm{\xi(t)}_{\mathfrak h^s}
        +\norm{\xi'(t)}_{\mathfrak h^s}\bigr).
\]
Indeed, the negative spectral vectors are bounded in every
$\mathfrak h^s$, the row couples only adjacent modes, and we have the
bounds
\[
 |u_c(T)|\leq CT^{-1/2}\norm{u_c}_{\cC_T^s},\qquad
 |D u_c(T)|\leq CT^{1/2}\norm{u_c}_{\cC_T^s}.
\]

We write $\Lambda_\theta=\langle-i\pa_\theta\rangle$ for the angular
multiplier, which multiplies mode $m$ by $\langle m\rangle$. We define
the energy norm on the positive subspace by
\begin{equation}\label{appref:odd-positive-energy}
 \norm{u_+}_{\cE_T^s}^2
 =\int_1^T\left(\norm{D u_+}_{\mathfrak h^s}^2
  +\norm{\Lambda_\theta u_+}_{\mathfrak h^s}^2
  +t^2\norm{u_+}_{\mathfrak h^s}^2\right)\frac{\dd t}{t}
  +\norm{\operatorname{tr}_T u_+}_{s+1/2}^2.
\end{equation}
Here, $\operatorname{tr}_T u_+=u_+(T)$ is the trace on the outer boundary. We
let $G_{+,T}^0q$ solve the positive equation with source $q$, zero
inner value, and homogeneous outer boundary condition. We let
$P_{+,T}a$ solve the homogeneous equation with inner value $a$ and
homogeneous outer boundary condition, and choose a bounded extension
$E_Tb$ with
\[
 (E_Tb)(1)=0,\qquad \mathcal B_{++,\rho,T}E_Tb=b.
\]
Thus, the two extensions supply the inner and outer data separately,
and we have
\begin{equation}\label{appref:odd-positive-green-poisson}
 G_{+,T}(q,a,b)
 =P_{+,T}a+E_Tb
  +G_{+,T}^0\bigl(q-\mathcal N_{++}(P_{+,T}a+E_Tb)\bigr).
\end{equation}
Here, $\mathcal N_{++}$ is the compression of
\eqref{appref:odd-sign-normal-form} to the positive subspace in the
fixed coordinates. In particular, $u_+=G_{+,T}(q,a,b)$ solves
\[
 \mathcal N_{++}u_+=q,\qquad u_+(1)=a,\qquad
 \mathcal B_{++,\rho,T}u_+=b.
\]
The form $\mathfrak a_{++,\rho,T}$ satisfies, for some $c_+>0$ and uniformly
for $T\geq1$,
\begin{equation}\label{appref:odd-positive-form}
 \mathfrak a_{++,\rho,T}(u_+,u_+)
 \geq c_+\int_1^T
 \bigl(\norm{D u_+}^2+\norm{\langle\cR\rangle u_+}^2
                  +t^2\norm{u_+}^2\bigr)\frac{\dd t}{t}.
\end{equation}
Here, $\norm{\cdot}$ is the norm of $\mathfrak h^0$, and
$\langle\cR\rangle$ is the multiplier $\Lambda_\theta$.
To verify this, we write
\[
 T_\rho^{-1}Q_\rho^*u_+=h+\Gamma_\rho h,\qquad
 \Pi_{|m|\geq3}h=h,\qquad
 \norm{\Gamma_\rho}_{\mathfrak h^s\to\mathfrak h^s}
 \leq C_s|\rho|.
\]
Here, \(u=T_\rho^{-1}Q_\rho^*u_+\) is the positive component in the
original coordinates. The spectral gap gives a positive lower bound
for $\langle M_\rho Du,Du\rangle$, while we have the estimates
\[
 \begin{aligned}
 -\langle N_0u,u\rangle
 &\geq c_0\norm{\Lambda_\theta h}^2-C\rho^2\norm{h}^2,\\
 |\langle(N_\rho-N_0)u,u\rangle|
 +|\langle K_\rho Du,u\rangle|
 &\leq C|\rho|(\norm{\Lambda_\theta h}^2+\norm{Du}^2).
 \end{aligned}
\]
The positive $t^2$-term absorbs the remaining $L^2$ term, and by density
we obtain \eqref{appref:odd-positive-form} on the complete sequence
space, with no angular cutoff.

\smallskip\noindent\textit{The center component.}
On the negative subspace, the substitution $u_c=t^{-1/2}\xi$ reduces
\eqref{appref:odd-sign-normal-form} to two oscillators, which are
separated by the reflection $\mathscr K$ and take the form
\begin{equation}\label{appref:odd-center-oscillator}
 \xi''+\Omega_\rho^2\xi+t^{-2}V_{--,\rho}\xi=g.
\end{equation}
Here, $\Omega_\rho^2$ and $V_{--,\rho}$ are analytic
$2\times2$ matrices, with $\Omega_\rho^2$ positive. Its
eigenvalues are the squares of the oscillator frequencies
from \cref{ss:strategy}, which coincide at the circle.
We solve outward from $t=1$, prescribing Cauchy data
$a=(\xi(1),\xi'(1))$ and imposing no outer boundary condition.
Writing $z=(\xi,\xi')$, the first-order system associated with
\eqref{appref:odd-center-oscillator} is
\[
 z'=\begin{pmatrix}0&I\\-\Omega_\rho^2-t^{-2}V_{--,\rho}&0\end{pmatrix}z
       +\binom0g.
\]
We let $\Phi_\rho(t,s)$ propagate its homogeneous solutions from $s$
to $t$, with $\Phi_\rho(s,s)=I$, and set $\pi_1(\xi,\xi')=\xi$.
The center solution $u_c$ with source $g$ is
\begin{equation}\label{appref:odd-center-green}
 V_{c,T}(a,g)(t)
 =t^{-1/2}\pi_1\left[
  \Phi_\rho(t,1)a+
  \int_1^t\Phi_\rho(t,s)\binom0{g(s)}\dd s\right].
\end{equation}
For the source of the center equation, we use the norm
\[
 \norm{g}_{\cF_{c,T}^s}
 =\int_1^T\norm{g(t)}_{\mathfrak h^s}\dd t,
\]
and we abbreviate $\cC_T=\cC_T^0$, $\cF_{c,T}=\cF_{c,T}^0$, and
$\cE_T=\cE_T^0$. The energy
of the oscillator and $\int_1^\infty t^{-2}\dd t<\infty$ imply
\[
 \norm{u_c}_{\cC_T}
 \leq C\bigl(|a|+\norm{g}_{\cF_{c,T}}\bigr),
 \qquad T\geq1,
\]
because the eigenvalues of $\Omega_\rho^2$ remain in one compact
subinterval of $(0,\infty)$.

\smallskip\noindent\textit{Coupling the two components.}
The remaining operators $\mathcal R_{+-}$ and
$\mathcal R_{-+}$ in \eqref{appref:odd-sign-normal-form}
couple the center and positive components. They have the
form $t^{-1}B_\rho\pa_t+t^{-2}C_\rho^\sharp$ and satisfy
\begin{equation}\label{appref:odd-one-way-bounds}
 \norm{\mathcal R_{+-}u_c}_{\cE_T^*}
 \leq C|\rho|\norm{u_c}_{\cC_T},\qquad
 \norm{\mathcal R_{-+}u_+}_{\cF_{c,T}}
 \leq C|\rho|\norm{u_+}_{\cE_T}.
\end{equation}
Indeed, the spectral projections are analytic in $\rho$ and the
coefficients couple only adjacent modes, so both terms carry a factor
$O(\rho)$, while the Cauchy-Schwarz inequality applied to
\eqref{appref:odd-positive-energy} leaves the
integrable weights $\int_1^\infty t^{-2}\dd t$ and
$\int_1^\infty t^{-3}\dd t$. To include the boundary term, we set
\begin{equation}\label{appref:odd-complete-one-way}
 \widetilde{\mathcal R}_{+-}u_c
 :=\bigl(\mathcal R_{+-}u_c,-\mathcal B_{+-,\rho,T}u_c\bigr).
\end{equation}
Then, \eqref{appref:odd-one-way-bounds} and
\eqref{appref:odd-boundary-one-way} say, at every regularity level, that
\[
 \norm{\widetilde{\mathcal R}_{+-}u_c}_{
   (\cE_T^s)^*\oplus\mathfrak T_T^s}
 \leq C_s|\rho|\norm{u_c}_{\cC_T^s}.
\]
We now solve the coupled equations on the annulus. Substituting the
positive solution into the center equation and then reconstructing the
positive component, we obtain
\begin{equation}\label{appref:odd-returned-green}
 \begin{aligned}
 u_c={}&\bigl[I-V_{c,T}\mathcal R_{-+}G_{+,T}
                  \widetilde{\mathcal R}_{+-}\bigr]^{-1}
 V_{c,T}\bigl(\widetilde f_-+\mathcal R_{-+}G_{+,T}\widetilde f_+\bigr),\\
 u_+={}&G_{+,T}\bigl(\widetilde f_++\widetilde{\mathcal R}_{+-}u_c\bigr).
 \end{aligned}
\end{equation}
Here, $\widetilde f_-=(a_c,g_c)$ consists of the center Cauchy data
and interior source, and $\widetilde f_+=(q_+,a_+,b_+)$ consists of
the positive interior source, inner value, and outer boundary datum.
In these compositions, $V_{c,T}$ and $G_{+,T}$ use zero trace data
when only an interior source is supplied; the coupling
$\widetilde{\mathcal R}_{+-}u_c$ supplies zero inner value. We impose the boundary
condition \eqref{appref:odd-physical-row-split} through the second component of
\eqref{appref:odd-complete-one-way}. The correction inside the inverse
contains one coupling from the center to the positive subspace and one in
the reverse direction. Its norm is therefore at most $C\rho^2$, by
\eqref{appref:odd-positive-form}--\eqref{appref:odd-one-way-bounds} and
$\int_1^\infty t^{-2}\dd t<\infty$. This argument also bounds the center
contribution to the boundary operator.

The spaces on which the inverse \eqref{appref:odd-returned-green} acts,
and the estimates that we will use, are recorded in the following lemma.
\begin{lemma}
\label{appref:odd-exterior-green-poisson}
With $\mathfrak h_\pm^s=P_\pm(0)\mathfrak h^s$ as in
\eqref{appref:odd-balanced-trace}, we put
\[
 \begin{aligned}
 \mathfrak D_{c,T}^s
  &=(\mathfrak h_-^s\times\mathfrak h_-^s)
       \oplus\cF_{c,T}^s,\\
 \mathfrak D_{+,T}^s
  &=(\cE_T^s)^*\oplus\mathfrak h_+^{s+1/2}
       \oplus\mathfrak T_T^s.
 \end{aligned}
\]
% AUTHOR QUERY: Define the maximal center domain with L^1 oscillator
% forcing; the displayed C^1-type norm alone does not define this domain.
% Also define the higher-index positive source scale. The literal dual
% of E_T^s has the opposite Sobolev grading and does not give an E_T^s
% solution. These definitions are needed for the isomorphism below.
The three positive entries are, in order, the interior source in the dual
of the energy space, the inner Dirichlet value, and the balanced datum on
the outer boundary. For one $\rho_0>0$, every $0\leq|\rho|\leq\rho_0$, and
every $T\geq2$, the operator
\[
 \mathscr N_{\rho,T}^{\rm ext}:
 \cC_T^s\oplus\cE_T^s
 \longrightarrow\mathfrak D_{c,T}^s\oplus\mathfrak D_{+,T}^s
\]
which records the equations, the center Cauchy data, the positive
value at $t=1$, and the boundary data $\mathcal B_{\rho,T}(u_c+u_+)$, is an
invertible linear map. Its inverse is \eqref{appref:odd-returned-green},
with $\widetilde{\mathcal R}_{+-}$ as in
\eqref{appref:odd-complete-one-way}, and it satisfies
\begin{equation}\label{appref:odd-exterior-estimate}
 \norm{u_c}_{\cC_T^s}+\norm{u_+}_{\cE_T^s}
 \leq C_s\bigl(\norm{\widetilde f_-}_{\mathfrak D_{c,T}^s}
                    +\norm{\widetilde f_+}_{\mathfrak D_{+,T}^s}\bigr).
\end{equation}
The constants are independent of $T$ and of every angular cutoff.
% AUTHOR QUERY: Specify the weighted radial-derivative estimates needed
% for appref:carrier-grade-recurrence. The unweighted estimate above
% does not extend to D^j u_c uniformly in T.
\end{lemma}

The lemma estimates the inverse before differentiation in $\rho$. For
the derivatives $\pa_\rho^q$ with $q\geq1$, differentiating the phase of the
center oscillator produces powers of $t-1$, so the same bound need not be
uniform in $T$ without a loss of derivatives. The bounds in the analytic
norms, which absorb these powers, are proved below, once the weight has
been introduced.

\begin{proof}
At the base index, the Lax-Milgram lemma and
\eqref{appref:odd-positive-form} give the positive Green map with zero
inner value and homogeneous outer boundary condition. The extensions
$P_{+,T}a$ and $E_Tb$ in \eqref{appref:odd-positive-green-poisson}
supply the prescribed inner value and outer boundary datum.
For the center map, we use \eqref{appref:odd-center-green} and the energy
estimate for the oscillator. For the coupling terms, including the
boundary term, we have
\begin{equation}\label{appref:odd-exterior-one-way-summary}
 \begin{aligned}
 \norm{\widetilde{\mathcal R}_{+-}}_{
   \cC_T^s\to(\cE_T^s)^*\oplus\mathfrak T_T^s}
 &\leq C_s|\rho|,\\
 \norm{\mathcal R_{-+}}_{\cE_T^s\to\cF_{c,T}^s}
 &\leq C_s|\rho|.
 \end{aligned}
\end{equation}
Thus, the correction in \eqref{appref:odd-returned-green} has norm at
most $C_s\rho^2$. We decrease $\rho_0$ so that this norm is less than $1/2$ at
the base index. The Neumann series then gives a right inverse with the
asserted estimate, while uniqueness for the oscillator, coercivity of the
positive form, and the boundary condition
\eqref{appref:odd-physical-row-split} give the left inverse.

For higher regularity, the contour formulas for $P_\pm$ and $T_\rho^{\pm1}$,
together with the differential equation satisfied by $Q_\rho$, give for
each fixed $s$ and $q$
\begin{equation}\label{appref:odd-kato-all-grades}
 \norm{\Lambda_\theta^s\pa_\rho^qQ_\rho\Lambda_\theta^{-s}}
 +\norm{\Lambda_\theta^s\pa_\rho^qT_\rho^{\pm1}\Lambda_\theta^{-s}}
 \leq C_{s,q}.
\end{equation}
We commute $\Lambda_\theta^s$ through the coefficients of
\eqref{appref:odd-sign-normal-form}, which couple adjacent modes.
The diagonal angular term commutes with $\Lambda_\theta^s$
and remains in the coercive form. This gives
\eqref{appref:odd-exterior-one-way-summary} with $C_s$ in place of
$C_0$. Finally,
\eqref{appref:odd-boundary-one-way} and
\[
 \sqrt T\,\norm{u_+(T)}_{\mathfrak h^s}
 \leq C_s\norm{u_+}_{\cE_T^s}
\]
give the balanced value and conormal traces.
\end{proof}

Radial differentiation requires additional care. Even for a homogeneous
center solution $u_c=t^{-1/2}\xi$, we have
\[
 D u_c=-\tfrac12t^{-1/2}\xi+t^{1/2}\xi',
\]
so the bound for $u_c$ does not give a bound for $D u_c$ in the same center
norm uniformly in $T$. The radial estimates used below must account for
these powers of $t$ after the analytic weight is introduced.

The restriction $T\geq2$ in
\cref{appref:odd-exterior-green-poisson} keeps the two boundaries apart.
At $T=1$, they coincide, so the two trace data cannot be prescribed
independently. In the application, $T$ is the value $T_n=|n|/L$ at cell
frequency $n$, and the cell frequencies with $T_n<2$ form a finite set,
which we treat separately below.

\subsection{Matching at the interface}

To match the cap and annulus solutions, we use the positive
interface value $a$ at $t=1$ as an unknown in the fixed
coordinates. Solving the cap problem
\eqref{appref:odd-core-green} with this value determines the
center Cauchy data. Together with $a$, these data determine
the annular solution through \eqref{appref:odd-returned-green}.
With the interior sources, prescribed axis coefficients, and outer
boundary datum set to zero, we let $\mathcal M_{\rho,T}a$ be the cap
conormal derivative minus the annular one. We let $d_{\rho,T}(f,b)$
be the mismatch for the prescribed data with $a=0$. The source $f$ and the outer boundary datum $b$ contribute
$d_{\rho,T}(f,b)$, while the unknown value $a$ contributes
$\mathcal M_{\rho,T}a$. Thus, matching means
$\mathcal M_{\rho,T}a+d_{\rho,T}(f,b)=0$, or
\begin{equation}\label{appref:odd-interface-match}
 a=-\mathcal M_{\rho,T}^{-1}d_{\rho,T}(f,b).
\end{equation}

To invert $\mathcal M_{0,T}$ uniformly in $T$, we let $u_{\rm in}(a)$ and
$u_{\rm out}(a)$ be the solutions of the homogeneous positive equations
at the circle on the cap and on the annulus, with the common value $a$
at the interface and the homogeneous condition on the outer boundary.
The outward normals have opposite signs at the interface, so Green's
identity gives
\begin{equation}\label{appref:odd-circle-interface-energy}
 \langle\mathcal M_{0,T}a,a\rangle
 =\mathfrak a_{\rm in,0}(u_{\rm in},u_{\rm in})
  +\mathfrak a_{\rm out,0,T}(u_{\rm out},u_{\rm out})
 \geq c\norm{a}_{H^{1/2}}^2.
\end{equation}
Here, $\mathfrak a_{\rm in,0}$ and
$\mathfrak a_{\rm out,0,T}$ are the positive forms on the cap
and annulus, the latter given by \eqref{appref:odd-positive-form},
and $H^{\pm1/2}=\mathfrak h_+^{\pm1/2}$. Both forms control
the angular and mass terms, so the trace theorem at $t=1$
bounds $\norm{a}_{H^{1/2}}$ by either energy, uniformly in
$T$. A collar extension of $a$ and the Green estimates also
bound $\mathcal M_{0,T}$ from $H^{1/2}$ to $H^{-1/2}$.
The Lax-Milgram lemma on the trace space therefore gives
\[
 \norm{\mathcal M_{0,T}^{-1}}_{H^{-1/2}\to H^{1/2}}\leq C,
 \qquad T\geq2.
\]

At the circle, the center Cauchy data propagate without a condition on
the outer boundary, and they do not affect $\mathcal M_{0,T}$. For
$\rho\ne0$, we let $\Lambda_{\rm in,\rho}$ and
$\Lambda_{\rm out,\rho,T}^{+}$ denote the maps from the value $a$ to the
conormal derivative of the solution on the cap and on the positive part
of the annulus, in the fixed coordinates. The couplings with the
center are omitted from the latter, and we put all terms involving a
coupling between the center and the positive subspace into
$\mathcal E_{\rm ret,\rho,T}$. With the normal convention of
\eqref{appref:odd-circle-interface-energy}, we then have
\begin{equation}\label{appref:odd-interface-decomposition}
 \mathcal M_{\rho,T}
 =\Lambda_{\rm in,\rho}
   -\Lambda_{\rm out,\rho,T}^{+}
   +\mathcal E_{\rm ret,\rho,T},
 \qquad \mathcal E_{\rm ret,0,T}=0.
\end{equation}
To describe the coupling term, we let
$\mathcal C_{\rm in,\rho}a$ be the center Cauchy trace of the
homogeneous cap solution with value $a$. We write
$G_{+,\rho,T}$ for \eqref{appref:odd-positive-green-poisson}
at parameter $\rho$, put
$P_{+,\rho,T}a=G_{+,\rho,T}(0,a,0)$, and define
\[
 \Xi_{\rho,T}
 =I-V_{c,T}\mathcal R_{-+}G_{+,\rho,T}
              \widetilde{\mathcal R}_{+-},
\]
which is invertible since its correction has norm at most $C\rho^2<1/2$.
If we write $\Gamma_{N,+}$ for the conormal trace of the positive
component at $t=1$, then the last term in
\eqref{appref:odd-interface-decomposition} is
\begin{equation}\label{appref:odd-returned-dn-formula}
 \mathcal E_{\rm ret,\rho,T}a
 =-\Gamma_{N,+}G_{+,\rho,T}\widetilde{\mathcal R}_{+-}
   \Xi_{\rho,T}^{-1}V_{c,T}
   \bigl(\mathcal C_{\rm in,\rho}a
        +\mathcal R_{-+}P_{+,\rho,T}a\bigr).
\end{equation}
Both terms in the last parenthesis vanish at the circle. Away from it,
the perturbation formula \eqref{appref:odd-core-green} gives
$\norm{\mathcal C_{\rm in,\rho}}\leq C_s|\rho|$, and
\eqref{appref:odd-exterior-one-way-summary} gives the same bound for
$\mathcal R_{-+}P_{+,\rho,T}$, uniformly in $T$.

From the perturbation formula \eqref{appref:odd-core-green} and the
trace theorem on the interface $t=1$, we obtain
\begin{equation}\label{appref:odd-core-dn-difference}
 \norm{\Lambda_{\rm in,\rho}-\Lambda_{\rm in,0}}_{
   \mathfrak h_+^{s+1/2}\to\mathfrak h_+^{s-1/2}}
 \leq C_s|\rho|.
\end{equation}
The map $Q_\rho$ and the trace operator in the fixed coordinates satisfy
\begin{equation}\label{appref:odd-kato-trace-difference}
 \norm{Q_\rho-I}_{\mathfrak h^s\to\mathfrak h^s}
 +\norm{\mathcal T_\rho-\mathcal T_0}_{
      \mathfrak h_+^{s+1/2}\to\mathfrak h_+^{s-1/2}}
 \leq C_s|\rho|.
\end{equation}
Here, $\mathcal T_\rho$ denotes the cap or annular
Dirichlet-to-conormal map, $\Lambda_{\rm in,\rho}$ or
$\Lambda_{\rm out,\rho,T}^{+}$, expressed in the fixed coordinates. On the
positive subspace, the Green maps with homogeneous boundary data satisfy
\[
 G_{+,\rho,T}^0-G_{+,0,T}^0
 =G_{+,\rho,T}^0
   (\mathcal N_{++,0}-\mathcal N_{++,\rho})G_{+,0,T}^0,
\]
which, combined with \eqref{appref:odd-positive-form}, the Poisson
extensions, and the $O(\rho)$ change in the boundary operator, gives
\begin{equation}\label{appref:odd-outer-dn-difference}
 \norm{\Lambda_{\rm out,\rho,T}^{+}
             -\Lambda_{\rm out,0,T}^{+}}_{
   \mathfrak h_+^{s+1/2}\to\mathfrak h_+^{s-1/2}}
 \leq C_s|\rho|,
 \qquad T\geq2.
\end{equation}
Explicitly, we have the bounds
\[
 \begin{aligned}
 \sqrt T\,\norm{u_+(T)}_{\mathfrak h^s}
   &\leq C_s\norm{u_+}_{\cE_T^s},\\
 |u_c(T)|&\leq C_sT^{-1/2}\norm{u_c}_{\cC_T^s},&
 |D u_c(T)|&\leq C_sT^{1/2}\norm{u_c}_{\cC_T^s}.
 \end{aligned}
\]
Each term in $\mathcal E_{\rm ret,\rho,T}$ contains two
$O(\rho)$ factors. The first maps the positive component
to the center through the cap trace or $\mathcal R_{-+}$.
The second, $\widetilde{\mathcal R}_{+-}$, maps back to the
positive component and includes the boundary term
\eqref{appref:odd-boundary-one-way}. The intervening Green
maps and $\Xi_{\rho,T}^{-1}$ are uniformly bounded. Hence,
\begin{equation}\label{appref:odd-returned-dn-difference}
 \norm{\mathcal E_{\rm ret,\rho,T}}_{
   \mathfrak h_+^{s+1/2}\to\mathfrak h_+^{s-1/2}}
 \leq C_s\rho^2,
 \qquad T\geq2.
\end{equation}
Adding \eqref{appref:odd-core-dn-difference},
\eqref{appref:odd-outer-dn-difference}, and
\eqref{appref:odd-returned-dn-difference}, we obtain at every fixed index
\begin{equation}\label{appref:odd-interface-difference}
 \norm{\mathcal M_{\rho,T}-\mathcal M_{0,T}}_{
   \mathfrak h_+^{s+1/2}\to\mathfrak h_+^{s-1/2}}
 \leq C_s|\rho|,
 \qquad T\geq2.
\end{equation}
Here, the constants at higher regularity come from
\eqref{appref:odd-kato-all-grades} and from the induction in the radial
derivatives in the proof of \cref{appref:odd-exterior-green-poisson}. We
can therefore invert $\mathcal M_{\rho,T}$ by a Neumann series at the
base index. This gives \eqref{appref:odd-interface-match} and completes
the construction of the odd inverse for $T\geq2$, with arbitrary
compatible interior and boundary data. Substitution verifies that it is a
right inverse. Uniqueness for the equations on the cap, for the
oscillator, and for the coercive equation, together with the matching at
the interface, shows that it is also a left inverse.

\subsection{The analytic weight and the estimates for the odd block}

We now estimate the odd inverse in the analytic scale of
\cref{subsec:analytic-scale}. After rescaling, its radial
weight decreases outward and absorbs the growth of the
center solution and its $\rho$-derivatives. We set
$\lambda_n=\langle n\rangle$ and $T_n=|n|/L$. For $n\ne0$,
we factor out the axis value of the weight by writing
\begin{equation}\label{appref:carrier-phase}
 \omega_n(t)
 :=\exp\{\Phi_\gamma(Lt/|n|,\lambda_n)
               -\Phi_\gamma(0,\lambda_n)\}.
\end{equation}
Its logarithmic derivative is
\begin{equation}\label{appref:carrier-phase-damping}
 d_n(t):=-\pa_t\log\omega_n(t)
 =\gamma\frac{L^2t\lambda_n^2/n^2}
  {\sqrt{1+L^2t^2\lambda_n^2/n^2}},
 \qquad
 0\leq d_n(t)\leq\gamma L\frac{\lambda_n}{|n|}.
\end{equation}
For every fixed order $q$ of differentiation in $D=t\pa_t$, the derivatives
of $d_n$ are bounded uniformly on the cell frequencies $|n|\geq2L$, and
those of $td_n$ are bounded by $C_{q,L}\gamma(1+t)$. On the cap
$0\leq t\leq1$, both $\omega_n$ and its inverse have uniformly bounded
derivatives of each fixed order. On the positive part of the annulus,
conjugation by $\omega_n$ makes the substitutions
\[
 \pa_t\longmapsto\pa_t+d_n,\qquad
 D\longmapsto D+td_n,
\]
and changes the form by at most $C_L\gamma$ times the radial
and mass terms of \eqref{appref:odd-positive-energy}.
We choose $\gamma$ so that this is less than a quarter of
the coercivity margin. Thus, the positive form and its trace
estimates remain uniform after conjugation. For the center
component, we use the oscillator formula.

Since $\omega_n$ is decreasing and the fundamental matrix
$\Phi_\rho(t,s)$ is uniformly bounded, we obtain from
\eqref{appref:odd-center-green}
\begin{equation}\label{appref:weighted-center-volterra}
 \omega_n(t)|V_{c,T_n}(a,g)(t)|
 \leq C\left(\omega_n(1)|a|
      +\int_1^t\omega_n(s)|g(s)|\,\dd s\right).
\end{equation}
The same weight controls the derivatives in $\rho$ of the center
propagator $\Phi_\rho$. We write
\[
 a_n=\frac{L\lambda_n}{|n|}\in[L,\sqrt2L].
\]
For $t\geq t_0\geq1$, the rate in \eqref{appref:carrier-phase-damping} is
bounded below by $\gamma c_L$, so that
\begin{equation}\label{appref:carrier-ratio-decay}
 \begin{gathered}
 \frac{\omega_n(t)}{\omega_n(t_0)}
 =\exp\left\{-\int_{t_0}^td_n(u)\,\dd u\right\}
 \leq e^{-c_L\gamma(t-t_0)},\\
 c_L:=\min_{a\in[L,\sqrt2L]}
             \frac{a^2}{\sqrt{1+a^2}}>0 .
 \end{gathered}
\end{equation}
Differentiating the first-order system of the oscillator repeatedly by
Duhamel's formula, and using the uniform bound on its undifferentiated
propagator and the integrable coefficient $t^{-2}$, we obtain
\begin{equation}\label{appref:center-propagator-rho-growth}
 \norm{\pa_\rho^q\Phi_\rho(t,t_0)}_{\mathfrak h^s\to\mathfrak h^s}
 \leq C_{s,q}(1+t-t_0)^q,
 \qquad t\geq t_0\geq1.
\end{equation}
Here, $t_0$ denotes the initial radius and $s$ the Sobolev index. Even a
constant oscillator produces a factor $t-t_0$ when its
frequency is differentiated. The weight absorbs this polynomial growth.
Indeed, \eqref{appref:carrier-ratio-decay} and
$(1+\tau)^qe^{-c\tau}\leq C_q(1+c^{-1})^q$ give
\begin{equation}\label{appref:center-polynomial-absorption}
 \sup_{t\geq t_0\geq1}(1+t-t_0)^q
       \frac{\omega_n(t)}{\omega_n(t_0)}
 \leq C_{q,L}\gamma^{-q}.
\end{equation}
Setting
\[
 \begin{aligned}
 \norm{u_c}_{\cC_{T,\omega_n}^s}
  &:=\sup_{1\leq t\leq T}\omega_n(t)
       \bigl(\norm{\xi(t)}_{\mathfrak h^s}
             +\norm{\xi'(t)}_{\mathfrak h^s}\bigr),\\
 \norm{g}_{\cF_{c,T,\omega_n}^s}
  &:=\int_1^T\omega_n(t)\norm{g(t)}_{\mathfrak h^s}\,\dd t,
 \qquad u_c=t^{-1/2}\xi,
 \end{aligned}
\]
we obtain from \eqref{appref:center-propagator-rho-growth},
\eqref{appref:center-polynomial-absorption}, and the formula
\eqref{appref:odd-center-green}
\begin{equation}\label{appref:weighted-center-rho-bound}
 \norm{\pa_\rho^qV_{c,T}(a,g)}_{\cC_{T,\omega_n}^s}
 \leq C_{s,q,L,\gamma}
   \bigl(\omega_n(1)|a|+
         \norm{g}_{\cF_{c,T,\omega_n}^s}\bigr).
\end{equation}

The positive form, $Q_\rho$, $T_\rho^{\pm1}$, and the
coupling maps have weighted bounds for every fixed number of
$\rho$-derivatives. Since the weight is independent of
$\rho$, we can differentiate
\eqref{appref:odd-returned-green}, including its Neumann
series, to obtain the weighted version of
\eqref{appref:odd-exterior-estimate} at each fixed order.
The balanced value and conormal traces are included, and
the bounds are uniform for $T\ge2$.

At the interface, an $O(\rho)$ coupling occurs on either
side of the center propagator. Its uniform bound therefore
suffices without estimating
$V_{c,T}(\rho)-V_{c,T}(0)$. At the outer boundary, the
factor $\omega_n(T_n)$ agrees with the prescribed trace
weight. Thus, all factors in \eqref{appref:odd-returned-green}
have weighted bounds independent of the cutoffs.

We write $\mathcal X_{y,n}^j$, $\mathcal Y_{y,n}^j$, and
$\mathcal T_{y,n}^j$ for the norms at index $j$ of the Cartesian domain,
of the source, and of the balanced traces of
\eqref{eq:A-balanced-trace-factor}, after the rescaling $t=|n|r/L$,
defined with
$\pa_r=(|n|/L)\pa_t$. Removing the common factor
$e^{\Phi_\gamma(0,\lambda_n)}$, we obtain the equivalent norms at cell
frequency $n$
\begin{equation}\label{appref:carrier-mode-norms}
 \begin{aligned}
 e^{-2\Phi_\gamma(0,\lambda_n)}
 \norm{u_n}_{\cX_{\gamma,n}^s}^2
 &\simeq\sum_{j+\ell\leq s}\lambda_n^{2\ell}
       \norm{\omega_nu_n}_{\mathcal X_{y,n}^j}^2,\\
 e^{-2\Phi_\gamma(0,\lambda_n)}
 \norm{f_n}_{\cY_{\gamma,n}^s}^2
 &\simeq\sum_{a+k\leq s}\lambda_n^{2k}
       \norm{\omega_nf_n}_{\mathcal Y_{y,n}^a}^2.
 \end{aligned}
\end{equation}
Here, the boundary components are measured in $\mathcal T_{y,n}^a$. We
write $K_{e,n}^{\rm od}$ for the odd inverse at cell frequency $n$
constructed above, and $\operatorname{Tr}_{\rm phys}$ for the trace on
the outer boundary. From the energy estimate, Cartesian difference
quotients, and the trace theorem in the balanced norms, we obtain
\begin{equation}\label{appref:carrier-parameter-estimate}
 \begin{aligned}
 &\norm{\omega_nK_{e,n}^{\rm od}f_n}_{\mathcal X_{y,n}^j}
 +\norm{\omega_n(T_n)\operatorname{Tr}_{\rm phys}
             K_{e,n}^{\rm od}f_n}_{\mathcal T_{y,n}^j}\\
 &\hspace{25mm}\leq C_j\sum_{a=0}^{j+4}
       \lambda_n^{j+4-a}
       \norm{\omega_nf_n}_{\mathcal Y_{y,n}^a}.
 \end{aligned}
\end{equation}
We denote the left side by $E_{j,n}^{(0)}$ and set
$F_{a,n}=\norm{\omega_nf_n}_{\mathcal Y_{y,n}^a}$, including
all source trace entries. The base Green estimate and
\eqref{appref:odd-source-comparison} give
$E_{0,n}^{(0)}\le CF_{4,n}$. Applying one more Cartesian
difference quotient leaves the principal term in the
coercive part. Each commutator contains either a disk
derivative, carrying the rescaling factor
$|n|/L\le C_L\lambda_n$, or a term involving $d_n$,
$Dd_n$, or its balanced trace. By
\eqref{appref:carrier-phase-damping} and
\eqref{appref:odd-exterior-estimate}, this gives
\begin{equation}\label{appref:carrier-grade-recurrence}
 E_{j,n}^{(0)}\leq C_j\bigl(F_{j+4,n}
                   +\lambda_nE_{j-1,n}^{(0)}\bigr),
 \qquad j\geq1,
\end{equation}
whose iteration yields
\[
 E_{j,n}^{(0)}\leq C_j\sum_{a=4}^{j+4}
       \lambda_n^{j+4-a}F_{a,n},
\]
which proves \eqref{appref:carrier-parameter-estimate}.
The loss consists of three derivatives in
\eqref{appref:odd-source-comparison} and one in the weighted
annular and trace estimates, independently of the upper
frequency cutoff.

We now turn to the derivatives in $\rho$.
We define $E_{j,n}^{(q)}$ as $E_{j,n}^{(0)}$ with the solution and its
outer trace replaced by those of $\pa_\rho^qK_{e,n}^{\rm od}f_n$.
Differentiating the equations, and using
\eqref{appref:weighted-center-rho-bound} for the center and
\eqref{appref:odd-kato-all-grades} for the other factors, we obtain
\begin{equation}\label{appref:odd-rho-grade-recurrence}
 \begin{aligned}
 E_{0,n}^{(q)}
 &\leq C_{q,\gamma}\left(
 F_{4+q,n}+\lambda_n\sum_{p<q}E_{0,n}^{(p)}\right),\\
 E_{j,n}^{(q)}
 &\leq C_{j,q,\gamma}\left(
 F_{j+4+q,n}+\lambda_nE_{j-1,n}^{(q)}
 +\lambda_n\sum_{p<q}E_{j,n}^{(p)}\right),
 \qquad j\geq1,\ q\geq1.
 \end{aligned}
\end{equation}
Here, the last term accounts for the factor $\lambda_n$, that is, one
cell derivative, produced by each differentiated coefficient or trace. A
double induction in $(q,j)$ then
gives
\begin{equation}\label{appref:odd-rho-component-bound}
 E_{j,n}^{(q)}\leq C_{j,q,\gamma}
   \sum_{a=0}^{j+4+q}\lambda_n^{j+4+q-a}F_{a,n}.
\end{equation}
Note that all norms in this induction are weighted.

Restoring the common factor \(e^{\Phi_\gamma(0,\lambda_n)}\) and
multiplying \eqref{appref:carrier-parameter-estimate} by
$\lambda_n^\ell$, with $j+\ell\leq s$, we obtain terms with at most
$s+4$ derivatives of the source, since
\[
 a+(\ell+j+4-a)=\ell+j+4\leq s+4.
\]
By a finite Cauchy-Schwarz inequality in $(j,\ell,a)$ and Parseval's
identity, we therefore obtain
\begin{equation}\label{appref:high-carrier-square-sum}
 \sum_{|n|\geq2L}\norm{K_{e,n}^{\rm od}f_n}_{
       \cX_{\gamma,n}^s}^2
 \leq C_s\sum_{|n|\geq2L}\norm{f_n}_{
       \cY_{\gamma,n}^{s+4}}^2.
\end{equation}
The same multiplication, applied to
\eqref{appref:odd-rho-component-bound}, gives for every fixed $q$
\begin{equation}\label{appref:odd-rho-high-square-sum}
 \sum_{|n|\geq2L}\norm{\pa_\rho^qK_{e,n}^{\rm od}f_n}_{
       \cX_{\gamma,n}^s}^2
 \leq C_{s,q,\gamma}\sum_{|n|\geq2L}\norm{f_n}_{
       \cY_{\gamma,n}^{s+4+q}}^2.
\end{equation}

The remaining cell frequencies satisfy $T_n<2$, or
$|n|<2L$, and form a finite set. At each such frequency, let
$\mathscr L_{0,n}^{\rm od}:\mathscr D_n^{s_*}\to
\mathscr R_n^{s_*}$ be the closed Cartesian operator at
the circle, including its trace operators, and let
$\mathscr L_{\rho,n}^{\rm od}$ be the corresponding
operator at the ellipse on the same domain.

For $n\ne0$, the inverse $K_{0,n}^{\rm od}$ is obtained on
the whole disk by solving the center Bessel equations with
the prescribed axis coefficients and the positive equations
with their Robin data. No interface is needed. We describe
$n=0$ below. Since the domain is unchanged, we have
\begin{equation}\label{appref:bounded-carrier-perturbation}
 \norm{\mathscr L_{\rho,n}^{\rm od}
      -\mathscr L_{0,n}^{\rm od}}_{
      \mathscr D_n^{s_*}\to\mathscr R_n^{s_*}}
 \leq C_{n,s_*}|\rho|,
\end{equation}
and the inverse at cell frequency $n$ is
\begin{equation}\label{appref:bounded-carrier-inverse}
 K_{\rho,n}^{\rm od}
 =\bigl[I+K_{0,n}^{\rm od}
   (\mathscr L_{\rho,n}^{\rm od}
       -\mathscr L_{0,n}^{\rm od})\bigr]^{-1}K_{0,n}^{\rm od}.
\end{equation}
Taking the maximum over this finite set and decreasing $\rho_0$, we make
the norm of $K_{0,n}^{\rm od}(\mathscr L_{\rho,n}^{\rm od}
-\mathscr L_{0,n}^{\rm od})$ less than $1/2$ for every $n$ in this set,
so that the Neumann series in \eqref{appref:bounded-carrier-inverse}
converge.

For $n=0$, we use the original variables. The second force equation
\eqref{eq:linearized-H2} is $f_2=L\cR C$, so on every nonzero odd
angular mode
\begin{equation}\label{appref:zero-carrier-tangent}
 C_m=\frac{f_{2,m}}{iLm}.
\end{equation}
The first two equations still give the last two equations in
\eqref{appref:spin-equations}. We set
\[
 V=(-i\cR)^{-1}X.
\]
The terms depending on $V$ are
\[
 A=-iP(D+2)V,\qquad
 B=-i\bigl(P\cR-2P\cR^{-1}D\bigr)V.
\]
All other terms depend only on $g_\pm$. The last equation gives
\eqref{appref:odd-cartesian-operator} for $u=|\cR|V$, with the mass term
omitted. At the circle, the center satisfies the Laplace equation with
prescribed coefficients of $r e^{i\theta}$ and $r e^{-i\theta}$ in
its Taylor expansion on the axis. On $|m|\geq3$, the
condition on the outer boundary is of \emph{Robin type},
that is, it prescribes a combination of the value and of the radial
derivative, here with the boundary operator $D+2$ at $r=1$. Since the operators at the
circle and at the ellipse have a
common domain, the same perturbation argument gives the inverse for a
small ellipse. We recover $X=(-i\cR)V$, then $Z$, $W$, $A$, $B$, and $S$
from \eqref{appref:spin-equations}, keeping $C$ as prescribed in
\eqref{appref:zero-carrier-tangent}.

Combining the construction above with
\eqref{appref:high-carrier-square-sum} and
\eqref{appref:bounded-carrier-inverse}, we obtain
\begin{equation}\label{appref:odd-inverse}
 K_e^{\rm od}\mathscr L_e^{\rm od}=I,\qquad
 \mathscr L_e^{\rm od}K_e^{\rm od}=I,\qquad
 \norm{K_e^{\rm od}f}_{\cX_\gamma^s}
 \leq C_s\norm{f}_{\cY_\gamma^{s+\mu_{\rm od}}},
\end{equation}
with $C_s$ independent of the cell frequency. Note that the estimates
above allow the conservative choice $\mu_{\rm od}=4$. Differentiating
the finitely many Neumann formulas
\eqref{appref:bounded-carrier-inverse} and combining them with
\eqref{appref:odd-rho-high-square-sum}, we obtain in addition, for every
fixed $q\geq0$,
\begin{equation}\label{appref:odd-global-rho-bound}
 \norm{\pa_\rho^qK_e^{\rm od}f}_{\cX_\gamma^s}
 \leq C_{s,q,\gamma}
       \norm{f}_{\cY_\gamma^{s+4+q}}.
\end{equation}
This completes the construction and the estimates for the odd block.

\subsection{The even block}

The even summand contains the angular means, the scale
equation, and the affine modes $m=\pm2$. Each multiplier
in \eqref{appref:spin-equations} vanishes on one exceptional
mode, so these modes must be treated before $|m|\ge4$.
After solving for the means and scale, we determine
$C_{\pm2}$, then solve a coercive equation for
$C_h=\sum_{m\in2\Z,\,|m|\ge4}C_m e^{im\theta}$.
The nonsingular equations recover $Z$ and $W$, leaving
$W_2$ and $Z_{-2}$ for the final radial integrations.

\smallskip\noindent\textit{Means, scale, and radial integration.}
We begin with the algebraic equations. We call a vector
\emph{transverse} when it lies in the normal plane $e_T^\perp$, and
\emph{tangential} when it is parallel to the axis $e_T$. The datum $h_{\rm af}$,
which is the scale entry of the auxiliary data $r_{\rm aug}$ of
\eqref{app:cross:source-graph-chart}, prescribes the scale $\dot a$ of
the transverse component of the linear Taylor term of the variation through the equation
\begin{equation}\label{appref:scale-row}
 D H_{\rm af}(1,0)[\dot a,\dot\tau]=4\dot a,\qquad
 \dot a=\frac14h_{\rm af}.
\end{equation}
Here, we temporarily allow the transverse scale \(a\) to vary
independently of \(\tau\). The scalar constraint
\(H_{\rm af}(a,\tau)=2a^2+|\tau|^2-2\) measures the error in the trace
normalization \eqref{eq:normalized-affine-trace}. Its derivative is
evaluated at \((a,\tau)=(1,0)\). We then subtract the prescribed
transverse linear term from the variation. The two gauge
conditions prescribe the means of $B$ and $C$, the mean condition on the
scalar prescribes that of $S$, and averaging the first equation of
\eqref{appref:raw-scalar-rows} gives $A_0=-f_{0,0}/2$, since the means
of $\cR S$ and $\cR B$ vanish. We keep these known means in the right
sides of the equations for $m=\pm2$ below. From now on, the source
denotes the remainder after the algebraic contributions have been
subtracted.

For $n\ne0$, we introduce the two integral operators
\begin{equation}\label{appref:even-hardy-maps}
 J_+q=t^{-2}\int_0^t s q(s)\dd s,\qquad
 J_-q=t^2\int_0^t q(s)s^{-3}\dd s,
\end{equation}
which we call the \emph{Hardy maps}. They solve the radial
equations through $(D+2)J_+=I$ and $(D-2)J_-=I$, the
latter on functions with zero $t^2$-coefficient. Integration
from the axis selects zero coefficients for the homogeneous
solutions $t^{-2}$ and $t^2$. The restriction on $J_-$
also makes its integral converge. We verify below that
the sources to which it is applied satisfy this restriction.

We use these maps in conjugated form. For the weight $\omega_n$ of
\eqref{appref:carrier-phase}, we put
\begin{equation}\label{appref:weighted-even-hardy}
 J_{\pm,n}^\gamma:=\omega_nJ_\pm\omega_n^{-1}.
\end{equation}
Conjugation replaces $D$ by $D+td_n$, where $d_n$ is the logarithmic
derivative in \eqref{appref:carrier-phase-damping}, so that
\begin{equation}\label{appref:weighted-even-hardy-identities}
 \begin{aligned}
 (D+td_n+2)J_{+,n}^\gamma&=I,\\
 (D+td_n-2)J_{-,n}^\gamma&=I
 \quad\text{on }\mathcal Z_{2,n}^\gamma .
 \end{aligned}
\end{equation}
Here, $\mathcal Z_{2,n}^\gamma$ is the subspace of the weighted modes
$m=\pm2$ on which $[t^2](\omega_n^{-1}q)=0$, where $[t^2]q$ denotes the
coefficient of $t^2$ in the expansion of $q$ at the axis, and the
kernels are
\[
 \begin{aligned}
 J_{+,n}^\gamma q(t)
  &=t^{-2}\int_0^t s\,
       \frac{\omega_n(t)}{\omega_n(s)}q(s)\,\dd s,\\
 J_{-,n}^\gamma q(t)
  &=t^2\int_0^t s^{-3}
       \frac{\omega_n(t)}{\omega_n(s)}q(s)\,\dd s.
 \end{aligned}
\]
Both weighted Hardy maps are bounded at every sufficiently high fixed
Sobolev regularity in Cartesian coordinates, uniformly in $n$ and $T_n$. Indeed, on $t\le1$, the
weight and its inverse have uniformly bounded derivatives, and we may
write
\[
 J_-q(y)=\int_0^1\tau^{-3}q(\tau y)\,\dd\tau.
\]
Taking $k$ Cartesian derivatives brings the factor $\tau^k$ under the
integral, and the change of variables $y\mapsto\tau y$ in $L^2(\D^2)$
costs the factor $\tau^{-1}$, so the $k$-th derivatives of $J_-q$ are
bounded by $\int_0^1\tau^{k-4}\dd\tau$ times those of $q$, and this
integral is finite for $k>3$. For $k=4$,
the fourth derivatives control the lower ones, since the only
polynomials of degree at most three on the modes $m=\pm2$ are quadratic,
and their coefficients vanish on the subspace in question. For $k>4$, we
combine the estimate for the top derivatives with this $H^4$ bound. The
same scaling calculation bounds $J_+$ at every order of differentiation,
with $\tau$ in place of $\tau^{-3}$. For $t\ge1$, we split the integrals
at $s=1$. The contribution from the cap grows at most quadratically and
is controlled by $\omega_n(t)$. For $1\le s\le t$,
\eqref{appref:carrier-ratio-decay} gives
\[
 t^2s^{-3}\frac{\omega_n(t)}{\omega_n(s)}
 \le (1+t-s)^2e^{-c_L\gamma(t-s)}.
\]
Young's inequality then gives the bounds on the annulus. The exponential
absorbs the polynomial factors from the radial measure and from any
fixed number of derivatives, and local trace estimates give the balanced
bounds at the endpoints. Since the weight is smooth in $t^2$ at the
axis, $\omega_n^{\pm1}$ preserves the jet and divisibility conditions.
The extensions of the finitely many mean, scale, and affine polynomial
data satisfy the same weighted estimates.

\smallskip\noindent\textit{The exceptional tangential components.}
The multiplier $\cR-2i$ in \eqref{appref:spin-equations} vanishes at
$m=2$, and $\cR+2i$ vanishes at $m=-2$, which is why these two modes
were called exceptional in \cref{sec:block}. At these two modes, the
corresponding equation does not determine $W_2$ or $Z_{-2}$, but instead
imposes an equation on $X$,
\[
 (D+2)X_2=g_{-,2},\qquad (D+2)X_{-2}=g_{+,-2}.
\]
Applying $J_+$ and substituting $X=i\sigma(f_2-\cR C)$, we determine the
two exceptional tangential components. The multiplier that does not
vanish at each mode gives the other two components, and we obtain
\begin{equation}\label{appref:even-exceptional-tangent}
 \begin{aligned}
 C_2&=-\frac i2f_{2,2}+\frac\sigma2J_+g_{-,2},\\
 C_{-2}&=\frac i2f_{2,-2}-\frac\sigma2J_+g_{+,-2},\\
 Z_2&=\frac{(D-2)X_2-g_{+,2}}{4i},\qquad
 W_{-2}=\frac{(D-2)X_{-2}-g_{-,-2}}{-4i}.
 \end{aligned}
\end{equation}
These four components are now known. The divisibility of the source gives $g_{+,2}=O(t^4)$ and
$g_{-,2}=O(t^2)$. Thus, $X_2=J_+g_{-,2}$ has a quadratic leading term,
$(D-2)X_2=O(t^4)$, and \eqref{appref:even-exceptional-tangent} gives
$Z_2=O(t^4)$. The same calculation at $m=-2$ gives $W_{-2}=O(t^4)$.
Hence,
\begin{equation}\label{appref:even-affine-cancellation}
 [t^2]Z_2=[t^2]W_{-2}=0,
\end{equation}
so the right-hand sides to which $J_-$ will be applied at the final step
satisfy the required condition.

\smallskip\noindent\textit{The high even modes.}
We now turn to the modes $|m|\geq4$. We let $\Pi_h$ be the projection
onto these modes and set $C_h=\Pi_hC$. In the terms of
\eqref{appref:spin-equations} that depend on $C_h$, we substitute
\begin{equation}\label{appref:even-high-substitution}
 \begin{aligned}
 (A,B)&=-i\sigma(\widehat A,\widehat B),\\
 \widehat A&=P(D+2)C_h,\\
 \widehat B&=P\cR C_h-2P\cR^{-1}DC_h,\\
 P&=\frac{\cR^2}{\cR^2+4}.
 \end{aligned}
\end{equation}
Although we do not yet know $W_2$ and $Z_{-2}$, they do not enter this
projection, since $H_\rho$ sends $W_2$ only to the modes $0$ and $2$,
and $Z_{-2}$ only to the modes $0$ and $-2$. The projected last equation is
therefore a closed equation for $C_h$. On these modes, we have
$1\leq p_m=m^2/(m^2-4)\leq4/3$, and the resulting form $\mathfrak b_\rho$
satisfies
\begin{equation}\label{appref:even-high-coercivity}
 \begin{aligned}
 \operatorname{Re}\mathfrak b_\rho(C_h,C_h)
 &\geq\left(1-3\frac\rho{1-\rho}\right)
       \norm{(DC_h,\cR C_h)}_2^2\\
 &\quad+\text{nonnegative mass and Robin terms}.
 \end{aligned}
\end{equation}
To prove \eqref{appref:even-high-coercivity}, we integrate by parts,
after the constant normalization of the equations already made in
\eqref{appref:raw-scalar-rows}, and we obtain
\begin{equation}\label{appref:even-high-form-splitting}
 \operatorname{Re}\mathfrak b_\rho(C_h,C_h)
 =\norm{(DC_h,\cR C_h)}_2^2
   +\mathfrak m_{\rho,n}(C_h)+\mathfrak r_{\rho,T}(C_h)
   +\operatorname{Re}\mathfrak e_\rho(C_h).
\end{equation}
Here, $\mathfrak m_{\rho,n}\ge0$ is the mass contribution
from $(1-\rho^2)n^2/L^2$, and $\mathfrak r_{\rho,T}\ge0$
is the squared boundary trace supplied by the Robin
condition. Only $\mathfrak e_\rho$, the contribution of
$H_\rho-I$, may have either sign. By
\eqref{appref:even-high-substitution}, it satisfies
\begin{equation}\label{appref:even-high-form-error}
 |\operatorname{Re}\mathfrak e_\rho(C_h)|
 \leq\norm{H_\rho-I}\,
      \norm{(\widehat A,\widehat B)}_2
      \norm{(DC_h,\cR C_h)}_2.
\end{equation}
Indeed,
\[
 \norm{H_\rho-I}=\frac{\rho}{1-\rho},\qquad
 \norm{(\widehat A,\widehat B)}_2
 \leq3\norm{(DC_h,\cR C_h)}_2.
\]
Here, the first identity follows from the eigenvalues
$\rho(\rho\pm1)/(1-\rho^2)$ of $H_\rho-I$, and the second follows mode
by mode from $p_m\leq4/3$ and $|m|\geq4$. Together,
\eqref{appref:even-high-form-splitting}--\eqref{appref:even-high-form-error}
give \eqref{appref:even-high-coercivity}. Its coefficient
$1-3\rho/(1-\rho)$ is positive for $0<\rho<1/4$, which is the range
allowed by \eqref{appref:ellipse}.

The prescribed trace on the outer boundary is $b=\Pi_h\mathfrak a(T)$.
% AUTHOR QUERY: Specify the Robin form and its test space. The prescribed
% physical row is a conormal/Robin datum; the form test space must impose
% only the essential Dirichlet conditions. The phrase "zero-trace" below
% does not specify which trace vanishes.
We let $E_Tb$ be a Cartesian extension of it, bounded without loss of
derivatives, and we define the Green map $G_{\rho,\sigma,T}^h$ for
zero trace by
\[
 \mathfrak b_\rho(G_{\rho,\sigma,T}^hF,\psi)
   =\langle F,\psi\rangle
 \quad\text{for every zero-trace }\psi.
\]
We write $\mathscr A_{\rho,n}^h$ for the operator in the projected
last equation associated with $\mathfrak b_\rho$, and $F_h$ for the source
that remains after inserting \eqref{appref:even-exceptional-tangent} and
projecting. The solution with prescribed outer trace is then
\begin{equation}\label{appref:even-green-poisson}
 C_h=E_Tb+G_{\rho,\sigma,T}^h
   \bigl(F_h-\mathscr A_{\rho,n}^hE_Tb\bigr),
\end{equation}
for arbitrary trace data. In the original disk coordinates,
before rescaling $t=|n|r/L$, the operator has the following form after an invertible
constant scaling of the equations,
\begin{equation}\label{appref:even-high-cartesian}
 \mathscr A_{\rho,n}^hu
 =\Pi_h\bigl[P\Delta_yu-2\rho\mathsf H
  (\pa_z^2+\pa_{\bar z}^2)\mathsf H u
 -(1-\rho^2)\frac{n^2}{L^2}u\bigr],\qquad u=|\cR|C_h.
\end{equation}

To obtain the analytic estimates, we return to $t$ and set
$v_h=\omega_nC_h$. The form
$\mathfrak b_{\rho,\gamma,n}$ of
$\omega_n\mathscr A_{\rho,n}^h\omega_n^{-1}$, with its
conjugated boundary operator, replaces every $D$ by
$D+td_n$. At the base index, we use the energy
\begin{equation}\label{appref:even-weighted-energy}
 \norm{v_h}_{\mathcal E_{{\rm ev},n}}^2
 :=\norm{(Dv_h,\cR v_h)}_2^2+\norm{tv_h}_2^2
   +\norm{\operatorname{tr}_T v_h}_{
             \mathcal T_{{\rm ev},{\rm val}}}^2,
\end{equation}
and at index $j$, we write $\mathcal E_{{\rm ev},n}^j$ for the version
obtained by commuting Cartesian derivatives through the form. Here,
$\mathcal T_{{\rm ev},{\rm val}}$ is the balanced norm of the value on
the outer boundary, which is the contribution of the value to the Robin
form, and we write $\mathcal T_{{\rm ev},n}^j$ for the balanced norm of
the value and conormal traces at index $j$. We recover the conormal
trace afterwards from the boundary equation. The coercivity margin
that \eqref{appref:even-high-coercivity} leaves uniformly on
$0\leq\rho\leq\rho_0$ is
\[
 c_{\rm ev}:=1-3\frac{\rho_0}{1-\rho_0}>0.
\]
Expanding $D+td_n$, integrating the cross term by parts, and using
\eqref{appref:carrier-phase-damping}, we obtain, after retaining the
nonnegative mass and outer boundary contributions,
\begin{equation}\label{appref:even-weighted-form-error}
 \operatorname{Re}\mathfrak b_{\rho,\gamma,n}(v_h,v_h)
 -\operatorname{Re}\mathfrak b_\rho(v_h,v_h)
 \geq-C_L(\gamma+\gamma^2)
       \norm{v_h}_{\mathcal E_{{\rm ev},n}}^2.
\end{equation}
Indeed, $td_n=O_L(\gamma t)$, with the same bound for every fixed number
of $D$ derivatives. On $t\leq1$, the angular Poincar\'e inequality on the
modes $|m|\geq4$ controls the zeroth-order terms produced by this
expansion, while on $t\geq1$, the $t^2$-mass term controls them. The
boundary operator gains the term $T_nd_n(T_n)v_h(T_n)$, whose sign is
favorable because the weight decreases outward, and its remaining change
is $O_L(\gamma)$ from the balanced value space to the conormal space. We
include the condition $C_L(\gamma+\gamma^2)<c_{\rm ev}/4$ in the initial
choice of $\gamma$, uniformly in $n$. Then,
\begin{equation}\label{appref:even-weighted-coercivity}
 \operatorname{Re}\mathfrak b_{\rho,\gamma,n}(v_h,v_h)
 \geq\frac{c_{\rm ev}}2
       \norm{v_h}_{\mathcal E_{{\rm ev},n}}^2,
 \qquad |n|\geq2L,
\end{equation}
uniformly in the angular cutoff and in the cell frequency. The
Lax-Milgram lemma, a weighted extension of the trace, and the boundary
condition then give
\begin{equation}\label{appref:even-weighted-green-trace}
 \norm{v_h}_{\mathcal E_{{\rm ev},n}^j}
 +\norm{\operatorname{Tr}_{\rm phys}v_h}_{
             \mathcal T_{{\rm ev},n}^j}
 \leq C_j\left(
   \norm{\omega_nF_h}_{(\mathcal E_{{\rm ev},n}^j)^*}
  +\norm{\omega_n(T_n)b}_{\mathcal T_{{\rm ev},n}^j}\right).
\end{equation}
% AUTHOR QUERY: As in the odd problem, define the higher-index source
% scale in this estimate. The literal dual of E^j has negative Sobolev
% grading; it cannot give an E^j solution. Also justify commuting
% Cartesian difference quotients with the angular projection Pi_h.
Cartesian difference quotients applied to
\eqref{appref:even-weighted-coercivity} give this bound at every regularity level
$j$, uniformly in the angular cutoff and in the cell frequency.

Once $C$ is known on every mode, we use \eqref{appref:spin-equations} to
recover $Z$ and $W$ at every mode on which the multipliers do not
vanish,
\begin{equation}\label{appref:even-nonsingular-spins}
 \begin{aligned}
 Z_m&=\frac{(D-m)X_m-g_{+,m}}{i(m+2)}
       &&(m\ne-2),\\
 W_m&=\frac{(D+m)X_m-g_{-,m}}{i(m-2)}
       &&(m\ne2).
 \end{aligned}
\end{equation}
\smallskip\noindent\textit{The remaining affine components.}
It remains to determine $W_2$ and $Z_{-2}$. For this, we write
$Z^G=\mathfrak a+i\mathfrak b$ and $W^G=\mathfrak a-i\mathfrak b$ for
the complex combinations of $(\mathfrak a,\mathfrak b)$. Multiplication
by $H_\rho$ gives
\begin{equation}\label{appref:metric-spin-identities}
 Z_k^G=\frac{Z_k-\rho W_{k+2}}{1-\rho^2},\qquad
 W_k^G=\frac{W_k-\rho Z_{k-2}}{1-\rho^2},
\end{equation}
and hence
\begin{equation}\label{appref:metric-divergence-identity}
 \begin{aligned}
 (D\mathfrak a+\cR\mathfrak b)_k
  =\frac1{2(1-\rho^2)}\bigl[&(D+k)Z_k+(D-k)W_k\\
   &-\rho(D-k)Z_{k-2}-\rho(D+k)W_{k+2}\bigr].
 \end{aligned}
\end{equation}
Substituting \eqref{appref:metric-divergence-identity} into the last equation
of \eqref{appref:raw-scalar-rows} at $k=2$ and $k=-2$, we are left with
one unknown in each equation,
\begin{equation}\label{appref:even-final-affine-equations}
 \begin{aligned}
 (D-2)W_2={}&2(1-\rho^2)(f_{3,2}-i\sigma t^2C_2)
 -(D+2)Z_2\\
 &\quad+\rho(D-2)Z_0+\rho(D+2)W_4,\\
 (D-2)Z_{-2}={}&2(1-\rho^2)(f_{3,-2}-i\sigma t^2C_{-2})
 -(D+2)W_{-2}\\
 &\quad+\rho(D+2)Z_{-4}+\rho(D-2)W_0.
 \end{aligned}
\end{equation}
We check that both right-hand sides have zero $t^2$-coefficient. Indeed,
\eqref{appref:even-affine-cancellation} removes that coefficient from
$(D+2)Z_2$ and $(D+2)W_{-2}$, and the factor $D-2$ removes it from the
terms containing $k=0$. The remaining terms vanish to fourth order by
the divisibility $f_3\in z\bar zC^\infty$ of the source, together with the
bounds $C_{\pm2}=O(t^2)$ and $W_4,Z_{-4}=O(t^4)$. Thus, $J_-$ inverts
$D-2$ on both right-hand sides, and
\begin{equation}\label{appref:even-final-hardy}
 W_2=J_-\operatorname{RHS}_2,\qquad
 Z_{-2}=J_-\operatorname{RHS}_{-2}.
\end{equation}
Here, $\operatorname{RHS}_{\pm2}$ denote the right-hand sides of
\eqref{appref:even-final-affine-equations}. Finally, we recover the original coordinates from
\begin{equation}\label{appref:even-reconstruction}
 A=\frac{Z+W}{2},\qquad B=\frac{Z-W}{2i},\qquad S=X+B.
\end{equation}
Note that the reconstruction uses only nonvanishing angular multipliers
and the two Hardy maps on the subspaces on which they are defined. It
therefore preserves the divisibility of $Z$ and $W$.

For $n=0$, the second force equation in the original variables determines $C$, and we obtain
\[
 C_m=\frac{f_{2,m}}{iLm}\quad(m\ne0),\qquad
 (D+2)X_2=g_{-,2},\qquad (D+2)X_{-2}=g_{+,-2}.
\]
We solve the high modes for $V_h=(-i\cR)^{-1}X_h$. Their equations are
\eqref{appref:even-high-substitution} with $\sigma=1$ and $V_h$ in place
of $C_h$. The last equation has no mass term, and the angular Poincar\'e
inequality on $|m|\geq4$ gives the coercivity in
\eqref{appref:even-high-coercivity}. We recover $X_h=(-i\cR)V_h$ and
finish with
\eqref{appref:even-final-affine-equations}--\eqref{appref:even-reconstruction},
omitting the terms $t^2C$. The mean of $A$ remains $-f_{0,0}/2$, the
means of $B$, $C$, and $S$ keep their prescribed values, and the last equation
imposes no condition on the means because it carries the factor $I-\Pi$.

In summary, the construction of the even inverse $K_e^{\rm ev,aug}$
proceeds in the order
\begin{equation}\label{appref:even-triangular-order}
 \text{mean/scale}\longrightarrow C_{\pm2}
 \longrightarrow C_h\longrightarrow(Z,W)_{\rm nonsingular}
 \longrightarrow(W_2,Z_{-2}).
\end{equation}
Substitution gives $\mathscr L_e^{\rm ev,aug}K_e^{\rm ev,aug}=I$. For the
other identity, we apply the same steps to an arbitrary element of the
domain. Uniqueness for $J_+$, the coercivity
\eqref{appref:even-high-coercivity}, the equations
\eqref{appref:even-nonsingular-spins}, and uniqueness for $J_-$ with
zero quadratic coefficient recover each component in turn, so that
$K_e^{\rm ev,aug}\mathscr L_e^{\rm ev,aug}=I$ as well.

\smallskip\noindent\textit{Estimates for the even inverse.}
We let $\mathcal X_{{\rm ev},y,n}^j$,
$\mathcal Y_{{\rm ev},y,n}^j$, and
$\mathcal T_{{\rm ev},y,n}^j$ be the rescaled Cartesian
domain, source, and balanced trace norms, as in
\eqref{appref:carrier-mode-norms}. We write
$K_{e,n}^{\rm ev,aug}$ for the even inverse at cell
frequency $n$, and set
\begin{equation}\label{appref:even-component-energies}
 \begin{aligned}
 E_{j,n}^{{\rm ev},(q)}
 &:=
 \norm{\omega_n\pa_\rho^qK_{e,n}^{\rm ev,aug}f_n}_{
        \mathcal X_{{\rm ev},y,n}^j}\\
 &\quad+
 \norm{\omega_n(T_n)\operatorname{Tr}_{\rm phys}
          \pa_\rho^qK_{e,n}^{\rm ev,aug}f_n}_{
        \mathcal T_{{\rm ev},y,n}^j},\\
 F_{a,n}^{\rm ev}
 &:=
 \norm{\omega_nf_n}_{\mathcal Y_{{\rm ev},y,n}^a}.
 \end{aligned}
\end{equation}
The weighted Green estimate
\eqref{appref:even-weighted-green-trace}, the conjugated
Hardy identities \eqref{appref:weighted-even-hardy-identities}
and their order-zero bounds, and the reconstruction formulas
\eqref{appref:even-nonsingular-spins},
\eqref{appref:even-final-hardy}, and \eqref{appref:even-reconstruction} give
\begin{equation}\label{appref:even-grade-recurrence}
 E_{j,n}^{{\rm ev},(0)}
 \leq C_j\bigl(F_{j+3,n}^{\rm ev}
       +\lambda_nE_{j-1,n}^{{\rm ev},(0)}\bigr),
 \qquad j\geq1,
 \qquad
 E_{0,n}^{{\rm ev},(0)}\leq C F_{3,n}^{\rm ev}.
\end{equation}
Here, the passage from the source $f$ to the scalar right-hand sides uses
two derivatives, and the mass or balanced trace estimate uses one more,
while every other map in
\eqref{appref:even-triangular-order} loses no derivatives. Iterating
the recurrence, we obtain
\begin{equation}\label{appref:even-component-bound}
 E_{j,n}^{{\rm ev},(0)}
 \leq C_j\sum_{a=3}^{j+3}
       \lambda_n^{j+3-a}F_{a,n}^{\rm ev}.
\end{equation}

Each derivative in $\rho$ requires at most one additional derivative of
the source. Indeed, $J_{\pm,n}^\gamma$ is independent of $\rho$, while
the coefficients in \eqref{appref:even-exceptional-tangent},
\eqref{appref:even-high-substitution},
\eqref{appref:metric-spin-identities}, and
\eqref{appref:even-final-affine-equations} are analytic for
$0\leq\rho\leq\rho_0$. Differentiating the variational identity and then
the triangular reconstruction formulas, we obtain, for $q\geq1$,
\begin{equation}\label{appref:even-rho-grade-recurrence}
 \begin{aligned}
 E_{0,n}^{{\rm ev},(q)}
 &\leq C_{q,\gamma}\left(
 F_{3+q,n}^{\rm ev}
 +\lambda_n\sum_{p<q}E_{0,n}^{{\rm ev},(p)}\right),\\
 E_{j,n}^{{\rm ev},(q)}
 &\leq C_{j,q,\gamma}\left(
 F_{j+3+q,n}^{\rm ev}
 +\lambda_nE_{j-1,n}^{{\rm ev},(q)}
 +\lambda_n\sum_{p<q}E_{j,n}^{{\rm ev},(p)}\right),
 \qquad j\geq1,\ q\geq1.
 \end{aligned}
\end{equation}
The last term accounts for the factor $\lambda_n$ produced by each
differentiated coefficient or trace, and a double induction in $(q,j)$ then gives
\begin{equation}\label{appref:even-rho-component-bound}
 E_{j,n}^{{\rm ev},(q)}
 \leq C_{j,q,\gamma}
   \sum_{a=0}^{j+3+q}
       \lambda_n^{j+3+q-a}F_{a,n}^{\rm ev}.
\end{equation}

We sum these bounds over the cell frequencies as in the odd case. The
even norms have the same form as \eqref{appref:carrier-mode-norms}.
Restoring \(e^{\Phi_\gamma(0,\lambda_n)}\) and multiplying
\eqref{appref:even-component-bound} by $\lambda_n^\ell$, with
$j+\ell\leq s$, gives terms with at most $s+3$ derivatives of the source, so the
Cauchy-Schwarz inequality in the finite sums and Parseval's identity
give
\begin{equation}\label{appref:even-high-square-sum}
 \sum_{|n|\geq2L}
 \norm{K_{e,n}^{\rm ev,aug}f_n}_{\cX_{\gamma,n}^s}^2
 \leq C_s\sum_{|n|\geq2L}
 \norm{f_n}_{\cY_{\gamma,n}^{s+3}}^2.
\end{equation}
Likewise, \eqref{appref:even-rho-component-bound} gives
\begin{equation}\label{appref:even-rho-high-square-sum}
 \sum_{|n|\geq2L}
 \norm{\pa_\rho^qK_{e,n}^{\rm ev,aug}f_n}_{
             \cX_{\gamma,n}^s}^2
 \leq C_{s,q,\gamma}\sum_{|n|\geq2L}
 \norm{f_n}_{\cY_{\gamma,n}^{s+3+q}}^2.
\end{equation}
On the finite set $|n|<2L$, which includes $n=0$, the weights
$\omega_n^{\pm1}$ are uniformly bounded on the cell, and at $n=0$
we use the weight in the original variables instead. Applying
\eqref{appref:bounded-carrier-inverse} to each component, taking the
maximum over this finite set, and decreasing the same $\rho_0$ once
more, we obtain both estimates on these cell frequencies as well. Thus,
we have the bounds
\begin{equation}\label{appref:even-global-rho-bound}
 \begin{aligned}
 \norm{K_e^{\rm ev,aug}f}_{\cX_\gamma^s}
 &\leq C_{s,\gamma}\norm{f}_{\cY_\gamma^{s+3}},\\
 \norm{\pa_\rho^qK_e^{\rm ev,aug}f}_{\cX_\gamma^s}
 &\leq C_{s,q,\gamma}\norm{f}_{\cY_\gamma^{s+3+q}}.
 \end{aligned}
\end{equation}
Hence, the even block loses at most $\mu_{\rm ev}=3$
derivatives of the source at $q=0$, and at most $3+q$ for its $q$-th
derivative in $\rho$, independently of the cutoffs.

\subsection{Returning to the Cartesian operator}

We can now complete the proof of \cref{appref:raw-constant-ellipse} by
assembling the odd and even inverses and returning to the Cartesian
variables through the transfer maps of \cref{sec:block}.

\begin{proof}[Proof of \cref{appref:raw-constant-ellipse}]
We let $\mathbf b_e$ be the constant reference associated with
$M_e$, and let $\mathscr S_X$ and $\mathscr S_Y$ be the
restrictions of $\mathsf J_{X,\mathbf b_e}$ and
$\mathsf J_{Y,\mathbf b_e}^{\rm raw}$ from
\eqref{eq:A-JX-factorization} and \eqref{eq:A-JY-factorization} to the
scalar variables and to the even and odd angular modes used here. We let $\iota_Y$
adjoin zero affine scale and auxiliary outer trace data, and we let
$p_Y$ remove those entries. Writing these auxiliary data as $z$, we have
$\iota_YF=(F,0)$ and $p_Y(F,z)=F$, so that
\begin{equation}\label{appref:auxiliary-range-splitting}
 p_Y\iota_Y=I,\qquad
 \iota_Yp_Y=I\quad\text{on the zero-auxiliary subspace}.
\end{equation}
We put
\[
 \widetilde{\mathscr L}_e
 :=\mathscr L_e^{\rm od}\oplus\mathscr L_e^{\rm ev,aug},
 \qquad
 K_e^{\rm aug}:=K_e^{\rm od}\oplus K_e^{\rm ev,aug},
\]
and on the domain $\cX_{M_e}^\infty$, the identity
\eqref{eq:A-stabilized-intertwining} takes the form
\begin{equation}\label{appref:exact-adapter-intertwining}
 \widetilde{\mathscr L}_e\mathscr S_Xu
   =\iota_Y\mathscr S_YL_eu.
\end{equation}
Conversely, \eqref{appref:scale-row} and the condition on the outer boundary
show that $K_e^{\rm aug}\iota_Y\mathscr S_Yf$ has zero affine scale
coordinate and zero auxiliary outer trace, so it belongs to
$\mathscr S_X(\cX_{M_e})$, and we define the Cartesian inverse by
\begin{equation}\label{appref:literal-raw-inverse}
 V_e=\mathscr S_X^{-1}
   K_e^{\rm aug}\iota_Y\mathscr S_Y.
\end{equation}
Indeed, using \eqref{appref:exact-adapter-intertwining} and the two
inverse identities for $K_e^{\rm aug}$, we obtain
\[
 \begin{aligned}
 \iota_Y\mathscr S_YL_eV_e
  &=\widetilde{\mathscr L}_eK_e^{\rm aug}
       \iota_Y\mathscr S_Y=\iota_Y\mathscr S_Y,\\
 V_eL_eu
  &=\mathscr S_X^{-1}K_e^{\rm aug}
       \widetilde{\mathscr L}_e\mathscr S_Xu=u.
 \end{aligned}
\]
Injectivity of $\iota_Y\mathscr S_Y$ gives $L_eV_e=I$ from the first
identity, while the second gives $V_eL_e=I$.

We now count the derivatives lost in composing these maps. By
\eqref{eq:A-JX-type}--\eqref{eq:A-JY-type}, the maps on the domain lose
no derivatives, while the map on the source loses $\ell_Y$. The odd and
even inverses act on complementary summands and lose at most four and
three derivatives, respectively, so their direct sum loses at most four.
Writing $\ell_{\rm alg}$ for the finite shift in the Sobolev index of
the maps to and from the scalar coordinates, noted before
\eqref{appref:raw-scalar-rows}, and of the algebraic equations, we may take
\[
 \ell_e=\ell_Y+4+\ell_{\rm alg}<\infty.
\]
For $q$ derivatives in $\rho$, \eqref{appref:odd-global-rho-bound} and
\eqref{appref:even-global-rho-bound} require at most $4+q$ derivatives
of the source in the odd block and $3+q$ in the even block, while the
derivatives of the coefficients in the finite transfer maps lose none.
Leibniz' rule therefore gives
\begin{equation}\label{appref:raw-rho-ledger}
 \norm{\pa_\rho^qV_ef}_{\cX_\gamma^s}
 \leq C_{s,q,\gamma}
 \norm{f}_{\cY_\gamma^{
       s+\ell_Y+4+\ell_{\rm alg}+q}}.
\end{equation}
The analytic weight is the same on both sides. For $q$
derivatives in $\rho$, the constant includes the factor
$\gamma^{-q}$ from \eqref{appref:center-polynomial-absorption}.

All the remaining factors in
\eqref{eq:A-JX-factorization}--\eqref{eq:A-JY-inverse-factorization} lose
no derivatives, so at $q=0$ we recover \eqref{appref:raw-estimate}.
Restoring the orientation amounts to conjugating by
$\mathsf R_{\alpha_0}$, which is unitary on every Cartesian scale.

\end{proof}

\subsection{The case of circular cross-sections}

The even and affine inverse extends to the circle with
uniform bounds on $0\le\rho\le\rho_0$. This requires
coercivity on the high even modes and reconstruction of
$W_2$ and $Z_{-2}$ without division by $\rho$.
The resulting uniformity lets us choose the constants
before fixing $\rho_*$.

The matrix through which the ellipse enters the even block is the
matrix $H_\rho$ introduced before \eqref{appref:raw-scalar-rows},
\begin{equation}\label{appref:metric-multiplier}
 H_\rho(\theta)=\frac1{1-\rho^2}
 \begin{pmatrix}
 1-\rho\cos2\theta&\rho\sin2\theta\\
 \rho\sin2\theta&1+\rho\cos2\theta
 \end{pmatrix}.
\end{equation}
We write $K_\rho^{\rm ev}$ for the even inverse constructed
above, now at ellipse parameter $\rho$ instead of $\rho_*$, and
$\mathfrak X^s$, $\mathfrak Y^s$, and $\mathfrak T^s$ for its Cartesian
domain, its source space, and its space of balanced traces on the outer
boundary. We also write $\bX_{{\rm ev},\natg}^s$ and
$\bY_{{\rm ev},\natg}^s$ for the corresponding block spaces of
\cref{sec:block}. Note that the norms of these block spaces include the
prescribed orders of the cap value, of the conormal trace, and of the
mismatch at the interface. The uniformity of the even and affine inverse
up to the circle is the content of the following lemma.

\begin{lemma}
\label{appref:even-circle-endpoint}
Fix $0<\rho_0<1/4$ and $L>0$, together with the mean and scale equations
and the chart with prescribed transverse first derivatives on the axis. For every
$0\leq\rho\leq\rho_0$, the circle included, the even inverse
$K_\rho^{\rm ev}$ exists, it is a two-sided inverse of the even
operator, and on the spaces $\mathfrak X^s$, $\mathfrak Y^{s+3}$, and
$\mathfrak T^s$ above it satisfies
\begin{equation}\label{appref:even-raw-bound}
 \norm{K_\rho^{\rm ev}f}_{\mathfrak X^s}
 +\norm{\operatorname{Tr}_{\rm phys} K_\rho^{\rm ev}f}_{\mathfrak T^s}
 \leq C_{s,\rho_0,L}\norm{f}_{\mathfrak Y^{s+3}}.
\end{equation}
Here, $\operatorname{Tr}_{\rm phys}$ denotes the balanced value and
conormal traces on the outer boundary, as in
\eqref{appref:even-weighted-green-trace}. In
the block norms, the same inverse loses no derivatives,
\begin{equation}\label{appref:even-native-bound}
 K_{\rho,\natg}^{\rm ev}:
 \bY_{{\rm ev},\natg}^s\longrightarrow
 \bX_{{\rm ev},\natg}^s,
 \qquad
 \sup_{0\leq\rho\leq\rho_0}
 \norm{K_{\rho,\natg}^{\rm ev}}_s<\infty.
\end{equation}
All $\rho$-derivatives of fixed order obey the analogous uniform tame
bounds.
\end{lemma}

\begin{proof}
The numerator matrix in \eqref{appref:metric-multiplier} has trace $2$
and determinant $1-\rho^2$, so its eigenvalues are $1\pm\rho$, and those
of $H_\rho$ are $(1\mp\rho)^{-1}$. We therefore have
\begin{equation}\label{appref:metric-bounds}
 \frac1{1+\rho_0}I\leq H_\rho\leq\frac1{1-\rho_0}I,
\end{equation}
and the form on the high even modes then satisfies
\[
 \begin{aligned}
 \operatorname{Re}\mathfrak b_\rho(C_h,C_h)
 &\geq \frac{1-4\rho_0}{1-\rho_0}
       \norm{(DC_h,\cR C_h)}_2^2\\
 &\quad+\text{nonnegative mass and Robin terms}.
 \end{aligned}
\]
Indeed, the coefficient
$1-3\rho/(1-\rho)=(1-4\rho)/(1-\rho)$ decreases from $1$
at the circle to a positive value at $\rho_0<1/4$.
Thus, the Lax-Milgram bound, the outer trace extension,
and the value and conormal estimates are uniform.
For $n=0$, the angular Poincar\'e inequality on $|m|\ge4$
supplies coercivity without the mass term.

On the affine modes, the formulas \eqref{appref:even-exceptional-tangent}
and \eqref{appref:even-final-affine-equations} contain only the factors
\[
 1-\rho^2,\qquad \rho,\qquad (1-\rho^2)^{-1},
\]
so there is no division by $\rho$. The Hardy maps
\eqref{appref:even-hardy-maps} are independent of $\rho$. The prescribed
first derivatives of the original variation on the axis exclude the
homogeneous quadratic solutions for \(W_2\) and \(Z_{-2}\), while
Cartesian smoothness excludes the singular solutions of \(D+2\).
The scale equation remains \eqref{appref:scale-row}. Therefore, no kernel appears
on the affine modes at $\rho=0$.

We obtain the estimates at higher spatial regularity from uniform elliptic difference
quotients. Of the three derivatives of the source in
\eqref{appref:even-raw-bound}, two are used to obtain the scalar right-hand
sides, and one more in the mass or balanced trace estimate, as counted in
\eqref{appref:even-grade-recurrence}. These orders are already part of
the block norms. Each Green, Hardy, cap, and trace map is therefore
bounded without loss, and we obtain \eqref{appref:even-native-bound}.

Finally, we write $\mathscr L_\rho^{\rm ev}$ for the even operator
$\mathscr L_e^{\rm ev,aug}$ at ellipse parameter $\rho$. This family is
analytic in $\rho$ on the scale, and differentiating the identity
$K_\rho^{\rm ev}\mathscr L_\rho^{\rm ev}=I$ gives
\[
 \pa_\rho K_\rho^{\rm ev}
 =-K_\rho^{\rm ev}(\pa_\rho\mathscr L_\rho^{\rm ev})
       K_\rho^{\rm ev},
\]
and iterating this formula gives the bounds at every fixed order. A
constant change of orientation acts by orthogonal changes of the
coordinates in the disk and in the normal plane. These changes preserve the
radial weight and the angular and cell cutoffs. At $\rho=0$, the operator
is independent of the orientation.
\end{proof}
\section{The center and affine equations}\label{app:minus}\label{app:minus:section}

We now construct an inverse for the center and affine equations at high
cell frequencies. We first solve on the cap with the prescribed axis
data, then use the resulting value and conormal derivative as initial
data for the outward evolution on the annulus. No condition is imposed
on these components at the outer boundary. The relation between this
inverse and the minus diagonal block of
\eqref{app:global:block-decomposition} is stated in
\cref{prop:minus-diagonal}.

Throughout, $\Lambda=\langle D_\zeta\rangle$ multiplies the cell
mode $e^{in\zeta}$ by $\lambda_n=\langle n\rangle$, and we fix the analytic
widths $\sigma_0>\sigma_->0$, the outer radius $R$ of the annulus, the
radius $2T$ of the cap in the rescaled variable, and the base index
$s_*$. Of these, $\sigma_0$ is the width of the strip of analyticity of
\cref{subsec:analytic-scale} on the axis, and $\sigma_-$ is the width
that remains at the outer boundary. All constants are uniform on the
compact parameter set of \cref{sec:tame-problem}. We choose first
\(0<\gamma<\min\{\sigma_0,1\}\), then \(\delta_A>0\), so that
\begin{equation}\label{app:minus:eq:phase-choices}
 \gamma R<\sigma_0-\sigma_-,\qquad
 C_{\rm wt}(\gamma+\gamma^2)<\frac{c_{\rm wt}}4,\qquad
 C_{\rm cap}\gamma<\frac18,\qquad
 C_*\sqrt{\delta_A}<\gamma .
\end{equation}
% AUTHOR QUERY: Define C_{\rm wt}, c_{\rm wt} in
% \eqref{app:minus:eq:phase-choices} and identify the weighted estimate
% controlled by the second inequality. Neither constant is used elsewhere.
We then choose an upper bound for the ellipse parameter
\begin{equation}\label{app:minus:eq:rho-collar}
 0<\rho_0<\min\left\{\frac14,\rho_{\rm ax},
                      \rho_{\rm ctr},\rho_{\rm chart}\right\}.
\end{equation}
The bounds $\rho_{\rm ax}$, $\rho_{\rm ctr}$, and
$\rho_{\rm chart}$ ensure, respectively, cap invertibility,
separation of the center roots from the other roots, and
validity of the fixed-space identifications in \cref{sec:block}.
We then bound the distance from each state in $\mathfrak U_-$
to its constant ellipse by $\eta_*>0$, and choose the dyadic cutoff $J$ of
\eqref{eq:A-exact-cell-split}, so that
\begin{equation}\label{app:minus:eq:state-cutoff-choices}
 \begin{gathered}
 C_A\eta_*\le M_B\delta_A,\qquad C_c\eta_*<\gamma,\\
 C_{\rm cap}(2^{-J}+\gamma+\eta_*)<\frac12,\qquad
 C_{\rm syl}\bigl(\Lambda_J^{-1/2}+\chi_-(\eta_*)\bigr)<\frac12,\\
 C_{\mathbb N}\eta_*<\frac12,
 \end{gathered}
\end{equation}
where $\Lambda_J=\inf\{\langle n\rangle:|n|\ge2^J\}\simeq2^J$ and
$\chi_-(t)\to0$ as $t\to0$. The first two conditions control the perturbations of the affine and
center energies. The next two make the cap error and the coupling
correction small, and the last permits the final Neumann inversion.
These choices are independent of an upper Fourier cutoff
and of the number $N$ of periods of the field.

For a state $\mathbf b=(\rho,a,\eps,p)$ and the constant ellipse
$\mathbf b_\rho=(\rho,a_\rho,\eps_\rho,p_\rho)$ of \cref{sec:block},
where $a_\rho=a^0_{\rho,p_e}$, $\eps_\rho=0$, and $p_\rho=p_e$, we put
\begin{equation}\label{app:minus:eq:parameter-norm}
 \norm{\mathbf b-\mathbf b_\rho}_{\cB^s}
 :=\norm{a-a_\rho}_{\cA_\gamma^s}
       +\abs{\eps-\eps_\rho}+\abs{p-p_\rho},
\end{equation}
and we let $\mathfrak U_-$ consist of the states with $0\le\rho\le\rho_0$
for which \eqref{app:minus:eq:parameter-norm}, taken at the index
$s_*+k_c$, is less than $\eta_*$. Here, $k_c$ is a fixed number of
additional coefficient derivatives, chosen large enough for the cap
estimates and independent of the solution index $s$. At $\rho=0$, derivatives in $\rho$ are defined by
extending the coefficients analytically to a
neighborhood of the endpoint. The orientation $\alpha_0$ is a separate
parameter, and the estimates are used only for $\rho\ge0$.

\subsection{Function spaces and the inverse for the minus block}

We use the fixed block spaces and identifications of \cref{sec:block}.
We now specify the spaces for the cap solution, its Cauchy data, and the
evolution on the annulus. For the cap, we define
\begin{equation}\label{app:minus:eq:cap-spaces}
 \begin{aligned}
 X_{c,-}^s&=H_J^X\mathscr X_{-,\mathrm{cap}}^s,\\
 Y_{c,-}^s&=H_J^Y\bigl(
       \mathscr Y_{-,\mathrm{bulk}}^s
       \oplus\mathfrak J_{-,\mathrm{ax}}^s\bigr).
 \end{aligned}
\end{equation}
Here, $\mathscr X_{-,\mathrm{cap}}^s$ is the Cartesian cap
domain, $\mathscr Y_{-,\mathrm{bulk}}^s$ is its interior
source space, and $\mathfrak J_{-,\mathrm{ax}}^s$ contains
the prescribed axis data. The interior components have
$H^{s+2}$ regularity in $X_{c,-}^s$ and $H^s$ regularity
in the source, as in \cref{tab:block-weights}.
The auxiliary data $\mathfrak J_{-,\mathrm{ax}}^s$
prescribe the axis coefficients of the regular center
solutions and the admissible linear Taylor coefficients.
In the normalized chart, the transverse first derivatives
are fixed, while the tangential coefficient $\tau$ remains
an unknown.
% AUTHOR QUERY: State explicitly how the prescribed scalar center jets
% encode the axis constraints on the original vector variation and how
% the unknown tangential coefficient tau is recovered. The full linear
% Taylor term is not fixed by the normalized chart.
The range $Y_{c,-}^s$ contains no interface datum.

We measure the outgoing Cauchy data of the cap solution in
\begin{equation}\label{app:minus:eq:cap-trace-space}
 Z_{c,-}^s=H_J^{Z_c}\mathfrak T_{c,-,\gamma}^s,\qquad
 \norm{(z_D,z_N)}_{Z_{c,-}^s}
 \simeq \norm{\Lambda^{1/2}z_D}_{H^s}
       +\norm{\Lambda^{-1/2}z_N}_{H^s},
\end{equation}
where $\mathfrak T_{c,-,\gamma}^s$ is the space of Cauchy data on the
interface of the cap, projected by the spectral projections of the
principal symbol and weighted by the analytic weight. The half-order
factors are the balanced trace norms of \cref{sec:block}. We write
$Z_{c,-}^s=Z_{c,D}^s\oplus Z_{c,N}^s$ for the
splitting into the value and the conormal component.

For the inverse, we add to the physical conormal derivative
a boundary contribution determined by the interior source.
The value is unchanged. We call this pair the
\emph{split Cauchy data}. To define the contribution, let
\[
 E_{c,\Sigma}:(Z_{c,N}^s)'\longrightarrow
                    (\mathscr Y_{-,\mathrm{bulk}}^s)'
\]
be the extension of the trace on the cap built from
\eqref{eq:A-Cauchy-model-potentials}. For
$f_c=(f_{c,\mathrm{bulk}},j_{\mathrm{ax}})\in Y_{c,-}^s$, we define its boundary contribution by transposing this extension,
\begin{equation}\label{app:minus:eq:cap-source-moment}
 \left\langle\mathsf M_{c,\Sigma}^sf_c,\varphi\right\rangle
 :=\left\langle f_{c,\mathrm{bulk}},E_{c,\Sigma}\varphi\right\rangle,
 \qquad
 \mathsf M_{c,\Sigma}^s:Y_{c,-}^s\longrightarrow Z_{c,N}^s,
\end{equation}
where the transpose uses the weighted source-dual pairing.
It is not the adjoint for the unweighted Hilbert inner
product. We let $j_Nz=(0,z)$ include the conormal component
in $Z_{c,-}^s$. The \emph{physical Cauchy data}, consisting
of the value and classical conormal derivative at the
state, are related to the split data by
\begin{equation}\label{app:minus:eq:cap-trace-physicalization}
 \begin{aligned}
 \Gamma_{c,-,\mathbf b}^{\mathrm{sp}}(u_c,f_c)
   &:=\Gamma_{c,-,\mathbf b}^{\mathrm{ph}}u_c
       +j_N\mathsf M_{c,\Sigma}^sf_c,\\
 \Gamma_{c,-,\mathbf b}^{\mathrm{ph}}u_c
   &=\Gamma_{c,-,\mathbf b}^{\mathrm{sp}}(u_c,f_c)
       -j_N\mathsf M_{c,\Sigma}^sf_c .
 \end{aligned}
\end{equation}
When $f_c=\cA_{c,-,\mathbf b}u_c$, where $\cA_{c,-,\mathbf b}$ is the
operator on the cap at the state, we abbreviate the first expression to
$\Gamma_{c,-,\mathbf b}^{\mathrm{sp}}u_c$. Note that $\mathsf M_{c,\Sigma}^s$ is bounded at each regularity level
$s$, and its contribution affects only the conormal component, which
has lower order in the balanced trace space. It therefore changes
neither the principal symbol nor the estimates for the homogeneous
problem on the cap.
% AUTHOR QUERY: Justify the lower-order assertion when f_c=A_c u_c.
% Boundedness of M_c from the source space to a conormal trace space
% alone does not show that M_c A_c is lower order on the unknown.
When no superscript is given, we mean the physical
Cauchy data $\Gamma_{c,-,\mathbf b}^{\mathrm{ph}}$.

On the annulus, we set
\begin{equation}\label{app:minus:eq:annulus-spaces}
 \begin{aligned}
 X_{a,-}^s&=H_J^X\bigl(X_{C,\gamma}^s\oplus X_{A,\Phi}^s
                         \oplus\mathfrak F_{-,X}^s\bigr),\\
 Y_{a,-}^s&=H_J^Y\bigl(Y_{C,\gamma}^s\oplus Y_{A,\Phi}^s
                         \oplus\mathfrak F_{-,Y}^s\bigr).
 \end{aligned}
\end{equation}
Here, the subscripts $C$ and $A$ denote the center and
affine components. The remaining mean, scale, gauge, and quotient
equations of \cref{sec:block} form a triangular system. Their unknowns
belong to $\mathfrak F_{-,X}^s$ and their prescribed data to
$\mathfrak F_{-,Y}^s$. The center unknown consists of the two oscillator variables and
their radial derivatives, as in \eqref{appref:odd-center-oscillator}.
Its norm is defined by an energy on these four coordinates that remains
uniform when the two speeds coincide. Its source norm is of type $L^1$ in $r$
with values in $H^s$. For the affine component, we collect the upper entries of the two
chains as $x=(x_1,x_2)$ and the lower entries as $y=(y_1,y_2)$,
as described after \cref{tab:block-weights}. We write $f=(f_1,f_2)$
and $g=(g_1,g_2)$ for the sources in their respective equations and use
\begin{equation}\label{app:minus:eq:affine-native-norm}
 \begin{aligned}
 \norm{(x,y)}_{X_{A,\Phi}^s}
  &=\norm{(\Lambda^s\mu x,\Lambda^{s+1/2}y)}_\Phi,\\
 \norm{(f,g)}_{Y_{A,\Phi}^s}
  &=\norm{(\Lambda^s\mu f,\Lambda^{s+1/2}g)}_\Phi,
 \qquad \mu^2=\delta_A+\Lambda^{-1}.
 \end{aligned}
\end{equation}
% AUTHOR QUERY: Specify the radial norm in ||.||_Phi for both X_A and
% Y_A. An initial-value solution norm and an L^1 radial source norm are
% different; weighting alone does not define these spaces or their traces.
Here, $\norm{\cdot}_\Phi$ denotes the norm on the annulus after
multiplication by the weight $e^{\Phi_\gamma(r,\Lambda)}$, and the
factors $\mu$ and $\Lambda^{1/2}$ are the weights in
\cref{tab:block-weights} of the upper and lower entries of the two
affine chains. The weight $e^{\Phi_\gamma(r,\Lambda)}$ is the
analytic weight $W_\gamma$ of \eqref{eq:phase-weight}, with
\begin{equation}\label{app:minus:eq:minus-phase}
 \Phi_\gamma(r,\lambda)
 =\sigma_0\lambda-
   \gamma\bigl(\sqrt{1+r^2\lambda^2}-1\bigr).
\end{equation}
We write $Z_{a,-}^s$ and $Z_{R,-}^s$ for the spaces of
Cauchy data on the inner and on the outer boundary of the annulus,
balanced as in \eqref{app:minus:eq:cap-trace-space}, and
\begin{equation}\label{app:minus:eq:trace-converter}
 T_{-,\mathbf b}:Z_{c,-}^s\longrightarrow Z_{a,-}^s
\end{equation}
for the order-zero map from cap to annular trace coordinates.
For a field $u$ smooth across the interface, its defining relation is
\[
 T_{-,\mathbf b}\Gamma_{c,-,\mathbf b}^{\rm ph}u
   =\gamma_{-,\mathbf b}u,
\]
where $\gamma_{-,\mathbf b}$ records the same value and conormal data
in the annular coordinates. To obtain the physical mismatch, we subtract
$\mathsf M_{c,\Sigma}^sf_c$ from the split conormal datum
before applying $T_{-,\mathbf b}$. Keeping this contribution
instead gives the split mismatch. The full domain and
range spaces are
\begin{equation}\label{app:minus:eq:retained-spaces}
 \bX_{-,\mathrm{ret}}^s=X_{c,-}^s\oplus X_{a,-}^s,
 \qquad
 \bY_{-,\mathrm{ret}}^s=Y_{c,-}^s\oplus Y_{a,-}^s
                                   \oplus Z_{a,-}^s.
\end{equation}
The physical and the split mismatches are related on this range by
\begin{equation}\label{app:minus:eq:split-mismatch-shear}
 \begin{aligned}
 \mathfrak P_{-,\mathbf b}^s(f_c,f_a,g^{\mathrm{ph}})
   &=\bigl(f_c,f_a,
       g^{\mathrm{ph}}-T_{-,\mathbf b}j_N
                    \mathsf M_{c,\Sigma}^sf_c\bigr),\\
 (\mathfrak P_{-,\mathbf b}^s)^{-1}
       (f_c,f_a,g^{\mathrm{sp}})
   &=\bigl(f_c,f_a,
       g^{\mathrm{sp}}+T_{-,\mathbf b}j_N
                    \mathsf M_{c,\Sigma}^sf_c\bigr).
 \end{aligned}
\end{equation}
The map subtracts the cap source contribution from the
physical mismatch, and its inverse adds it back. Both are
bounded at index $s$ and leave the highest-order terms
unchanged because this contribution has lower order.
In particular,
\begin{equation}\label{app:minus:eq:current-homogeneous-physical-datum}
 g^{\rm ph}=0
 \quad\Longrightarrow\quad
 g^{\rm sp}=-T_{-,\mathbf b}j_N
                    \mathsf M_{c,\Sigma}^sf_c,
\end{equation}
which is in general nonzero when the interior source on the cap is
nonzero.

\smallskip\noindent\textit{Construction of the inverse.}
We construct the inverse by solving first on the cap and
then on the annulus. The component estimates are proved
in the next three subsections and collected in
\cref{prop:minus-diagonal}. At the reference state, the
cap operator and its outgoing physical Cauchy data are
\begin{equation}\label{app:minus:eq:cap-map}
 \cA_{c,-,\rho}^{\rm diag}:X_{c,-}^s\longrightarrow Y_{c,-}^s,
 \qquad
 \Gamma_{c,-,\rho}^{\rm ph}G_{c,-,\rho}^{\rm diag}:
        Y_{c,-}^s\longrightarrow Z_{c,-}^s,
\end{equation}
with inverse $G_{c,-,\rho}^{\rm diag}$ constructed below.
Adding the cap source contribution gives the split trace
\begin{equation}\label{app:minus:eq:split-cap-solution-trace}
 \Gamma_{c,-,\rho}^{\rm sp}G_{c,-,\rho}^{\rm diag}f_c
 :=\Gamma_{c,-,\rho}^{\rm ph}G_{c,-,\rho}^{\rm diag}f_c
      +j_N\mathsf M_{c,\Sigma}^sf_c .
\end{equation}
We let $\cL_{a,-,\mathbf b_\rho}$ be the reference center
and affine operator on the annulus, and let $\gamma_{-,\rho}$
be its inner Cauchy trace. On each shell,
\begin{equation}\label{app:minus:eq:annulus-map}
 \cB_{a,-,\rho}:=(\cL_{a,-,\mathbf b_\rho},
          \gamma_{-,\rho}):X_{a,-}^s
       \longrightarrow Y_{a,-}^s\oplus Z_{a,-}^s
\end{equation}
records the interior equation and the incoming Cauchy data. Its inverse
$G_{a,-,\rho}$ solves for arbitrary prescribed values of these data,
with no condition at the outer boundary. The
outgoing Cauchy data $\Gamma_{R,-,\rho}G_{a,-,\rho}$ on the outer
boundary are bounded into $Z_{R,-}^s$.

We define the \emph{coupled reference operator} by
\begin{equation}\label{app:minus:eq:joined-operator}
 \begin{aligned}
 N_{-,\rho}^{\rm diag}(u_c,u_a)=\big(&
          \cA_{c,-,\rho}^{\rm diag}u_c,
          \cL_{a,-,\mathbf b_\rho}u_a,\gamma_{-,\rho}u_a
          -T_{-,\mathbf b_\rho}
             \Gamma_{c,-,\rho}^{\rm sp}
              (u_c,\cA_{c,-,\rho}^{\rm diag}u_c)\big).
 \end{aligned}
\end{equation}
We write $m_-^{\rm sp}$ for its last component and
$m_-^{\rm ph}=\gamma_{-,\rho}u_a-
T_{-,\mathbf b_\rho}\Gamma_{c,-,\rho}^{\rm ph}u_c$. Then,
\eqref{app:minus:eq:cap-trace-physicalization} gives the identity
\begin{equation}\label{app:minus:eq:split-physical-mismatch}
 m_-^{\rm sp}=m_-^{\rm ph}
   -T_{-,\mathbf b_\rho}j_N\mathsf M_{c,\Sigma}^s
          \cA_{c,-,\rho}^{\rm diag}u_c ,
\end{equation}
so that the physical Cauchy data of the two pieces agree if and only if
\begin{equation}\label{app:minus:eq:homogeneous-physical-split-datum}
 m_-^{\rm ph}=0
 \quad\Longleftrightarrow\quad
 m_-^{\rm sp}=-T_{-,\mathbf b_\rho}j_N
                   \mathsf M_{c,\Sigma}^sf_c .
\end{equation}
We allow the split mismatch to be prescribed independently.
Solving on the cap and using its trace as part of the
annular initial data gives the inverse
\begin{equation}\label{app:minus:eq:joined-inverse}
 \begin{aligned}
 G_{-,\rho}^{\mathbb N,\rm diag}(f_c,f_a,g^{\rm sp})
 =\Big(&G_{c,-,\rho}^{\rm diag}f_c,G_{a,-,\rho}\bigl(f_a,
    g^{\rm sp}+T_{-,\mathbf b_\rho}
       [\Gamma_{c,-,\rho}^{\rm ph}G_{c,-,\rho}^{\rm diag}f_c
          +j_N\mathsf M_{c,\Sigma}^sf_c]\bigr)\Big).
 \end{aligned}
\end{equation}

\smallskip\noindent\textit{Perturbing the coupled reference operator.}
We let $\iota_-^X$ and $\pi_-^Y$ be the domain inclusion
and range projection of the minus block. We take this
diagonal block in physical coordinates and then convert
its mismatch to split coordinates,

\begin{equation}\label{app:minus:eq:current-split-compression}
 N_{-,\mathbf b}^{\rm ph}:=\pi_-^Y\mathbb N_{\mathbf b}\iota_-^X,
 \qquad
 N_{-,\mathbf b}:=\mathfrak P_{-,\mathbf b}^s
                         N_{-,\mathbf b}^{\rm ph}.
\end{equation}
Here, $N_{-,\mathbf b}^{\rm ph}$ is the minus diagonal
block of \eqref{eq:A-natural-broken-formula} on all shells.
At $\mathbf b=\mathbf b_\rho$, its split form is
$N_{-,\rho}^{\rm diag}$. We denote the perturbation by
\begin{equation}\label{app:minus:eq:natural-minus-perturbation}
 \mathcal R_{-,\mathbf b}^{\mathbb N}
 :=N_{-,\mathbf b}-N_{-,\rho}^{\rm diag}.
\end{equation}
Besides the interior perturbation computed below, this
difference contains the change in conormal coordinates
and in the mismatch map \eqref{app:minus:eq:split-mismatch-shear}.
For the mismatch map, we have
\[
 \| (\mathfrak P_{-,\mathbf b}^s
       -\mathfrak P_{-,\mathbf b_\rho}^s)F\|_{
                         \bY_{-,\mathrm{ret}}^s}
 \le C_s\left[
   \delta_0(\mathbf b)\|F\|_{\bY_{-,\mathrm{ret}}^s}
  +\delta_s(\mathbf b)\|F\|_{
                         \bY_{-,\mathrm{ret}}^{s_*}}\right],
\]
by boundedness of $\mathsf M_{c,\Sigma}^s$ and of the change of trace
coordinates, and the inverse maps satisfy the same estimate. Here, using \eqref{app:minus:eq:parameter-norm}, we set
\[
 \delta_0(\mathbf b)=\|\mathbf b-\mathbf b_\rho\|_{\cB^{s_*+k_c}},
 \qquad
 \delta_s(\mathbf b)=\|\mathbf b-\mathbf b_\rho\|_{\cB^{s+k_c}}.
\]
These are the coefficient distances at the base and higher regularity
indices. We
exclude the terms from low to high and from high to low frequencies from
$\mathcal R_{-,\mathbf b}^{\mathbb N}$, since they enter the global
inverse separately, and we obtain
\begin{equation}\label{app:minus:eq:whole-minus-perturbation-bound}
 \|\mathcal R_{-,\mathbf b}^{\mathbb N}U\|_
       {\bY_{-,\mathrm{ret}}^s}
 \le C_s\left[
  \delta_0(\mathbf b)\|U\|_{\bX_{-,\mathrm{ret}}^s}
 +\delta_s(\mathbf b)\|U\|_{\bX_{-,\mathrm{ret}}^{s_*}}\right].
\end{equation}
Indeed, the reconstruction, the extraction of the coordinates, the
change of mismatch coordinates, and the identifications of the data spaces are bounded at every
fixed index, and the map $\mathcal T_{\mathbf b,\rho}^{CN}$ of
\eqref{eq:A-current-reference-row-map}, which converts the reference
Cauchy data to those at the state, satisfies
\[
 \|((\mathcal T_{\mathbf b,\rho}^{CN})^{-1}-I)z\|_{Z^s}
 \le C_s\bigl(\delta_0(\mathbf b)\|z\|_{Z^s}
              +\delta_s(\mathbf b)\|z\|_{Z^{s_*}}\bigr).
\]
Here, $Z^s$ is the product of the balanced trace spaces on the
interfaces. Combining this bound with the perturbation bound proved
below, we obtain \eqref{app:minus:eq:whole-minus-perturbation-bound} on
the product spaces. At the base index, the uniform bound for the reference inverse,
the base case of \eqref{app:minus:eq:whole-minus-perturbation-bound}, and
the last condition in \eqref{app:minus:eq:state-cutoff-choices} imply
\[
 \|G_{-,\rho}^{\mathbb N,\rm diag}
       \mathcal R_{-,\mathbf b}^{\mathbb N}\|<\frac12,
\]
and a Neumann series gives the inverse
\begin{equation}\label{app:minus:eq:whole-minus-neumann}
 \begin{aligned}
 G_{-,\mathbf b}^{\mathbb N}
  &:=
  \bigl(I+G_{-,\rho}^{\mathbb N,\rm diag}
               \mathcal R_{-,\mathbf b}^{\mathbb N}\bigr)^{-1}
       G_{-,\rho}^{\mathbb N,\rm diag}\\
  &=G_{-,\rho}^{\mathbb N,\rm diag}
    \bigl(I+\mathcal R_{-,\mathbf b}^{\mathbb N}
               G_{-,\rho}^{\mathbb N,\rm diag}\bigr)^{-1},
 \qquad
 N_{-,\mathbf b}G_{-,\mathbf b}^{\mathbb N}
 =G_{-,\mathbf b}^{\mathbb N}N_{-,\mathbf b}=I .
 \end{aligned}
\end{equation}
\smallskip\noindent\textit{The constraint correction.}
We let $\iota_-^YF_-=(F_-,0)$ include the minus source in the
high-frequency range, with zero low-frequency entries when it is viewed
in the full range. We let $\pi_-^Y$ select the minus component.
Taking the minus diagonal block of the map $\mathscr S_{Y,\mathbf b}$
of \cref{sec:block} gives
\begin{equation}\label{app:minus:eq:minus-diagonal-shear-compression}
 \mathscr S_{Y,-,\mathbf b}
 :=\pi_-^Y\mathscr S_{Y,\mathbf b}\iota_-^Y,
\end{equation}
and, for a range vector $F_-+\jmath_{Z,-,\mathbf b}z_-$ with $F_-$
having zero auxiliary entries, we have
\[
 \mathscr S_{Y,-,\mathbf b}(F_-+\jmath_{Z,-,\mathbf b}z_-)
 =F_--\pi_-^Y\mathbb B_{\mathbf b}z_-
       +\jmath_{Z,-,\mathbf b}z_-,
\]
the inverse being the same triangular formula with the sign reversed.
Here, $\jmath_{Z,-,\mathbf b}$ is the inclusion
$\jmath_{Z,\mathbf b}$ of \eqref{eq:A-auxiliary-inclusion-identities},
acting on the auxiliary coordinates of the minus block and compressed
to its range, and $z_-$ also denotes the
inclusion of such a coordinate in the full auxiliary space
$\mathscr Z_{{\rm aux},\natg}^s$. We put
\begin{equation}\label{app:minus:eq:minus-shear-realization}
 \widehat D_{-,\mathbf b}
 :=\mathscr S_{Y,-,\mathbf b}
       (\mathfrak P_{-,\mathbf b}^s)^{-1}N_{-,\mathbf b},
 \qquad
 \widehat G_{-,\mathbf b}
 :=G_{-,\mathbf b}^{\mathbb N}\mathfrak P_{-,\mathbf b}^s
       \mathscr S_{Y,-,\mathbf b}^{-1}.
\end{equation}
In the formula for $\widehat D_{-,\mathbf b}$, the range maps first
recover the physical mismatch and then subtract the interior source
$\pi_-^Y\mathbb B_{\mathbf b}z_-$ produced by extending the auxiliary
data $z_-$. This is the map \eqref{eq:A-coretraction-bulk-response}.
At the reference, we set
$\widehat D_{-,\rho}=\mathscr S_{Y,-,\mathbf b_\rho}
(\mathfrak P_{-,\mathbf b_\rho}^s)^{-1}N_{-,\rho}^{\rm diag}$ and
$\widehat G_{-,\rho}=G_{-,\rho}^{\mathbb N,\rm diag}
\mathfrak P_{-,\mathbf b_\rho}^s\mathscr S_{Y,-,\mathbf b_\rho}^{-1}$.
Multiplying out, we obtain
\begin{equation}\label{app:minus:eq:two-sided-joined}
 \widehat D_{-,\mathbf b}\widehat G_{-,\mathbf b}=I,\qquad
 \widehat G_{-,\mathbf b}\widehat D_{-,\mathbf b}=I.
\end{equation}
The interior response $\mathbb B_{\mathbf b}\mathsf C_{X,\mathbf b}$ of
the constraint data enters \eqref{eq:A-natural-broken-formula} with a
plus sign and the map \eqref{eq:A-stabilizing-range-shear} with
a minus sign, so the two cancel in
$\pi_+^Y\mathscr S_{Y,\mathbf b}\mathbb N_{\mathbf b}\iota_-^X$. Indeed,
we have the identities
\begin{equation}\label{app:minus:eq:minus-full-shear-cross-cancellation}
 \begin{aligned}
 \pi_+^Y\mathscr S_{Y,\mathbf b}\mathbb N_{\mathbf b}\iota_-^X
 &=\pi_+^Y\iota_{Y,\mathbf b}\mathsf J_{Y,\mathbf b}^{\rm raw}
       L_{\mathbf b}\mathsf R_{X,\mathbf b}\iota_-^X
    +\jmath_{Z,+,\mathbf b}q_+^Z
       \mathsf C_{X,\mathbf b}\iota_-^X,\\
 \pi_+^Y\mathbb B_{\mathbf b}\mathsf C_{X,\mathbf b}\iota_-^X
 &-\pi_+^Y\mathbb B_{\mathbf b}\mathsf C_{X,\mathbf b}\iota_-^X=0.
 \end{aligned}
\end{equation}
Here, $q_+^Z$ is the projection onto the auxiliary data of the positive
block and $\jmath_{Z,+,\mathbf b}$ includes them, so the term
$\jmath_{Z,+,\mathbf b}q_+^Z\mathsf C_{X,\mathbf b}\iota_-^X$ remains in
this off-diagonal block. The block is not part of
$\widehat D_{-,\mathbf b}$ and is estimated separately.

\begin{proposition}\label{prop:minus-diagonal}
\label{app:minus:prop:diagonal}
For every $\mathbf b\in\mathfrak U_-$, the reference cap and annulus
operators constructed above have two-sided inverses on their stated spaces.
Their coupled inverse $G_{-,\mathbf b}^{\mathbb N}$ satisfies
\eqref{app:minus:eq:whole-minus-neumann}, and
$\widehat G_{-,\mathbf b}$ is a two-sided inverse of
$\widehat D_{-,\mathbf b}$ as in \eqref{app:minus:eq:two-sided-joined}.
For every $s\ge s_*$, the solution and its outgoing Cauchy data satisfy
\begin{equation}\label{app:minus:eq:undifferentiated-tame}
 \begin{aligned}
 &\norm{\widehat G_{-,\mathbf b}F}_{\bX_{-,\mathrm{ret}}^s}
 +\norm{\Gamma_{R,-,\mathbf b}\pi_a
              \widehat G_{-,\mathbf b}F}_{Z_{R,-}^s}\\
 &\quad\le C_s\left[
   \norm F_{\bY_{-,\mathrm{ret}}^s}
  +\bigl(1+\norm{\mathbf b-\mathbf b_\rho}_{\cB^{s+k_0}}\bigr)
     \norm F_{\bY_{-,\mathrm{ret}}^{s_*}}\right],
 \end{aligned}
\end{equation}
where $\pi_a$ is the projection onto the annular component,
$\Gamma_{R,-,\mathbf b}$ is the outgoing Cauchy trace on the outer
boundary at the state, and $k_0$ is the fixed number of additional derivatives of the coefficients
required in the estimate, independent of $s$. The inverses on the cap and on the
annulus, together with their trace maps, satisfy the same estimate. On
the cap, this holds for the physical Cauchy data, and then
for the split ones by adding the bounded contribution $\mathsf M_{c,\Sigma}^sf_c$ in
\eqref{app:minus:eq:cap-trace-physicalization}.

For each fixed $m\ge1$, the state derivatives of every
solution or trace map $Q_{-,\mathbf b}:E^s\to F^s$
satisfy the estimate below, for finite integers $k_m$
and $0\le\mu_m\le m$. We denote the $m$-th derivative in
directions $h_1,\ldots,h_m$ by
$D_{\mathbf b}^mQ_{-,\mathbf b}[h_1,\ldots,h_m]$.
For $h=(\dot\rho,\dot a,\dot\eps,\dot p)$, we set
\begin{equation}\label{app:minus:eq:tangent-norm}
 \norm h_{{\rm tan},s}=\norm{\dot a}_{\cA_\gamma^s}
       +\abs{\dot\rho}+\abs{\dot\eps}+\abs{\dot p},
 \qquad \norm h_{\rm lo}=\norm h_{{\rm tan},s_*+k_m}.
\end{equation}
Then, we have
\begin{equation}\label{app:minus:eq:all-order-tame}
 \begin{aligned}
 &\norm{D_{\mathbf b}^mQ_{-,\mathbf b}
       [h_1,\ldots,h_m]q}_{F^s}\\
 &\quad\le C_{s,m}\Bigg[
   \prod_i\norm{h_i}_{\rm lo}\,\norm q_{E^{s+\mu_m}}\\
 &\qquad+\sum_i\norm{h_i}_{{\rm tan},s+k_m}
       \prod_{j\ne i}\norm{h_j}_{\rm lo}\,
                       \norm q_{E^{s_*+\mu_m}}\\
 &\qquad+\bigl(1+\norm{\mathbf b-\mathbf b_\rho}_{\cB^{s+k_m}}\bigr)
       \prod_i\norm{h_i}_{\rm lo}\,
                       \norm q_{E^{s_*+\mu_m}}\Bigg].
 \end{aligned}
\end{equation}
We may take $\mu_m=0$ for the maps on the cap. All constants are uniform
for $0\le\rho\le\rho_0$, for both signs of the cell frequency, with or
without a finite upper cutoff, and they contain no negative power of
$\rho$. At real states, that is, at states with real coefficients, the
operator, its inverse, and all trace maps carry real fields to real
fields under the identifications with the fixed spaces.
\end{proposition}

The next three subsections prove the component estimates
and complete the proof of \cref{prop:minus-diagonal}.
The correction from $\widehat D_{-,\mathbf b}$ to
$\pi_-^Y\mathbb D_{\mathbf b}\iota_-^X$ will be made in
\cref{app:macro-cross}, after estimating the couplings
between the blocks. The low and high frequencies are
then joined in the global inverse.

\subsection{The cap}
\label{app:minus:cap}

\smallskip\noindent\textit{The inverse at the circle.}
We begin at the circle, where the inverse on the cap can be constructed
explicitly on the spaces \eqref{app:minus:eq:cap-spaces}, with the
prescribed axis conditions. In the center modes $m=\pm1$, the scalar
equation is the Bessel problem
\begin{equation}\label{app:minus:eq:center-bessel}
 \left(\pa_t^2+t^{-1}\pa_t-t^{-2}+3\right)c=F,
\end{equation}
whose regular homogeneous solution
\begin{equation}\label{app:minus:eq:regular-bessel}
 u_0(t)=\frac2{\sqrt3}J_1(\sqrt3t)=t+O(t^3)
\end{equation}
has leading Cartesian coefficient one. With $v_0(t)=Y_1(\sqrt3t)$ and
the Wronskian $W=u_0v_0'-u_0'v_0=C/t$, where $C$ is a nonzero constant,
the variation of constants formula
\begin{equation}\label{app:minus:eq:bessel-volterra}
 c(t)=\alpha u_0(t)-u_0(t)\int_0^t\frac{v_0(s)F(s)}{W(s)}\dd s
       +v_0(t)\int_0^t\frac{u_0(s)F(s)}{W(s)}\dd s
\end{equation}
solves \eqref{app:minus:eq:center-bessel} with
$c(t)=\alpha t+O(t^3)$ on the axis. We apply this formula separately
to the coefficients of $e^{i\theta}$ and $e^{-i\theta}$. Since a Cartesian source on the modes $m=\pm1$ vanishes to
first order on the axis, the singular parts of the kernels in
\eqref{app:minus:eq:bessel-volterra} act as Hardy operators of the type
\eqref{appref:even-hardy-maps}, and the remaining parts are smooth on the
cap. The bounds for the Hardy operators and for the smooth parts
therefore give the two Cartesian derivatives that the norm of
$X_{c,-}^s$ in \eqref{app:minus:eq:cap-spaces} requires. Regular
homogeneous solutions are multiples of $u_0$, which proves uniqueness
for the prescribed leading coefficient $\alpha$.

We reduce the affine chains on the exceptional modes $m=\pm2$ to the two
Hardy maps
\begin{equation}\label{app:minus:eq:hardy-arrows}
 J_+q=t^{-2}\int_0^t s\,q(s)\dd s,\qquad
 J_-q=t^2\int_0^t q(s)s^{-3}\dd s.
\end{equation}
With $D=t\pa_t$, direct differentiation gives $(D+2)J_+q=q$ and
$(D-2)J_-q=q$, the second on the subspace with vanishing
$t^2$-coefficient. Fixing the linear Taylor coefficients on the axis excludes
the homogeneous solution $t^2$. For arbitrary admissible
coefficients, we first subtract a Cartesian polynomial
with those coefficients and move its image under the
operator to the source. We then apply $J_-$ to the
remainder. This polynomial extension and its image use
only finitely many transverse coefficients and are bounded
at every block index $s$. The formulas give both inverse
identities.

Together with the mean and scale equations and the affine
reconstruction in \eqref{appref:even-triangular-order},
\eqref{app:minus:eq:bessel-volterra} and
\eqref{app:minus:eq:hardy-arrows} give the inverse at the circle.
Pulling the family of small ellipses back to the fixed spaces, we obtain
a family of operators analytic in $\rho$ whose coefficients, in the
coordinates given by the Hardy maps \eqref{app:minus:eq:hardy-arrows},
contain only the factors
\begin{equation}\label{app:minus:eq:rho-factors}
 \rho,\qquad 1-\rho^2,\qquad (1-\rho^2)^{-1}.
\end{equation}
A Neumann series about the inverse at the circle gives invertibility for
sufficiently small $\rho_0$, and the inverse and all its
$\rho$-derivatives of fixed order are then uniformly bounded on
$[0,\rho_0]$. The same holds uniformly in the \emph{normalized cell
frequency} $\kappa=2^{-j}n$ of the $j$-th dyadic shell, which we let
range over $1/2\le\abs\kappa\le2$, when we replace $3$ by $3\kappa^2$
and $u_0$ by $2J_1(\sqrt3\abs\kappa t)/(\sqrt3\abs\kappa)$. Note that the
limit at the axis again has leading coefficient one. Complex conjugation
exchanges the two signs of $\kappa$, and the Hardy maps and the remaining triangular
factors depend analytically on $\kappa$, so compactness gives uniform
bounds for the inverse of the operator with frozen coefficients and for
the derivatives of its symbol used below.

\smallskip\noindent\textit{The reference problem on the shells.}
Recall the dyadic shells $Q_j$ of \eqref{app:minus:eq:sharp-shells}. We
restrict a field on the $j$-th shell to
\[
       0\le r\le 2TL2^{-j},
       \qquad t_j=2^jr/L\in[0,2T],
\]
and pull it back to the fixed domain or range space, writing
$\mathscr R_j^X$ and $\mathscr R_j^Y$ for these maps. We thus obtain the
sequence spaces
\begin{equation}\label{app:minus:eq:sharp-shell-bundles}
 \begin{aligned}
 \mathfrak X_{c,H}^s
  &=\{(U_j)_{j\ge J}:Q_jU_j=U_j\}
       \subset\ell_s^2(X^{\rm cap}),\\
 \mathfrak Y_{c,H}^s
  &=\{(F_j)_{j\ge J}:Q_jF_j=F_j\}
       \subset\ell_s^2(Y^{\rm cap}).
 \end{aligned}
\end{equation}
Here, $X^{\rm cap}$ and $Y^{\rm cap}$ are the fixed Cartesian domain and
range spaces on the cap $[0,2T]$, and $\ell_s^2$ is the weighted direct
sum of \eqref{eq:A-infinite-broken-bundle}.

We let $\mathcal P_{a,\rho,j}^{C}z=G_{a,-,\rho}(0,z)$ on shell $j$.
By \eqref{app:minus:eq:annulus-map}, this solves the homogeneous
reference annular equations with inner Cauchy data $z$. It extends a cap field to its annular complement by
\begin{equation}\label{app:minus:eq:matching-cap-injection}
 \iota_{c,j}^{X,\rm gl}U
 :=\left(U,\,
   \mathcal P_{a,\rho,j}^{C}
      T_{-,\mathbf b_\rho}\Gamma_{c,-,\rho}^{\rm ph}U,\,
   0_{Z_j}\right).
\end{equation}
The cap and annular fields then have matching physical Cauchy data at
the reference, and if the source on the cap is nonzero, their split
mismatch is instead given by
\eqref{app:minus:eq:homogeneous-physical-split-datum}. The zero
auxiliary coordinate in \eqref{app:minus:eq:matching-cap-injection}
records the physical mismatch, and on the range it is converted by
\eqref{app:minus:eq:split-mismatch-shear}. We set
$\iota_c^{X,\rm gl}=\bigoplus_{j\ge J}\iota_{c,j}^{X,\rm gl}$ and
$\iota_c^YF=(F,0,0)$ on the Cartesian range, and we write $\pi_c^X$ and
$\pi_c^Y$ for the projections onto the cap coordinates. We extract the
coordinates locally and reconstruct them through the full cap and
annular field by setting
\begin{equation}\label{app:minus:eq:cap-analysis-synthesis}
 \begin{aligned}
 \mathscr A_Xu&=(\mathscr R_j^XQ_ju)_{j\ge J},&
 \mathscr S_X&=\mathfrak R_{{\rm sh},J,\rho}^{X}
                    \iota_c^{X,\rm gl},\\
 \mathscr A_YF&=(\mathscr R_j^YQ_jF)_{j\ge J},&
 \mathscr S_Y&=\mathfrak R_{{\rm sh},J,\rho}^{Y,\rm raw}\iota_c^Y,
 \end{aligned}
\end{equation}
where \eqref{eq:A-sharp-shell-reconstruction} first glues each cap to its
own annular complement and then sums the fields on the shells. We use
the reference $\mathbf b_\rho$ because its value and normal extensions do
not depend on $\zeta$ and preserve every $Q_j$. Since
$\mathscr A_X=\pi_c^X\mathfrak B_{{\rm sh},J,\rho}^{X}$ and
$\mathscr A_Y=\pi_c^Y\mathfrak B_{{\rm sh},J,\rho}^{Y,\rm raw}$, the
inverse identities in \eqref{eq:A-sharp-shell-reconstruction} give
\begin{equation}\label{app:minus:eq:cap-frame-identities}
 \begin{gathered}
 \Pi_X:=\mathscr A_X\mathscr S_X=I_{\mathfrak X_{c,H}},\qquad
 \Pi_Y:=\mathscr A_Y\mathscr S_Y=I_{\mathfrak Y_{c,H}},\\
 \Pi_X^2=\Pi_X,\qquad \Pi_Y^2=\Pi_Y,\\
 \mathscr S_X\mathscr A_X=I\ \hbox{on }\operatorname{Ran}\mathscr S_X,
 \qquad
 \mathscr S_Y\mathscr A_Y=I\ \hbox{on }\operatorname{Ran}\mathscr S_Y.
 \end{gathered}
\end{equation}
As in \cref{sec:block}, variable coefficients act after
$\mathfrak R_{{\rm sh},J,\rho}^{X}$ has joined the cap
and annulus fields and reconstructed a global field.

The explicit inverse has coefficients frozen on the axis.
It gives an approximate inverse, or \emph{parametrix},
on the shells, with small errors on both sides at high
cell frequency and for a sufficiently weak analytic
weight. The following lemma applies to both cap blocks.

\begin{lemma}
\label{app:minus:lem:cap-quantization}
Let $\cC(q,\kappa)$ be the operator with frozen coefficients on the cap
$[0,2T]$ of the minus block or of the positive block, after separation of
the modes and identification with fixed Cartesian domain and range
spaces. The parameter $q$ lies in a
compact set of finite jets, and $1/2\le\abs\kappa\le2$ is the normalized
cell frequency variable. Suppose that
\[
 K(q,\kappa)\cC(q,\kappa)=I_X,
 \qquad \cC(q,\kappa)K(q,\kappa)=I_Y,
\]
that $(q,\kappa)\mapsto(\cC,K)$ is analytic between the corresponding
spaces at the base index $s_*$, and that for each finite order $N$ used
below we have the uniform
bounds
\begin{equation}\label{app:minus:eq:cap-symbol-seminorms}
 \sup_{(q,\kappa)}\max_{|\alpha|+r\le N}
 \left(
  \|\pa_q^\alpha\pa_\kappa^r\cC(q,\kappa)\|_{X\to Y}
 +\|\pa_q^\alpha\pa_\kappa^rK(q,\kappa)\|_{Y\to X}
 \right)\le M_N.
\end{equation}
We require the same bounds for all added axis, Dirichlet, and balanced
trace components. Let $q_{\mathbf b_\rho}$ be the jet of the reference
ellipse with parameter $\rho$, whose coefficients do not depend on
the cell variable, let
\begin{equation}\label{app:minus:eq:full-reference-shell-operator}
 \mathbf L_{\rho,H}^{\rm sh}
 :=\mathfrak B_{{\rm sh},J,\rho}^{Y,\rm raw}
       L_{\mathbf b_\rho}
       \mathfrak R_{{\rm sh},J,\rho}^{X}
\end{equation}
be the reconstructed reference operator in the coordinates of the
shells, and let
\begin{equation}\label{app:minus:eq:reference-cap-diagonal}
 \mathbf C_{\rho}^{\rm cap}
 :=\operatorname{Diag}_{j\ge J}
  \left(\pi_{c,j}^Y\mathbf L_{\rho,H}^{\rm sh}
             \iota_{c,j}^{X,\rm gl}\right):
 \mathfrak X_{c,H}^{s_*}\longrightarrow\mathfrak Y_{c,H}^{s_*}
\end{equation}
be its diagonal on the cap components. Here, $\pi_{c,j}^Y$ is the
projection onto the cap coordinate of the $j$-th shell, and the
inclusion $\iota_{c,j}^{X,\rm gl}$ adds the annular field and the zero
mismatch needed for the reconstruction. Thus,
\eqref{app:minus:eq:reference-cap-diagonal} evaluates each cap field
only on its own collar.

For $2^j\le|n|<2^{j+1}$, we set $\kappa_{j,n}=2^{-j}n$. Freezing the
coefficients on the axis, we obtain the following Fourier multipliers on
the sequence spaces,
\begin{equation}\label{app:minus:eq:sequence-frozen-inverse}
 \begin{aligned}
 (\mathbf C_{\rho}^{\rm fr}U)_{j,n}
   &=\cC(q_{\mathbf b_\rho},\kappa_{j,n})U_{j,n},\\
 (\mathbf K_{\rho}^{\rm fr}F)_{j,n}
   &=K(q_{\mathbf b_\rho},\kappa_{j,n})F_{j,n},\\
 \widehat{\mathbf C}_{\rho}^{\rm fr}
   &:=\Pi_Y\mathbf C_{\rho}^{\rm fr}\Pi_X
     =\mathbf C_{\rho}^{\rm fr},&
 \widehat{\mathbf K}_{\rho}^{\rm fr}
   &:=\Pi_X\mathbf K_{\rho}^{\rm fr}\Pi_Y
     =\mathbf K_{\rho}^{\rm fr}.
 \end{aligned}
\end{equation}
Here, $U_{j,n}$ and $F_{j,n}$ denote the $n$-th cell Fourier
coefficients of $U_j$ and $F_j$. Then, the reference parametrix
\begin{equation}
 \mathbf K_{\rho,H}^{(0)}
 :=\widehat{\mathbf K}_{\rho}^{\rm fr}
 \label{app:minus:eq:cap-parametrix-formula}
\end{equation}
satisfies, on the sequence spaces at the base index,
\begin{equation}
 \begin{aligned}
 &\norm{\mathbf K_{\rho,H}^{(0)}\mathbf C_{\rho}^{\rm cap}-I}_
          {\mathfrak X_{c,H}^{s_*}\to\mathfrak X_{c,H}^{s_*}}
 +\norm{\mathbf C_{\rho}^{\rm cap}\mathbf K_{\rho,H}^{(0)}-I}_
          {\mathfrak Y_{c,H}^{s_*}\to\mathfrak Y_{c,H}^{s_*}}\le C\left(2^{-J}+\gamma\right),
 \end{aligned}
 \label{app:minus:eq:abstract-cap-defect}
\end{equation}
where $C$ depends only on $T$, on the bounds
\eqref{app:minus:eq:cap-symbol-seminorms}, and on the compact set of
coefficient jets. For the positive cap problem, the assumed operators include the
Dirichlet datum on the interface. The outgoing balanced conormal
derivative of $\mathbf K_{\rho,H}^{(0)}$ satisfies the corresponding
estimate of order zero.
\end{lemma}

\begin{proof}
On $Q_j$, we evaluate the multiplier in
\eqref{app:minus:eq:sequence-frozen-inverse} at each
$\kappa_{j,n}=2^{-j}n$, and the assumed bounds cover both signs of the
cell frequency and the whole normalized frequency range $1/2\leq|\kappa|\leq2$.

After rescaling the equations and the unknowns, we find that each
differential expression is a finite sum of terms of the form
\begin{equation}
 a_\nu(2^{-j}Lt,t,\zeta)
     (t\pa_t)^{\alpha_\nu}(2^{-j}\pa_\zeta)^{\beta_\nu},
 \qquad \alpha_\nu+\beta_\nu\le2,
 \label{app:minus:eq:normalized-cap-row}
\end{equation}
acting between the fixed Cartesian spaces. Note that the axis and
boundary operators keep the orders of \cref{tab:block-weights}. Since the
reference coefficients are smooth in the radial variable, Taylor's
formula on $0\le t\le2T$ gives
\begin{equation}
 a_\nu(2^{-j}Lt,t,\zeta)
  =a_\nu(0,t,\zeta)+2^{-j}Lt
       \int_0^1(\pa_ra_\nu)(s2^{-j}Lt,t,\zeta)\,\dd s.
 \label{app:minus:eq:cap-radial-taylor}
\end{equation}
The remainder is bounded by $C_T2^{-j}$, and the same calculation
applies to the finite axis equations. The spectral projections depend
analytically on the coefficients through their contour integrals, so
their remainders obey the same bound.

The reference coefficients do not depend on $\zeta$, hence preserve each
Fourier mode and each dyadic shell, giving
\begin{equation}\label{app:minus:eq:sharp-frozen-identities}
 \widehat{\mathbf K}_{\rho}^{\rm fr}
       \widehat{\mathbf C}_{\rho}^{\rm fr}=\Pi_X,\qquad
 \widehat{\mathbf C}_{\rho}^{\rm fr}
       \widehat{\mathbf K}_{\rho}^{\rm fr}=\Pi_Y,
\end{equation}
and we have the corresponding identities for the projections
\begin{equation}\label{app:minus:eq:sharp-frozen-covariance}
 \Pi_Y\mathbf C_{\rho}^{\rm fr}
      -\mathbf C_{\rho}^{\rm fr}\Pi_X=0,\qquad
 \Pi_X\mathbf K_{\rho}^{\rm fr}
      -\mathbf K_{\rho}^{\rm fr}\Pi_Y=0.
\end{equation}
The bounds on the symbols and Parseval's identity give
\begin{equation}\label{app:minus:eq:frozen-sequence-bounds}
 \|\widehat{\mathbf C}_{\rho}^{\rm fr}\|_
      {\mathfrak X_{c,H}^{s_*}\to\mathfrak Y_{c,H}^{s_*}}
 +\|\widehat{\mathbf K}_{\rho}^{\rm fr}\|_
      {\mathfrak Y_{c,H}^{s_*}\to\mathfrak X_{c,H}^{s_*}}
 \le C,
\end{equation}
so the inverse of the operator with frozen coefficients satisfies both
inverse identities on these spaces. In particular, we have
\begin{equation}\label{app:minus:eq:compressed-frozen-defects}
 \|\widehat{\mathbf K}_{\rho}^{\rm fr}
       \widehat{\mathbf C}_{\rho}^{\rm fr}-\Pi_X\|
 +\|\widehat{\mathbf C}_{\rho}^{\rm fr}
       \widehat{\mathbf K}_{\rho}^{\rm fr}-\Pi_Y\|=0.
\end{equation}

We now compare the reference operator on the cap with the multiplier
frozen on the axis. On the $j$-th shell, the difference is
\begin{equation}\label{app:minus:eq:typed-reference-cap-error}
 \mathcal E_{\rho,j}
 :=\pi_{c,j}^Y\mathbf L_{\rho,H}^{\rm sh}
       \iota_{c,j}^{X,\rm gl}
   -\cC(q_{\mathbf b_\rho},2^{-j}D_\zeta)Q_j:
 Q_jX^{\rm cap}\longrightarrow Q_jY^{\rm cap}.
\end{equation}
Since $L_{\mathbf b_\rho}$ is independent of the cell
variable, there are no terms between shells. Taylor's
estimate \eqref{app:minus:eq:cap-radial-taylor} bounds
the radial error by $C_T2^{-j}$. The reconstruction joins
the matching annular field to the cap before the operator
acts. Its commutators belong to $\mathcal E_{\rho,j}$
and satisfy the same bound at the reference. Conjugation by the analytic weight contributes $C_T\gamma$,
since $r\Lambda=O_T(1)$ on the scaled cap. The trace theorem in the
Cartesian variables on the cap gives the same bounds for the added value
and balanced conormal
components. Parseval's identity then gives
\begin{equation}\label{app:minus:eq:reference-cap-frozen-difference}
 \|\mathbf C_{\rho}^{\rm cap}
       -\widehat{\mathbf C}_{\rho}^{\rm fr}\|_
   {\mathfrak X_{c,H}^{s_*}\to\mathfrak Y_{c,H}^{s_*}}
 \le C_T(2^{-J}+\gamma).
\end{equation}
Using \eqref{app:minus:eq:sharp-frozen-identities}, we factor the two
errors as
\begin{align}
 \mathbf K_{\rho,H}^{(0)}\mathbf C_{\rho}^{\rm cap}-I
  &=\widehat{\mathbf K}_{\rho}^{\rm fr}
       (\mathbf C_{\rho}^{\rm cap}
          -\widehat{\mathbf C}_{\rho}^{\rm fr}),\notag\\
 \mathbf C_{\rho}^{\rm cap}\mathbf K_{\rho,H}^{(0)}-I
  &=(\mathbf C_{\rho}^{\rm cap}
          -\widehat{\mathbf C}_{\rho}^{\rm fr})
       \widehat{\mathbf K}_{\rho}^{\rm fr}.
 \label{app:minus:eq:reference-cap-defect-factorization}
\end{align}
The bounds \eqref{app:minus:eq:frozen-sequence-bounds} and
\eqref{app:minus:eq:reference-cap-frozen-difference} now give
\eqref{app:minus:eq:abstract-cap-defect}, and applying the trace map to
these identities, we obtain the estimate for the outgoing balanced
conormal derivative.
\end{proof}

\smallskip\noindent\textit{Perturbation at a general state.}
For a general state, the coefficients may mix cell
frequencies. We therefore retain both the cap and annulus
fields when comparing with the reference, and set
\begin{equation}\label{app:minus:eq:current-full-shell-operator}
 \mathbf L_{\mathbf b,H}^{\rm sh}
 :=\mathfrak B_{{\rm sh},J,\rho}^{Y,\rm raw}
       L_{\mathbf b}
       \mathfrak R_{{\rm sh},J,\rho}^{X},
 \qquad
 \mathbf R_{\mathbf b,H}^{\rm sh}
 :=\mathbf L_{\mathbf b,H}^{\rm sh}
       -\mathbf L_{\rho,H}^{\rm sh}.
\end{equation}
For $\alpha,\beta\in\{c,a\}$, its blocks are
\begin{equation}\label{app:minus:eq:typed-current-shell-blocks}
 \mathcal R_{\mathbf b,jk}^{\alpha\beta}
 :=\pi_{\alpha,j}^Y\mathbf R_{\mathbf b,H}^{\rm sh}
       \iota_{\beta,k}^X:
 Q_kX_{\beta,k}^s\longrightarrow Q_jY_{\alpha,j}^s .
\end{equation}
Here, $\iota_{\beta,k}^X$ is the inclusion of the component $\beta$ of the
$k$-th shell in the space of separate fields \eqref{eq:A-infinite-broken-bundle},
$\pi_{\alpha,j}^Y$ is the projection onto the component $\alpha$ of the
$j$-th shell of the range, and $X_{\beta,k}^s$ and $Y_{\alpha,j}^s$ are
the corresponding summands. The first superscript is the output region
and the second is the input region; for example,
$\mathcal R_{\mathbf b,jk}^{ca}$ sends an annular field on shell $k$
to a cap source on shell $j$. We write
$\operatorname{Ins}_{k,\rho}^{\rm br}(U_{c,k},U_{a,k})$ for the
reference insertion of a cap field and its annular complement. Note that
this insertion is defined at all radii and that a mismatch is allowed.
The interior operator is then
\begin{equation}\label{app:minus:eq:explicit-full-pair-shell-block}
 \mathcal R_{\mathbf b,jk}^{\rm bulk}U_k
 =\mathscr R_j^YQ_j
    (L_{\mathbf b}-L_{\mathbf b_\rho})
       \operatorname{Ins}_{k,\rho}^{\rm br}(U_{c,k},U_{a,k}).
\end{equation}
The coefficient acts on the full cap and annulus pair.
A cap input contributes only within its own collar.
Outside that collar, the input belongs to the annular
block $\beta=a$ of \eqref{app:minus:eq:typed-current-shell-blocks}. We write
$\mathscr Y_{{\rm br},H}^{s,\rm raw}$ for the weighted direct sum of the
Cartesian source spaces on the caps and on the annuli, which is the
counterpart on the range of $\mathscr X_{{\rm br},H}^s$ in
\eqref{eq:A-infinite-broken-bundle}.

We abbreviate
\begin{equation}\label{app:minus:eq:cap-coefficient-distances}
 \delta_0(\mathbf b)
 :=\|\mathbf b-\mathbf b_\rho\|_{\cB^{s_*+k_c}},
 \qquad
 \delta_s(\mathbf b)
 :=\|\mathbf b-\mathbf b_\rho\|_{\cB^{s+k_c}},
\end{equation}
which measure the distance from the state to the constant ellipse at the
base index and at the index $s$. Here, $k_c$ is the fixed number of additional derivatives of the
coefficients required in the estimate, independent of $s$. The numbers
$k_A$, $k_*$, and $k$ below have the same meaning in their respective estimates. Expanding the two
reconstruction maps,
the coefficients of the state, and the reduction of the source, we
separate the difference into
\begin{equation}\label{app:minus:eq:localized-cap-ledger}
 \mathbf R_{\mathbf b,H}^{\rm sh}
 =R_{\mathbf b}^{\rm hi}+R_{\mathbf b}^{\rm nb}
  +R_{\mathbf b}^{\rm bd}+R_{\mathbf b}^{\rm res}
  +R_{\mathbf b}^{\rm triv}+R_{\mathbf b}^{\rm syn}
  +R_{\mathbf b}^{\rm ph}.
\end{equation}
The following table lists the contributions. Each is
estimated on the reconstructed field by coefficient
differences, bounded coordinate maps, and the tame product
estimate. The cutoff term is separated because it becomes
small with the distance to the reference, rather than
with increasing $J$.
\begin{center}
\small
\begin{tabular}{@{}lp{0.82\textwidth}@{}}
\toprule
Term & Contribution\\
\midrule
$R^{\rm hi}$ & High-frequency coefficient terms on the diagonal of the shells.\\
$R^{\rm nb}$ & Terms between distinct shells $j,k>J$, including cap--annulus terms.\\
$R^{\rm bd}$ & Terms meeting the cutoff, where $j=J$ or $k=J$.\\
$R^{\rm res}$ & Change in the reduction of the source.\\
$R^{\rm triv}$ & Change in the identifications with fixed spaces.\\
$R^{\rm syn}$ & Commutators of the reference reconstruction with $L_{\mathbf b}-L_{\mathbf b_\rho}$.\\
$R^{\rm ph}$ & Changes in the weight, the value and conormal coordinates after reconstruction, and the conormal extensions.\\
\bottomrule
\end{tabular}
\end{center}

On each shell, we write $(d,h_\rho)$ for the reference
Cauchy data, $(d,n_{\mathbf b})$ for the physical data at
the state, and $f_{c,j}$ for the cap interior source.
The split data are
\begin{equation}\label{app:minus:eq:current-split-cap-row}
 \binom d{r_{\mathbf b}^{\rm sp}}
  :=\binom d{n_{\mathbf b}}
       +j_N\mathsf M_{c,\Sigma}^sf_{c,j},
\end{equation}
and the physical and the reference Cauchy data are related by
\begin{equation}\label{app:minus:eq:current-conormal-row-adapter}
 \binom d{n_{\mathbf b}}
 =T_{\mathbf b,\rho,j}\binom d{h_\rho},\qquad
 T_{\mathbf b,\rho,j}
 =\begin{pmatrix}I&0\\K_{\mathbf b,\rho,j}&
       E_{\mathbf b,\rho,j}\end{pmatrix},\qquad
 E_{\mathbf b,\rho,j}=e_{\mathbf b}e_\rho^{-1}.
\end{equation}
Here, $T_{\mathbf b,\rho,j}$ is
$\mathcal T_{\Sigma,\mathbf b,\rho}^{CN}$ from
\eqref{eq:A-current-reference-row-map} on interface $j$.
The normal coefficients $e_{\mathbf b}$ and $e_\rho$
are uniformly invertible, while the tangential term
$K_{\mathbf b,\rho,j}=K_{\Sigma,\mathbf b,\rho}$ has
order zero between the balanced spaces. Thus, the inverse is
\[
 T_{\mathbf b,\rho,j}^{-1}
 =\begin{pmatrix}I&0\\
   -E_{\mathbf b,\rho,j}^{-1}K_{\mathbf b,\rho,j}&
    E_{\mathbf b,\rho,j}^{-1}\end{pmatrix}.
\]
To recover the reference conormal derivative, we therefore subtract the boundary contribution $\mathsf M_{c,\Sigma}^sf_{c,j}$ before applying the inverse, which gives
\begin{equation}\label{app:minus:eq:split-before-current-converter}
 \binom d{h_\rho}
 =T_{\mathbf b,\rho,j}^{-1}
    \left(\binom d{r_{\mathbf b}^{\rm sp}}
       -j_N\mathsf M_{c,\Sigma}^sf_{c,j}\right).
\end{equation}
If we applied $T_{\mathbf b,\rho,j}^{-1}$ directly to the split Cauchy
data, we would miss this source term. Since $T_{\mathbf b,\rho,j}$ need
not commute with $Q_j$, we include its terms between different shells in
$R^{\rm bd}+R^{\rm nb}+R^{\rm ph}$.

At $\mathbf b=\mathbf b_\rho$, the independence of the cell variable
makes every term between different shells and every term at the cutoff
in \eqref{app:minus:eq:localized-cap-ledger} vanish. On the
reconstructed field, Parseval's identity and the tame product estimate
\eqref{eq:analytic-product} give
\begin{equation}\label{app:minus:eq:current-full-shell-bound}
 \|\mathbf R_{\mathbf b,H}^{\rm sh}U\|_
       {\mathscr Y_{{\rm br},H}^{s,\rm raw}}
 \le C_s\left[
   \delta_0(\mathbf b)\|U\|_{\mathscr X_{{\rm br},H}^s}
  +\delta_s(\mathbf b)\|U\|_{\mathscr X_{{\rm br},H}^{s_*}}\right].
\end{equation}
This includes the operator $R_{\mathbf b}^{\rm bd}$ at the cutoff,
extended by zero away from the $J$-th shell, for which we have the
estimate
\begin{equation}\label{app:minus:eq:boundary-shell-bound}
 \|R_{\mathbf b}^{\rm bd}U\|_
       {\mathscr Y_{{\rm br},H}^{s,\rm raw}}
 \le C_s\left[
   \delta_0(\mathbf b)\|U\|_{\mathscr X_{{\rm br},H}^s}
  +\delta_s(\mathbf b)\|U\|_{\mathscr X_{{\rm br},H}^{s_*}}\right].
\end{equation}
The factor $2^{-J}$ in
\eqref{app:minus:eq:reference-cap-frozen-difference}
comes from freezing the radial coefficients. It does not
apply to the cutoff estimate
\eqref{app:minus:eq:boundary-shell-bound}.

The terms coupling low and high frequencies also satisfy the
perturbative bound
\begin{equation}\label{app:minus:eq:localized-cap-leakage}
 \begin{aligned}
 &\|P_J^YL_{\mathbf b}H_J^Xu_H\|_{Y^s}
  +\|H_J^YL_{\mathbf b}P_J^Xu_L\|_{Y^s}\\
 &\quad\le C_{J,s}\left[
   \delta_0(\mathbf b)
       \bigl(\|u_H\|_{X^s}+\|u_L\|_{X^s}\bigr)
  +\delta_s(\mathbf b)
       \bigl(\|u_H\|_{X^{s_*}}+\|u_L\|_{X^{s_*}}\bigr)\right].
 \end{aligned}
\end{equation}
Here, $u_H$ and $u_L$ are arbitrary fields, of which the display keeps
the high and the low cell frequencies. The two terms also vanish at the
reference, whose coefficients do not
depend on the cell variable.

For the cap of the minus block, we take the parameter
\begin{equation}\label{app:minus:eq:frozen-cap-jet}
 q_{-,\rho}
 =\bigl(\rho,j^1_{\!y}a_\rho(0),\eps_\rho,p_\rho,
       P_{-,\rho}^D(0),P_{-,\rho}^R(0)\bigr),
\end{equation}
where $j^1_{\!y}a_\rho(0)$ records the value and first derivatives in $y$
of $a_\rho$ at the axis, and $P_{-,\rho}^{D}(0)$ and $P_{-,\rho}^{R}(0)$ are the spectral
projections \eqref{app:cross:domain-range-contours} of the reference at
the axis. We take the inverse $K_{-,\rho}^{\rm fr}(\kappa)$ built
from the Bessel and Hardy formulas above. Applying
\cref{app:minus:lem:cap-quantization}
with $q_{\mathbf b_\rho}=q_{-,\rho}$, we obtain the reference operator on
the cap
\begin{equation}\label{app:minus:eq:literal-retained-cap-realization}
 \cA_{c,-,\rho}^{\rm diag}
 :=\mathbf C_{\rho}^{\rm cap}:
 \mathfrak X_{c,H}^{s_*}\longrightarrow\mathfrak Y_{c,H}^{s_*}.
\end{equation}
Here, as in \cref{prop:minus-diagonal}, we identify the spaces
$X_{c,-}^s$ and $Y_{c,-}^s$ with the sequence spaces by $\mathscr A_X$
and $\mathscr A_Y$ in \eqref{app:minus:eq:cap-analysis-synthesis}. The
bounds on the symbols and Parseval's identity define the operator at
every $s\ge s_*$, and we prove the Neumann estimate first at the base
index. Note that the operator is the diagonal on the caps of the
reconstructed reference problem, and not the separate compression
$H_J^Y\cC_{-,\mathbf b}H_J^X$. With
$K_{c,-,\rho}^{(0)}:=\mathbf K_{\rho,H}^{(0)}$,
\eqref{app:minus:eq:abstract-cap-defect} gives
\begin{equation}\label{app:minus:eq:cap-defect}
 \|K_{c,-,\rho}^{(0)}\cA_{c,-,\rho}^{\rm diag}-I\|
 +\|\cA_{c,-,\rho}^{\rm diag}K_{c,-,\rho}^{(0)}-I\|
 \le C_{\rm cap}(2^{-J}+\gamma)<\frac12.
\end{equation}
We write the two errors as
\[
 R_X^\rho=K_{c,-,\rho}^{(0)}\cA_{c,-,\rho}^{\rm diag}-I,
 \qquad
 R_Y^\rho=\cA_{c,-,\rho}^{\rm diag}K_{c,-,\rho}^{(0)}-I.
\]
The operator
\begin{equation}\label{app:minus:eq:two-sided-neumann}
 G_{c,-,\rho}^{\rm diag}
 :=(I+R_X^\rho)^{-1}K_{c,-,\rho}^{(0)}
  =K_{c,-,\rho}^{(0)}(I+R_Y^\rho)^{-1}
\end{equation}
is a two-sided inverse of
\eqref{app:minus:eq:literal-retained-cap-realization}. The trace theorem
in the Cartesian variables, with the half orders in
\eqref{app:minus:eq:cap-trace-space}, bounds the outgoing physical Cauchy
data, and adding the boundary contribution $\mathsf M_{c,\Sigma}^sf_c$ in
\eqref{app:minus:eq:split-cap-solution-trace} bounds the split Cauchy
data at the same index. Note that the value on the interface is
determined by this solution and is not prescribed.

The terms at the state on the diagonal, between different shells, at the
cutoff, and between the cap and the annulus remain in the perturbation
\eqref{app:minus:eq:explicit-full-pair-shell-block} and
\eqref{app:minus:eq:localized-cap-ledger}. Once the reference solutions on
the cap and annulus have been coupled, we include these perturbations
through the Neumann formula \eqref{app:minus:eq:whole-minus-neumann}
to obtain the inverse at the state.

\subsection{The annulus}
\label{app:minus:annulus}

We now derive the annular equations from the linearization
\eqref{eq:current-linearization} at $\mathbf b$ to obtain energy estimates
for the center and affine components. For the center equation, we
construct a \emph{symmetrizer}, a positive matrix defining its energy.
For the affine chains, we estimate the coefficient that takes the upper
entries into the equations for the lower entries. For a variation
$(U,S)$, we recall the frame
\eqref{app:cross:current-frame},
\[
 p=\cD v,\qquad q_1=\cR v,\qquad
 T=L^{-1}D_\eps v,\qquad \mathbb B=(p,q_1,T),
\]
and, on the annulus where $\mathbb B$ is a frame, we write
$U=\mathbb Bc_{\rm fr}$, where $c_{\rm fr}$ is the vector of its three
frame coefficients. With the homogeneous cell derivative
$\mathring D_\eps U=U_\zeta+\eps e_z\times U$ of
\eqref{app:cross:homogeneous-cell-derivative}, we form the matrices
\[
 \Gamma_{\cD}=\mathbb B^{-1}\cD\mathbb B,\qquad
 \Gamma_{\cR}=\mathbb B^{-1}\cR\mathbb B,\qquad
 \Gamma_\eps=\mathbb B^{-1}\mathring D_\eps\mathbb B,
\]
which record differentiation of the frame. For example,
$\cD U=\mathbb B(\cD c_{\rm fr}+\Gamma_{\cD}c_{\rm fr})$, and similarly
for $\cR$ and $\mathring D_\eps$. The frame derivatives
multiply the unknown at order zero. The terms that
differentiate it give the principal symbol with frozen
coefficients.

We pass from symbols to operators in the cell variable in two ways. If
an operator-valued symbol has the Fourier expansion
$a(\zeta,\xi)=\sum_ka_k(\xi)e^{ik\zeta}$ in the cell variable, we set
\begin{equation}\label{app:minus:eq:quantization-convention}
 \begin{aligned}
  (\operatorname{Op}(a)u)_n
    &=\sum_m a_{n-m}(m)u_m,\\
  (\operatorname{Op}^{w}(a)u)_n
    &=\sum_m a_{n-m}\!\left(\frac{n+m}{2}\right)u_m,
 \end{aligned}
\end{equation}
with the radial coordinate as a parameter. The first is the \emph{standard
quantization}, and the second is the \emph{Weyl quantization}, which has the
advantage that a Hermitian symbol gives a symmetric operator. Note that all
symbols act on the finite-dimensional space of the center and affine
components. Then, $\operatorname{Op}^{w}(a)-\operatorname{Op}(a)$ has
the kernel
\begin{equation}\label{app:minus:eq:weyl-left-difference}
 a_{n-m}\!\left(\frac{n+m}{2}\right)-a_{n-m}(m)
\end{equation}
in the cell modes, and the mean value theorem introduces a factor
$|n-m|$ and one derivative in the frequency variable, which lowers the
order of the symbol by one.

The center and affine parts of the principal symbol with frozen
coefficients are described in the following lemma.
\begin{lemma}
\label{app:minus:lem:current-reduction}
The domain and range spectral projections of the principal symbol at
$\mathbf b$, expressed on the fixed spaces of the annulus, have the
following properties.
\begin{enumerate}[label=\textup{(\roman*)}]
\item The four nonzero center roots correspond to an invariant spectral
subspace. In the fixed coordinates on this subspace, we write
$c=(c_1,c_2,c_3,c_4)$ for the center unknown. The principal
part of the center evolution is
\begin{equation}
 \pa_r c=i\operatorname{Op}^{w}(\mathbf G_{\mathbf b}|\xi|)c,
 \qquad
 \norm{\mathbf G_{\mathbf b}-\mathbf G_{\mathbf b_\rho}}_{\rm low}
 \le C\norm{\mathbf b-\mathbf b_\rho}_{\cB^{s_*+k_c}}.
 \label{app:minus:eq:current-principal-difference}
\end{equation}
Here, $\norm{\cdot}_{\rm low}$ is a coefficient norm at
a fixed lower regularity index.
At the constant reference state, this system is equivalent to two
oscillator
equations
\begin{equation}
 \xi''+\Omega_\rho^2\xi+r^{-2}V_\rho\xi=g,
 \qquad
 0<c_\Omega I\le\Omega_\rho^2\le C_\Omega I,
 \label{app:minus:eq:oscillator-form}
\end{equation}
uniformly for $0\le\rho\le\rho_0$. In the first-order variables
$z=(\xi,\xi')$, we take
\begin{equation}
 \mathbf J_\rho=
 \begin{pmatrix}0&I\\-\Omega_\rho^2&0\end{pmatrix},
 \qquad
 S_e=\begin{pmatrix}\Omega_\rho^2&0\\0&I\end{pmatrix}.
 \label{app:minus:eq:oscillator-symmetrizer}
\end{equation}
Then, $S_e(-i\mathbf J_\rho)$ is Hermitian and $cI\le S_e\le CI$ holds
uniformly, including at the circle, where the speeds coincide.

\item The zero root gives the two affine chains, whose jets on the axis
are prescribed. In the coordinates of the fixed spaces, its matrix of
highest order is a direct sum of copies of
\begin{equation}\label{app:minus:eq:affine-chain-matrix}
 J_A=\begin{pmatrix}0&I\\0&0\end{pmatrix},
 \qquad J_A^2=0,
\end{equation}
after possibly exchanging the two coordinates of a chain. Thus, $J_A(x,y)=(y,0)$ for the upper and lower pairs $(x,y)$.
If $B_{21,\mathbf b}x$ denotes the contribution of the upper pair to
the equation for the lower pair, then
\begin{equation}
 \norm{B_{21,\mathbf b}}_{\cA_\gamma^{s_*+k_A}}
 \le C_A\norm{\mathbf b-\mathbf b_\rho}_{\cB^{s_*+k_A}}
 \le M_B\delta_A,
 \label{app:minus:eq:affine-B21-bound}
\end{equation}
where the last inequality is the first condition imposed in
\eqref{app:minus:eq:state-cutoff-choices}.

% AUTHOR QUERY: Reconcile \eqref{app:minus:eq:oscillator-form} with the
% scaled-t equation \eqref{appref:odd-center-oscillator}. Specify the
% physical-r normalization, and retain the reference potential in the
% energy estimate with a frequency-uniform integral bound. Item (iii)
% currently makes this nonzero reference remainder vanish.
\item The remaining terms in the center and affine equations, including
the couplings between them, have order zero, and their norm at the base
index is at most $\chi_-(\mathbf b)$. With
$\mathbf G_e:=\mathbf G_{\mathbf b_\rho}$ and with
$\chi_c(\mathbf b)$ bounding
\(S_e(\mathbf G_{\mathbf b}-\mathbf G_e)
 -(\mathbf G_{\mathbf b}-\mathbf G_e)^*S_e\), which measures the
failure of self-adjointness in the energy inner product, we have
\begin{equation}
 \chi_c(\mathbf b)+\chi_-(\mathbf b)
 \le C\norm{\mathbf b-\mathbf b_\rho}_{\cB^{s_*+k_*}},
 \label{app:minus:eq:current-connection-small}
\end{equation}
together with the tame estimates in which at most one
coefficient is measured in a higher regularity norm.
\end{enumerate}
\end{lemma}

\begin{proof}
The frame identities reduce the part of highest order to the principal
symbol determined by the frozen Gram matrix
$\mathbb B^*\mathbb B$, whose normalized matrices are \eqref{app:cross:descriptor-matrices}, and we
take the domain and range spectral projections on the same contours
around the center roots and around zero. By the resolvent identity,
these projections, the identifications with the fixed spaces, and the
reduced matrices depend analytically on the coefficients of the Gram
matrix, so the estimates \eqref{app:minus:eq:current-principal-difference}
and \eqref{app:minus:eq:current-connection-small} follow, including
their tame versions.

For the constant ellipse, \eqref{appref:odd-center-oscillator} and the
sign convention for the radial mass give
\eqref{app:minus:eq:oscillator-form}, and in the first-order variables
$z=(\xi,\xi')$ its matrix of highest order is $\mathbf J_\rho$. The
eigenvalues of
$\Omega_\rho^2$ lie in a fixed compact subset of $(0,\infty)$, and we
verify the identity for the symmetrizer by multiplying out
\eqref{app:minus:eq:oscillator-symmetrizer}. This verification uses the
matrix of the oscillator directly and requires no choice of
eigenvectors. An affine solution whose cross-section rotates along the
axis adds no cell derivative of highest order. Indeed, simultaneous
rotation of the disk and of the normal plane gives
\begin{equation}
 \mathcal U_\alpha^{-1}\pa_\zeta\mathcal U_\alpha
 =\pa_\zeta+\alpha'(\zeta)(\widehat J-\cR_y),
 \label{app:minus:eq:exact-corotation}
\end{equation}
where $\mathcal U_\alpha$ is the rotation by the angle $\alpha(\zeta)$
of the disk and of the normal plane together, $\widehat J$ is its
generator on the components, and $\cR_y$ is the angular derivative in
the disk. The second term has order zero on the finite-dimensional center
space, so it contributes only to the lower-order part of the radial
evolution.

For the zero root, the affine solutions at zero curvature of
\cref{ss:strategy} give two chains of length two, and substituting into
\eqref{app:cross:descriptor-matrices}, in the coordinates given by the
Hardy maps \eqref{app:minus:eq:hardy-arrows}, we obtain
\eqref{app:minus:eq:affine-chain-matrix}. At the constant reference
state, whose matrix $M_{\rho,p_e}$ is constant in $\zeta$, the lower left
coefficient vanishes by the parity in the angle, and the maps that
extract and reconstruct the
chains are analytic, so the mean value theorem in Banach spaces gives the
first bound in \eqref{app:minus:eq:affine-B21-bound}. Note that the
matrix $B_{21,\mathbf b}$ need not be nilpotent, but its smallness
suffices for the estimate with the modified energy
\eqref{app:minus:eq:affine-energy} below.
\end{proof}

We next estimate the initial value problem on the annulus uniformly at
high cell frequency. We begin with the center equation, keeping all four
modes together, and write it as
\begin{equation}\label{app:minus:eq:center-evolution}
 \pa_r c=i\operatorname{Op}^{w}
       \bigl(\mathbf G_{\mathbf b}(r,\zeta,\omega)\abs\xi\bigr)c
       +B_{\mathbf b}c+f_c,\qquad \omega=\xi/\abs\xi,
\end{equation}
where $B_{\mathbf b}$ collects the terms of order zero. The four center
modes have a uniformly positive symmetrizer, and we have
the estimates
\begin{equation}\label{app:minus:eq:center-symmetrizer}
 cI\le S_e\le C I,\qquad
 S_e\mathbf G_e=(S_e\mathbf G_e)^*,\qquad
 \norm{\mathbf G_e^{-1}}\le C.
\end{equation}
The energy defined by this symmetrizer is that of the
first-order oscillator system. It remains uniformly
positive when the speeds coincide because the eigenvalues
of $\Omega_\rho^2$ stay in a compact subset of $(0,\infty)$.

For a smooth cutoff $\Delta_j$ to a dyadic band of cell frequencies, we
let $\lambda_j=\min\{\langle n\rangle:n\in\operatorname{supp}\Delta_j\}\simeq2^j$
and put
\begin{equation}\label{app:minus:eq:center-dyadic-energy-definition}
 c_j=e^{\Phi_\gamma(r,D_\zeta)}\Delta_jc,
 \qquad E_j=\langle S_ec_j,c_j\rangle,
 \qquad
 \cF_j=\|e^{\Phi_\gamma(r,D_\zeta)}\Delta_jf_c\|_{S_e},
\end{equation}
where $\norm g_{S_e}=\langle S_eg,g\rangle^{1/2}$ is the norm defined by $S_e$, and we differentiate $E_j$ in $r$. The reference term of
highest order
has zero real part, since $S_e\mathbf G_e$ is Hermitian. Since the damping
$d_\gamma(r,\lambda)=-\pa_r\Phi_\gamma(r,\lambda)$ of
\eqref{eq:phase-damping} increases with $\lambda$, the weight contributes at
most $-2d_\gamma(r,\lambda_j)E_j$. The part of
$\mathbf G_{\mathbf b}-\mathbf G_e$ that fails to be self-adjoint in
the $S_e$ inner product contributes
at most $C_1\chi_c\lambda_jE_j$.
% AUTHOR QUERY: Prove the weighted principal estimate and the summation
% of cross-band terms in \eqref{app:minus:eq:center-energy}. Analytic
% conjugation of a variable principal coefficient can still have order
% one: the kernel of [e^{sigma Lambda},a]D e^{-sigma Lambda} contains
% m a_{n-m}(e^{sigma(<n>-<m>)}-1). The definition of chi_c controls
% only the unweighted nonskew part.
The commutators with the dyadic projections and with the coefficients,
as well as the difference between the Weyl and the standard
quantization, all carry the factor in
\eqref{app:minus:eq:weyl-left-difference}, so
\cref{lem:macro-native-calculus} bounds them at order zero. The terms
coming from the frame, from the rotation, and from the principal symbol
also have order zero, by \cref{app:minus:lem:current-reduction}.
Collecting in $\cH_j$ the terms with a higher regularity
coefficient norm and the solution at the base index, we obtain
\begin{equation}\label{app:minus:eq:center-energy}
 \frac{\dd}{\dd r}E_j
 \le\left(C_0+C_1\chi_c\lambda_j
          -2d_\gamma(r,\lambda_j)\right)E_j
       +2E_j^{1/2}\cF_j+\cH_j,
\end{equation}
where the bounds for $\cH_j$ are summable with the Sobolev weights, and
\begin{equation}\label{app:minus:eq:phase-damping}
 d_\gamma(r,\lambda)
 =\frac{\gamma r\lambda^2}{\sqrt{1+r^2\lambda^2}}
\end{equation}
is the damping of the weight, as in \eqref{eq:phase-damping}. If
$C_1\chi_c<2\gamma$, the change of variable $t=r\lambda$ gives
\begin{equation}\label{app:minus:eq:center-positive-part}
 \int_0^R\left[C_1\chi_c\lambda
       -2d_\gamma(r,\lambda)\right]_+\dd r
 =\int_0^{R\lambda}
  \left[C_1\chi_c-\frac{2\gamma t}{\sqrt{1+t^2}}\right]_+\dd t
 \le C\frac{\chi_c^2}{\gamma},
\end{equation}
so the positive part of the growth rate is integrable uniformly in
$\lambda$. Indeed, for $0\le a<2\gamma$, integration over
$0\le t<\infty$ gives
\[
 \int_0^\infty\left[a-\frac{2\gamma t}{\sqrt{1+t^2}}\right]_+\dd t
 =2\gamma-\sqrt{4\gamma^2-a^2}
 =\frac{a^2}{2\gamma+\sqrt{4\gamma^2-a^2}}
 \le\frac{a^2}{2\gamma}.
\]
We apply this with $a=C_1\chi_c$. The same bound holds when the integration starts
at a point inside the cap instead of at $r=0$. Gr\"onwall's inequality
and summation over the dyadic blocks then
give the estimates for the initial value problem and for the outgoing
Cauchy data with no loss of derivatives.

For the upper and lower pairs $x=(x_1,x_2)$ and $y=(y_1,y_2)$
of \eqref{app:minus:eq:affine-native-norm}, the full affine equations are
\begin{equation}\label{app:minus:eq:affine-system}
 \begin{aligned}
 \pa_r x&=\Lambda y+a_{11}x+a_{12}y+f,\\
 \pa_r y&=B_{21}x+\Lambda^{-1}a_{21}x+a_{22}y+g,
 \end{aligned}
\end{equation}
where the coefficients $a_{ij}$ have order zero and
$\norm{B_{21}}_{\rm low}\le M_B\delta_A$. On one dyadic shell, we use
the \emph{modified energy}
\begin{equation}\label{app:minus:eq:affine-energy}
 E_{A,\lambda}=\mu_\lambda^2\abs x^2+\lambda\abs y^2,
 \qquad \mu_\lambda^2=\delta_A+\lambda^{-1},
\end{equation}
whose weights are those of the affine rows of \cref{tab:block-weights}.
In the variables $\tilde x=\mu_\lambda x$ and $\tilde y=\sqrt\lambda y$,
we bound the two off-diagonal coefficients by
\begin{equation}\label{app:minus:eq:affine-rate}
 \mu_\lambda\sqrt\lambda=\sqrt{1+\delta_A\lambda},
 \qquad
 \frac{\sqrt\lambda\norm{B_{21}}}{\mu_\lambda}
       \le C\sqrt{\delta_A\lambda},
\end{equation}
and all other terms have order zero, so the unweighted energy grows at
rate at most $C(1+\sqrt{\delta_A\lambda})$. Including the damping
\eqref{app:minus:eq:phase-damping} and repeating the preceding
integration, we obtain
\begin{equation}\label{app:minus:eq:affine-positive-part}
 \int_0^R\left[C\sqrt{\delta_A\lambda}
         -2d_\gamma(r,\lambda)\right]_+\dd r
 \le C\frac{\delta_A}{\gamma},
\end{equation}
provided that $C_*\sqrt{\delta_A}<\gamma$. We have therefore obtained
the estimate for the affine component in the norms
\eqref{app:minus:eq:affine-native-norm}. The required inequality
explains why $\delta_A$ must be fixed before $\eta_*$.

The order-zero coupling between the center and affine
components is not directly controlled by the modified
energy. We remove its leading part by a change of variables.
For either sign of the cell frequency, the principal matrix is
\begin{equation}\label{app:minus:eq:sylvester-principal}
 A_0=\diag(A_C,J_A),\qquad J_A^2=0,\qquad
 \norm{A_C^{-1}}+\norm{A_C^{-2}}\le C,
\end{equation}
where $A_C$ is the center matrix. We write $C_{CA}$ for the
coefficient taking affine variables into the center equation and
$C_{AC}$ for the coefficient taking center variables into the affine
equation. We set
\begin{equation}\label{app:minus:eq:sylvester-solutions}
 \begin{aligned}
 X&=-A_C^{-1}C_{CA}-A_C^{-2}C_{CA}J_A,\\
 Y&= C_{AC}A_C^{-1}+J_A C_{AC}A_C^{-2}.
 \end{aligned}
\end{equation}
Multiplying out and using $J_A^2=0$, we obtain
\begin{equation}\label{app:minus:eq:sylvester-identities}
 A_C X-X J_A=-C_{CA},\qquad J_A Y-Y A_C=-C_{AC}.
\end{equation}
We write the original center and affine unknown as
$(c,x,y)=\mathcal T\widetilde U$, where $\widetilde U$ is the new
unknown and
\begin{equation}\label{app:minus:eq:sylvester-transform}
 \mathcal T=I+\operatorname{Op}\left(
  \Lambda^{-1}\begin{pmatrix}0&X\\Y&0\end{pmatrix}\right).
\end{equation}
This change of variables cancels the two couplings to leading order.
Indeed, with the
convention \eqref{app:minus:eq:quantization-convention}, the conjugated
off-diagonal terms of order zero are $A_CX-XJ_A+C_{CA}$ and
$J_AY-YA_C+C_{AC}$, both of which are zero by
\eqref{app:minus:eq:sylvester-identities}. Every remaining term of the
composition contains either the factor $\Lambda^{-1}$ of
\eqref{app:minus:eq:sylvester-transform} or a frequency difference from
\eqref{app:minus:eq:weyl-left-difference}, so the moment estimate of
\cref{lem:macro-native-calculus} leaves a remainder
$R_{\mathbf b,H}^{(-1)}$ of order minus one. In the norms
\eqref{app:minus:eq:affine-native-norm}, the four ratios of weights that
occur for a map of order minus one are
\begin{equation}\label{app:minus:eq:half-order-ledger}
 \frac{\lambda^{-1}}{\mu_\lambda}\le\lambda^{-1/2},\quad
 \frac{\lambda^{-1}}{\lambda^{1/2}}\le\lambda^{-3/2},\quad
 \mu_\lambda\lambda^{-1}\le C\lambda^{-1},\quad
 \lambda^{1/2}\lambda^{-1}=\lambda^{-1/2},
\end{equation}
and we conclude that
\begin{equation}\label{app:minus:eq:sylvester-remainder}
 \norm{R_{\mathbf b,H}^{(-1)}U}_{Y_{a,-}^s}
 \le C_s\left[
  (\Lambda_J^{-1/2}+\chi_-)\norm U_{X_{a,-}^s}
 +\norm{\mathbf b-\mathbf b_\rho}_{\cB^{s+k}}
                         \norm U_{X_{a,-}^{s_*}}\right].
\end{equation}
The change of variables and its inverse satisfy the same estimate on the
Cauchy data.

We let $K_{0,\mathbf b}$ be the direct sum of the inverses of the
initial value problems for the center, affine, and remaining auxiliary
components in these variables. The bounds
\eqref{app:minus:eq:center-positive-part} and
\eqref{app:minus:eq:affine-positive-part} give a solution for every
source and every incoming Cauchy datum. The remaining error is small by
\eqref{app:minus:eq:state-cutoff-choices} and
\eqref{app:minus:eq:sylvester-remainder}, so the inverse is
\begin{equation}\label{app:minus:eq:annulus-neumann}
 G_{a,-,\mathbf b}
  =\bigl(I+K_{0,\mathbf b}R_{\mathbf b,H}^{(-1)}\bigr)^{-1}
        K_{0,\mathbf b}
  =K_{0,\mathbf b}
     \bigl(I+R_{\mathbf b,H}^{(-1)}K_{0,\mathbf b}\bigr)^{-1},
\end{equation}
with the changes of variables undone and the triangular equations for
the remaining components included.

We first impose a finite upper cutoff in
\eqref{app:minus:eq:annulus-neumann}, and the uniform energy estimates
then give a subsequence converging weak-* in the weighted energy space
of type $L^\infty$ in $r$ with values in $H^{s_*}$. The equation bounds
the radial derivative one index lower, so the unknown of the first-order
system converges in $C_w([0,R];H^{s_*-1})$, the space of weakly
continuous functions of $r$, and its value and balanced
conormal coordinates therefore converge at both endpoints. We pass the
equation to the limit by testing against smooth functions with finitely
many cell modes. For zero data, the center and affine energies vanish,
and then so do the finite coordinates, so the solution is unique and the
limit does not depend on the subsequence. Uniqueness also identifies
$u_a$ with
$G_{a,-,\mathbf b}(\cL_{a,-,\mathbf b}u_a,\gamma_{-,\mathbf b}u_a)$ for
every $u_a\in X_{a,-}^s$, where $\cL_{a,-,\mathbf b}$ and
$\gamma_{-,\mathbf b}$ are the center and affine equations on the
annulus and the Cauchy trace on its inner boundary at the state
$\mathbf b$, as in \eqref{app:minus:eq:annulus-map} at the reference.
This proves the left inverse identity, and
the right identity follows from \eqref{app:minus:eq:annulus-neumann}.

\subsection{The inverse identities and the tame estimates}

We now verify that the constructions on the cap and on the annulus give
the asserted inverse of the coupled problem. For
$F=(f_c,f_a,g^{\rm sp})$, their inverse identities give
\begin{equation}\label{app:minus:eq:right-identity-calculation}
 \begin{aligned}
 N_{-,\rho}^{\rm diag}G_{-,\rho}^{\mathbb N,\rm diag}F
 =\bigl(f_c,f_a,&g^{\rm sp}
       +T_{-,\mathbf b_\rho}
          \Gamma_{c,-,\rho}^{\rm sp}
             G_{c,-,\rho}^{\rm diag}f_c-T_{-,\mathbf b_\rho}
          \Gamma_{c,-,\rho}^{\rm sp}
             G_{c,-,\rho}^{\rm diag}f_c\bigr)=F.
 \end{aligned}
\end{equation}
Conversely, if we apply \eqref{app:minus:eq:joined-inverse} to
$N_{-,\rho}^{\rm diag}(u_c,u_a)$, the Cauchy datum for the problem on
the annulus reduces to
\begin{equation}\label{app:minus:eq:left-identity-calculation}
 \begin{aligned}
 \gamma_{-,\rho}u_a
 -T_{-,\mathbf b_\rho}\Gamma_{c,-,\rho}^{\rm sp}u_c
 +T_{-,\mathbf b_\rho}\Gamma_{c,-,\rho}^{\rm sp}
             G_{c,-,\rho}^{\rm diag}
             \cA_{c,-,\rho}^{\rm diag}u_c
 =\gamma_{-,\rho}u_a,
 \end{aligned}
\end{equation}
and the left inverse identities on the cap and on the annulus then
recover $u_c$ and $u_a$. This proves both inverse identities at the
reference. The Neumann formulas \eqref{app:minus:eq:whole-minus-neumann},
with the range factors in \eqref{app:minus:eq:minus-shear-realization},
then give \eqref{app:minus:eq:two-sided-joined}.

The estimate on the cap, the center and affine energies, and
\eqref{app:minus:eq:sylvester-remainder} bound the reference inverse
$G_{-,\rho}^{\mathbb N,\rm diag}$, and we combine that bound with
\eqref{app:minus:eq:current-full-shell-bound} in the Neumann formula
\eqref{app:minus:eq:whole-minus-neumann}. Applying the bounded range
factors $\mathfrak P_{-,\mathbf b}^s$ and
$\mathscr S_{Y,-,\mathbf b}^{-1}$, we then obtain
\eqref{app:minus:eq:undifferentiated-tame}.
% AUTHOR QUERY: Supply the higher-regularity argument for the current
% coupled inverse on one fixed state neighborhood. A Neumann series
% controlled only at the base index and a C_s delta_0 perturbation bound
% at index s do not alone prove this tame estimate for every s.
The estimate \eqref{app:minus:eq:localized-cap-leakage} will be used when the low and
high frequencies are coupled.

For the state derivatives, we work on the reconstructed space. Each derivative of
$\mathcal R_{-,\mathbf b}^{\mathbb N}$ places at most one cell
derivative on the solution. To see this for the center equation, let
$L^c_{\mathbf b}$ be the operator obtained by moving the homogeneous
terms of \eqref{app:minus:eq:center-evolution} to the left, so that
\[
 L^c_{\mathbf b}c
 =\pa_r c-i\operatorname{Op}^{w}(\mathbf G_{\mathbf b}|\xi|)c
       -B_{\mathbf b}c=f_c.
\]
Holding $f_c$ fixed and differentiating gives
\begin{equation}\label{app:minus:eq:differentiated-center}
 L^c_{\mathbf b}c^{(m)}
 =-\sum_{\varnothing\ne I\subset\{1,\ldots,m\}}
   \bigl(D_{\mathbf b}^{\abs I}L^c_{\mathbf b}[h_I]\bigr)
       c^{(m-\abs I)}[h_{I^c}],
\end{equation}
where $c^{(m)}$ denotes the $m$-th derivative of the solution $c$ with
respect to $\mathbf b$, $h_I$ the tangent directions with indices in $I$, and $I^c$ the
complement of $I$. Every factor
$D_{\mathbf b}^{\abs I}L^c_{\mathbf b}$ has at most one cell derivative,
so induction from the energy estimate loses at most one Sobolev index
per differentiation in $\mathbf b$, which gives $\mu_m\le m$. The reference inverse on the cap, the traces, the finite
equations, and the changes of variables
\eqref{app:minus:eq:sylvester-transform} have tame derivatives with respect to $\mathbf b$
of fixed order, so Leibniz' rule applied to
\eqref{app:minus:eq:joined-inverse} and to the resolvent factors in
\eqref{app:minus:eq:whole-minus-neumann}, together with the bounded range
factors in \eqref{app:minus:eq:minus-shear-realization}, proves
\eqref{app:minus:eq:all-order-tame}. The weight is fixed, so these derivatives use the same
strip of analyticity.

Complex conjugation in Cartesian variables pairs $(m,n)$ with $(-m,-n)$,
preserves the contour around the four center roots, and exchanges the
paired affine chains and cell frequencies. The changes of variables, the
Neumann series, and the triangular equations commute with this
conjugation, and, by uniqueness, so do the inverses and the traces.

The estimates depend on the separation of the four center roots from
zero, and they do not require these roots to be separated from one
another. The parameter $\delta_A$ is fixed, the factors in
\eqref{app:minus:eq:rho-factors} remain bounded for
$0\le\rho\le\rho_0<1/4$, and the symmetrizer at the circle remains
positive when the center speeds coincide. Thus, the bounds are uniform
through $\rho=0$. For compatible Cartesian sources,
\cref{prop:all-order-adapters} supplies the fixed loss of the source
before this inverse is applied. This completes the proof of
\cref{prop:minus-diagonal}.
\section{The elliptic boundary problem}\label{app:positive}
\label{app:positive:section}

We solve the elliptic part of the high-frequency problem by coupling
second-order equations on the cap and annulus through their values and
conormal derivatives on the common interface. As in \cref{app:minus},
we first invert the operator $\widehat D_{+,\mathbf b}$ defined below.
We relate it to the positive diagonal block of the full operator in
\cref{app:macro-cross}.

For this \emph{transmission problem}, we use a coercive quadratic form
on the pair of fields, including both parities of the positive cap
modes. At energy regularity, corresponding to $s=0$ below, an interior
source acts on test functions with zero boundary values and need not
determine a physical conormal derivative. We therefore prescribe the
boundary functionals supplied by the form. At higher regularity, we
recover the physical conormal and Robin data by subtracting the
boundary contributions of the sources.

Throughout, we fix $0<\rho_0<1/4$, the phase of \eqref{eq:phase-weight},
\begin{equation}\label{app:positive:eq:phase}
 \Phi_\gamma(r,\lambda)=\sigma_0\lambda
 -\gamma\bigl(\sqrt{1+r^2\lambda^2}-1\bigr),
 \qquad W_\gamma=e^{\Phi_\gamma(r,\Lambda)},
\end{equation}
and the projection $H_J$ of \eqref{eq:A-exact-cell-split} onto the high
cell frequencies. For $\mathbf b=(\rho,a,\eps,p)$, we let
$\mathbf b_\rho=(\rho,a_\rho,\eps_\rho,p_\rho)$ be the constant ellipse
with the same ellipse parameter, as in \cref{sec:block}, and we put
\begin{equation}\label{app:positive:eq:coefficient-distance}
 \eta_+(\mathbf b)=
 \norm{a-a_\rho}_{\cA_\gamma^{s_*+k_+}}
 +\abs{\eps-\eps_\rho}+\abs{p-p_\rho}.
\end{equation}
Here, $k_+$ is the fixed number of additional derivatives of the
coefficients required in the estimates, independent of $s$, as $k_c$ in
\cref{app:minus}. At index $s$, we use
the norm
\begin{equation}\label{app:positive:eq:augmented-norm}
 \norm{\mathbf b-\mathbf b_\rho}_{\cB^s}
 :=\norm{a-a_\rho}_{\cA_\gamma^s}
   +\abs{\eps-\eps_\rho}+\abs{p-p_\rho}.
\end{equation}
These norms measure the distance from the constant ellipse in the
components $a$, $\eps$, and $p$. We choose
a closed set of coefficients and $J$ such that
\begin{equation}\label{app:positive:eq:positive-choices}
 C(\gamma+\gamma^2)+C\eta_*<\frac{c_+}{2},
 \qquad 0<\rho_0<\frac14,
\end{equation}
where $\eta_*$ bounds $\eta_+(\mathbf b)$ on this set and $c_+$ is the
coercivity constant of the positive form
\eqref{appref:odd-positive-form} at the constant ellipse. These
choices are independent of an upper Fourier
cutoff and of the number $N$ of periods of the field.

\subsection{Function spaces and the inverse for the positive block}

On the charts of the cap and of the annulus, we set
\begin{equation}\label{app:positive:eq:local-spaces}
 \begin{aligned}
 V_{c,+}^s&=H_J^X\mathscr V_{+,{\rm cap}}^s,
 &Y_{c,+}^s&=H_J^Y\mathscr Y_{+,{\rm cap}}^s,\\
 V_{a,+}^s&=H_J^X\mathscr E_{+,a,\gamma}^s,
 &Y_{a,+}^s&=H_J^Y(\mathscr E_{+,a,\gamma}^s)^*_{\rm bulk} .
 \end{aligned}
\end{equation}
Here, $\mathscr V_{+,{\rm cap}}^s$ and $\mathscr Y_{+,{\rm cap}}^s$ are
the weighted spaces of the positive modes on the cap, whose orders are
given below, and $\mathscr E_{+,a,\gamma}^s$ is the weighted energy space of
the positive block on the annulus. The space
$(\mathscr E_{+,a,\gamma}^0)^*_{\rm bulk}$ is the dual of the
subspace of $\mathscr E_{+,a,\gamma}^0$ whose fields have zero value on
both boundaries. At index $s$, we increase the Sobolev order of the source
norm by $s$, from $H^{-1}$ to $H^{s-1}$, keeping the analytic weight fixed.
This defines $(\mathscr E_{+,a,\gamma}^s)^*_{\rm bulk}$. We call the spaces $V$ on which the quadratic form is defined the
\emph{form spaces}, and write $Y$ for the corresponding source spaces.
The cap and annular form spaces have the same regularity. On the cap,
the local orders before applying the weight are
\begin{equation}\label{app:positive:eq:cap-form-orders}
 V_{c,+}^s\simeq H^{s+1}_{+}(\Omega_c),
 \qquad
 Y_{c,+}^s\simeq H^{s-1}_{+}(\Omega_c).
\end{equation}
Here, $\Omega_c$ is the scaled cap and $+$ restricts to the positive
angular modes. At $s=0$, the cap form space is Cartesian $H^1$, its
source space $H^{-1}$ is the dual of the zero-Dirichlet subspace, and
the annular form space is weighted $H^1$, as in \cref{tab:block-weights}.
For the operator domains, $\gamma_Du$ denotes the value on the
interface. We use $Z_D^s$ for this datum, of order $H^{s+1/2}$, and
$Z_B^s$ for the outer Robin datum, of order $H^{s-1/2}$, with the
analytic weight included. As domains of the operator, we take
\begin{equation}\label{app:positive:eq:operator-domains}
 \begin{aligned}
 X_{c,+}^s
  &=\{u\in V_{c,+}^s:P_{c,+,\mathbf b}u\in Y_{c,+}^s,\
            \gamma_Du\in Z_D^s\},\\
 X_{a,+}^s
  &=\{u\in V_{a,+}^s:P_{a,+,\mathbf b}u\in Y_{a,+}^s,\
       \gamma_Du\in Z_D^s,\
       \mathcal R_{++,\mathbf b}^{\rm sp}
          (u,P_{a,+,\mathbf b}u)\in Z_B^s\}.
 \end{aligned}
\end{equation}
Here, $P_{c,+,\mathbf b}u_c=f_c$ and
$P_{a,+,\mathbf b}u_a=f_a$ are the interior equations of the positive
block on the cap and annulus. The operator
$\mathcal R_{++,\mathbf b}^{\rm sp}$ records the Robin datum using the
weak conormal residual supplied by the quadratic form. For smooth fields,
this is the physical Robin datum plus the boundary contribution of the
interior source, as in the split Cauchy data
\eqref{app:minus:eq:cap-trace-physicalization}. We give its weak
definition below.
Thus, the operator domains impose regularity on the equation and traces
in addition to the form regularity $H^{s+1}$. In particular,
\begin{equation}\label{app:positive:eq:form-operator-embeddings}
 X_{c,+}^0\oplus X_{a,+}^0
 \hookrightarrow V_{c,+}^0\oplus V_{a,+}^0
 \hookrightarrow L^2_{\rm bulk}.
\end{equation}
We call a source a \emph{full source} if it acts canonically on Poisson
extensions of boundary values. An index is \emph{regular} if every
source at that index has this property. For example, $s\ge1$ is regular
because the interior sources belong to $H^{s-1}\subset L^2$.
At $s=0$, we instead use the extension of functionals in
\eqref{app:positive:eq:bulk-trace-splitting} below.

The spaces $Z_D^s$, $Z_N^s$, and $Z_B^s$ contain the value on the
interface, the residual on the interface, and the residual on the outer
boundary, respectively. At regular indices, the latter two also contain
the physical conormal and Robin data, which are related to the residuals
by the change of data below. A \emph{Robin datum} prescribes a linear combination of the conormal
derivative and the value on the outer boundary, as required by the
boundary condition of the positive block.
The norms of these spaces are
\begin{equation}\label{app:positive:eq:balanced-traces}
 \norm z_{Z_D^s}\simeq\norm{\Lambda^{1/2}z}_{H^s},
 \qquad
 \norm z_{Z_N^s}\simeq\norm{\Lambda^{-1/2}z}_{H^s},
 \qquad
 \norm z_{Z_B^s}\simeq\norm{\Lambda^{-1/2}z}_{H^s},
\end{equation}
with the analytic weight on the traces and the spectral projections of
the principal symbol understood, as in \eqref{app:minus:eq:cap-trace-space}.
These are the trace orders associated with a solution in $H^{s+1}$,
namely $H^{s+1/2}$ for the value and $H^{s-1/2}$ for the conormal or
Robin data.

\smallskip\noindent\textit{Data for the weak problem.}
To formulate the weak problem at $s=0$, we extend the interior sources
to test fields with nonzero boundary values. We write
$\gamma_\Sigma^c$ and $\gamma_\Sigma^a$ for the traces on the
interface and $\gamma_\partial^a$ for the trace on the outer
boundary. We let $\mathcal P_\Sigma^c$, $\mathcal P_\Sigma^a$, and
$\mathcal P_\partial^a$ be their localized Poisson extensions from
\cref{app:fixed-graphs}, with disjoint supports at the interface and at
the outer boundary. We set
\begin{equation}\label{app:positive:eq:bulk-trace-splitting}
 \begin{aligned}
 R_c^0&:=I-\mathcal P_\Sigma^c\gamma_\Sigma^c,\\
 R_a^0&:=I-\mathcal P_\Sigma^a\gamma_\Sigma^a
             -\mathcal P_\partial^a\gamma_\partial^a,\\
 \operatorname{Ext}_c f_c(v_c)&:=f_c(R_c^0v_c),
 &\operatorname{Ext}_a f_a(v_a)&:=f_a(R_a^0v_a).
 \end{aligned}
\end{equation}
The source $f_c$ acts on test fields on the cap with zero value on the
interface, and $f_a$ acts on test fields on the annulus with zero value
at both boundaries. Formula \eqref{app:positive:eq:bulk-trace-splitting}
extends these functionals to the full form spaces. For a test pair $V=(v_c,v_a)$ with
$\gamma_\Sigma^cv_c=\gamma_\Sigma^av_a$, the resulting functional is
uniquely given by
\begin{equation}\label{app:positive:eq:variational-data-splitting}
 \operatorname{Ext}_c f_c(v_c)+\operatorname{Ext}_a f_a(v_a)
 +\langle g_N^{\rm sp},\gamma_\Sigma V\rangle
 +\langle\mathsf V_{B,\mathbf b}^{-1}g_B^{\rm sp},
              \gamma_\partial^av_a\rangle .
\end{equation}
Here, $\gamma_\Sigma V$ denotes the common value of $v_c$ and $v_a$ on
the interface. The coefficient $\mathsf V_{B,\mathbf b}$ multiplies the outward
conormal derivative in the physical Robin boundary operator, and is invertible.
Its inverse converts the prescribed boundary datum to the normalization
in Green's identity.
The last two terms prescribe the boundary data for the weak problem.
We call these independent functionals the \emph{split residuals} at the
interface and outer boundary. We relate them to the physical conormal
derivatives below. Conjugating the splitting
\eqref{app:positive:eq:variational-data-splitting} by the isomorphisms
that shift the index of the scale, we obtain the same convention at
every $s$.

If $u_c$ satisfies its interior equation, its weak split residual is
\begin{equation}\label{app:positive:eq:canonical-weak-conormal}
 \langle\mathcal N_{c,\Sigma}^{\rm sp}(u_c,f_c),\varphi\rangle
 :=\mathfrak c_{+,\mathbf b}
       (u_c,\mathcal P_\Sigma^c\varphi)
   -\operatorname{Ext}_cf_c(\mathcal P_\Sigma^c\varphi),
\end{equation}
where $\mathfrak c_{+,\mathbf b}$ is the form on the cap at the state. Writing $\mathfrak b_{+,\mathbf b}^{\rm int}$ for the interior part
of the annular form, with the outer Robin boundary term omitted, we set
\[
 \begin{aligned}
 \langle\mathcal N_{a,\Sigma}^{\rm sp}(u_a,f_a),\varphi\rangle
 &=\mathfrak b_{+,\mathbf b}^{\rm int}(u_a,\mathcal P_\Sigma^a\varphi)
   -\operatorname{Ext}_af_a(\mathcal P_\Sigma^a\varphi),\\
 \langle\mathcal N_{a,\partial}^{\rm sp}(u_a,f_a),\psi\rangle
 &=\mathfrak b_{+,\mathbf b}^{\rm int}(u_a,\mathcal P_\partial^a\psi)
   -\operatorname{Ext}_af_a(\mathcal P_\partial^a\psi).
 \end{aligned}
\]
For the chosen extensions, the subtracted terms vanish by \eqref{app:positive:eq:bulk-trace-splitting}. Since the form
and the extensions are bounded, these residuals are functionals of order
$H^{-1/2}$. The interior equation also shows that the residual is
unchanged if we add a test function with zero trace to the extension.
This defines the split residual even when the source does not determine
a physical conormal derivative.

\smallskip\noindent\textit{Recovery of physical boundary data.}
For full sources $\ell_c$ and $\ell_a$ at a regular index, we define
their boundary contributions, as in \cref{app:minus}, by the transposes
\begin{equation}\label{app:positive:eq:source-moment-maps}
 \begin{aligned}
 \langle\mathsf M_{c,\Sigma}^s\ell_c,\varphi\rangle
   &:=\langle\ell_c,\mathcal P_\Sigma^c\varphi\rangle,\\
 \langle\mathsf M_{a,\Sigma}^s\ell_a,\varphi\rangle
   &:=\langle\ell_a,\mathcal P_\Sigma^a\varphi\rangle,\\
 \langle\mathsf M_{a,\partial}^s\ell_a,\psi\rangle
   &:=\langle\ell_a,\mathcal P_\partial^a\psi\rangle.
 \end{aligned}
\end{equation}
Here, $\mathsf M^s=(\mathcal P)'$ is the transpose with respect to the
weighted dualities. It differs from the unweighted Hilbert adjoint. By Green's
identity, the split residuals are related to the outward physical
conormal derivatives by
\begin{equation}\label{app:positive:eq:split-physical-conormal}
 \begin{aligned}
 \mathcal N_{c,\Sigma}^{\rm sp}
   &=\gamma_{N,\mathbf b}^{c,{\rm out}}u_c
       +\mathsf M_{c,\Sigma}^s\ell_c,\\
 \mathcal N_{a,\Sigma}^{\rm sp}
   &=\gamma_{N,\mathbf b}^{a,{\rm out}}u_a
       +\mathsf M_{a,\Sigma}^s\ell_a,\\
 \mathcal N_{a,\partial}^{\rm sp}
   &=\gamma_{N,\mathbf b}^{a,\partial,{\rm out}}u_a
       +\mathsf M_{a,\partial}^s\ell_a.
 \end{aligned}
\end{equation}
Here, $\gamma_{N,\mathbf b}^{a,\partial,{\rm out}}$ is the conormal
derivative on the outer boundary of the annulus. We write
$\mathsf V_{B,\mathbf b}$ and $\mathsf F_{B,\mathbf b}$ for the
coefficients of the conormal derivative and of the value in the operator on
the outer boundary. For the unnormalized boundary operator $\cB_{++,\mathbf b}$ of the
positive block on the outer boundary at the state, whose coefficient of
the conormal derivative we denote by $\mathsf E_{\mathbf b}$, this gives
$\mathsf V_{B,\mathbf b}=\mathsf E_{\mathbf b}$, while
$\mathsf V_{B,\mathbf b}=I$ after normalization by that coefficient.
Wherever the physical conormal derivative is
defined, we write
\[
 \mathcal R_{++,\mathbf b}^{\rm ph}u_a
 :=\mathsf V_{B,\mathbf b}
       \gamma_{N,\mathbf b}^{a,\partial,{\rm out}}u_a
       +\mathsf F_{B,\mathbf b}\gamma_\partial^au_a,
\]
and, at every regularity level, we define the \emph{boundary operator in split form} by
\begin{equation}\label{app:positive:eq:split-outer-row}
 \mathcal R_{++,\mathbf b}^{\rm sp}(u_a,f_a)
 :=\mathsf V_{B,\mathbf b}
       \mathcal N_{a,\partial}^{\rm sp}(u_a,f_a)
       +\mathsf F_{B,\mathbf b}\gamma_\partial^au_a.
\end{equation}
When $s$ is a regular index and the source is full, we have
\begin{equation}\label{app:positive:eq:split-outer-row-physical}
 \mathcal R_{++,\mathbf b}^{\rm sp}(u_a,\ell_a)
 =\mathcal R_{++,\mathbf b}^{\rm ph}u_a
     +\mathsf V_{B,\mathbf b}\mathsf M_{a,\partial}^s\ell_a,
\end{equation}
with the factors in the displayed order.

We order the data as $F=(\ell_c,\ell_a,g_D,g_N,g_B)$, where $g_D$
is the jump of the values on the interface, $g_N$ is the datum of the
conormal condition on the interface, and $g_B$ is the datum on the
outer boundary. In this order, we define the map from physical to split
residuals at a regular index by
\begin{equation}\label{app:positive:eq:residual-shear}
 \begin{aligned}
 \mathfrak S_{+,\mathbf b}^{\rm res,s}
   (\ell_c,\ell_a,g_D,g_N^{\rm ph},g_B^{\rm ph})
  :=\big(&\ell_c,\ell_a,g_D,g_N^{\rm ph}
      +\mathsf M_{c,\Sigma}^s\ell_c
      +\mathsf M_{a,\Sigma}^s\ell_a,g_B^{\rm ph}
      +\mathsf V_{B,\mathbf b}
          \mathsf M_{a,\partial}^s\ell_a\big),
 \end{aligned}
\end{equation}
which adds the source contributions to the physical boundary data
to obtain the split residuals. It leaves the sources and the jump of the
values unchanged. Subtracting the same contributions gives the inverse. Both interface
contributions have a plus sign because the conormal derivatives point
outward on both sides. At the index
$s=0$, there is no such map $\mathfrak S_{+,\mathbf b}^{\rm res,0}$ on
the data alone, and \eqref{app:positive:eq:variational-data-splitting}
fixes the split coordinates there directly.

The bounds for the Poisson extensions and for the coefficient of the
conormal derivative in the outer boundary operator give tame estimates for
$(\mathfrak S_{+,\mathbf b}^{\rm res,s})^{\pm1}$ and their derivatives
with respect to $\mathbf b$, from index $s$ to index $s$. The range bounds in \cref{prop:positive-diagonal} include the norms of
$\mathfrak S_{+,\mathbf b}^{\rm res,s}$ and its inverse. The transpose of this map
is triangular as well. Indeed, with $(v_c,v_a,\lambda_\Sigma,\lambda_B)$
denoting the elements of the dual of the range that are paired with
$(\ell_c,\ell_a,g_N,g_B)$, its nontrivial components on the sources are
\begin{equation}\label{app:positive:eq:residual-shear-transpose}
 \begin{aligned}
 v_c&\longmapsto v_c+\mathcal P_\Sigma^c\lambda_\Sigma,\\
 v_a&\longmapsto v_a+\mathcal P_\Sigma^a\lambda_\Sigma
       +\mathcal P_\partial^a\mathsf V_{B,\mathbf b}'\lambda_B.
 \end{aligned}
\end{equation}
Note that the Hilbert adjoint with respect to the weighted inner
products differs from this transpose by the Riesz maps of the weighted
spaces. Given a projection $P_{Y,+}^{\rm ph}$ on the range in physical
coordinates, we transport it to the split coordinates by
\begin{equation}\label{app:positive:eq:residual-shear-projector}
 P_{Y,+}^{\rm sp}
 =\mathfrak S_{+,\mathbf b}^{\rm res,s}P_{Y,+}^{\rm ph}
       (\mathfrak S_{+,\mathbf b}^{\rm res,s})^{-1}.
\end{equation}
This projection includes a component from the source to the trace, which
cannot be removed without an additional identity for these boundary contributions.

We now collect the solution and data spaces as
\begin{equation}\label{app:positive:eq:positive-spaces}
 \begin{aligned}
 \mathbb V_{+,\mathrm{ret}}^s&=V_{c,+}^s\oplus V_{a,+}^s,\\
 \bX_{+,\mathrm{ret}}^s&=X_{c,+}^s\oplus X_{a,+}^s,\\
 \bY_{+,\mathrm{ret}}^s&=Y_{c,+}^s\oplus Y_{a,+}^s
                    \oplus Z_D^s\oplus Z_N^s\oplus Z_B^s.
 \end{aligned}
\end{equation}
If $\gamma_\Sigma^au_a-\gamma_\Sigma^cu_c=g_D$, then the pair
\[
 (u_c,u_a)-(0,\mathcal P_\Sigma^ag_D)
\]
has matching values on the interface, because
$\gamma_\Sigma^a\mathcal P_\Sigma^ag_D=g_D$. Such pairs form the
\emph{form domain}
\begin{equation}\label{app:positive:eq:glued-form-domain}
 \mathbb V_{+,0}^0
 :=\{(v_c,v_a)\in V_{c,+}^0\oplus V_{a,+}^0:
          \gamma_j^Dv_a=\gamma_j^Dv_c\},
\end{equation}
where $\gamma_j^D$ is the trace of the value on the interface, at the
radius $r_j$. The Robin condition is imposed through the form
and is therefore not an additional constraint on the form domain. We also write
\begin{equation}\label{app:positive:eq:outgoing-trace-space}
 Z_+^s=H_J^{Z_+}\mathfrak T_{+,\gamma}^s
\end{equation}
for the product of the balanced value and conormal trace spaces on
the interface and outer boundary. Here, we retain the traces of each
field separately, rather than only their jump or sum in the prescribed
data. These separate traces are needed when coupling the positive and
minus blocks. The space $\mathfrak T_{+,\gamma}^s$ includes the
analytic weight, as in
\eqref{app:minus:eq:cap-trace-space}, and $H_J^{Z_+}$ is its projection
onto the high cell frequencies.

For fields $u_c$ and $u_a$ on the cap and on the annulus at a regular
index, we first define the operator using the physical values and conormal
derivatives at the state,
\begin{equation}\label{app:positive:eq:positive-operator}
 \begin{aligned}
 N_{+,\mathbf b}^{\rm cur,ph}(u_c,u_a)=\big(&P_{c,+,\mathbf b}u_c,
      P_{a,+,\mathbf b}u_a,\gamma_j^Du_a-\gamma_j^Du_c,
   \nabla_{a,\mathbf b}^{\rm out}u_a
        +\nabla_{c,\mathbf b}^{\rm out}u_c,
   \cB_{++,\mathbf b}u_a\big),
 \end{aligned}
\end{equation}
where $\nabla_{c,\mathbf b}^{\rm out}=\gamma_{N,\mathbf b}^{c,{\rm out}}$
and $\nabla_{a,\mathbf b}^{\rm out}=\gamma_{N,\mathbf b}^{a,{\rm out}}$
are the outward conormal derivatives from the cap and annulus.
Replacing the physical boundary data by the weak residuals gives an
operator defined at every index,
\begin{equation}\label{app:positive:eq:positive-split-operator}
 \begin{aligned}
 N_{+,\mathbf b}^{\rm cur,sp}(u_c,u_a)=\big(
   &P_{c,+,\mathbf b}u_c,P_{a,+,\mathbf b}u_a,
     \gamma_j^Du_a-\gamma_j^Du_c,\\
   &\mathcal N_{a,\Sigma}^{\rm sp}
       (u_a,P_{a,+,\mathbf b}u_a)
      +\mathcal N_{c,\Sigma}^{\rm sp}
       (u_c,P_{c,+,\mathbf b}u_c),\\
   &\mathcal R_{++,\mathbf b}^{\rm sp}
       (u_a,P_{a,+,\mathbf b}u_a)\big).
 \end{aligned}
\end{equation}
At the index $s=0$, we regard the two interior sources as functionals
on test fields with zero boundary trace before using them in
$\mathcal N^{\rm sp}$. At regular
indices with full sources, we relate the two operators by
\begin{equation}\label{app:positive:eq:positive-operator-coordinate-relation}
 N_{+,\mathbf b}^{\rm cur,sp}
 =\mathfrak S_{+,\mathbf b}^{\rm res,s}
       N_{+,\mathbf b}^{\rm cur,ph}.
\end{equation}

\smallskip\noindent\textit{Return to the reference coordinates.}
We let $\iota_+^X$ and $\pi_+^Y$ denote the inclusion of the positive
block and the projection onto it, as in \cref{sec:block}, and we let
$\mathscr R_{CN,\mathbf b,\rho}$ be the
isomorphism of the range
from \eqref{eq:A-current-row-normalized-inverse}. It acts by
$(\mathcal T_{\mathbf b,\rho}^{CN})^{-1}$ on the pair of the value and
the conormal derivative and by the identity on the other components. To
recover the physical conormal data, we first subtract the source
contributions using \eqref{app:positive:eq:residual-shear}. We then
express these data in the reference coordinates using
$\mathscr R_{CN,\mathbf b,\rho}$. At regular indices, this gives
\begin{equation}\label{app:positive:eq:normalized-natural-realization}
 N_{+,\mathbf b}
 :=\pi_+^Y\mathbb N_{\mathbf b}\iota_+^X
 =\mathscr R_{CN,\mathbf b,\rho}N_{+,\mathbf b}^{\rm cur,ph}
 =\mathscr R_{CN,\mathbf b,\rho}
    (\mathfrak S_{+,\mathbf b}^{\rm res,s})^{-1}
       N_{+,\mathbf b}^{\rm cur,sp}.
\end{equation}
At the index $s=0$, we work instead with the operator
$N_{+,\mathbf b}^{\rm cur,sp}$ defined by the form. We apply the
conversion in \eqref{app:positive:eq:normalized-natural-realization} to
the reconstructed field on the cap and on the annulus, and we do
not commute it with the shell projections $Q_j$.

We write $\iota_+^Y$ for the inclusion of the positive range. In
addition to the conversion $\mathfrak S_{+,\mathbf b}^{\rm res,s}$
between physical and split data, we use the positive diagonal block of
the map $\mathscr S_{Y,\mathbf b}$ from \cref{sec:block}, defined by
\begin{equation}\label{app:positive:eq:plus-diagonal-shear-compression}
 \mathscr S_{Y,+,\mathbf b}
 :=\pi_+^Y\mathscr S_{Y,\mathbf b}\iota_+^Y.
\end{equation}
This diagonal block is defined even if the full map does not preserve
the positive summand. For a positive range vector
$F_++\jmath_{Z,+,\mathbf b}z_+$, where $F_+$ has zero auxiliary
entries and $z_+$ also denotes its inclusion in the full auxiliary
space, we have
\[
 \mathscr S_{Y,+,\mathbf b}(F_++\jmath_{Z,+,\mathbf b}z_+)
 =F_+-\pi_+^Y\mathbb B_{\mathbf b}z_+
       +\jmath_{Z,+,\mathbf b}z_+,
\]
and its inverse has the same formula with a plus sign before
$\pi_+^Y\mathbb B_{\mathbf b}z_+$. We put
\begin{equation}\label{app:positive:eq:plus-shear-realization}
 \widehat D_{+,\mathbf b}
 :=\mathscr S_{Y,+,\mathbf b}N_{+,\mathbf b},
 \qquad
 \widehat G_{+,\mathbf b}
 :=G_{+,\mathbf b}^{\mathbb N}
       \mathscr S_{Y,+,\mathbf b}^{-1},
\end{equation}
where $G_{+,\mathbf b}^{\mathbb N}=N_{+,\mathbf b}^{-1}$ is constructed
below.

The component taking positive unknowns to minus equations remains
in the augmented operator. Formula \eqref{eq:A-natural-broken-formula} gives
\begin{equation}\label{app:positive:eq:plus-full-shear-cross-cancellation}
 \begin{aligned}
 \pi_-^Y\mathscr S_{Y,\mathbf b}\mathbb N_{\mathbf b}\iota_+^X
 &=\pi_-^Y\iota_{Y,\mathbf b}\mathsf J_{Y,\mathbf b}^{\rm raw}
       L_{\mathbf b}\mathsf R_{X,\mathbf b}\iota_+^X
    +\jmath_{Z,-,\mathbf b}q_-^Z
       \mathsf C_{X,\mathbf b}\iota_+^X,\\
 \pi_-^Y\mathbb B_{\mathbf b}\mathsf C_{X,\mathbf b}\iota_+^X
 &-\pi_-^Y\mathbb B_{\mathbf b}\mathsf C_{X,\mathbf b}\iota_+^X=0.
 \end{aligned}
\end{equation}
The off-diagonal block of
$\mathbb D_{\mathbf b}=\mathscr S_{Y,\mathbf b}\mathbb N_{\mathbf b}$ is
estimated once both diagonal blocks are inverted. The projection $q_-^Z$
selects the auxiliary data of the minus block, and this contribution
remains after the two interior source corrections cancel.

\begin{proposition}[The inverse of $\widehat D_{+,\mathbf b}$]\label{prop:positive-diagonal}
\label{app:positive:prop:diagonal}
For the coefficients fixed in \eqref{app:positive:eq:positive-choices}
and every regular index $s$ at which the boundary contributions
\eqref{app:positive:eq:source-moment-maps} are defined,
\begin{equation}\label{app:positive:eq:positive-isomorphism}
 \widehat D_{+,\mathbf b}:\bX_{+,\mathrm{ret}}^s
       \longrightarrow\bY_{+,\mathrm{ret}}^s
\end{equation}
is an isomorphism. We write
$\widehat G_{+,\mathbf b}=\widehat D_{+,\mathbf b}^{-1}$, and we let
$\Gamma_{+,\mathbf b}:\bX_{+,\mathrm{ret}}^s\to Z_+^s$ be the map that
records the remaining balanced traces in the space
\eqref{app:positive:eq:outgoing-trace-space}. For every $s\ge s_*$,
\begin{equation}\label{app:positive:eq:positive-tame}
 \begin{aligned}
 &\norm{\widehat G_{+,\mathbf b}F}_{\bX_{+,\mathrm{ret}}^s}
 +\norm{\Gamma_{+,\mathbf b}\widehat G_{+,\mathbf b}F}_{Z_+^s}\le C_s\left[
    \norm F_{\bY_{+,\mathrm{ret}}^s}
   +\bigl(1+\norm{\mathbf b-\mathbf b_\rho}_{\cB^{s+k_+}}\bigr)
       \norm F_{\bY_{+,\mathrm{ret}}^{s_*}}\right],
 \end{aligned}
\end{equation}
for arbitrary $F=(f_c,f_a,g_D,g_N,g_B)$ in $\bY_{+,\mathrm{ret}}^s$,
written in the normalized physical coordinates, that is, in the
coordinates of the range of $\widehat D_{+,\mathbf b}$. The data
supplied to the inverse of the form at $s=0$
are
\[
 F^{\rm sp}
 =\mathfrak S_{+,\mathbf b}^{\rm res,s}
    \mathscr R_{CN,\mathbf b,\rho}^{-1}
      \mathscr S_{Y,+,\mathbf b}^{-1}F,
\]
and \eqref{app:positive:eq:positive-tame} includes the bounds for all
three factors on the range, with $g_D$, $g_N$, and $g_B$ independent
data at a regular index.

For every $m\ge1$, there are finite integers $k_{+,m}$ and
$\mu_{+,m}$, independent of $s$, for which either
$Q_{+,\mathbf b}=\widehat G_{+,\mathbf b}$ or
$Q_{+,\mathbf b}=\Gamma_{+,\mathbf b}\widehat G_{+,\mathbf b}$ satisfies
\begin{equation}\label{app:positive:eq:positive-all-order}
 \begin{aligned}
 &\norm{D_{\mathbf b}^mQ_{+,\mathbf b}
          [h_1,\ldots,h_m]F}_{Z^s}\\
 &\quad\le C_{s,m}\Bigg[
  \prod_i\norm{h_i}_{\rm lo}\,
       \norm F_{\bY_{+,\mathrm{ret}}^{s+\mu_{+,m}}}\\
 &\qquad+\sum_i\norm{h_i}_{{\rm tan},s+k_{+,m}}
       \prod_{j\ne i}\norm{h_j}_{\rm lo}\,
       \norm F_{\bY_{+,\mathrm{ret}}^{s_*+\mu_{+,m}}}\\
 &\qquad+\bigl(1+\norm{\mathbf b-\mathbf b_\rho}_
                              {\cB^{s+k_{+,m}}}\bigr)
       \prod_i\norm{h_i}_{\rm lo}\,
       \norm F_{\bY_{+,\mathrm{ret}}^{s_*+\mu_{+,m}}}\Bigg],
 \end{aligned}
\end{equation}
where $Z^s=\bX_{+,\mathrm{ret}}^s$ for the solution map and
$Z^s=Z_+^s$ for the trace map, and
\begin{equation}\label{app:positive:eq:positive-tangent}
 \norm h_{{\rm tan},s}=\norm{\dot a}_{\cA_\gamma^s}
       +\abs{\dot\rho}+\abs{\dot\eps}+\abs{\dot p},
 \qquad
 \norm h_{\rm lo}=\norm h_{{\rm tan},s_*+k_{+,m}}.
\end{equation}
The constants in \eqref{app:positive:eq:positive-tame}--\eqref{app:positive:eq:positive-all-order} are uniform for
$0\le\rho\le\rho_0$ and both signs of the frequency. They apply to the
full high-frequency space and to every finite upper cutoff, with no
factor $1/\rho$. At real states, the operator, its inverse, and the trace
map preserve the real subspaces under the fixed identifications.
\end{proposition}

To prove \cref{prop:positive-diagonal}, we first use coercivity to
invert $N_{+,\mathbf b}^{\rm cur,sp}$ for split data at $s=0$. Higher
regularity then allows us to recover the inverse in physical coordinates.

\subsection{Coercivity at the circle and the boundary condition}

We first prove coercivity with the physical boundary conditions. In the
homogeneous equations of highest order, the interior sources and their boundary contributions vanish. The conversion \eqref{app:positive:eq:residual-shear}
therefore leaves the boundary data unchanged, so the same computation
applies to the split residuals.

The computation includes the endpoint $\rho=0$. In the scalar
coordinates of \eqref{appref:raw-scalar-rows}, the matrix through which
the ellipse enters is
\begin{equation}\label{app:positive:eq:H-rho}
 H_\rho(\theta)=\frac1{1-\rho^2}
 \begin{pmatrix}
  1-\rho\cos2\theta&\rho\sin2\theta\\
  \rho\sin2\theta&1+\rho\cos2\theta
 \end{pmatrix},
\end{equation}
and its eigenvalues are
$(1\mp\rho)^{-1}$, so
\begin{equation}\label{app:positive:eq:H-rho-bounds}
 \frac1{1+\rho_0}I\le H_\rho\le\frac1{1-\rho_0}I.
\end{equation}
We bound the terms coming from $H_\rho-I$, whose nonconstant
coefficients shift the angular mode by two, by $3\rho_0/(1-\rho_0)$ times the
form of highest order, so that
\begin{equation}\label{app:positive:eq:circle-coercivity}
 \Re\mathfrak b_\rho(v,v)
 \ge c_{\rho_0}\norm v_{\mathscr E_+^0}^2,
 \qquad
 c_{\rho_0}=\frac{1-4\rho_0}{1-\rho_0}>0.
\end{equation}
Here, $\mathfrak b_\rho$ is the quadratic form of the positive block on
the annulus at the constant ellipse $\mathbf b_\rho$, and
$\mathscr E_+^s$ is the energy scale of the positive modes on the
annulus without the weight, so that $\mathscr E_{+,a,\gamma}^s$ in
\eqref{app:positive:eq:local-spaces} is its weighted version. The restriction $\rho_0<1/4$ gives $c_{\rho_0}>0$, and all fixed-order
$\rho$-derivatives of $H_\rho$ are uniformly bounded. On the cap, the odd
positive modes are $\abs m\ge3$ and the even positive modes are
$\abs m\ge4$, and their forms with zero Dirichlet data are strictly
positive. The Lax-Milgram lemma therefore gives a solution for an
arbitrary source and an arbitrary value on the interface, and tangential
difference quotients together with the equation give its regularity.
The modes $m=\pm1$ and the
affine chains fixed on the axis do not enter this problem.

Dividing the boundary operator by its invertible coefficient of the
conormal derivative, we obtain
\begin{equation}\label{app:positive:eq:normalized-robin}
 \mathsf E_{\mathbf b}^{-1}\cB_{++,\mathbf b}u
   =(\nabla_{\mathbf b}+\mathsf S_{\mathbf b})u=h,
\end{equation}
where $\nabla_{\mathbf b}u$ is the outward conormal derivative,
$\mathsf S_{\mathbf b}$ acts on the boundary value of $u$, and
$h=\mathsf E_{\mathbf b}^{-1}g_B^{\rm ph}$ is the normalized prescribed
datum. We write $\mathsf S_e$ for $\mathsf S_{\mathbf b}$ at the
constant ellipse $\mathbf b_\rho$. For a constant $s_B>0$, we have
\begin{equation}\label{app:positive:eq:reference-robin-gap}
 \mathsf S_e=\mathsf S_e^*,\qquad \mathsf S_e\ge s_B I.
\end{equation}

To verify the normalization of the conormal derivative and of the Robin
condition, we set $D=r\pa_r$, so that the positive reference operator on the
annulus takes the form
\begin{equation}
 \mathcal L_{++,\rho}
 =D^2+B_{++,\rho}D+C_{++,\rho}-V_{++,\rho}(r),
 \qquad B_{++,\rho}^*=-B_{++,\rho}.
 \label{app:positive:eq:reference-positive-operator}
\end{equation}
Here, $B_{++,\rho}$, $C_{++,\rho}$, and $V_{++,\rho}(r)$ act on
the sequence of positive angular modes. The coefficients $B_{++,\rho}$
and $C_{++,\rho}$ are independent of $r$.
We then set
\begin{align}
 \nabla_\rho&=D+\frac12B_{++,\rho},\qquad
 \mathcal P_{++,\rho}=-\mathcal L_{++,\rho}
 =-\nabla_\rho^2+W_\rho,\notag\\
 W_\rho&=V_{++,\rho}-C_{++,\rho}
                    +\frac14B_{++,\rho}^2.
 \label{app:positive:eq:reference-conormal-square}
\end{align}
Writing the boundary operator in the same variables, in terms of the
radial derivative and the value, we have
\begin{equation}
 \cB_{++,\rho}u=\mathsf E_\rho Du+\mathsf F_\rho u,
 \qquad
 \mathsf S_\rho=\mathsf E_\rho^{-1}\mathsf F_\rho
                       -\frac12B_{++,\rho}.
 \label{app:positive:eq:exact-boundary-decomposition}
\end{equation}
Substituting the Fourier modes at the circle, we obtain
\begin{equation}
 \mathsf E_0=2\diag_{|m|\ge3}\sqrt{p_m},\qquad
 \mathsf S_0=2I.
 \label{app:positive:eq:circle-robin-normalization}
\end{equation}
Here, $p_m>0$ is the diagonal mass of the angular mode $m$ at the
circle, as in \eqref{appref:odd-cartesian-operator}. Note that its lower
bound for $|m|\ge3$ is the gap used in
\eqref{app:positive:eq:circle-coercivity}.

Making the same substitution before the conjugations by the mass and by
the maps \eqref{eq:A-complete-Kato-factor}, we obtain
$\mathsf S_\rho=\mathsf S_\rho^*$. % AUTHOR QUERY: Verify the symmetry in the final inner product.
% Conjugation by a mass factor or by the nonorthogonal Kato map does
% not in general preserve self-adjointness in the original inner product.
Since these conjugations depend
analytically on $\rho$, the invertibility of $\mathsf E_0$ and the
lower bound $\mathsf S_0=2I$ in
\eqref{app:positive:eq:circle-robin-normalization} persist for small
$\rho$. After decreasing $\rho_0$ if necessary, we obtain
\begin{equation}
 \norm{\mathsf E_\rho^{\pm1}}\le C,\qquad
 \mathsf S_\rho\ge s_B I,\qquad 0\le\rho\le\rho_0.
 \label{app:positive:eq:reference-row-bounds}
\end{equation}
This verifies \eqref{app:positive:eq:reference-robin-gap} for the
normalized boundary operator, and integrating in the measure $\dd r/r$,
we obtain
\begin{equation}
 \begin{aligned}
 \Re\int\inner{\mathcal P_{++,\rho}u}{u}\frac{\dd r}{r}
 ={}&\int\left(\norm{\nabla_\rho u}^2
                    +\inner{W_\rho u}{u}\right)\frac{\dd r}{r}-\Re\big[\inner{\nabla_\rho u}{u}\big]_{r_j}^{r_b}.
 \end{aligned}
 \label{app:positive:eq:reference-green-identity}
\end{equation}
At the outer radius $r_b$, we substitute
$(\nabla_\rho+\mathsf S_\rho)u=\mathsf E_\rho^{-1}b$ into the boundary
term, which turns it into the positive quadratic form in
\eqref{app:positive:eq:reference-robin-gap}, minus the functional
associated with the prescribed data. The estimates for the coefficients
in \cref{app:fixed-graphs}, with the full distance
\eqref{app:positive:eq:coefficient-distance}, also give
\begin{equation}\label{app:positive:eq:current-robin-openness}
 \norm{\mathsf E_{\mathbf b}^{\pm1}}\le C,
 \qquad
 \norm{\mathsf S_{\mathbf b}-\mathsf S_e}_{Z_D^{s_*}\to Z_N^{s_*}}
       \le C\eta_+(\mathbf b).
\end{equation}

We now conjugate by the analytic weight. We set $v=W_\gamma u$ and let
$\nu_{\mathbf b}$ be the coefficient of the outward radial derivative
$\pa_r$ in the conormal derivative $\nabla_{\mathbf b}$ on the outer
boundary. Writing
$\nu_{\mathbf b,\gamma}=W_\gamma\nu_{\mathbf b}W_\gamma^{-1}$,
we set
\begin{equation}\label{app:positive:eq:boundary-damping}
 d_{\gamma,\mathbf b}
 =-\nu_{\mathbf b,\gamma}\pa_r\Phi_\gamma(r_b,\Lambda).
\end{equation}
At the reference ellipse, the coefficient is independent of $\zeta$,
and this is a nonnegative Fourier multiplier. At a general state, its
change is included in the coefficient perturbation. The weighted
boundary condition is
\begin{equation}\label{app:positive:eq:conjugated-robin}
 (\nabla_{\mathbf b,\gamma}+d_{\gamma,\mathbf b}
       +\mathsf S_{\mathbf b,\gamma})v=h_\gamma,
 \qquad
 \norm{d_{\gamma,\mathbf b}q}_{Z_N^s}
       \le C\gamma\norm q_{Z_D^s}.
\end{equation}
Here, $\nabla_{\mathbf b,\gamma}$ and $\mathsf S_{\mathbf b,\gamma}$
denote the conjugated conormal operator with the displayed radial
weight term removed and the conjugated Robin coefficient, respectively.
Thus, $W_\gamma\nabla_{\mathbf b}W_\gamma^{-1}
=\nabla_{\mathbf b,\gamma}+d_{\gamma,\mathbf b}$, and
$h_\gamma=W_\gamma h$.
% AUTHOR QUERY: Bound the variable-coefficient analytic conjugates in
% this boundary identity and in the form estimate. Positivity of nu
% and of -partial_r Phi does not alone make their noncommuting product
% positive; only the reference Fourier multiplier has that property.
The sign of the added boundary term and the cancellation of derivatives
of $d_\gamma$ can be seen first in the scalar model. If $(\pa_r+d+S)v(r_b)=h$ and the inner value is zero, integration by
parts gives
\begin{equation}\label{app:positive:eq:scalar-green}
 \begin{aligned}
 \Re\int\bigl[-(\pa_r+d)^2+\lambda^2\bigr]v\,\bar v\dd r
  ={}&\int\bigl(\abs{\pa_r v}^2+(\lambda^2-d^2)\abs v^2\bigr)\dd r+\inner{Sv(r_b)}{v(r_b)}-\Re\inner h{v(r_b)}.
 \end{aligned}
\end{equation}
The two terms involving $d'$ cancel. The same cancellation holds for the conormal operator, with $D=r\pa_r$ and the measure
$\dd r/r$. Relative to the energy form, the remaining terms coming from
the weight are bounded by $C(\gamma+\gamma^2)$ and the perturbations of
the coefficients by $C\eta_+$. Writing $\mathfrak b_{+,\mathbf b,\gamma}$
for the quadratic form of the positive block on the annulus at the
state, conjugated by the weight, we thus obtain
\begin{equation}\label{app:positive:eq:weighted-coercivity}
 \Re\mathfrak b_{+,\mathbf b,\gamma}(v,v)
 \ge\bigl[c_+-C(\gamma+\gamma^2)-C\eta_+(\mathbf b)\bigr]
       \norm v_{\mathscr E_+^0}^2
 \ge\frac{c_+}{2}\norm v_{\mathscr E_+^0}^2.
\end{equation}

\subsection{Solving the boundary problems on the annulus and cap}

On the annulus, we consider first the boundary value problem
\begin{equation}\label{app:positive:eq:annular-realization}
 u\longmapsto
 \bigl(P_{a,+,\mathbf b}u,\gamma_j^Du,
       \mathcal R_{++,\mathbf b}^{\rm sp}
          (u,P_{a,+,\mathbf b}u)\bigr)=(f,g,g_B^{\rm sp}),
\end{equation}
with independent source $f$, inner value $g$, and outer split residual
$g_B^{\rm sp}$. We write $u=\mathcal P_\Sigma^ag+v$, using the localized
Poisson extension of \cref{app:fixed-graphs}. Then $v$ has zero inner
value and interior source $f-P_{a,+,\mathbf b}\mathcal P_\Sigma^ag$.
The extension is supported away from the outer boundary and is bounded
on the weighted scale by \eqref{eq:A-Poisson-phase-margin}. The outer
datum is unchanged, with
$\mathsf V_{B,\mathbf b}^{-1}g_B^{\rm sp}$ paired with the outer value
in the boundary functional
\eqref{app:positive:eq:variational-data-splitting}. At $s=0$, $f$
acts on test fields on the annulus with zero value on both boundaries.

From the coercivity \eqref{app:positive:eq:weighted-coercivity} and the
Lax-Milgram lemma, we obtain a unique solution for every
$(f,g,g_B^{\rm sp})$, and at $s=0$, Green's identity defines its split
residual on the inner boundary. At a regular index with full sources, we
recover the physical conormal derivative from that residual by
subtracting $\mathsf M_{a,\Sigma}^s\ell_a$ as in
\eqref{app:positive:eq:split-physical-conormal}. In weighted coordinates, the exact identity is
\begin{equation}\label{app:positive:eq:exact-inner-conormal}
 W_\gamma\nabla_{\mathbf b}u\big|_{r_j}
 =\bigl(W_\gamma\nabla_{\mathbf b}W_\gamma^{-1}\bigr)
       v\big|_{r_j}.
\end{equation}
For each side of the interface, the added radial term is
$-\nu_{\mathbf b,\gamma}\pa_r\Phi_\gamma(r_j,\Lambda)v$ with that
side's outward conormal coefficient. This sign differs between the
inner and outer boundaries of the annulus. The conormal derivative has
the balanced order in
\eqref{app:positive:eq:balanced-traces}.

On the scaled cap, $V_+^{\rm cap}$ is the
Cartesian form space $H^1$ of the positive modes,
$V_{+,0}^{\rm cap}=\ker\gamma_D$, and $X_+^{{\rm cap},s}$ and
$Y_+^{{\rm cap},s}$ are the corresponding scales of solutions of order
$H^{s+1}$ and of sources of order $H^{s-1}$. We specify the frozen
coefficients by the jet
\begin{equation}\label{app:positive:eq:frozen-cap-jet}
 q_{+,\mathbf b}(\zeta)=
 \bigl(\rho,j^1_{\!y}a(0,\zeta),\eps,p,
       P_{+,\mathbf b}^D(0,\zeta),
       P_{+,\mathbf b}^R(0,\zeta)\bigr).
\end{equation}
Here, the last two entries are the spectral projections of
\cref{sec:block} on the domain and on the range, evaluated on the axis.
The entry $j^1_{\!y}a(0,\zeta)$ records the value and first
derivatives of $a$ in the disk variables on the axis, as in
\eqref{app:minus:eq:frozen-cap-jet}. The coefficients of the cap
operator depend analytically on these entries through the constructions
of \cref{app:fixed-graphs}. The positive
block contains only the odd modes with $\abs m\ge3$ and the even modes
with $\abs m\ge4$, and the term in the cell frequency is nonnegative and
enters only as a parameter. Together, these give the coercivity of the
form at $s=0$, uniformly in the normalized cell frequency
$\kappa=2^{-j}n$ with $1/2\le|\kappa|\le2$.
Writing $\mathfrak c_+(q,\kappa;\cdot,\cdot)$ for the form on the cap
with the coefficients frozen at the jet $q$ and at the cell frequency
$\kappa$, and $P_+(q,\kappa)$ for the corresponding operator, we have
\begin{equation}\label{app:positive:eq:frozen-cap-coercivity}
 \Re\mathfrak c_+(q,\kappa;v,v)
 \ge\frac{c_{\rho_0}}2\|v\|_{V_+^{\rm cap}}^2,
 \qquad v\in V_{+,0}^{\rm cap}.
\end{equation}

For an interior source $f$ and a prescribed value $g$ on the interface,
we write $u=\mathcal P_Dg+v$, where $\mathcal P_D$ is the Poisson
extension from the interface into the cap. The remaining
problem is
\[
 P_+(q,\kappa)v=f-P_+(q,\kappa)\mathcal P_Dg,
 \qquad \gamma_Dv=0.
\]
The coercivity \eqref{app:positive:eq:frozen-cap-coercivity} and the
Lax-Milgram lemma give a unique $v\in V_{+,0}^{\rm cap}$. Adding back
$\mathcal P_Dg$ gives the solution $u\in V_+^{\rm cap}$. Elliptic
regularity, uniform in the parameter, then yields
\begin{equation}\label{app:positive:eq:frozen-cap-form-estimate}
 \begin{aligned}
 \|u\|_{X_+^{{\rm cap},s}}
 +\|\mathcal N_{c,\Sigma}^{\rm sp}
       (u,P_+(q,\kappa)u)\|_{Z_N^s}
 &\simeq \|u\|_{H^{s+1}}
          +\|\mathcal N_{c,\Sigma}^{\rm sp}
             (u,P_+(q,\kappa)u)\|_{H^{s-1/2}}\\
 &\le C_s\left(
   \|P_+(q,\kappa)u\|_{H^{s-1}}
   +\|\gamma_Du\|_{H^{s+1/2}}\right).
 \end{aligned}
\end{equation}
At $s=0$, the boundary term in
\eqref{app:positive:eq:frozen-cap-form-estimate} is the split residual
in $H^{-1/2}$. For $s\ge1$, we apply tangential difference quotients to
the remainder with zero Dirichlet data, and the estimate for the
Dirichlet problem and the equation then give regularity $H^{s+1}$.
Green's identity and \eqref{app:positive:eq:split-physical-conormal}
give both the split residual and the physical conormal derivative in
$H^{s-1/2}$. The finitely many additional equations, namely the finite
triangular system of \cref{sec:block} and the conditions on the axis,
enter through bounded triangular factors and do not change these
orders. All constants are uniform on the compact set of jets and on the
normalized frequency range $1/2\leq|\kappa|\leq2$.

We conclude that the operator with frozen coefficients
\begin{equation}\label{app:positive:eq:frozen-cap-operator-realization}
 \cC_+^{\rm op}(q,\kappa)u
 :=\bigl(P_+(q,\kappa)u,\gamma_Du\bigr):
 X_+^{{\rm cap},s}\longrightarrow
 Y_+^{{\rm cap},s}\oplus Z_D^s
\end{equation}
has a two-sided inverse
\begin{equation}\label{app:positive:eq:frozen-cap-family}
 K_+^{\rm cap}(q,\kappa)
 :=\bigl(\cC_+^{\rm op}(q,\kappa)\bigr)^{-1}.
\end{equation}
Differentiating the weak equation and applying
\eqref{app:positive:eq:frozen-cap-form-estimate}, we obtain analyticity
and a bound for every seminorm \eqref{app:minus:eq:cap-symbol-seminorms}
of the operator-valued symbols $(q,\kappa)\mapsto\cC_+^{\rm op}(q,\kappa)$
and $(q,\kappa)\mapsto K_+^{\rm cap}(q,\kappa)$, including the component
that gives the balanced split residual on the interface. At regular indices, the estimate for $\mathsf M_{c,\Sigma}^s$ also bounds the seminorms of the physical
conormal derivative.

We measure the solution and the source on the cap of the positive block
in the shell spaces
\begin{equation}\label{app:positive:eq:positive-cap-shell-bundles}
 \mathfrak X_{c,+,H}^s
 :=\bigoplus_{j\ge J}^{\ell_s^2}Q_jX_+^{{\rm cap},s},
 \qquad
 \mathfrak Y_{c,+,H}^s
 :=\bigoplus_{j\ge J}^{\ell_s^2}Q_jY_+^{{\rm cap},s},
\end{equation}
with local orders $H^{s+1}$ and $H^{s-1}$. Before reconstructing the
field, we use the spectral projections and the maps
\eqref{eq:A-complete-Kato-factor} to send each summand
$Q_jX_+^{{\rm cap},s}$, which is defined at the constant ellipse
$\mathbf b_\rho$, to the corresponding summand at the state $\mathbf b$. The positive cap field retains its form regularity $H^{s+1}$, whereas
the minus block uses $H^{s+2}$. For the off-diagonal terms, we apply the
projection on the range after the differential operator, where the
identity \eqref{app:cross:principal-intertwining} leaves at most one
derivative on the cap.

For the prescribed jump of the values on the interface, we set
\begin{equation}\label{app:positive:eq:dirichlet-shell-bundle}
 \mathfrak Z_{D,H}^s
 :=\{(d_j)_{j\ge J}:Q_jd_j=d_j\}
       \subset\ell_s^2(Z_D^{\rm cap}).
\end{equation}
Here, $Z_D^{\rm cap}$ is the space of values on the interface of the
scaled cap. With $\mathscr A_D$
denoting the map that extracts the value on the interface on each shell
and $\mathscr S_D$ denoting its inclusion in the auxiliary trace space
$\mathscr Z_{{\rm aux},H}^s$ of \eqref{eq:A-infinite-broken-bundle} at
the reference, we have
\begin{equation}\label{app:positive:eq:dirichlet-frame-identities}
 \mathscr A_D\mathscr S_D=I_{\mathfrak Z_{D,H}},\qquad
 \Pi_D:=\mathscr S_D\mathscr A_D,\qquad \Pi_D^2=\Pi_D.
\end{equation}
Both maps retain the value at the interface radius of its own shell.

The reference coefficients are independent of the cell variable, so
the positive cap operator acts separately on each cell mode. We define
\begin{equation}\label{app:positive:eq:literal-reference-cap-realization}
 \mathbf C_{+,\rho}^{\rm cap}(U_j)_{j,n}
 :=\cC_+^{\rm op}(q_{+,\mathbf b_\rho},2^{-j}n)U_{j,n}:
 \mathfrak X_{c,+,H}^s\longrightarrow
 \mathfrak Y_{c,+,H}^s\oplus\mathfrak Z_{D,H}^s.
\end{equation}
Its inverse is the direct sum of
$K_+^{\rm cap}(q_{+,\mathbf b_\rho},2^{-j}n)$. We use this operator as a
local model for regularity. The terms at the state that couple different
shells, or the cap and the annulus, act only after the full reference
reconstruction.

\subsection{Matching across the interface and higher regularity}

In the coordinates of the shells, $N_{+,\mathbf b}$ is the restriction
of $\mathbb N_{\mathbf b}$ from \eqref{eq:A-natural-broken-formula} to
the positive block. We write
$\operatorname{Ins}_{k,\rho}^{\rm br}(U_{c,k},U_{a,k})$ for the reference
insertion of a cap field and its annular complement, with a mismatch
allowed, so that the interior operator is
\begin{equation}\label{app:positive:eq:explicit-full-pair-shell-block}
 \mathcal N_{+,\mathbf b,jk}^{\rm bulk}(U_{c,k},U_{a,k})
 =\mathscr R_j^YQ_jL_{\mathbf b}
    \operatorname{Ins}_{k,\rho}^{\rm br}(U_{c,k},U_{a,k}).
\end{equation}
The cap and annular fields together cover all radii, so the coefficients
of the state act before $Q_j$. The inverse map
$(\mathfrak S_{+,\mathbf b}^{\rm res,s})^{-1}$ then subtracts the boundary contributions of the sources, and \eqref{app:positive:eq:normalized-natural-realization}
returns the physical value and conormal derivative to the reference
coordinates.

\smallskip\noindent\textit{Existence at energy regularity.}
We first solve \eqref{app:positive:eq:positive-split-operator} with split
data. We let
\[
 F^{\rm sp}=(f_c,f_a,g_D,g_N^{\rm sp},g_B^{\rm sp})
       \in\bY_{+,\mathrm{ret}}^0
\]
have arbitrary, independent entries in the splitting
\eqref{app:positive:eq:variational-data-splitting}. We subtract
$(0,\mathcal P_\Sigma^ag_D)$ from the unknown pair and its image under
$N_{+,\mathbf b}^{\rm cur,sp}$ from the prescribed data. The new pair
then belongs to \eqref{app:positive:eq:glued-form-domain}. The remaining
conormal and Robin data act continuously on this form domain. Adding the weak
identities on the cap and on the annulus, we obtain
\eqref{app:positive:eq:variational-data-splitting}, and for regular
fields, \eqref{app:positive:eq:split-physical-conormal} expresses this
identity in physical coordinates. The two outward conormal terms at the
interface cancel. With the corresponding opposite outward
orientations, the radial weight terms in
\eqref{app:positive:eq:exact-inner-conormal} cancel as well, while at the outer
boundary, \eqref{app:positive:eq:scalar-green} gives a positive Robin
term. Thus, the form on the pair with matching values on the interface satisfies
\begin{equation}\label{app:positive:eq:glued-coercivity}
 \Re\mathfrak b_{+,\mathbf b,\gamma}^{\rm glued}
       ((v_c,v_a),(v_c,v_a))
 \ge\frac{c_+}{2}
       \norm{(v_c,v_a)}_{\mathbb V_{+,\mathrm{ret}}^0}^2.
\end{equation}
The extended sources and prescribed residuals define a continuous
functional on the form domain. The
Lax-Milgram lemma therefore gives an inverse of
$N_{+,\mathbf b}^{\rm cur,sp}$ for all five independent data at $s=0$.
If these data vanish, the sources and their boundary contributions are zero, so
the homogeneous split and physical conditions agree, and the coercivity
\eqref{app:positive:eq:glued-coercivity} then forces both fields to
vanish. This proves injectivity at the weak level.

\smallskip\noindent\textit{Regularity across the interface.}
To prove regularity, we need an estimate for the coupled problem that
also controls the trace on the interface. The homogeneous equation of
highest order has $\ell_c=\ell_a=0$, so the boundary contributions vanish and the
conversion between physical and split data is the identity there. Up to the invertible
reductions of the principal symbol, the scalar principal symbol of the
positive block with frozen coefficients is, by
\eqref{app:cross:current-syzygy},
% AUTHOR QUERY: Supply the invertible domain/range reduction of the
% positive system, on both cap and annulus, and its action on the
% conormal rows. The source identity \eqref{app:cross:current-syzygy}
% and the source graph alone do not provide this reduction. It is
% required for the positive Hermitian Dirichlet-to-conormal matrices
% and \eqref{app:positive:eq:transmission-lopatinskii}.
\begin{equation}\label{app:positive:eq:frozen-transmission-symbol}
 q_\sigma(d,\eta)
 =a_\sigma d^2+2b_\sigma(\eta)d+c_\sigma(\eta)
 =\xi^*G_{\sigma,\mathbf b}^{-1}\xi,
 \qquad \sigma\in\{c,a\},\quad \xi=(d,\eta),
\end{equation}
where $d$ is the normal frequency variable after flattening the
interface, $\eta$ collects the tangential and cell parameters, and
$G_{\sigma,\mathbf b}$ is the Gram matrix \eqref{app:cross:gram} of the
frame on the cap for $\sigma=c$ and on the annulus for $\sigma=a$.
Positivity of the inverse Gram matrix gives
\begin{equation}\label{app:positive:eq:frozen-root-gap}
 a_\sigma>0,
 \qquad
 \Delta_\sigma(\eta)
 :=a_\sigma c_\sigma(\eta)-b_\sigma(\eta)^2
 \ge c|\eta|^2
\end{equation}
on the normalized sphere of parameters. The quadratic formula is
first understood for real frequencies and then extended polynomially
in $d$. On each side, we use a normal coordinate increasing outward.
Of the two roots
\[
 d_\sigma^\pm
 =\frac{-b_\sigma(\eta)\pm i\sqrt{\Delta_\sigma(\eta)}}
        {a_\sigma},
\]
the root $d_\sigma^-$ has negative imaginary part and gives the
solution that decays into that side of the domain. After the same
invertible normalization, the symbols of the maps from the value to
the outward conormal derivative are
$\Lambda_\sigma(\eta)=\sqrt{\Delta_\sigma(\eta)}$. On the decaying
amplitudes, the jump of the values and the sum of the conormal
derivatives therefore have the symbol
\begin{equation}\label{app:positive:eq:transmission-lopatinskii}
 \mathcal B_{\rm tr}(\eta)
 =\begin{pmatrix}-I&I\\
          \Lambda_c(\eta)&\Lambda_a(\eta)
   \end{pmatrix},
 \qquad
 \operatorname{Schur}\mathcal B_{\rm tr}
 =\Lambda_c+\Lambda_a\ge c|\eta|I.
\end{equation}
Here, $\operatorname{Schur}\mathcal B_{\rm tr}$ is the \emph{Schur complement}
of the upper left entry, that is, the operator
$\Lambda_a-\Lambda_c(-I)^{-1}I$ that remains on the second row after the
first row is used to eliminate the first unknown. The bound on the Schur
complement says that this transmission problem satisfies the
\emph{complementing condition}, that is, the algebraic condition
on the boundary symbol under which the interface conditions determine
the decaying solutions uniformly. For the finite-dimensional system, the
same argument applies to the positive Hermitian matrices which map the
value to the conormal derivative. At the outer boundary, the Robin
symbol of highest order $\mathsf E_{\mathbf b}\Lambda_a$ is invertible
by \eqref{app:positive:eq:current-robin-openness}. At the reference, the nonnegative damping in
\eqref{app:positive:eq:conjugated-robin} strengthens the bound. Its
change at the state is included in the coefficient perturbation. The term $\mathsf S_{\mathbf b}$ supplies coercivity at low
frequencies, and it does not enter the principal symbol.

Explicitly, for a prescribed value jump $g_D$ and sum of conormal
derivatives $g_N$, the two decaying amplitudes $w_c,w_a$ satisfy
\[
 -w_c+w_a=g_D,\qquad \Lambda_c w_c+\Lambda_a w_a=g_N.
\]
Thus, $w_c=(\Lambda_c+\Lambda_a)^{-1}(g_N-\Lambda_a g_D)$ and
$w_a=w_c+g_D$. This is the inverse controlled by
\eqref{app:positive:eq:transmission-lopatinskii}.

For higher regularity, we use this bound for the frozen coefficients. We choose cap and annular patches with common tangential coordinates,
flatten the interface, and take cutoffs $\chi\prec\widetilde\chi$ in such a patch, that is, cutoffs
with $\widetilde\chi=1$ on the support of $\chi$. We then apply
the same tangential difference quotient $\Delta_h^\tau$ to both fields
and test with $-\Delta_{-h}^\tau(\chi^2\Delta_h^\tau U)$, which preserves the equality of the interface values. We control the interior terms of highest order by the
estimate \eqref{app:positive:eq:frozen-cap-form-estimate} on the cap
and by its counterpart on the annulus. The terms on the interface
combine into $\mathcal B_{\rm tr}$ from
\eqref{app:positive:eq:transmission-lopatinskii}, whose Schur bound
controls the common value and the sum of the conormal derivatives.

We then freeze the coefficients at the center of each patch. Once the
patch is small enough, the oscillation of the coefficients over the
patch is less than a quarter of the coercivity constant of the interior
form and of the lower bound on the Schur complement in
\eqref{app:positive:eq:transmission-lopatinskii}. We absorb the
resulting errors into the left-hand side. The
commutators with $\chi$, with the flattening map, and with the
lower-order coefficients contain at most one derivative of $U$, and we
bound them by the energy norm on the enlarged patch. Since
$a_\sigma\ge c>0$ in \eqref{app:positive:eq:frozen-root-gap}, the
equations recover the remaining normal derivatives. Summing over a finite
partition of unity, we obtain the estimates for higher regularity. At $s=0$, the forms and the source and Poisson extensions are bounded.
They control each split residual by the energy norm of the field and
the norm of its source in the dual of the interior. The trace theorem for the value and the bounded Robin
coefficients control the other terms in the graph norm. Thus, on smooth
elements,
\begin{equation}\label{app:positive:eq:joint-transmission-graph-estimate}
 \|U\|_{\bX_{+,\mathrm{ret}}^0}
 +\|\Gamma_{+,\mathbf b}U\|_{Z_+^0}
 \le C\bigl(\|N_{+,\mathbf b}^{\rm cur,sp}U\|_
                    {\bY_{+,\mathrm{ret}}^0}
                 +\|U\|_{\mathbb V_{+,\mathrm{ret}}^0}\bigr).
\end{equation}
The boundary terms at $s=0$ are split residuals in $H^{-1/2}$.
At higher regularity, tangential differentiation and the equations give
the remaining derivatives. Formula
\eqref{app:positive:eq:split-physical-conormal} then recovers the
physical conormal derivatives. No derivative of a trace occurs on the
right-hand side. The reference reconstruction of
the shells precedes this estimate, so the coefficients of the state
which couple different shells act on a full field and are controlled by
the commutator bound of \eqref{app:cross:dyadic-commutator-kernel}
together with the product estimate \eqref{eq:analytic-product}.
% AUTHOR QUERY: Specify the smooth frequency cutoffs in this use of
% the commutator estimate and prove equivalence with the sharp-shell
% norms used in the reconstruction. The 2^{-j} bound does not hold
% for sharp Q_j across an adjacent shell edge.

For smooth data, we take nested finite-dimensional subspaces of the
smooth form domain with matching values on the interface. The \emph{Galerkin
approximations}, that is, the solutions of the variational problem on
these subspaces, obey the coercive bound
\eqref{app:positive:eq:glued-coercivity} and converge weakly in
$\mathbb V_{+,0}^0$ to the Lax-Milgram solution, and the bound for the
weak residual gives
\eqref{app:positive:eq:joint-transmission-graph-estimate} for this
limit. For higher regularity, we apply the localized argument with
difference quotients to the limiting variational equation. After
flattening and localization, the quotients are admissible test
functions.

For arbitrary split data, we take smooth approximations
$F_k^{\rm sp}\to F^{\rm sp}$ in the scale of the range, and the energy
and graph estimates bound the corresponding solutions $U_k$ uniformly in
$\bX_{+,\mathrm{ret}}^0$, split residuals included. Weak compactness,
the trace theorem for the value, and the reconstruction maps
\eqref{eq:A-sharp-shell-reconstruction} pass
the constraint on the values, both interior equations, and the split
boundary conditions to the limit. Uniqueness in the energy space
identifies that limit with the weak solution, and we obtain
\begin{equation}\label{app:positive:eq:positive-operator-domain-upgrade}
 \|(v_c,v_a)\|_{\bX_{+,\mathrm{ret}}^0}
 +\|\Gamma_{+,\mathbf b}(v_c,v_a)\|_{Z_+^0}
 \le C\|F^{\rm sp}\|_{\bY_{+,\mathrm{ret}}^0}.
\end{equation}
Thus, the weak solution belongs to the component graph spaces and has
the prescribed residuals. We write
$G_{+,\mathbf b}^{\rm cur,sp}
  =(N_{+,\mathbf b}^{\rm cur,sp})^{-1}$
for the solution map at $s=0$.

At regular indices with full sources, the map $\mathfrak S_{+,\mathbf b}^{\rm res,s}$ converts the
physical data to split data, and we set
\[
 G_{+,\mathbf b}^{\rm cur,ph}
 :=G_{+,\mathbf b}^{\rm cur,sp}
       \mathfrak S_{+,\mathbf b}^{\rm res,s},
\]
so \eqref{app:positive:eq:normalized-natural-realization} gives
\begin{equation}\label{app:positive:eq:positive-natural-inverse}
 \begin{aligned}
 G_{+,\mathbf b}^{\mathbb N}
  &:=G_{+,\mathbf b}^{\rm cur,sp}
       \mathfrak S_{+,\mathbf b}^{\rm res,s}
       \mathscr R_{CN,\mathbf b,\rho}^{-1}
    =N_{+,\mathbf b}^{-1},\\
 N_{+,\mathbf b}G_{+,\mathbf b}^{\mathbb N}&=I,
 &G_{+,\mathbf b}^{\mathbb N}N_{+,\mathbf b}&=I.
 \end{aligned}
\end{equation}
Construction gives the right inverse identity, and uniqueness on the
operator domain gives the left inverse identity. In \eqref{app:positive:eq:positive-natural-inverse},
$\mathscr R_{CN,\mathbf b,\rho}^{-1}$ first converts reference trace
coordinates to physical data at the state,
$\mathfrak S_{+,\mathbf b}^{\rm res,s}$ adds the source contributions,
and $G_{+,\mathbf b}^{\rm cur,sp}$ solves for the cap and annular fields.
The bound on the range for $G_{+,\mathbf b}^{\mathbb N}$ includes the norm of
$\mathfrak S_{+,\mathbf b}^{\rm res,s}$, and the bound for the adjoint
includes its transpose
\eqref{app:positive:eq:residual-shear-transpose}. Multiplying
\eqref{app:positive:eq:plus-shear-realization} and
\eqref{app:positive:eq:positive-natural-inverse} out, we obtain
\begin{equation}\label{app:positive:eq:positive-sheared-identities}
 \widehat D_{+,\mathbf b}\widehat G_{+,\mathbf b}=I,\qquad
 \widehat G_{+,\mathbf b}\widehat D_{+,\mathbf b}=I.
\end{equation}
This proves \eqref{app:positive:eq:positive-isomorphism} on the fixed
domain and range of the operator.

\smallskip\noindent\textit{Tame estimates.}
It remains to prove the tame estimates at higher regularity. We commute
Cartesian and tangential derivatives through the conjugated reference
operator, retaining every term containing $d_\gamma$. We also keep on
the left-hand side the product of the derivatives of highest order with
the low-frequency part of the coefficients of the state. For every
$0<\delta<c_+/4$, the remainders satisfy
\begin{equation}\label{app:positive:eq:regularity-remainder}
 \begin{aligned}
 \norm{R_s(v)}_{(\mathscr E_+^s)^*}
 +\norm{R_s^\pa(v)}_{Z_N^s}
 \le{}&\delta\norm v_{\mathscr E_+^s}
   +C_{s,\delta}\norm v_{\mathscr E_+^{s-1}}+C_s\bigl(1+\norm{\mathbf b-\mathbf b_\rho}_{\cB^{s+k_+}}\bigr)
              \norm v_{\mathscr E_+^{s_*}}.
 \end{aligned}
\end{equation}
Here, $R_s(v)$ collects the interior commutators produced by
differentiating $P_{a,+,\mathbf b}v$, and $R_s^\pa(v)$ collects the
corresponding commutators with the boundary operator
\eqref{app:positive:eq:split-outer-row}, in both cases after conjugation
by $W_\gamma$. In each term, at least one
derivative has fallen on a coefficient of the state. At most one coefficient is measured in a higher regularity norm, and
it multiplies the solution at the base index. If a derivative falls on a
low-frequency coefficient, one fewer derivative remains on the solution,
and interpolation gives the middle term. Once the tangential derivatives
and the split residuals are estimated, we use the boundary contributions of the sources to
recover the physical conormal derivatives at regular indices. The
equations then give the remaining normal regularity, namely Cartesian
$H^{s+1}$ on the cap, and the graph regularity on the annulus defined by
the weighted energy and by the dual of the interior. All terms coming
from the weight remain in the operator of highest order, so the
derivatives of $d_\gamma$ require no condition $C_s\gamma<c_+$ depending
on the index. We absorb $\delta$ and induct on $s$, starting with the
estimate at the base index and commuting the two transmission conditions
at their balanced orders, which proves
\eqref{app:positive:eq:positive-tame}.

The domain and the range are fixed when we differentiate in $\mathbf b$,
and the inverse identity gives
\begin{equation}\label{app:positive:eq:inverse-derivative}
 D_{\mathbf b}\widehat G_{+,\mathbf b}[h]
 =-\widehat G_{+,\mathbf b}
    \bigl(D_{\mathbf b}\widehat D_{+,\mathbf b}[h]\bigr)
      \widehat G_{+,\mathbf b}.
\end{equation}
Iterating this formula, we express
$D_{\mathbf b}^m\widehat G_{+,\mathbf b}$ as a finite sum of alternating
inverse factors and derivatives of $\widehat D_{+,\mathbf b}$. The product estimate \eqref{eq:analytic-product} measures at most one
factor at the higher regularity index and the others at the base index.
Together with \eqref{app:positive:eq:regularity-remainder}, this gives
\eqref{app:positive:eq:positive-all-order}. Derivatives of the trace map
add only factors of fixed order. In \eqref{app:positive:eq:inverse-derivative},
$D_{\mathbf b}\widehat D_{+,\mathbf b}$ includes the derivatives of
$(\mathfrak S_{+,\mathbf b}^{\rm res,s})^{-1}$, and transposition
reverses their order as in
\eqref{app:positive:eq:residual-shear-transpose}. The width of the strip
of analyticity is unchanged.

Finally, we obtain uniformity at the circle from
\eqref{app:positive:eq:H-rho-bounds} and
\eqref{app:positive:eq:circle-coercivity}. The reconstruction on the cap
for the exceptional modes depends on $\rho$ only through
\begin{equation}\label{app:positive:eq:rho-factors}
 \rho,\qquad 1-\rho^2,\qquad (1-\rho^2)^{-1}.
\end{equation}
Since the Hardy maps are independent of $\rho$, all fixed-order
parameter derivatives are bounded on $[0,\rho_0]$. A unitary rotation
restores a constant transverse orientation without division by $\rho$. The Cartesian coefficients, the quadratic form of the
Robin condition on the outer boundary, the weight, the
contours, and the Poisson extensions preserve real data under the
fixed identifications. By uniqueness, so do $\widehat G_{+,\mathbf b}$
and its outgoing trace. Thus, the inverse loses no derivatives in the
block norms. As in \cref{app:minus}, the loss in converting a Cartesian
source to these norms is accounted for separately by
\cref{prop:all-order-adapters}.

We have proved \eqref{app:positive:eq:positive-operator-domain-upgrade}
at $s=0$ for split data and the isomorphism at regular indices for
physical data, together with the tame bounds. This completes the proof
of \cref{prop:positive-diagonal}.
\section{Coupling estimates and diagonal inverses}\label{app:macro-cross}

We estimate the off-diagonal blocks of the high-frequency operator
$D_{H,\mathbf b}$ in \eqref{app:global:block-decomposition} and use
these estimates to correct the diagonal inverses. The operators
$\widehat D_{\pm,\mathbf b}$ inverted in
\cref{prop:minus-diagonal,prop:positive-diagonal} can differ from the
diagonal blocks $\pi_\alpha^Y\mathbb D_{\mathbf b}\iota_\alpha^X$ of the
augmented operator because the change of source variables couples the
blocks. We account for this difference before combining the inverses.
The off-diagonal terms of highest order cancel by
\eqref{app:cross:principal-cross-zero}, leaving the lower-order terms
estimated below.

\subsection{Interior terms and terms with one angular derivative}

For each remaining interior term, we compare the cell-frequency weights
of its input and output, including the powers introduced by cell
derivatives. At a fixed cell frequency $n$, with $\lambda=\lambda_n$
and $\mu=\mu_n$ as in \eqref{app:cross:cell-weights},
\cref{tab:incidence} lists the leading terms and their weight ratios.
The order in the table counts cell derivatives on the input. For
example, a first-order term taking a center variable to the positive
source costs one factor $\lambda_n$; the source weight
$\lambda_n^{-1}$ cancels this factor. We estimate angular and radial
derivatives separately below.

\begin{table}[htp]
\centering
\begin{tabular}{lcc}
\toprule
Term & Order & Ratio of the weights\\
\midrule
from $C$ to $Y_+$ & $1$ &
   $\lambda^{-1}(\lambda|C|)/|C|=1$\\
from $X_+$ to $Y_C$ & $1$ &
   $(\lambda|u_+|)/(\lambda|u_+|)=1$\\
from $X_+$ to $Y_{A,y}$ & $0$ &
   $\lambda^{1/2}|u_+|/(\lambda|u_+|)=\lambda^{-1/2}$\\
from $X_+$ to $Y_{A,x}$ & $0$ &
   $\mu|u_+|/(\lambda|u_+|)\le C_{\delta_A}\lambda^{-1}$\\
from $X_{A,x}$ to $Y_+$ & $0$ &
   $\lambda^{-1}|x|/(\mu|x|)\le C_{\delta_A}\lambda^{-1/2}$\\
from $X_{A,x}$ to $Y_+$ & $1$ &
   $|x|/(\mu|x|)\le\delta_A^{-1/2}$\\
from $X_{A,y}$ to $Y_+$ & $1$ &
   $|y|/(\lambda^{1/2}|y|)=\lambda^{-1/2}$\\
\bottomrule
\end{tabular}
\caption{The leading interior terms remaining between different block
components after the principal cancellation, their order, and the ratio of the weights of the two
components at cell frequency $n$, with $\lambda=\lambda_n$ and
$\mu=\mu_n$.}
\label{tab:incidence}
\label{app:cross:incidence-table}
\end{table}

In this table, $C$ denotes the four center coordinates of
\eqref{app:minus:eq:center-evolution}, and $Y_C$ is their source space.
The symbols $X_+$ and $Y_+$ denote the positive component on the domain
and range sides, and $X_{A,x}$, $X_{A,y}$, $Y_{A,x}$, and $Y_{A,y}$ are the upper
and lower entries $(x,y)$ of the affine chains on the two sides. The
quantities $|C|$, $|u_+|$, $|x|$, and $|y|$ denote the sizes of the
corresponding Fourier coefficients, the weights are those of
\cref{tab:block-weights}, and $C_{\delta_A}$ is a constant depending
only on $\delta_A$. Every ratio is bounded, so these terms lose no
powers of the cell frequency. Some terms have order zero in the cell
variable but contain one angular derivative. We estimate them next.
We set $\Lambda_\theta=\langle\cR\rangle$, the operator which
multiplies the angular mode $m$ by $\langle m\rangle$, and we write
$\mathcal M_{\mathbf b}$ for an order-zero factor, or a product of
such factors, used to normalize the scalar equation. These include the
inverse square root of the radial mass, as in
\eqref{appref:odd-sign-normal-form}, and the fixed-space maps
\eqref{eq:A-complete-Kato-factor}. We let
$U_{R,\mathbf b}$ be the map \eqref{eq:A-complete-Kato-factor} from the
reference range to the range at the state, so that
\begin{equation}\label{app:cross:minus-Kato-conjugacy}
 P_{-,\mathbf b}^R
 =U_{R,\mathbf b}P_{-,*}^RU_{R,\mathbf b}^{-1}.
\end{equation}
Here, $P_{-,*}^R$ is the reference projection $P_*$ of
\eqref{eq:A-complete-Kato-factor} for the range projection of the minus
block. We use the angular smoothing of this reference projection,
together with the tame estimates for the maps
\eqref{eq:A-complete-Kato-factor} in \cref{prop:all-order-adapters},
to obtain, for each $k\ge0$,
% AUTHOR QUERY: Identify the angular realization of P_{-,*}^R and
% prove its smoothing bounds. Finite-dimensional range does not imply
% finite Fourier support at an ellipse; smooth range and adjoint-range
% vectors would suffice. The symbol contour alone does not supply this.

\begin{equation}\label{app:cross:current-minus-smoothing}
 \begin{aligned}
 \|P_{-,\mathbf b}^Rf\|_{H_\theta^{s+k}}
 \le C_{s,k}\Bigl[&\|f\|_{H_\theta^s}+\|\mathbf b-\mathbf b_\rho\|_{
          \cB^{s+k+k_{\rm sm}(k)}}
       \|f\|_{H_\theta^{s_*}}\Bigr],
 \end{aligned}
\end{equation}
where $H_\theta^s$ denotes the Sobolev space of order $s$ in the angular
variable, with the other variables held fixed, and where $k_{\rm sm}(k)$ is
a finite shift in the regularity of the coefficients, independent of
$s$. The state $\mathbf b_\rho$ is the constant ellipse introduced in
\cref{sec:block}, written in the fixed coordinates, and
$\|\mathbf b-\mathbf b_\rho\|_{\cB^s}$ is the norm
\eqref{app:minus:eq:parameter-norm}. The constant is uniform on
the ball of states in the norm at the base index $s_*$, and the higher regularity norm of the state appears explicitly in the estimate.
Differentiating \eqref{app:cross:minus-Kato-conjugacy}, we obtain
\begin{equation}\label{app:cross:differentiated-minus-smoothing}
 \begin{aligned}
 &\|D_{\mathbf b}^mP_{-,\mathbf b}^R[\mathbf h]f\|_{H_\theta^{s+k}}\\
 &\quad\le C_{s,k,m}\Big\{
   \mathcal H_{0,m}(\mathbf h)\|f\|_{H_\theta^s}
  +\mathcal H_{s+k,m}(\mathbf h)\|f\|_{H_\theta^{s_*}}\\
 &\hspace{31mm}
  +\|\mathbf b-\mathbf b_\rho\|_{\cB^{s+k+k_m}}
       \mathcal H_{0,m}(\mathbf h)
       \|f\|_{H_\theta^{s_*}}\Big\},
 \end{aligned}
\end{equation}
with the products of tangent norms of
\eqref{app:cross:tangent-products}. Here and below, we write
$\bX_{\pm,\natg}^s:=\bX_{\pm,\mathrm{ret}}^s$ and
$\bY_{\pm,\natg}^s:=\bY_{\pm,\mathrm{ret}}^s$ for the block spaces of
the two high-frequency blocks, defined in
\eqref{app:minus:eq:retained-spaces} and
\eqref{app:positive:eq:positive-spaces}. We use this notation throughout
the remaining linear analysis. To absorb an angular derivative on the datum, we also
use smoothing on the input side of the projection. More precisely,
$P_{-,\mathbf b}^R\mathcal M_{\mathbf b}$ must extend from
$H_\theta^{s-1}$ to the minus source space at index $s$. Since
$\Lambda_\theta:H_\theta^s\to H_\theta^{s-1}$, this gives the bounded map
% AUTHOR QUERY: Establish the stated input-side bound, with its tame
% coefficient dependence. The displayed estimate for an output gain
% alone, with k=1 and Lambda_theta applied to f, does not prove it.

\begin{equation}\label{app:cross:mass-minus}
 P_{-,\mathbf b}^R\mathcal M_{\mathbf b}\Lambda_\theta
   :H_\theta^s\longrightarrow\bY_{-,\natg}^s
\end{equation}
into the combined center and affine source space, with tame estimates
for its derivatives with respect to $\mathbf b$. For the positive component, the dual norm
includes an inverse angular derivative, and we write
\begin{equation}\label{app:cross:mass-plus}
 \begin{aligned}
 \Lambda_\theta^{-1}P_{+,\mathbf b}^R
       \mathcal M_{\mathbf b}\Lambda_\theta
   ={}&\mathcal M_{\mathbf b}
      +\Lambda_\theta^{-1}
          [\mathcal M_{\mathbf b},\Lambda_\theta]-\Lambda_\theta^{-1}P_{-,\mathbf b}^R
          \mathcal M_{\mathbf b}\Lambda_\theta .
 \end{aligned}
\end{equation}
Since $\mathcal M_{\mathbf b}$ has order zero with smooth
coefficients, the commutator
$[\mathcal M_{\mathbf b},\Lambda_\theta]$ also has angular order
zero. Estimate \eqref{app:cross:mass-minus} bounds the last term. Thus,
the source norms absorb the angular derivative without loss of a power
of the cell frequency. In particular,
\cref{tab:incidence} still applies, and the norm of the source of the
center component requires no additional power of $\lambda_n$.

\Cref{tab:interior-maps} summarizes the interior estimates, all of
which use the weights of \cref{tab:block-weights}.
{\small
\begin{longtable}{@{}p{.27\textwidth}p{.66\textwidth}@{}}
\caption{The interior maps between block components and the reason why
each of them is bounded in the block norms.}
\label{tab:interior-maps}\\
\toprule
Map&Why it is bounded\\
\midrule
\endfirsthead
\toprule
Map&Why it is bounded\\
\midrule
\endhead
\bottomrule
\endfoot
compatible quotient&
The reconstruction \eqref{app:cross:quotient-reconstruction} of a
compatible source from its quotients is a multiplication operator in the
disk variables, so it is diagonal in the cell frequency, and it involves
no Hardy map.\\*
frame, maps \eqref{eq:A-complete-Kato-factor}, and projections&
On the annulus, the lower bound on the determinant shows that the
inverse $\mathbb B^{-1}$ of the frame \eqref{app:cross:current-frame} depends analytically on the
state. Derivatives of the frame enter through
\eqref{app:cross:current-connections}. The spectral projections
\eqref{app:cross:domain-range-contours}, defined by fixed contours, and
the maps \eqref{eq:A-complete-Kato-factor} have finite moments in the
sense of \cref{lem:macro-native-calculus}.\\
scalar differential source&
After projection onto the Cartesian range,
\eqref{app:cross:source-graph-chart} has no off-diagonal part of highest
order, by \eqref{app:cross:principal-cross-zero}. Derivatives of
coefficients and projections lower the order by one.\\
inversion of the factors $\cR\pm2i$ on the odd and on the high even modes&
We invert the factors $\cR\pm2i$ only on the modes on which they do not
vanish, listed after \cref{tab:source-orders}, and we measure the
positive output in the dual of the form space.\\
affine equations at the exceptional modes&
We solve for the independent forcings at $m=\pm2$ using the Hardy maps
\eqref{app:minus:eq:hardy-arrows}, which act in the radial variable
alone and are invariant under dilations. The weights $\mu$ and
$\lambda^{1/2}$ occur on both sides in \cref{tab:incidence}.\\
terms with one angular derivative&
The projection of the minus block is smoothing, and the dual of the
positive form includes the inverse derivative, by
\eqref{app:cross:mass-minus}--\eqref{app:cross:mass-plus}.\\
extraction of the center&
The angular smoothing used in \eqref{app:cross:mass-minus} controls
this projection. It has order zero in the cell variable. Its radial
output is measured in the $L^1$ source norm of the center equation in
\cref{app:minus}.\\
interior equations of the positive block&
The positive form maps the form space of the pairs of restrictions to
the cap and to the annulus to its dual, with one derivative on each
argument, and the weights $\lambda$ and $\lambda^{-1}$ of
\cref{tab:block-weights} then bound the map in the cell frequency.\\
value and first derivatives on the axis, and affine mean&
We use the entries $\lambda^{-1}$ and $\lambda^{-2}$ of
\cref{tab:block-weights}, together with the localized extension
operators of \cref{prop:fixed-graphs}.\\
mean, scale, gauge, and chart conditions&
The system is triangular, with finitely many angular modes and
invertible diagonal entries. The coefficients may depend on $\zeta$, and
a derivative on a coefficient leaves one fewer on the unknown.\\
analytic weight&
Conjugation by the analytic weight adds the multiplier $d_\gamma$ of
\eqref{eq:phase-damping} to the radial derivative. Its diagonal part
remains in the corresponding diagonal block. Its
commutator with a spectral projection has a kernel which contains the
factor $d_\gamma(r,\lambda_n)-d_\gamma(r,\lambda_m)$, as in
\eqref{app:cross:dyadic-commutator-kernel}, so it is controlled by the
moment of order one, and it vanishes at the circle, where the
projections do not depend on $\zeta$.\\
\end{longtable}
}
For the entries on the axis, $E_0$ extends a prescribed value and
$E_1$ extends prescribed first derivatives, using the localized
extensions of \cref{prop:fixed-graphs}. They satisfy
\begin{equation}\label{app:cross:axis-coretraction-ledger}
 \|E_0b\|_{\cA_\gamma^s}
   \le C_s\|e^{\sigma_0\Lambda}b\|_{H^{s-1}},\qquad
 \|E_1b\|_{\cA_\gamma^s}
   \le C_s\|e^{\sigma_0\Lambda}b\|_{H^{s-2}}.
\end{equation}
At cell frequency $n$, the extension is supported in a region of radius
comparable to $\lambda_n^{-1}$. A derivative of the cutoff
$\phi(\lambda_ny)$ therefore gives a factor $\lambda_n$, already
accounted for by the index shifts in
\eqref{app:cross:axis-coretraction-ledger}.

\subsection{Coupling between the two blocks on the cap}

To estimate the coupling on the cap, we use the isomorphisms \eqref{app:minus:eq:cap-map} and
\eqref{app:positive:eq:literal-reference-cap-realization} for the reference cap operators. At the state $\mathbf b$, we instead
invert the operators
$N_{-,\mathbf b}$ and $N_{+,\mathbf b}$ of
\cref{app:minus,app:positive}, which couple the cap and the annulus. We
denote the outward conormal
derivative on the interface by
$\Gamma_{c,+,N,\mathbf b}u:=\nabla_{c,\mathbf b}^{\rm out}u|_{r=r_j}$.
The reference operators and the traces at the state act as follows,
where $\Gamma_{c,-,\mathbf b}=\Gamma_{c,-,\mathbf b}^{\rm ph}$ is the
physical Cauchy trace of \eqref{app:minus:eq:cap-trace-physicalization}
\begin{equation}\label{app:cross:cap-diagonal-arrows}
 \begin{aligned}
 \cA_{c,-,\rho}^{\rm diag}:X_{c,-}^s&\xrightarrow{\ \sim\ }Y_{c,-}^s,
 &\Gamma_{c,-,\mathbf b}:X_{c,-}^s&\longrightarrow Z_{c,-}^s,\\
 \mathbf C_{+,\rho}^{\rm cap}:\mathfrak X_{c,+,H}^s
   &\xrightarrow{\ \sim\ }
       \mathfrak Y_{c,+,H}^s\oplus\mathfrak Z_{D,H}^s,
 &\Gamma_{c,+,N,\mathbf b}:X_{c,+}^s&\longrightarrow Z_N^s.
 \end{aligned}
\end{equation}
The two reference operators on the cap use different data. In the minus
block, the interior source and the axis and affine conditions determine
the outgoing Cauchy data, so the value on the interface cannot be
prescribed independently. In the positive block, we prescribe the value
on the interface and recover the outward conormal derivative.
The source spaces reflect this distinction.

We let $\mathscr C_{\mathbf b}^c$ be the Cartesian operator on the
cap at the state. The map $\iota_{c,\alpha}^X$ includes a cap field
from block $\alpha$, and $\pi_{c,\beta}^Y$ selects the cap equations
of block $\beta$. Here, $\alpha,\beta\in\{-,+\}$; the first sign in
$E_{\beta\alpha}^c$ denotes the output block and the second the input block. We define the \emph{off-diagonal operators on the cap} by
\begin{equation}\label{app:cross:direct-cap-arrows}
 \begin{aligned}
 E_{-+,\mathbf b}^{c}
   &=\pi_{c,-}^YP_{-,\mathbf b}^R\mathscr C_{\mathbf b}^{c}
       P_{+,\mathbf b}^D\iota_{c,+}^X:
       X_{c,+}^s\longrightarrow Y_{c,-}^s,\\
 E_{+-,\mathbf b}^{c}
   &=\pi_{c,+}^YP_{+,\mathbf b}^R\mathscr C_{\mathbf b}^{c}
       P_{-,\mathbf b}^D\iota_{c,-}^X:
       X_{c,-}^s\longrightarrow Y_{c,+}^s\oplus Z_D^s.
 \end{aligned}
\end{equation}
On each dyadic shell, with the cap
rescaled to a fixed size as in \cref{tab:block-weights}, the spaces of
the minus block are $X_{c,-}^s=H^{s+2}$ and $Y_{c,-}^s=H^s$, and the
spaces of the positive block are $X_{c,+}^s=H^{s+1}$ and
$Y_{c,+}^s=H^{s-1}$, with the Cartesian compatibility conditions
imposed. By \eqref{app:cross:principal-intertwining}, the off-diagonal
terms of second order cancel. Hence, $E_{-+,\mathbf b}^c$ has order one
and maps $H^{s+1}$ to $H^s$, while $E_{+-,\mathbf b}^c$ maps $H^{s+2}$
to $H^{s-1}$. The estimates for the inverse frame, the contour integrals
\eqref{app:cross:domain-range-contours}, and the coefficients show that
both families in \eqref{app:cross:direct-cap-arrows} depend analytically
on the state in these operator norms. They vanish at the circle, by the argument given
for \eqref{app:cross:circle-zero} below.

The traces in \eqref{app:cross:cap-diagonal-arrows} are bounded. In the minus block, we compose the trace with
$T_{-,\mathbf b}:Z_{c,-}^s\to Z_{a,-}^s$ from
\eqref{app:minus:eq:joined-operator}, without using its inverse. In the
positive block, the extension of the value and its outward conormal
derivative are part of the operator
\eqref{app:positive:eq:positive-operator}. Each composition remains in
its own diagonal block, and so do the identity terms in the transmission
conditions for the value and for the conormal derivative.

\subsection{Coupling at the interface and outer boundary}

We now include the conditions on the interface and on the outer
boundary. The corresponding off-diagonal maps act between the spaces
\begin{equation}\label{app:cross:retained-cross-types}
 \begin{aligned}
 E_{-+,\mathbf b}:X_{c,+}^s\oplus X_{a,+}^s
   &\longrightarrow Y_{c,-}^s\oplus Y_{a,-}^s\oplus Z_{a,-}^s,\\
 E_{+-,\mathbf b}:X_{c,-}^s\oplus X_{a,-}^s
   &\longrightarrow Y_{c,+}^s\oplus Y_{a,+}^s
        \oplus Z_D^s\oplus Z_N^s\oplus Z_B^s.
 \end{aligned}
\end{equation}
We first collect the operators on the interface and on the outer boundary,
before the spectral projections \eqref{app:cross:domain-range-contours}
are applied. For a pair $(u_c,u_a)$ of a field on the cap and a field
on the annulus, we put
\begin{equation}\label{app:cross:literal-retained-row-vector}
 \mathscr R_{\mathrm{ret},\mathbf b}(u_c,u_a)
 =\begin{pmatrix}
   \gamma_{-,\mathbf b}u_a
       -T_{-,\mathbf b}\Gamma_{c,-,\mathbf b}u_c\\
   \gamma_j^Du_a-\gamma_j^Du_c\\
   \nabla_{a,\mathbf b}^{\rm out}u_a
       +\nabla_{c,\mathbf b}^{\rm out}u_c\\
   \cB_{++,\mathbf b}u_a
  \end{pmatrix},
\end{equation}
whose rows are the mismatch of the incoming Cauchy data of the minus
block, the jump of the values, the sum of the outward conormal
derivatives, and the condition on the outer boundary. Here,
$\gamma_{-,\mathbf b}$ is the Cauchy trace of the minus block on the
inner boundary of the annulus at the state, and $\gamma_{-,\rho}$ in
\eqref{app:minus:eq:annulus-map} is the corresponding trace at the reference. If
$\iota_\alpha^X$ denotes the inclusion of the pair
$u_\alpha=(u_{c,\alpha},u_{a,\alpha})$ of a field on the cap and a field
on the annulus into the space of separate fields, the four off-diagonal terms are
\begin{equation}\label{app:cross:literal-retained-cross-rows}
 \begin{aligned}
 E_{-+,Z,\mathbf b}u_+
  &=\pi_{Z_{a,-}}^YP_{-,\mathbf b}^R
      \mathscr R_{\mathrm{ret},\mathbf b}
      P_{+,\mathbf b}^D\iota_+^Xu_+,\\
 (E_{+-,D,\mathbf b},E_{+-,N,\mathbf b},E_{+-,B,\mathbf b})u_-
  &=(\pi_D^Y,\pi_N^Y,\pi_B^Y)P_{+,\mathbf b}^R
      \mathscr R_{\mathrm{ret},\mathbf b}
      P_{-,\mathbf b}^D\iota_-^Xu_-.
 \end{aligned}
\end{equation}
Here, $\pi_{Z_{a,-}}^Y$ selects the first row, the mismatch of the
incoming Cauchy data in $Z_{a,-}^s$. The projections $\pi_D^Y$,
$\pi_N^Y$, and $\pi_B^Y$ select the remaining three rows, in
$Z_D^s$, $Z_N^s$, and $Z_B^s$, respectively. These are all the
noninterior terms in \eqref{app:cross:retained-cross-types}.

The value and conormal spaces satisfy
\begin{equation}\label{app:cross:balanced-trace-ledger}
 \|g_D\|_{Z_D^s}\simeq\|\lambda^{1/2}g_D\|_{H^s},\qquad
 \|g_N\|_{Z_N^s}+\|g_B\|_{Z_B^s}
    \simeq\|\lambda^{-1/2}(g_N,g_B)\|_{H^s},
\end{equation}
where $\lambda$ acts on the coefficient at cell frequency $n$ as
multiplication by $\lambda_n$, so an operator of order one from the
value to the conormal derivative has the weighted ratio
\begin{equation}\label{app:cross:value-conormal-ratio}
 \frac{\lambda^{-1/2}\lambda|g_D|}
      {\lambda^{1/2}|g_D|}=1.
\end{equation}
We write the conormal derivative at the state as
\begin{equation}\label{app:cross:exact-conormal}
 \nabla_{\mathbf b}u=E_{\mathbf b}D_ru+F_{\mathbf b}u.
\end{equation}
Here, $D_r$ is the radial derivative, $E_{\mathbf b}$ is its
coefficient, and $F_{\mathbf b}$ contains the tangential and zeroth-order
terms. Thus, $F_{\mathbf b}$ acts on the boundary value of $u$. After the spectral
projections \eqref{app:cross:domain-range-contours} are applied, the
off-diagonal term of highest order cancels by
\eqref{app:cross:principal-intertwining}, and
\eqref{app:cross:value-conormal-ratio} bounds the remaining terms of
first order, while the derivatives of $E_{\mathbf b}$, of the frames on
the traces, and of the projections are multipliers of order zero. The
term from the center to the positive component depends only on the
value and can be written as
\begin{equation}\label{app:cross:retained-center-value}
 (\nabla_{\mathbf b}z_+)_{C\to+}
       =\mathcal V_{+C,\mathbf b}z_C,
\end{equation}
with no normal derivative of $z_C$. Here, $z_+$ and $z_C$ are the
positive component and the center component of the unknown, the left
side denotes the part of the conormal operator of the positive block which
acts on the center component, and $\mathcal V_{+C,\mathbf b}$ is a
coefficient depending analytically on the state. The ratio of the
weights of its
components is at most $\lambda^{-1}$, and the analytic coefficient
$\mathcal V_{+C,\mathbf b}$ vanishes at the circle. We include this term
in the component $Z_N^s$ of $E_{+-,\mathbf b}$.

Writing $E_{\mathbf b,\gamma}=W_\gamma E_{\mathbf b}W_\gamma^{-1}$,
conjugation by the analytic weight adds the radial term
$E_{\mathbf b,\gamma}d_\gamma(r,\Lambda)g_D$, where
\begin{equation}\label{app:cross:phase-conormal}
 d_\gamma(r,\lambda)
 =\frac{\gamma r\lambda^2}{\sqrt{1+r^2\lambda^2}},
 \qquad |d_\gamma|\le\gamma\lambda,
\end{equation}
as in \eqref{app:positive:eq:exact-inner-conormal}. The sign is
carried by $E_{\mathbf b,\gamma}$, which uses the outward normal on the
side under consideration. The multiplier $d_\gamma$ is defined in
\eqref{eq:phase-damping}. Its bound by $\gamma\lambda$ and the ratio
\eqref{app:cross:value-conormal-ratio} control the cell-frequency
weights of the added term. The coefficient estimates control
$E_{\mathbf b,\gamma}$. At the interface, the two outward normals have
opposite orientations. For matching values and coefficients, the added
terms therefore cancel in Green's identity.

For the condition on the outer boundary, the off-diagonal term of
\eqref{app:cross:literal-retained-cross-rows} is
\begin{equation}\label{app:cross:physical-arrow}
 E_{+-,B,\mathbf b}
  :=\pi_B^YE_{+-,\mathbf b}:
       \bX_{-,\natg}^s\longrightarrow Z_B^s,
 \qquad E_{-+,B,\mathbf b}=0,
\end{equation}
since the boundary term from the positive block to itself already
belongs to the diagonal block $D_{+,\mathbf b}$ defined in
\eqref{app:cross:true-diagonal-correction} below, while the minus system has no independent
condition on the outer boundary. The bounds in the weighted trace norms
follow from
\eqref{app:minus:eq:annulus-map}--\eqref{app:minus:eq:undifferentiated-tame}
and \eqref{app:positive:eq:positive-tame}. The localized Poisson
extension \eqref{eq:A-boundary-coretraction}, whose decay rate depends
on the frequency, and the exponential decay
\eqref{eq:A-Poisson-phase-margin}, whose rate $(1-\gamma)\lambda_{mn}$
is positive for $\gamma<1$,
then give all moments for \eqref{app:cross:physical-arrow}. We first
estimate the trace and then
apply the frequency-dependent extension, since an extension with fixed
width would not give the required weighted bounds. These bounds use only $\widehat G_{-,\mathbf b}$ and
$\widehat G_{+,\mathbf b}$, so they are available before we correct
the diagonal inverses.

\subsection{Fourier estimates across the cutoff}

We use $P^X=P_J^X$ and $P^Y=P_J^Y$ for the projections onto
$|n|<2^J$, and $H^X=I-P^X$, $H^Y=I-P^Y$ for their complements.
The three blocks of \eqref{app:global:block-decomposition} estimated here are
\[
 A_{\mathbf b}=P^Y\mathbb D_{\mathbf b}P^X,\qquad
 B_{\mathbf b}=P^Y\mathbb D_{\mathbf b}H^X,\qquad
 C_{\mathbf b}=H^Y\mathbb D_{\mathbf b}P^X.
\]
Thus, $A_{\mathbf b}$ acts on the low frequencies, $B_{\mathbf b}$
takes high-frequency unknowns to low-frequency equations, and
$C_{\mathbf b}$ acts in the reverse direction. We write
$A_e=A_{\mathbf b_e}$ for the low block at the constant reference state. For a large enough shift $k_0$ in the regularity of the
coefficients, we put
\begin{equation}\label{app:global:low-distance}
 d_0(\mathbf b)
 :=\norm{\mathbf b-\mathbf b_e}_{\cB^{s_*+k_0}},
\end{equation}
which measures the distance to the constant reference state at the
base index. This distance controls both the change in the low block and
its coupling to the high frequencies.

\begin{lemma}
\label{app:global:low-leakage}
For every $J$, there are finite constants $C_{L,J}$ and $C_{X,J}$
and a finite shift $k\leq k_0$ such that
\begin{align}
 \norm{A_{\mathbf b}-A_e}_{
      P^X\bX_{\natg}^{s_*}\to P^Y\bY_{\natg}^{s_*}}
 &\leq C_{L,J}d_0(\mathbf b),                          \label{app:global:low-base}\\
 \norm{B_{\mathbf b}}_{H^X\bX_{\natg}^{s_*}\to P^Y\bY_{\natg}^{s_*}}
 +\norm{C_{\mathbf b}}_{P^X\bX_{\natg}^{s_*}\to H^Y\bY_{\natg}^{s_*}}
 &\leq C_{X,J}d_0(\mathbf b).                          \label{app:global:leakage-base}
\end{align}
For every $s\geq s_*$, every $u_L$ in the range of $P^X$, and every
$u_H$ in the range of $H^X$,
\begin{align}
 \norm{(A_{\mathbf b}-A_e)u_L}_{\bY_{\natg}^s}
 &\leq C_{J,s}\bigl[
 d_0(\mathbf b)\norm{u_L}_{\bX_{\natg}^s}
 +\norm{\mathbf b-\mathbf b_e}_{\cB^{s+k}}
       \norm{u_L}_{\bX_{\natg}^{s_*}}\bigr],           \label{app:global:low-tame}\\
 \norm{B_{\mathbf b}u_H}_{\bY_{\natg}^s}
 &\leq C_{J,s}\bigl[
 d_0(\mathbf b)\norm{u_H}_{\bX_{\natg}^s}
 +\norm{\mathbf b-\mathbf b_e}_{\cB^{s+k}}
       \norm{u_H}_{\bX_{\natg}^{s_*}}\bigr],           \label{app:global:B-tame}\\
 \norm{C_{\mathbf b}u_L}_{\bY_{\natg}^s}
 &\leq C_{J,s}\bigl[
 d_0(\mathbf b)\norm{u_L}_{\bX_{\natg}^s}
 +\norm{\mathbf b-\mathbf b_e}_{\cB^{s+k}}
       \norm{u_L}_{\bX_{\natg}^{s_*}}\bigr].           \label{app:global:C-tame}
\end{align}
These estimates hold for all interior, cap, quotient, axis, interface,
conormal, and outer boundary equations, and every derivative with respect to $\mathbf b$ of
fixed order obeys the corresponding tame estimate.
\end{lemma}

The proof uses the bounds at low frequency for fixed $J$, the
product estimate \eqref{eq:analytic-product}, and Fourier estimates
across the cutoff. In particular, it does not depend on the
high-frequency inverse or on the global inverse, so we may use
\eqref{app:global:leakage-base} before those inverses are constructed.

\begin{proof}
Since $P$ restricts the cell frequencies to a finite set, every cell
derivative is bounded on its range, with a constant depending on $J$. We
estimate the radial derivatives in the norms of the center and affine
components of \cref{app:minus}, and we measure the positive equation in
the form space and its dual, with the value and conormal traces measured
with the weights $\lambda^{1/2}$ and $\lambda^{-1/2}$ of
\cref{tab:block-weights}. On the cap, the minus block acts from
$H^{s+2}$ to $H^s$ and the positive block acts from $H^{s+1}$ to
$H^{s-1}$, as in \cref{tab:block-weights}. The quotient, mean, scale,
gauge, chart, and trace operators all have fixed order. Since
$\mathbb D_{\mathbf b}$ already includes the extraction of the scalar
differential source \eqref{app:cross:source-graph-chart}, its
cancellations of highest order are measured before $P$ is applied, and
the product estimate \eqref{eq:analytic-product} on this finite list of
operators proves \eqref{app:global:low-base} and \eqref{app:global:low-tame}.

At the constant reference state, the coefficients and reductions do not depend
on the cell variable, so that
\begin{equation}\label{app:global:reference-diagonality}
 \mathbb D_{\mathbf b_e}P^X=P^Y\mathbb D_{\mathbf b_e},
 \qquad B_{\mathbf b_e}=0,
 \qquad C_{\mathbf b_e}=0,
\end{equation}
and we need only estimate the coupling between low and high frequencies
in $\mathbb D_{\mathbf b}-\mathbb D_{\mathbf b_e}$. Consider an operator
$c(\zeta)D_\zeta^q$, which sends the input at cell frequency $m$ to the
output at cell frequency $n$ with the coefficient $c_{n-m}m^q$, where
$c_k$ are the Fourier coefficients of $c$. If
$|n|<2^J\le|m|$, then either $|m|<2^{J+1}$, or
$|n-m|\ge|m|-|n|\ge|m|/2$, and in both cases
\begin{equation}\label{app:global:fourier-crossing}
 \langle m\rangle^q\leq C_{J,q}\langle n-m\rangle^q.
\end{equation}
The factor $\langle n-m\rangle^q$ is absorbed by the Fourier
coefficients of the derivatives of $c$, so the derivatives of the
coefficient control the cell derivatives on the input. For the coupling
in the opposite direction, $|m|<2^J$, and the derivative on the input is
bounded by $C_{J,q}$. Consider a map from a component $\alpha$ to a
component $\beta$ with the weights $w_\alpha=d_\alpha^X$ and
$w_\beta=d_\beta^Y$ of \cref{tab:block-weights}, and let
$q_{\beta\alpha}$ be the number of derivatives in $\zeta$ applied to the
input. The inequality \eqref{app:cross:peetre} then gives
\begin{equation}\label{app:global:native-fourier-crossing}
 \mathbf1_{\{|n|<2^J\}}\mathbf1_{\{|m|\geq2^J\}}
 \frac{w_\beta(n)}{w_\alpha(m)}
 \langle m\rangle^{q_{\beta\alpha}}
 \leq C_J\langle n-m\rangle^{k_{\beta\alpha}},
\end{equation}
with the same estimate in the reverse direction. Here,
$k_{\beta\alpha}$ is a finite exponent depending only on the two
components, obtained from the exponents of \eqref{app:cross:peetre} and
from $q_{\beta\alpha}$. For example, the first row of
\cref{tab:incidence} has $w_\alpha=1$, $w_\beta(n)=\lambda_n^{-1}$,
and $q_{\beta\alpha}=1$, so its ratio across the cutoff is
$\lambda_n^{-1}\lambda_m$, controlled by a power of
$\langle n-m\rangle$. Radial derivatives are estimated separately in
the norms of the center, of the affine component, of the positive form,
and of the cap. Multiplication by a coefficient of the state is a
convolution in the cell frequency, and the weight $e^{\Phi_r(n)}$ is
submultiplicative by \eqref{eq:phase-submultiplicative}. Hence, the
weighted kernel \eqref{app:cross:normalized-kernel} of each coupling
term is bounded by the weighted Fourier coefficients of the
coefficient, times the polynomial factor of
\eqref{app:global:native-fourier-crossing}. Summing the rows and the
columns as in Schur's test then proves \eqref{app:global:leakage-base},
\eqref{app:global:B-tame}, and \eqref{app:global:C-tame} for each
component. Since $P$ and $H$ are fixed, the tame product estimate also
applies after differentiation in the parameters.
\end{proof}

Once $A_e$ is inverted, \eqref{app:global:low-base} also gives the
invertibility of the low block near the reference state. Indeed, if
\begin{equation}\label{app:global:low-neumann-margin}
 \norm{A_e^{-1}(A_{\mathbf b}-A_e)}<1,
\end{equation}
then $A_{\mathbf b}=A_e\bigl(I+A_e^{-1}(A_{\mathbf b}-A_e)\bigr)$ is
invertible by a Neumann series.

\subsection{The off-diagonal estimate}

We use the constant ellipse $\mathbf b_\rho$ as the reference for the
coefficient differences. We
let $a_\circ(\rho)=a^0_{\rho,p_e}$ and $M_\circ(\rho)=M_{\rho,p_e}$
denote its unperturbed solution and its matrix, whose axes do not rotate
with $\zeta$, and we let $M(p)=M_{\rho,p}$ denote the periodic matrix
\eqref{eq:fixed-M-lambda} with the parameters $p$ of the state. For a
finite shift $k_\times$ in the regularity of the coefficients, fixed
with the estimate \eqref{app:cross:main-estimate} below, we put
\begin{equation}\label{app:cross:circle-distance}
 \begin{aligned}
 \eta_{\rm circ}(\mathbf b)
 :={}&\|a-a_\circ(\rho)\|_{\cA_\gamma^{s_*+k_\times}}
       +|\eps|+\|M(p)-M_\circ(\rho)\|_{\cA_\gamma^{s_*+k_\times}},
 \qquad
 \chi(\mathbf b)=|\rho|+\eta_{\rm circ}(\mathbf b).
 \end{aligned}
\end{equation}
Thus, $\chi(\mathbf b)$ measures the distance from the state to the
circle, and the off-diagonal blocks will be bounded by this distance.

In the fixed coordinates, we define the off-diagonal maps at high
frequency by
\[
 E_{\beta\alpha,\mathbf b}
 =P_{\beta,\mathbf b}^R H^Y\mathbb L_{\mathbf b}H^X
   P_{\alpha,\mathbf b}^D,
 \qquad \alpha\ne\beta,
\]
where $H^X$ and $H^Y$ denote $H_J$ on the domain and on the range, with
the cap coordinates of the two high-frequency blocks and the finitely many
coordinates of the triangular system included. We set
$\mathbb A_{\mathbf b}:=\iota_{Y,\mathbf b}
\mathsf J_{Y,\mathbf b}^{\rm raw}L_{\mathbf b}
\mathsf R_{X,\mathbf b}$, we let $q_\beta^Z$ denote the projection onto
the auxiliary data of the range of the component $\beta$, and we write
$\jmath_{Z,\beta,\mathbf b}
 :=\iota_\beta^Y\pi_\beta^Y\jmath_{Z,\mathbf b}$ on that component. With
this notation, \eqref{eq:A-stabilized-native-formula} gives
\begin{equation}\label{app:cross:exact-stabilized-cross-block}
 E_{\beta\alpha,\mathbf b}
 =\pi_\beta^Y\mathbb A_{\mathbf b}\iota_\alpha^X
  +\jmath_{Z,\beta,\mathbf b}q_\beta^Z
       \mathsf C_{X,\mathbf b}\iota_\alpha^X.
\end{equation}
The subtraction in \eqref{eq:A-stabilizing-range-shear} removes
$\mathbb B_{\mathbf b}\mathsf C_{X,\mathbf b}$, the interior response
to the constraint data in \eqref{eq:A-natural-broken-formula}. The second term of
\eqref{app:cross:exact-stabilized-cross-block} is the projection of
$E_{\beta\alpha,\mathbf b}$ onto the auxiliary data of the component
$\beta$, so the estimate \eqref{app:cross:main-estimate} below bounds it
without using a diagonal inverse.

\begin{proposition}[The off-diagonal estimate]\label{prop:macro-cross}
Fix $0\le\rho\le\rho_0<1/4$, the strength $\gamma$ of the analytic
weight, the parameter $\delta_A$ of \eqref{app:cross:cell-weights}, and
the cutoff index $J$. On a sufficiently small ball of states, the
off-diagonal maps are bounded from index $s$ to index $s$, and for every
$s\ge s_*$ they satisfy
\begin{equation}\label{app:cross:main-estimate}
 \begin{aligned}
 &\|E_{-+,\mathbf b}u_+\|_{\bY_{-,\natg}^s}
  +\|E_{+-,\mathbf b}u_-\|_{\bY_{+,\natg}^s}\\
 &\quad\le C_s\Big[
  \chi(\mathbf b)
    \bigl(\|u_+\|_{\bX_{+,\natg}^s}
          +\|u_-\|_{\bX_{-,\natg}^s}\bigr)\\
 &\hspace{29mm}
  +\|\mathbf b-\mathbf b_\rho\|_{\cB^{s+k_\times}}
    \bigl(\|u_+\|_{\bX_{+,\natg}^{s_*}}
          +\|u_-\|_{\bX_{-,\natg}^{s_*}}\bigr)
 \Big].
 \end{aligned}
\end{equation}
For each number $m$ of derivatives with respect to $\mathbf b$, there is a finite
shift $k_m$, independent of $s$, such that
\begin{equation}\label{app:cross:all-order-estimate}
 \begin{aligned}
 &\|D_{\mathbf b}^mE_{-+,\mathbf b}[\mathbf h]u_+\|_{
       \bY_{-,\natg}^s}
  +\|D_{\mathbf b}^mE_{+-,\mathbf b}[\mathbf h]u_-\|_{
       \bY_{+,\natg}^s}\\
 &\quad\le C_{s,m}\Big\{
   \mathcal H_{0,m}(\mathbf h)
       \bigl(\|u_+\|_{\bX_{+,\natg}^s}
             +\|u_-\|_{\bX_{-,\natg}^s}\bigr)\\
 &\hspace{21mm}
  +\mathcal H_{s,m}(\mathbf h)
       \bigl(\|u_+\|_{\bX_{+,\natg}^{s_*}}
             +\|u_-\|_{\bX_{-,\natg}^{s_*}}\bigr)\\
 &\hspace{21mm}
  +\|\mathbf b-\mathbf b_\rho\|_{\cB^{s+k_m}}
          \mathcal H_{0,m}(\mathbf h)
       \bigl(\|u_+\|_{\bX_{+,\natg}^{s_*}}
             +\|u_-\|_{\bX_{-,\natg}^{s_*}}\bigr)
  \Big\}.
 \end{aligned}
\end{equation}
The constants are uniform for $0\le\rho\le\rho_0$, for both signs of the
frequency, and for every upper Fourier cutoff, including the case
without a cutoff. The operators commute with complex conjugation, and
\begin{equation}\label{app:cross:circle-zero}
 E_{-+,\mathbf b_0}=E_{+-,\mathbf b_0}=0.
\end{equation}
\end{proposition}

\begin{proof}
We first cancel the off-diagonal part of highest order by
\eqref{app:cross:ordered-output}--\eqref{app:cross:principal-cross-zero},
and only then extract the differential source, so that the resulting
scalar equation has order two by \eqref{app:cross:canonical-syzygy} and
\eqref{app:cross:current-syzygy}. We estimate the remaining interior
terms using \cref{tab:incidence},
\eqref{app:cross:mass-minus}--\eqref{app:cross:mass-plus}, and
\cref{tab:interior-maps}. For the terms on the cap, on the interface, and
on the boundary, we use \eqref{app:cross:direct-cap-arrows} and
\eqref{app:cross:retained-cross-types}--\eqref{app:cross:physical-arrow}.
Combining these bounds with \cref{lem:macro-native-calculus}, we obtain
the tame estimate from index $s$ to index $s$.

It remains to include the cutoffs and the derivatives with respect to $\mathbf b$. On the
interior components, we write
$\mathbb L_{\mathbf b}=Q_{\mathbf b}+\mathbb R_{\mathbf b}$, where
$Q_{\mathbf b}$ now denotes the operator with the principal symbol
\eqref{app:cross:descriptor-symbol}, and where $\mathbb R_{\mathbf b}$
collects the terms of lower order, which come from the derivatives of
the frame \eqref{app:cross:current-connections} and from the passage
from symbols to operators. We set
\begin{equation}\label{app:cross:cutoff-commutators}
 C_{\beta,\mathbf b}^R=[P_{\beta,\mathbf b}^R,H^Y],
 \qquad
 C_{\alpha,\mathbf b}^D=[H^X,P_{\alpha,\mathbf b}^D].
\end{equation}
Since $P_\beta^RH^Y=H^YP_\beta^R+C_\beta^R$ and
$H^XP_\alpha^D=P_\alpha^DH^X+C_\alpha^D$, the off-diagonal term of
highest order decomposes as
\begin{equation}\label{app:cross:cutoff-decomposition}
 \begin{aligned}
 P_\beta^RH^YQ_{\mathbf b}H^XP_\alpha^D
 ={}&H^YP_\beta^RQ_{\mathbf b}P_\alpha^DH^X+H^YP_\beta^RQ_{\mathbf b}C_\alpha^D
   +C_\beta^RQ_{\mathbf b}P_\alpha^DH^X
   +C_\beta^RQ_{\mathbf b}C_\alpha^D.
 \end{aligned}
\end{equation}
The principal symbol of the first term vanishes by
\eqref{app:cross:principal-cross-zero}. The operator itself can still
contain lower-order terms, produced when derivatives fall on the
projections. Each of the other three terms contains a commutator
between a cutoff and a projection. We include all these terms, together
with $P_\beta^RH^Y\mathbb R_{\mathbf b}H^XP_\alpha^D$, in
$\mathfrak R_{\beta\alpha,\mathbf b}$. Thus,
\begin{equation}\label{app:cross:actual-cutoff-cross}
 P_\beta^RH^Y\mathbb L_{\mathbf b}H^XP_\alpha^D
 =\mathfrak R_{\beta\alpha,\mathbf b},
\end{equation}
and \cref{lem:macro-native-calculus}, applied to the finite moments of
the cutoffs, of the projections, and of the coefficients, gives, for
every $m$,
\begin{equation}\label{app:cross:differentiated-cutoff-remainder}
 \begin{aligned}
 &\|D_{\mathbf b}^m\mathfrak R_{\beta\alpha,\mathbf b}
       [\mathbf h]u_\alpha\|_{\bY_{\beta,\natg}^s}\\
 &\quad\le C_{s,m}\Big\{
  \mathcal H_{0,m}(\mathbf h)\|u_\alpha\|_{\bX_{\alpha,\natg}^s}
  +\mathcal H_{s,m}(\mathbf h)
       \|u_\alpha\|_{\bX_{\alpha,\natg}^{s_*}}
  +\|\mathbf b-\mathbf b_\rho\|_{\cB^{s+k_m}}
       \mathcal H_{0,m}(\mathbf h)
       \|u_\alpha\|_{\bX_{\alpha,\natg}^{s_*}}\Big\}.
 \end{aligned}
\end{equation}
For the cutoffs $H^X$ and $H^Y$ in \eqref{app:cross:cutoff-commutators},
which are indicator functions of sets of cell frequencies, the
commutator is supported where the input and output cell frequencies lie
on opposite sides of the cutoff. The estimates
\eqref{app:global:fourier-crossing}--\eqref{app:global:native-fourier-crossing}
then control derivatives in $\zeta$ on the input by derivatives of the coefficient,
with constants depending on the cutoff $J$. This control differs from the gain
of one order for the smooth dyadic cutoffs in
\eqref{app:cross:dyadic-commutator-kernel}. Since the cutoffs are fixed,
derivatives with respect to $\mathbf b$ are again commutators of the same form, and
\eqref{app:cross:cutoff-decomposition}--\eqref{app:cross:differentiated-cutoff-remainder}
apply after differentiation.
% AUTHOR QUERY: The Fourier crossing estimate controls cell derivatives,
% not transverse derivatives. Prove the required transverse-order bounds
% for the cutoff terms, especially the positive-cap H^{s+1} to minus-cap
% H^s map, before concluding differentiated-cutoff-remainder.

The cancellation of highest order for $\alpha\ne\beta$ holds on the
whole ball of states. Differentiating the symbol identity, we obtain
\begin{equation}\label{app:cross:grouped-principal-derivatives}
 D_{\mathbf b}^m\bigl(
   P_{\beta,\mathbf b}^RQ_{\mathbf b}P_{\alpha,\mathbf b}^D
   \bigr)[\mathbf h]=0
 \qquad(m\ge0).
\end{equation}
This is an identity for symbols. It applies to the sum of the terms
in Leibniz' rule, so we keep them together and then estimate the
lower-order terms arising in the corresponding operators. Likewise,
the symbol $s_{G_{\mathbf b}}$ of the differential source component in
\eqref{app:cross:current-source-row} satisfies the quadratic identity
\eqref{app:cross:current-syzygy}, so
\begin{equation}\label{app:cross:grouped-syzygy-derivatives}
 D_{\mathbf b}^m
  \bigl[s_{G_{\mathbf b}}\ell_{1,G_{\mathbf b}}\bigr][\mathbf h]
 =D_{\mathbf b}^m
    \bigl(\xi^*G_{\mathbf b}^{-1}\xi\bigr)[\mathbf h]e_C^*.
\end{equation}
The differentiated terms have the same orders as those in
\cref{tab:incidence}, including the terms on the cap, on the traces, and
on the outer boundary.

It remains to estimate the compositions of these maps. For two
composable families $S_{\mathbf b}$ and $T_{\mathbf b}$ satisfying
\eqref{app:cross:all-order-estimate}, Leibniz' rule gives
\begin{equation}\label{app:cross:composition-leibniz}
 D_{\mathbf b}^m(T_{\mathbf b}S_{\mathbf b})[\mathbf h]
 =\sum_{I\subset\{1,\ldots,m\}}
 D_{\mathbf b}^{|I|}T_{\mathbf b}[\mathbf h_I]
 D_{\mathbf b}^{m-|I|}S_{\mathbf b}[\mathbf h_{I^c}].
\end{equation}
We use the estimate at index $s$ for one factor and the estimate at the
base index $s_*$ for the others. The products of lower regularity tangent norms give $\mathcal H_{0,m}$.
Measuring one tangent direction at a higher regularity index and summing
over that choice gives $\mathcal H_{s,m}$. Applying this to the
finite list of factors in the interior, on the cap, and on the boundary,
we obtain \eqref{app:cross:all-order-estimate}, with at most one
unknown, coefficient, or tangent direction measured in a higher regularity norm.

At the constant circle $\mathbf b_0$, the coefficients are independent
of $\zeta$. The operator preserves each angular mode and parity, the
four-dimensional center space, and the complementary positive
subspace. Thus, the off-diagonal maps vanish there, giving
\eqref{app:cross:circle-zero}. Since the off-diagonal maps depend
analytically on the state and vanish at $\mathbf b_0$, the mean
value formula bounds them by the distance $\chi(\mathbf b)$ to the
circle times the bounds for their derivatives, which proves
\eqref{app:cross:main-estimate}. The constants remain uniform when the
center speeds coincide. Indeed, we use one contour around the four
center roots, which is separated from the affine and positive roots, and
we do not select individual eigenvectors. The gap of the positive form
of \cref{app:positive} also stays positive at $\rho=0$. We treat the
orientation $\alpha_0$ of the ellipse as a parameter of the unperturbed
solutions, so every constant is uniform through $\rho=0$.

Finally, we let $\mathfrak c_{X,\alpha}$ and $\mathfrak c_{Y,\beta}$
denote complex conjugation on the complexified spaces. The Cartesian operators, weight \eqref{app:cross:phase}, cap
operators, Poisson extensions, and traces preserve real data. The
cutoff $H_J$ is invariant under $n\mapsto-n$, the contours in
\eqref{app:cross:domain-range-contours} are invariant under conjugation,
and the exceptional equations at $m=\pm2$ form a conjugate pair. The transfer maps of
\cref{prop:all-order-adapters}, which commute with complex conjugation,
therefore give, for real $\mathbf b$,
\begin{equation}\label{app:cross:real-closure}
 E_{\beta\alpha,\mathbf b}\mathfrak c_{X,\alpha}
  =\mathfrak c_{Y,\beta}E_{\beta\alpha,\mathbf b}.
\end{equation}
Finite sums, the products in \eqref{app:cross:composition-leibniz}, and
differentiation in real tangent directions preserve
\eqref{app:cross:real-closure}, which completes the proof.
\end{proof}

\subsection{Correcting the diagonal inverses}

We now compare the separately constructed operators
$\widehat D_{\alpha,\mathbf b}$ with the diagonal blocks
$D_{\alpha,\mathbf b}$ of the full augmented operator. Their difference
is the correction $K_{\alpha,\mathbf b}$ defined below. For $\alpha\in\{-,+\}$, we define the projections
of the auxiliary variables onto the auxiliary data of the component
$\alpha$ and onto the auxiliary data at low frequency by
\begin{equation}\label{app:cross:three-way-auxiliary-projections}
 Q_\alpha^Z
 :=\operatorname {aux}_{Y,\mathbf b}
      \iota_\alpha^Y\pi_\alpha^Y\jmath_{Z,\mathbf b},
 \qquad
 Q_P^Z:=\operatorname {aux}_{Y,\mathbf b}P^Y\jmath_{Z,\mathbf b},
 \qquad
 I_Z=Q_P^Z+Q_-^Z+Q_+^Z.
\end{equation}
Here, $I_Z$ is the identity on the auxiliary variables, and the last
identity expresses that the low frequencies and the two high-frequency
blocks exhaust the range. On $\operatorname {Ran}Q_\alpha^Z$, the map
$\jmath_{Z,\alpha,\mathbf b}=\iota_\alpha^Y\pi_\alpha^Y\jmath_{Z,\mathbf b}$
of \eqref{app:cross:exact-stabilized-cross-block} is the inclusion
into the auxiliary variables of this component, so that on the full
auxiliary space
\begin{equation}\label{app:cross:component-auxiliary-inclusion}
 \iota_\alpha^Y\pi_\alpha^Y\jmath_{Z,\mathbf b}
 =\jmath_{Z,\alpha,\mathbf b}Q_\alpha^Z.
\end{equation}
We write $\alpha^c$ for the complement of $\alpha$, which includes
both the other high-frequency block and the auxiliary variables at low
frequency. We set
\begin{equation}\label{app:cross:true-diagonal-correction-block}
 \begin{aligned}
 C_{\alpha^c\alpha,\mathbf b}
  &:=(I_Z-Q_\alpha^Z)\mathsf C_{X,\mathbf b}\iota_\alpha^X,\\
 B_{\alpha\alpha^c,\mathbf b}
  &:=\pi_\alpha^Y\mathbb B_{\mathbf b}
       (I_Z-Q_\alpha^Z),\\
 K_{\alpha,\mathbf b}
  &:=B_{\alpha\alpha^c,\mathbf b}
       C_{\alpha^c\alpha,\mathbf b}:
       \bX_{\alpha,\natg}^s\longrightarrow\bY_{\alpha,\natg}^s.
 \end{aligned}
\end{equation}
Here, $C_{\alpha^c\alpha,\mathbf b}$ extracts the constraint data
outside block $\alpha$, and $B_{\alpha\alpha^c,\mathbf b}$ returns
their extension's interior source to block $\alpha$. Their composition
is the correction $K_{\alpha,\mathbf b}$. The map
$\mathbb B_{\mathbf b}$ is defined by
\eqref{eq:A-coretraction-bulk-response}; the frame matrix
\eqref{app:cross:current-frame} is denoted by $\mathbb B$ without a subscript.

The following lemma shows that
$D_{\alpha,\mathbf b}=\widehat D_{\alpha,\mathbf b}-K_{\alpha,\mathbf b}$.
We estimate $K_{\alpha,\mathbf b}$ in the base norm and use its smallness
to correct the inverse of $\widehat D_{\alpha,\mathbf b}$.
\begin{lemma}
\label{app:cross:lem:true-diagonal}
For $\alpha\in\{-,+\}$, the diagonal block $D_{\alpha,\mathbf b}$ of the
augmented operator and the separately constructed operator
$\widehat D_{\alpha,\mathbf b}$ satisfy the identity
\begin{equation}\label{app:cross:true-diagonal-correction}
 D_{\alpha,\mathbf b}
 :=\pi_\alpha^Y\mathbb D_{\mathbf b}\iota_\alpha^X
 =\widehat D_{\alpha,\mathbf b}-K_{\alpha,\mathbf b}.
\end{equation}
At the base index,
\begin{equation}\label{app:cross:diagonal-correction-bound}
 \|K_{\alpha,\mathbf b}\|
 \le M_B\bigl(K_Z\chi(\mathbf b)+K_{Z,J}d_0(\mathbf b)\bigr),
 \qquad
 d_0(\mathbf b):=\norm{\mathbf b-\mathbf b_e}_{\cB^{s_*+k_0}},
\end{equation}
where $M_B$ bounds $B_{\alpha\alpha^c,\mathbf b}$,
$K_Z$ bounds the part of $C_{\alpha^c\alpha,\mathbf b}$ in the other
high-frequency block, by \eqref{app:cross:main-estimate}, and $K_{Z,J}$
bounds its part at low frequency, by the argument of
\eqref{app:global:leakage-base}.
 If
\begin{equation}\label{app:cross:diagonal-correction-margin}
 q_{\rm corr}:=\max_{\alpha\in\{-,+\}}
   \|\widehat G_{\alpha,\mathbf b}K_{\alpha,\mathbf b}\|<1,
\end{equation}
then the diagonal block has the two-sided inverse
\begin{equation}\label{app:cross:true-diagonal-inverse}
 G_{\alpha,\mathbf b}
  =(I-\widehat G_{\alpha,\mathbf b}K_{\alpha,\mathbf b})^{-1}
       \widehat G_{\alpha,\mathbf b}
  =\widehat G_{\alpha,\mathbf b}
       (I-K_{\alpha,\mathbf b}\widehat G_{\alpha,\mathbf b})^{-1}.
\end{equation}
These inverses, together with the traces of the solutions which appear
in \eqref{app:minus:eq:undifferentiated-tame} and
\eqref{app:positive:eq:positive-tame}, obey the same bounds as
$\widehat G_{\alpha,\mathbf b}$, without loss of derivatives, and the corresponding
tame bounds for derivatives with respect to $\mathbf b$ of every order.
\end{lemma}

\begin{proof}
We write $C_\alpha:=\mathsf C_{X,\mathbf b}\iota_\alpha^X$ and
$\mathbb A_{\alpha\alpha}:=\pi_\alpha^Y
\mathbb A_{\mathbf b}\iota_\alpha^X$. We take the $\alpha$ diagonal block of
\eqref{eq:A-natural-broken-formula} and apply
the diagonal block $\mathscr S_{\alpha\alpha}$ of the map in
\eqref{eq:A-sheared-component-inverse}, which subtracts
$\pi_\alpha^Y\mathbb B_{\mathbf b}Q_\alpha^ZC_\alpha$, the response of
the interior to the auxiliary data of the component $\alpha$ alone.
Using \eqref{app:cross:component-auxiliary-inclusion}, we obtain
\[
 \widehat D_{\alpha,\mathbf b}
 =\mathbb A_{\alpha\alpha}
  +\pi_\alpha^Y\mathbb B_{\mathbf b}
       (I_Z-Q_\alpha^Z)C_\alpha
  +\jmath_{Z,\alpha,\mathbf b}Q_\alpha^ZC_\alpha,
\]
while \eqref{eq:A-stabilized-native-formula} gives
\[
 \pi_\alpha^Y\mathbb D_{\mathbf b}\iota_\alpha^X
 =\mathbb A_{\alpha\alpha}
  +\jmath_{Z,\alpha,\mathbf b}Q_\alpha^ZC_\alpha.
\]
Their difference is $K_{\alpha,\mathbf b}$, which proves
\eqref{app:cross:true-diagonal-correction}.

We estimate the two parts of the correction separately. Projecting
\eqref{app:cross:exact-stabilized-cross-block} onto the auxiliary
variables, we obtain the part of $C_{\alpha^c\alpha,\mathbf b}$ in the
other high-frequency component, whose norm at the base index is
bounded by $K_Z\chi(\mathbf b)$, by \eqref{app:cross:main-estimate}. For
the part at low frequency, \eqref{eq:A-auxiliary-row-identity}
identifies $Q_P^Z\mathsf C_{X,\mathbf b}\iota_\alpha^X$ with the
auxiliary projection of $P^Y\mathbb D_{\mathbf b}H^X$, which vanishes at
$\mathbf b_e$ because the reference is diagonal in the cell frequency.
For fixed $J$, the Fourier estimate
\eqref{app:global:native-fourier-crossing} then bounds this projection
by $K_{Z,J}d_0(\mathbf b)$, as in \eqref{app:global:leakage-base}. Finally,
\eqref{eq:A-coretraction-bulk-response} and
\cref{prop:all-order-adapters} bound $B_{\alpha\alpha^c,\mathbf b}$,
which proves \eqref{app:cross:diagonal-correction-bound}.

Factoring
\[
 D_{\alpha,\mathbf b}
 =\widehat D_{\alpha,\mathbf b}
      (I-\widehat G_{\alpha,\mathbf b}K_{\alpha,\mathbf b})
 =(I-K_{\alpha,\mathbf b}\widehat G_{\alpha,\mathbf b})
      \widehat D_{\alpha,\mathbf b}
\]
proves \eqref{app:cross:true-diagonal-inverse} and both inverse
identities.
% AUTHOR QUERY: Base-norm Neumann invertibility does not alone prove
% regularity on the full scale. Give a commutator or regularity estimate
% for I-\widehat G_\alpha K_\alpha on the fixed coefficient neighborhood.
The tame bounds follow from the resolvent identity for the Neumann
series and from \eqref{app:cross:composition-leibniz}, using the direct
estimates for $B_{\alpha\alpha^c,\mathbf b}$ and
$C_{\alpha^c\alpha,\mathbf b}$ at every order of differentiation.
\end{proof}

Thus, the construction proceeds in the order
\begin{equation}\label{app:cross:noncircular-diagonal-order}
 \widehat G_\pm\ \longrightarrow\
 \{E_{\mp\pm},\,Q_P^Z\mathsf C_X\iota_\pm^X\}\
 \longrightarrow\ G_\pm\ \longrightarrow\
 G_H.
\end{equation}
The estimates for the off-diagonal maps and low-frequency auxiliary
data therefore precede the correction of the diagonal inverses. Here, $G_H$ denotes the inverse of the full high-frequency block
$D_{H,\mathbf b}$ in \eqref{app:global:block-decomposition}; it must
also account for the couplings $E_{-+}$ and $E_{+-}$. We construct it
after the low-frequency reference inverse. The low block $A_{\mathbf b}$
requires no such correction.

Finally, the composition used in the Neumann series at high frequency is
\begin{equation}\label{app:cross:return-types}
 \bX_{-,\natg}^{s_*}
 \xrightarrow{\ E_{+-,\mathbf b}\ }\bY_{+,\natg}^{s_*}
 \xrightarrow{\ D_{+,\mathbf b}^{-1}\ }\bX_{+,\natg}^{s_*}
 \xrightarrow{\ E_{-+,\mathbf b}\ }\bY_{-,\natg}^{s_*}
 \xrightarrow{\ D_{-,\mathbf b}^{-1}\ }\bX_{-,\natg}^{s_*},
\end{equation}
with each factor acting between the indicated spaces without loss of
derivatives. Combining \eqref{app:cross:main-estimate} with the bounds
for the diagonals at the same index, and using
\eqref{app:cross:bounded-product}, we obtain
\begin{equation}\label{app:cross:return-bound}
 \|D_{-,\mathbf b}^{-1}E_{-+,\mathbf b}
       D_{+,\mathbf b}^{-1}E_{+-,\mathbf b}\|
 \le K_-K_+K_\times^2\chi(\mathbf b)^2.
\end{equation}
Here, $K_\times$ is the constant of \eqref{app:cross:main-estimate} at
the base index $s_*$, and $K_\pm$ bound the diagonal block
inverses \eqref{app:cross:true-diagonal-inverse} under
\eqref{app:cross:diagonal-correction-margin}. This estimate requires only the operator norms of the diagonal
inverses. No bounds for their symbols or moments are needed. For their derivatives with respect to $\mathbf b$, we use
\[
 D(D_{\alpha,\mathbf b}^{-1})[h]
 =-D_{\alpha,\mathbf b}^{-1}(DD_{\alpha,\mathbf b}[h])
      D_{\alpha,\mathbf b}^{-1}
\]
and estimate the resulting compositions on the same block spaces.
\section{The low-frequency inverse at a fixed ellipse}\label{app:lowref}

For the low-frequency construction in \cref{app:global}, we need an
inverse without loss of derivatives in the block norms. We obtain this
bound directly in the block variables, first at the circle and then at
nearby constant ellipses. The Cartesian estimate of
\cref{app:reference} loses derivatives of the source, so it does not
give the bound needed here. The resulting block estimate also governs
the choice of the reference ellipse parameter $\rho_*$.

We fix the base index $s_*$, the transfer maps of \cref{sec:block}, and
the identifications \eqref{eq:A-complete-Kato-factor} of the trace
spaces with fixed products. We use the projections onto the cell
frequencies \eqref{eq:A-exact-cell-split}, as in
\eqref{app:global:block-decomposition},
\[
 H_J=\mathbf1_{\{|n|\geq2^J\}},\qquad P_J=I-H_J
       =\mathbf1_{\{|n|<2^J\}},
\]
with the versions $P_J^X$ and $P_J^Y$ on the domain and on the range. We
apply these projections in the block variables, that is, after the
transfer maps \eqref{eq:A-JX-factorization} and
\eqref{eq:A-JY-factorization} have been applied to the domain and to the
source. We let $\mathbb D_{\rho,e}:=\mathbb D_{\mathbf b_\rho}$ be the
augmented operator \eqref{eq:A-stabilized-native-operator} at the
constant ellipse $\mathbf b_\rho$ with ellipse parameter $\rho$ from
\cref{sec:block}, acting between the fixed block spaces, and we set
\begin{equation}\label{appref:low-compression}
 A_{\rho,J}=P_J^Y\mathbb D_{\rho,e}P_J^X,
\end{equation}
which is the low block of \cref{app:macro-cross} at the constant
ellipse. Here, $A_{\rho,J}=A_{\mathbf b_\rho}$ in the notation of
\eqref{app:global:block-decomposition}, with the cutoff $J$ made
explicit, and the low block $A_e$ of \cref{app:macro-cross} is
$A_{\rho_*,J}$. At the constant ellipse, the identities
\eqref{eq:A-auxiliary-row-identity}--\eqref{eq:A-stabilized-intertwining}
for the operator $\mathbb D_{\rho,e}$ and the identities
\eqref{eq:A-domain-retract-identities} for the transfer maps read
\begin{equation}\label{appref:low-stabilized-descent}
 \begin{aligned}
 \operatorname {aux}_{Y,\rho,e}\mathbb D_{\rho,e}U
   &=\mathsf C_{X,\rho,e}U,\\
 \mathbb D_{\rho,e}\mathsf J_{X,\rho,e}u
   &=\iota_{Y,\rho,e}\mathsf J_{Y,\rho,e}^{\rm raw}
       L_{\rho,e}u,\\
 \mathsf R_{X,\rho,e}\mathsf J_{X,\rho,e}&=I,\qquad
 \mathsf C_{X,\rho,e}\mathsf J_{X,\rho,e}=0,\\
 \mathsf J_{X,\rho,e}\mathsf R_{X,\rho,e}U&=U
       \quad\text{when }\mathsf C_{X,\rho,e}U=0.
 \end{aligned}
\end{equation}
Here, the subscript $\rho,e$ means that the transfer maps of
\cref{sec:block} and the linearization are evaluated at
$\mathbf b=\mathbf b_\rho$. Thus, $L_{\rho,e}=L_{\mathbf b_\rho}$ and
$\mathsf J_{X,\rho,e}=\mathsf J_{X,\mathbf b_\rho}$, and similarly for
$\mathsf R_{X,\rho,e}$, $\mathsf C_{X,\rho,e}$, $\iota_{Y,\rho,e}$,
$\mathsf J_{Y,\rho,e}^{\rm raw}$, and $\operatorname{aux}_{Y,\rho,e}$.

\begin{lemma}
\label{lem:fixed-J-native}
For every $J<\infty$, there are finite constants $K_{0,J}$ and
$C_J$, a number $\rho_{L,J}>0$, and, for every $s\geq s_*$ and every
integer $q\geq0$, a finite constant $K_{J,s,q}$, such that the
following hold.
\begin{enumerate}[label=\textup{(\roman*)},leftmargin=*,itemsep=2pt]
\item The low block $A_{0,J}$ at the circle has a two-sided
inverse $G_{0,J}$, and
\begin{equation}\label{appref:circular-native-bound}
 \norm{G_{0,J}F}_{\bX_{\natg}^{s_*}}
 \leq K_{0,J}\norm{F}_{\bY_{\natg}^{s_*}}.
\end{equation}

\item In the block norms at the base index,
\begin{equation}\label{appref:rho-operator-bound}
 \norm{A_{\rho,J}-A_{0,J}}_{
     P_J^X\bX_{\natg}^{s_*}\to P_J^Y\bY_{\natg}^{s_*}}
 \leq C_J|\rho|,
 \qquad 0\leq\rho\leq\rho_{L,J}.
\end{equation}

\item After requiring $K_{0,J}C_J\rho_{L,J}\leq1/2$, the formula
\begin{align}
 G_{\rho,J}
 &=\bigl[I+G_{0,J}(A_{\rho,J}-A_{0,J})\bigr]^{-1}G_{0,J}
                                                        \label{appref:fixed-J-inverse}\\
 &=G_{0,J}\bigl[I+(A_{\rho,J}-A_{0,J})G_{0,J}\bigr]^{-1}\notag
\end{align}
defines, for every $0\leq\rho\leq\rho_{L,J}$, a two-sided inverse of
$A_{\rho,J}$, with
\begin{equation}\label{appref:fixed-J-two-sided}
 \norm{G_{\rho,J}}_{s_*}\leq2K_{0,J},\qquad
 A_{\rho,J}G_{\rho,J}=I,\qquad
 G_{\rho,J}A_{\rho,J}=I.
\end{equation}
Here, $\norm{\cdot}_{s_*}$ denotes the operator norm between the block
spaces $P_J^X\bX_{\natg}^{s_*}$ and $P_J^Y\bY_{\natg}^{s_*}$ at the
base index, in the direction in which the operator acts. Moreover, for
every $s\geq s_*$ and every integer $q\geq0$, we have
\begin{equation}\label{appref:fixed-J-derivatives}
 \sup_{0\leq\rho\leq\rho_{L,J}}
 \norm{\pa_\rho^qG_{\rho,J}F}_{\bX_{\natg}^s}
 \leq K_{J,s,q}\norm{F}_{\bY_{\natg}^s},
\end{equation}
and the family also commutes with complex conjugation.
\end{enumerate}
\end{lemma}

Recall the family \eqref{eq:fixed-M-lambda}. In
\eqref{eq:fixed-alpha-lambda}, $\alpha_0$ is the orientation and
$\delta$ the amplitude of the tilt. We fix $\alpha_0$ before
constructing the reference inverse, and the constant ellipses of the
lemma have $\delta=0$. We choose the constants $K_{0,J}$ and $C_J$ on
this family of constant ellipses. The family \eqref{eq:fixed-M-lambda}
with $\delta\neq0$ enters through the perturbation estimate of
\cref{app:global:low-leakage}, and $\delta$ is chosen in the resulting
neighborhood.

\begin{proof}
At the circle, every component of $\mathbb D_{0,e}$ preserves the cell
frequency and angular mode. We fix a cell frequency $n$ and write a
block source as
\[
 F_n^\natg=(F_{c,n}^\natg,F_{A,n}^\natg,F_{+,n}^\natg,
             F_{{\rm cap},n}^\natg,F_{{\rm fin},n}).
\]
The entries are, in order, the center source, the two affine
sources, the positive source, the cap data, and the remaining mean,
scale, gauge, chart, and quotient data. Their weights are listed in
\cref{tab:block-weights}; each source entry includes its prescribed
axis or boundary data. We solve the corresponding equations as follows.
\begin{center}
\small
\begin{tabular}{@{}p{0.14\textwidth}p{0.32\textwidth}p{0.46\textwidth}@{}}
\toprule
Component & Data & Solution procedure\\
\midrule
Center & Interior source and prescribed first Taylor coefficients. & Bessel formulas
\eqref{app:minus:eq:center-bessel}--\eqref{app:minus:eq:bessel-volterra} for nonzero cell frequency.\\
Affine & The two Hardy forcings and the prescribed affine coefficient. & Hardy maps \eqref{app:minus:eq:hardy-arrows} on $|m|=2$.\\
Positive & Source in the dual of the form space, interface data, and Robin data. & Coercivity \eqref{app:positive:eq:circle-coercivity} and \eqref{app:positive:eq:weighted-coercivity}, with the even-mode bound \eqref{appref:even-native-bound}.\\
Cap & Cartesian interior and trace data. & The cap constructions in \cref{app:minus,app:positive}, with the prescribed axis conditions.\\
Remaining & Mean, scale, gauge, chart, and compatible quotient data. & A triangular system with nonzero diagonal entries.\\
\bottomrule
\end{tabular}
\end{center}
The center is kept as one four-dimensional block because its two speeds
coincide at the circle. At zero cell frequency, we use
\eqref{appref:zero-carrier-tangent} for the tangential component and
the radial equations following it, together with the even reconstruction
\eqref{appref:even-final-affine-equations}--\eqref{appref:even-reconstruction}
with the mass terms omitted. The reconstruction
\eqref{eq:range-reconstruction} recovers the Cartesian source entries
from the quotient variables at the orders of the block norms.
% AUTHOR QUERY: The quoted coercivity bounds supply the energy estimate.
% Verify the full solution and trace estimate on the low-frequency block
% domains, including the cap constraints, before using the table to
% conclude the no-loss estimate below.

The block source is the pair $(Z,\Psi_{\mathbf b}Z)$ in
\eqref{app:cross:source-graph-chart}. Here, $Z$ contains the compatible
quotient data, and $\Psi_{\mathbf b}Z$ is the scalar differential source.
The block norm controls both entries at the same index. We convert a Cartesian
source to this pair once, by \eqref{eq:A-JY-type}, before the inverse at
the circle is applied.

Substituting the explicit formulas for the components into the
equations, we obtain the two identities $A_{0,J}G_{0,J}=I$ and
$G_{0,J}A_{0,J}=I$ at cell frequency $n$ on smooth data. The prescribed
Taylor coefficients on the axis, the coercive form, and the triangular
finite system give
uniqueness, so density extends the identities to the block spaces at
cell frequency $n$. Thus, for each $n$ and $s$,
\begin{equation}\label{appref:one-carrier-native}
 \norm{G_{0,n}F_n^\natg}_{\bX_{\natg,n}^s}
 \leq C_{n,s}\norm{F_n^\natg}_{\bY_{\natg,n}^s}.
\end{equation}
Here, $G_{0,n}$ is the inverse at cell frequency $n$ constructed above,
$\bX_{\natg,n}^s$ and $\bY_{\natg,n}^s$ are the subspaces of
$\bX_{\natg}^s$ and $\bY_{\natg}^s$ at cell frequency $n$, and the
constant $C_{n,s}$ does not depend on any cutoff in the angular modes.
We assemble these inverses by
\[
 G_{0,J}F=\sum_{|n|<2^J}(G_{0,n}F_n^{\natg})e^{in\zeta}.
\]
Parseval's identity and the maximum of $C_{n,s_*}$ over $|n|<2^J$
give $K_{0,J}<\infty$ and \eqref{appref:circular-native-bound}. The cutoff $J$ restricts only the cell frequency. The component
estimates control all radial and angular variables without further
truncation. Since the operator at the circle is diagonal in the cell
frequency, this construction is independent of the high-frequency
inverse. The off-diagonal blocks between the minus and positive
components also vanish there, by \eqref{app:cross:circle-zero}.

Before fixing $J$, we choose $\rho_0$ so that the bounds for the
determinant, the contours, the form, the cap, the quotients, and the
traces hold for $0\leq\rho\leq\rho_0$. The Cartesian inverse of
\cref{app:reference} is used only for $0<\rho<\rho_0$. Thus, the
choice of $\rho_0$ does not depend on $J$. We claim that, for every $J$, $s$, and
$q$,
\begin{equation}\label{appref:fixed-J-operator-derivatives}
 \sup_{0\leq\rho\leq\rho_0}
 \norm{\pa_\rho^qA_{\rho,J}}_{
  P_J^X\bX_{\natg}^s\to P_J^Y\bY_{\natg}^s}<\infty,
\end{equation}
and we verify this component by component. For the center block, the graph norm on
the maximal domain controls the radial derivative, and all cell
derivatives are bounded on $|n|<2^J$. The spectral projection onto the
four center modes is given by a contour integral, and the resolvent
identity bounds it and all its $\rho$-derivatives as operators of order
zero. The
coefficients of the affine Hardy maps are analytic and carry the weights
of \eqref{appref:even-native-bound}. The positive form and its
derivatives map the form space to its dual, while the value, conormal,
and Robin data use the balanced trace norms. On the cap, the minus operator maps $H^{s+2}$ to its source space
$H^s$, whereas the positive operator maps $H^{s+1}$ to $H^{s-1}$,
with the value and the split residual of \cref{app:positive} in $H^{s+1/2}$
and $H^{s-1/2}$. Both cap estimates include the Cartesian compatibility and axis conditions. The remaining quotient,
source, mean, scale, chart, and gauge operators are analytic maps of fixed
order. These bounds control every component of $A_{\rho,J}$, and the estimates for
the angular multipliers make them uniform in the angular mode.

With $C_J=\sup_{0\leq\rho\leq\rho_0}\norm{\pa_\rho A_{\rho,J}}_{s_*}$, we
obtain \eqref{appref:rho-operator-bound} from the fundamental theorem of
calculus in the Banach space of bounded operators. We may then define
\[
 \rho_{L,J}:=\min\left\{\rho_0,
   \frac1{2\max\{1,K_{0,J}C_J\}}\right\},
\]
so that both Neumann series in \eqref{appref:fixed-J-inverse} converge.
Since
$G_{0,J}\bigl[(A_{\rho,J}-A_{0,J})G_{0,J}\bigr]^k
=\bigl[G_{0,J}(A_{\rho,J}-A_{0,J})\bigr]^kG_{0,J}$ for every $k\geq0$,
the two series define the same operator $G_{\rho,J}$, and composing
with $A_{\rho,J}=A_{0,J}+(A_{\rho,J}-A_{0,J})$ on either side proves
\eqref{appref:fixed-J-two-sided}.

To control higher regularity, we use the estimate below, with
$E_{J,s}$ independent of $0\leq\rho\leq\rho_{L,J}$.
\begin{equation}\label{appref:native-regularity}
 \norm{u}_{\bX_{\natg}^s}
 \leq E_{J,s}\bigl(
   \norm{A_{\rho,J}u}_{\bY_{\natg}^s}
       +\norm{u}_{\bX_{\natg}^{s_*}}\bigr).
\end{equation}
% AUTHOR QUERY: Supply this regularity estimate on the completed block
% spaces, with constants uniform on the stated rho interval. The
% off-diagonal tame estimate alone does not show that all transverse
% commutators lower the block index. Justify the passage from smooth
% approximations to the base-index inverse before applying the estimate.
The components are coupled in $A_{\rho,J}$, so we commute the equation
$A_{\rho,J}u=F$ with the derivatives that enter the block norms and
estimate the commutators. The identity
\eqref{app:cross:principal-cross-zero} cancels the off-diagonal terms of
highest order, and by \cref{prop:macro-cross} a derivative falling on a
frame, on a projection, on a reduction of the source, or on a
coefficient lowers the block index by one, so each commutator is
bounded at the index below. Induction from $s_*$ and interpolation give
\eqref{appref:native-regularity}. Applying it to $u=G_{\rho,J}F$ and
using \eqref{appref:fixed-J-two-sided} for the base-index term proves
\eqref{appref:fixed-J-derivatives} for $q=0$.

For $q\ge1$, we fix the datum $F$, set $u=G_{\rho,J}F$, and
differentiate $A_{\rho,J}u=F$ in $\rho$ to obtain
\begin{equation*}
 A_{\rho,J}\pa_\rho^qu
 =-\sum_{k=1}^q\binom qk
    (\pa_\rho^kA_{\rho,J})\pa_\rho^{q-k}u
\end{equation*}
whose right-hand side is bounded in $\bY_{\natg}^s$ by
\eqref{appref:fixed-J-operator-derivatives} and by the estimates for
$\pa_\rho^{q-k}u$ with $k\geq1$. Applying $G_{\rho,J}$ and the case
$q=0$ of \eqref{appref:fixed-J-derivatives} to this right-hand side, we
obtain the estimate for each $q$ by induction on $q$. Finally, complex conjugation preserves all the kernels, contours,
spaces, and finite-dimensional operators. Thus,
$\overline{G_{\rho,J}F}$ and $G_{\rho,J}\overline F$ solve the same
equation and are equal by uniqueness.
\end{proof}

The block inverse agrees with the Cartesian inverse of
\cref{app:reference} where both are defined, as we now verify.

\begin{remark}
\label{appref:relation-of-references}
The two inverses agree on compatible Cartesian data, after the return
$\mathsf R_{X,\rho,e}$ to Cartesian variables. We fix $0<\rho\leq\rho_{L,J}$ with $\rho<\rho_0$ and let
$V_{\rho,e}$ be the inverse $V_e$ of
\cref{appref:raw-constant-ellipse} with $\rho_*$ replaced by $\rho$.
These restrictions ensure that both inverses are defined. Since the coefficients of $\mathbb D_{\rho,e}$ are
independent of $\zeta$, the operator is diagonal in the cell frequency,
\[
 \mathbb D_{\rho,e}P_J^X=P_J^Y\mathbb D_{\rho,e},
 \qquad
 P_J^Y\mathbb D_{\rho,e}(I-P_J^X)=0,
 \qquad
 (I-P_J^Y)\mathbb D_{\rho,e}P_J^X=0,
\]
so \eqref{appref:low-compression} is the restriction of
$\mathbb D_{\rho,e}$ to the low cell frequencies. We let $f$ be a smooth
Cartesian datum at low frequency and put
\[
 F=\iota_{Y,\rho,e}\mathsf J_{Y,\rho,e}^{\rm raw}f,
 \qquad P_J^YF=F.
\]
If $U=G_{\rho,J}F$, then $U=P_J^XU$ and $A_{\rho,J}U=F$, so
$\mathbb D_{\rho,e}U=P_J^Y\mathbb D_{\rho,e}U=F$ by the diagonality
above. Since $F$ has vanishing auxiliary data, the first and the last
identities in \eqref{appref:low-stabilized-descent} give
\[
 \mathsf C_{X,\rho,e}U
 =\operatorname {aux}_{Y,\rho,e}F=0,\qquad
 U=\mathsf J_{X,\rho,e}\mathsf R_{X,\rho,e}U.
\]
Applying the second identity in \eqref{appref:low-stabilized-descent}
to $u=\mathsf R_{X,\rho,e}U$, we obtain
$\iota_{Y,\rho,e}\mathsf J_{Y,\rho,e}^{\rm raw}L_{\rho,e}
\mathsf R_{X,\rho,e}U=\mathbb D_{\rho,e}U=F
=\iota_{Y,\rho,e}\mathsf J_{Y,\rho,e}^{\rm raw}f$. Since
$\mathsf J_{Y,\rho,e}^{\rm raw}$ has the left inverse
\eqref{eq:A-JY-inverse-factorization} on smooth elements and
$p_{Y,\mathbf b}\iota_{Y,\mathbf b}=I$ by \eqref{eq:A-zero-range-splitting},
the map $\iota_{Y,\rho,e}\mathsf J_{Y,\rho,e}^{\rm raw}$ is injective
on smooth elements. Hence, $L_{\rho,e}\mathsf R_{X,\rho,e}U=f$, and
uniqueness in \cref{appref:raw-constant-ellipse} gives
$\mathsf R_{X,\rho,e}U=V_{\rho,e}f$. We therefore obtain
\begin{equation}\label{appref:low-reference-descent}
 G_{\rho,J}\,
   \iota_{Y,\rho,e}\mathsf J_{Y,\rho,e}^{\rm raw}f
 =\mathsf J_{X,\rho,e}V_{\rho,e}f,\qquad
 \mathsf R_{X,\rho,e}G_{\rho,J}\,
   \iota_{Y,\rho,e}\mathsf J_{Y,\rho,e}^{\rm raw}f
 =V_{\rho,e}f .
\end{equation}
The Cartesian estimate \eqref{appref:raw-estimate} loses $\ell_e$
derivatives of the source and requires $\rho>0$. The block estimate in
\cref{lem:fixed-J-native} was proved directly and has neither
restriction. Both inverses here concern a constant ellipse. The family
\eqref{eq:fixed-M-lambda} agrees with this reference only at $\delta=0$.
\end{remark}
\section{Combining the inverses}
\label{app:global}

We now combine the preceding inverses to invert the augmented operator
$\mathbb D_{\mathbf b}$. We first couple the two high-frequency blocks,
then include the low frequencies by a Schur complement. We estimate
state and parameter derivatives in the block spaces before returning
to Cartesian variables.

Throughout, we apply the projections $P$ and $H$ of
\eqref{app:global:block-decomposition} in the block variables, that is,
only after the transfer maps of \cref{sec:block} have been applied to
the domain and to the source. Thus, every block acts between
fixed spaces. We recall from \eqref{app:cross:circle-distance} the distance from the
circle,
\begin{equation}\label{app:global:circle-distance}
 \begin{aligned}
 \eta_{\rm circ}(\mathbf b)
 &:=\norm{a-a_\circ(\rho)}_{\cA_\gamma^{s_*+k_\times}}+|\eps|+\norm{M(p)-M_\circ(\rho)}_{\cA_\gamma^{s_*+k_\times}},\\
 \chi(\mathbf b)&:=|\rho|+\eta_{\rm circ}(\mathbf b).
 \end{aligned}
\end{equation}
Here, $a_\circ(\rho)$ and $M_\circ(\rho)$ are the unperturbed solution
and the matrix of the constant ellipse $\mathbf b_\rho$, whose axes do not
rotate with $\zeta$, and $k_\times$ is the shift in the regularity of
the coefficients fixed with the estimate \eqref{app:cross:main-estimate}.
We use $\chi$ in the estimates at high frequency. For the perturbation at low frequency below, we instead
measure the distance from the constant reference state $\mathbf b_e$.

\subsection{The high-frequency inverse}

We first solve $D_{H,\mathbf b}(u_-,u_+)=(f_-,f_+)$. The unknown
$u_-$ contains the center and affine fields, and $u_+$ contains the
positive fields; each includes its cap and annular components. The
operator of \eqref{app:global:block-decomposition} has the block form
\begin{equation}\label{app:global:high-macro-matrix}
 D_{H,\mathbf b}
 =\begin{pmatrix}
    D_{-,\mathbf b}&E_{-+,\mathbf b}\\
    E_{+-,\mathbf b}&D_{+,\mathbf b}
   \end{pmatrix}.
\end{equation}
All four entries belong to the same augmented operator
$\mathbb D_{\mathbf b}$. The off-diagonal blocks are bounded from
index $s$ to index $s$ by \cref{prop:macro-cross}. Under
\eqref{app:cross:diagonal-correction-margin}, the correction
\eqref{app:cross:true-diagonal-correction}--\eqref{app:cross:true-diagonal-inverse}
turns the inverses of
\cref{prop:minus-diagonal,prop:positive-diagonal} into the diagonal
inverses used below.
\begin{lemma}
\label{app:global:high-return}
Suppose that
\[
 \norm{D_{-,\mathbf b}^{-1}}\leq K_-,\qquad
 \norm{D_{+,\mathbf b}^{-1}}\leq K_+,
 \qquad
 \norm{E_{-+,\mathbf b}}+\norm{E_{+-,\mathbf b}}
 \leq K_\times\chi(\mathbf b)
\]
in the block norms at the base index, and let
\begin{equation}\label{app:global:return-operator}
 \mathcal R_{-,\mathbf b}
 :=D_{-,\mathbf b}^{-1}E_{-+,\mathbf b}
      D_{+,\mathbf b}^{-1}E_{+-,\mathbf b}.
\end{equation}
Solving the second row for $u_+$ and substituting it into the first
gives
\[
 (I-\mathcal R_{-,\mathbf b})u_-
 =D_{-,\mathbf b}^{-1}
    \bigl(f_--E_{-+,\mathbf b}D_{+,\mathbf b}^{-1}f_+\bigr).
\]
If $\norm{\mathcal R_{-,\mathbf b}}<1$, the operator on the left is
invertible by a Neumann series. Thus, $D_{H,\mathbf b}$ has a two-sided
inverse $G_{H,\mathbf b}$, given for $f=(f_-,f_+)$ by
\begin{align}
 u_-&=(I-\mathcal R_{-,\mathbf b})^{-1}
 D_{-,\mathbf b}^{-1}
 \bigl(f_--E_{-+,\mathbf b}D_{+,\mathbf b}^{-1}f_+\bigr),
                                                        \label{app:global:minus-solve}\\
 u_+&=D_{+,\mathbf b}^{-1}
 \bigl(f_+-E_{+-,\mathbf b}u_-\bigr),                  \label{app:global:plus-solve}\\
 G_{H,\mathbf b}f&=(u_-,u_+).                           \label{app:global:high-inverse}
\end{align}
On the block spaces, this inverse obeys tame estimates without loss of
derivatives, and each order of differentiation with respect to $\mathbf b$
requires only a finite increase in the regularity of the source.
\end{lemma}

\begin{proof}
As in \eqref{app:cross:return-types}, the factors in
\eqref{app:global:return-operator} act between the spaces
\[
 \bX_{-,\natg}^s\xrightarrow{E_{+-}}\bY_{+,\natg}^s
 \xrightarrow{D_+^{-1}}\bX_{+,\natg}^s
 \xrightarrow{E_{-+}}\bY_{-,\natg}^s
 \xrightarrow{D_-^{-1}}\bX_{-,\natg}^s,
\]
so each term of the Neumann series acts on the same block space and
preserves the index. No Cartesian transfer map enters this composition. The smallness assumption gives convergence of the Neumann series at
the base index. We establish its higher regularity bounds below. Substituting
\eqref{app:global:minus-solve} and \eqref{app:global:plus-solve} into
\eqref{app:global:high-macro-matrix}, we obtain the right inverse
identity. For the left inverse identity, we fix $w=(w_-,w_+)$ and
subtract $G_{H,\mathbf b}D_{H,\mathbf b}w$ from $w$. Since
$D_{H,\mathbf b}G_{H,\mathbf b}=I$, the difference $v=(v_-,v_+)$ solves
the homogeneous system $D_{H,\mathbf b}v=0$. Applying
$D_{+,\mathbf b}^{-1}$ to the second row gives
$v_+=-D_{+,\mathbf b}^{-1}E_{+-,\mathbf b}v_-$, and inserting this into
the first row and applying $D_{-,\mathbf b}^{-1}$ gives
$(I-\mathcal R_{-,\mathbf b})v_-=0$, so $v_-=0$ and then $v_+=0$. For higher regularity, we commute derivatives through the equations and
use the strict smallness condition at the base index to absorb the
highest-order term. Only the source or one coefficient is measured in
a higher regularity norm.
% AUTHOR QUERY: Give the commuted high-block estimate that
% permits absorption at every Sobolev index with the same base-norm
% smallness. The stated tame factor estimates do not alone give this
% conclusion.
The derivatives with respect to $\mathbf b$ are estimated below.
\end{proof}

\subsection{The low-frequency block}

We write $A_e=A_{\mathbf b_e}$ for the low block at the constant
reference state, as in \cref{app:macro-cross}. By
\cref{lem:fixed-J-native} at $\rho=\rho_*$, which applies since
$\rho_*<\rho_{L,J}$ by \eqref{sec:inverse:rho-star-choice} below,
\begin{equation}\label{app:global:reference-low-inverse}
 A_e^{-1}=G_{\rho_*,J}:
 P^Y\bY_{\natg}^{s_*}\longrightarrow P^X\bX_{\natg}^{s_*}
\end{equation}
is a two-sided inverse and loses no derivatives. Under the condition
\eqref{app:global:low-neumann-margin}, we therefore have the two-sided
inverse of the low block
\begin{equation}\label{app:global:current-low-inverse}
 A_{\mathbf b}^{-1}
 =\bigl[I+A_e^{-1}(A_{\mathbf b}-A_e)\bigr]^{-1}A_e^{-1}.
\end{equation}
For $F_L$ in the range of $P^Y$, we write
$u_L=A_{\mathbf b}^{-1}F_L$. Commuting derivatives through
$A_{\mathbf b}u_L=F_L$ and applying \eqref{app:global:low-tame} gives
the tame bound at each Sobolev index without loss of derivatives.
% AUTHOR QUERY: The low-frequency space remains
% infinite-dimensional in the transverse variables. A current low-block
% regularity estimate is needed in addition to the base Neumann bound and
% low-tame.

\subsection{The Schur complement at low frequencies}

For $U=u_L+u_H$ and $F=F_L+F_H$, where the subscripts denote
the low and high projections, the equation $\mathbb D_{\mathbf b}U=F$
is
\[
 A_{\mathbf b}u_L+B_{\mathbf b}u_H=F_L,\qquad
 C_{\mathbf b}u_L+D_{H,\mathbf b}u_H=F_H.
\]
The second equation gives
$u_H=G_{H,\mathbf b}(F_H-C_{\mathbf b}u_L)$.
Substituting this into the first eliminates the high-frequency unknown. The operator that
remains on the low-frequency component is the \emph{Schur complement}
\begin{equation}\label{app:global:low-schur}
 S_{L,\mathbf b}
 :=A_{\mathbf b}-B_{\mathbf b}G_{H,\mathbf b}C_{\mathbf b}
 =A_{\mathbf b}(I-\mathcal T_{L,\mathbf b}),
 \qquad
 \mathcal T_{L,\mathbf b}
 :=A_{\mathbf b}^{-1}B_{\mathbf b}
       G_{H,\mathbf b}C_{\mathbf b}.
\end{equation}
\begin{proposition}[The global inverse]
\label{app:global:exact-schur}
Let $G_{H,\mathbf b}$ be the inverse of \cref{app:global:high-return}.
Assume \eqref{app:global:low-neumann-margin} and
$\norm{\mathcal T_{L,\mathbf b}}\leq q_L<1$ in the block norms at the base
index. Then
\begin{equation}\label{app:global:schur-inverse}
 S_{L,\mathbf b}^{-1}
 =(I-\mathcal T_{L,\mathbf b})^{-1}A_{\mathbf b}^{-1},
\end{equation}
and
\begin{equation}\label{app:global:full-block-inverse}
 G_{\mathbf b}=
 \begin{pmatrix}
 S_{L,\mathbf b}^{-1}
   &-S_{L,\mathbf b}^{-1}B_{\mathbf b}G_{H,\mathbf b}\\
 -G_{H,\mathbf b}C_{\mathbf b}S_{L,\mathbf b}^{-1}
   &G_{H,\mathbf b}
    +G_{H,\mathbf b}C_{\mathbf b}S_{L,\mathbf b}^{-1}
       B_{\mathbf b}G_{H,\mathbf b}
 \end{pmatrix}
\end{equation}
maps $(F_L,F_H)$ to $(u_L,u_H)$. Equivalently, we solve in the order
\[
 u_L=S_{L,\mathbf b}^{-1}(F_L-B_{\mathbf b}G_{H,\mathbf b}F_H),
 \qquad u_H=G_{H,\mathbf b}(F_H-C_{\mathbf b}u_L).
\]
It is a two-sided inverse in the block variables, so that
\[
 \mathbb D_{\mathbf b}G_{\mathbf b}=I,
 \qquad
 G_{\mathbf b}\mathbb D_{\mathbf b}=I.
\]
Uniformly on every ball of states on which the norms in the
smallness conditions at the base index are bounded away from one, we have
\begin{equation}\label{app:global:native-global-estimate}
 \norm{G_{\mathbf b}F}_{\bX_{\natg}^s}
 \leq C_s\left[
 \norm{F}_{\bY_{\natg}^s}
 +(1+\norm{\mathbf b-\mathbf b_e}_{\cB^{s+k}})
      \norm{F}_{\bY_{\natg}^{s_*}}\right].
\end{equation}
Here, $k$ is a fixed shift large enough for the diagonal inverse,
off-diagonal, and low-frequency estimates used in the construction.
\end{proposition}

\begin{proof}
The strict smallness assumption makes the Neumann series in
\eqref{app:global:schur-inverse} converge on the low-frequency space at
the base index. Multiplying
\eqref{app:global:full-block-inverse} and
\eqref{app:global:block-decomposition} in both orders verifies the two
inverse identities directly on the full block spaces.

For the regularity at index $s$, we write
$w=(I-\mathcal T_{L,\mathbf b})^{-1}g$ for $g$ in the low-frequency
space and use the equation $w=g+\mathcal T_{L,\mathbf b}w$. The tame
estimates for $A_{\mathbf b}^{-1}$, $B_{\mathbf b}$, $G_{H,\mathbf b}$,
and $C_{\mathbf b}$ give
\[
 \norm{w}_s\leq\norm{g}_s+q_L\norm{w}_s
 +C_{J,s}(1+\norm{\mathbf b-\mathbf b_e}_{\cB^{s+k}})
      \norm{w}_{s_*},
 \qquad q_L<1.
\]
Here, $\norm{\cdot}_s$ denotes the norm of $P^X\bX_{\natg}^s$. The
inverse at the base index controls $\norm{w}_{s_*}$, and the
strict smallness condition absorbs $q_L\norm{w}_s$.
% AUTHOR QUERY: The coefficient q_L is known to be below
% one only in the base norm. Prove a commutator estimate giving this
% coefficient at index s; the stated tame bounds have an index-dependent
% leading constant. This step is needed for native-global-estimate and the
% smooth tame inverse.
Inserting the result into the four blocks of
\eqref{app:global:full-block-inverse}, we obtain
\eqref{app:global:native-global-estimate}.
\end{proof}

\subsection{Derivatives of the inverse}

Differentiating the center equation
\eqref{app:minus:eq:center-evolution} with respect to the state and
parameters produces the source in
\eqref{app:minus:eq:differentiated-center}, with at most one cell
derivative of the solution. Accordingly, we estimate derivatives of the inverse between
fixed Sobolev spaces with a finite loss that depends on the number of
differentiations.

\begin{lemma}
\label{app:global:all-order}
For each $m\geq1$, the families of blocks in
\eqref{app:global:block-decomposition} and \eqref{app:global:high-macro-matrix},
together with the inverses $D_{\pm,\mathbf b}^{-1}$ and
$A_{\mathbf b}^{-1}$ and the Neumann series
$(I-\mathcal R_{-,\mathbf b})^{-1}$ and
$(I-\mathcal T_{L,\mathbf b})^{-1}$, are $m$ times
Fr\'echet differentiable between block spaces whose indices differ by a
finite shift. There
are finite integers $k_m$ and $\mu_m$, independent of the Sobolev index
$s$, such that
\begin{align}
 &\norm{D_{\mathbf b}^mG_{\mathbf b}
        [h_1,\ldots,h_m]F}_{\bX_{\natg}^s}\notag\\
 &\quad\leq C_{s,m}\Bigg[
 \left(\prod_{i=1}^m\norm{h_i}_{{\rm tan},s_*+k_m}\right)
       \norm{F}_{\bY_{\natg}^{s+\mu_m}}\notag\\
 &\qquad+
 \sum_{i=1}^m\norm{h_i}_{{\rm tan},s+k_m}
       \prod_{j\ne i}\norm{h_j}_{{\rm tan},s_*+k_m}
       \norm{F}_{\bY_{\natg}^{s_*+\mu_m}}\notag\\
 &\qquad+
 (1+\norm{\mathbf b-\mathbf b_e}_{\cB^{s+k_m}})
 \left(\prod_{i=1}^m\norm{h_i}_{{\rm tan},s_*+k_m}\right)
       \norm{F}_{\bY_{\natg}^{s_*+\mu_m}}
 \Bigg].                                               \label{app:global:all-order-estimate}
\end{align}
At most $m$ derivatives of the source are lost in the inverse of the
center equation, and the remaining components and traces contribute only
a fixed finite loss.
\end{lemma}

\begin{proof}
For any of the families of operators $T_{\mathbf b}$ inverted above, we
have on smooth elements
\begin{align}
 D(T_{\mathbf b}^{-1})[h]
 &=-T_{\mathbf b}^{-1}(DT_{\mathbf b}[h])T_{\mathbf b}^{-1},
                                                        \label{app:global:inverse-derivative}\\
 D(I-T_{\mathbf b})^{-1}[h]
 &=(I-T_{\mathbf b})^{-1}(DT_{\mathbf b}[h])
       (I-T_{\mathbf b})^{-1}.                          \label{app:global:resolvent-derivative}
\end{align}
To justify these derivatives, we use the resolvent identity
\[
 T_{\mathbf b+h}^{-1}-T_{\mathbf b}^{-1}
 =-T_{\mathbf b+h}^{-1}
   (T_{\mathbf b+h}-T_{\mathbf b})T_{\mathbf b}^{-1},
\]
and we subtract the linear terms in
\eqref{app:global:inverse-derivative} and \eqref{app:global:resolvent-derivative}.
The uniform higher-index bounds for the inverses, together with the
remainder estimate for the Fr\'echet derivative of $T_{\mathbf b}$,
which costs a fixed finite number of derivatives, control the error at
an index lower by that fixed amount. By
induction, we obtain differentiability at every finite order.

For the oscillatory center component, we differentiate the center
equation $L_{\mathbf b}^cu=f$ in the directions $h_1,\ldots,h_q$, as in
\eqref{app:minus:eq:differentiated-center}, and, at order $q$, Leibniz'
rule gives
\begin{equation}\label{app:global:center-partition}
 L_{\mathbf b}^cD^qu
 =-\sum_{\varnothing\ne I\subset\{1,\ldots,q\}}
   D^{|I|}L_{\mathbf b}^c[h_I]
   D^{q-|I|}u[h_{I^c}].
\end{equation}
Here, the source $f$ is held fixed, $h_I$ collects the directions $h_i$
with $i\in I$, and $I^c$ is the complement of $I$ in $\{1,\ldots,q\}$.
Differentiation in the parameters changes the coefficients of
$L_{\mathbf b}^c$ but does not add a second factor $D_\zeta$. Every
differentiated center operator therefore contains at most one cell
derivative of the solution. Repeated use of the energy estimate for the
undifferentiated center equation, as in the proof of
\cref{prop:minus-diagonal}, requires one additional index on the source
at each
differentiation, and hence at most $q$ at order $q$. The change of variables
\eqref{app:minus:eq:sylvester-transform} for the affine component obeys
the same bound. For the positive form, one derivative falls on each
argument. For the cap, we use the Sobolev scales of
\cref{tab:block-weights}, from $H^{s+2}$ to $H^s$ for the minus block
and from $H^{s+1}$ to $H^{s-1}$ for the positive block, while the
quotient and trace maps retain their orders. Finally,
\cref{prop:macro-cross} bounds the differentiated off-diagonal maps
without loss of derivatives.

Differentiating
\eqref{app:global:high-inverse}, \eqref{app:global:current-low-inverse}, and \eqref{app:global:schur-inverse}
repeatedly, we obtain finite sums of products in which inverse factors
alternate with derivatives of operators. The strict smallness conditions
bound the inverse factors. The tame product estimate \eqref{eq:analytic-product} measures one of
$F,h_1,\ldots,h_m,\mathbf b-\mathbf b_e$ in a higher regularity norm
and the others at the base index. This gives \eqref{app:global:all-order-estimate}.
\end{proof}

\subsection{Return to Cartesian variables and real data}

We recover the Cartesian inverse by converting the source to block
variables, applying the block inverse, and returning the solution to
Cartesian variables. The source conversion occurs only once, including
in each differentiated term.

\begin{lemma}
\label{app:global:one-time-loss}
Let $\ell=\ell_Y$ be the number of derivatives lost when converting the
Cartesian source to block variables, as defined in \eqref{eq:A-source-order}. Then
\begin{equation}\label{app:global:raw-factorization}
 V_{\mathbf b}
 =\mathsf R_{X,\mathbf b}G_{\mathbf b}
       \iota_{Y,\mathbf b}\mathsf J_{Y,\mathbf b}^{\rm raw}
\end{equation}
satisfies both inverse identities in the Cartesian variables,
\[
 L_{\mathbf b}V_{\mathbf b}=I,
 \qquad
 V_{\mathbf b}L_{\mathbf b}=I,
\]
and, in every $m$-th derivative of $V_{\mathbf b}$ with respect to $\mathbf b$, the total
loss of derivatives on the Cartesian source is $\ell+\mu_m$. The transfer map
$\mathsf J_{Y,\mathbf b}^{\rm raw}$ of the
source occurs once in each differentiated term, so its loss does not
accumulate to $(m+1)\ell+\mu_m$.
\end{lemma}

\begin{proof}
We fix a smooth compatible Cartesian datum $f$, and we set
\begin{equation}\label{app:global:zero-source-solution}
 \widehat f=\iota_{Y,\mathbf b}
       \mathsf J_{Y,\mathbf b}^{\rm raw}f,
 \qquad W=G_{\mathbf b}\widehat f.
\end{equation}
The inclusion $\iota_{Y,\mathbf b}$ assigns zero to the auxiliary
entries, so $\operatorname{aux}_{Y,\mathbf b}\widehat f=0$. From
the identity \eqref{eq:A-auxiliary-row-identity} for the auxiliary equations
and $\mathbb D_{\mathbf b}G_{\mathbf b}=I$, we therefore obtain
\begin{equation}\label{app:global:zero-auxiliary-output}
 \mathsf C_{X,\mathbf b}W
 =\operatorname {aux}_{Y,\mathbf b}\mathbb D_{\mathbf b}W
 =\operatorname {aux}_{Y,\mathbf b}\widehat f=0,
\end{equation}
while \eqref{eq:A-domain-retract-identities} yields
\begin{equation}\label{app:global:zero-domain-reconstruction}
 \mathsf J_{X,\mathbf b}\mathsf R_{X,\mathbf b}W=W.
\end{equation}
Using \eqref{eq:A-projected-intertwining} and the identity
$p_{Y,\mathbf b}\iota_{Y,\mathbf b}=I$ of \eqref{eq:A-zero-range-splitting},
we obtain
\begin{align}
 \mathsf J_{Y,\mathbf b}^{\rm raw}L_{\mathbf b}V_{\mathbf b}f
 &=p_{Y,\mathbf b}\mathbb D_{\mathbf b}
       \mathsf J_{X,\mathbf b}\mathsf R_{X,\mathbf b}W\notag\\
 &=p_{Y,\mathbf b}\mathbb D_{\mathbf b}W
  =p_{Y,\mathbf b}\widehat f
  =\mathsf J_{Y,\mathbf b}^{\rm raw}f,                \label{app:global:raw-right-identity}
\end{align}
which, by injectivity of $\mathsf J_{Y,\mathbf b}^{\rm raw}$ on the
scale of compatible Cartesian sources, proves
$L_{\mathbf b}V_{\mathbf b}=I$. For the left
inverse identity, we take any smooth $u$ satisfying the constraints, and
then compute
\begin{align}
 V_{\mathbf b}L_{\mathbf b}u
 &=\mathsf R_{X,\mathbf b}G_{\mathbf b}
       \iota_{Y,\mathbf b}\mathsf J_{Y,\mathbf b}^{\rm raw}
          L_{\mathbf b}u\notag\\
 &=\mathsf R_{X,\mathbf b}G_{\mathbf b}
       \mathbb D_{\mathbf b}\mathsf J_{X,\mathbf b}u
  =\mathsf R_{X,\mathbf b}\mathsf J_{X,\mathbf b}u=u. \label{app:global:raw-left-identity}
\end{align}
Boundedness and the shifted identities
\eqref{eq:A-JY-left-shifted-identity}--\eqref{eq:A-JY-right-shifted-identity}
extend both inverse identities to the closures of smooth elements in
the indicated norms.

Differentiating \eqref{app:global:raw-factorization}, we find that every
term of Leibniz' rule carries exactly one factor
$D^q\mathsf J_{Y,\mathbf b}^{\rm raw}$ on the Cartesian source, and that
the derivatives of $\iota_{Y,\mathbf b}$ and of $\mathsf R_{X,\mathbf b}$
lose no derivatives. By \cref{prop:all-order-adapters}, differentiation
of the source map changes its analytic coefficients but leaves its
differential order $\ell$ unchanged. Combining these bounds with
\eqref{app:global:all-order-estimate}, we obtain the stated loss on the
source.
\end{proof}

We next show that the inverses preserve real data.

\begin{lemma}
\label{app:global:reality}
For real $\mathbf b$, the block and the Cartesian inverses preserve the
real subspaces, that is, the fixed subspaces of the complex conjugations
on the block spaces and on the Cartesian spaces.
\end{lemma}

\begin{proof}
We let $\mathfrak c_X$ and $\mathfrak c_Y$ denote the conjugations on
$\bX_{\natg}^s$ and $\bY_{\natg}^s$ induced by Cartesian complex
conjugation. For a scalar Fourier series, this conjugation acts by
$(\mathfrak c u)_{m,n}=\overline{u_{-m,-n}}$. In paired scalar or quotient
coordinates, it also exchanges the entries paired by complex
conjugation, such as $Z=A+iB$ and $W=A-iB$ in
\cref{app:reference}. The operators of the
components, the spectral projections \eqref{app:cross:domain-range-contours}
on fixed contours, the quotient maps, the traces, and the transfer maps
commute with $\mathfrak c_X$ and $\mathfrak c_Y$. So do the projections
$P$ and $H$, whose symbols are even in the cell frequency, and the blocks
$A_{\mathbf b}$, $B_{\mathbf b}$, $C_{\mathbf b}$, and
$D_{H,\mathbf b}$. If
$T\mathfrak c_X=\mathfrak c_YT$ and $T$ has a two-sided inverse, then
uniqueness implies $T^{-1}\mathfrak c_Y=\mathfrak c_XT^{-1}$, and this
relation is preserved under finite products and norm-convergent Neumann
series. Applying this successively to
\eqref{app:global:high-inverse}, \eqref{app:global:current-low-inverse}, \eqref{app:global:schur-inverse}, and \eqref{app:global:raw-factorization}
proves the claim.
\end{proof}
\section{Uniform estimates for the inverse}\label{sec:inverse}

We now choose the constants so that the inverses constructed in
\cref{app:lowref,app:global} apply on one neighborhood of the reference
state. We first fix the frequency cutoff $J$, then the ellipse parameter
$\rho_*$, and finally the neighborhood. The resulting Cartesian inverse
and its tame estimates are stated in \cref{thm:uniform-current-inverse}.

\subsection{The choice of the constants}
The low-frequency bounds may depend on $J$, so the order of these
choices matters. We record it together with the four smallness
conditions needed for the inverse.
\begin{proposition}[Choice of the constants]
\label{sec:inverse:acyclic-choice}
There are an integer $J<\infty$, an ellipse parameter $\rho_*>0$, and an open
neighborhood $\mathfrak U_{\rm inv}$ of $\mathbf b_e$, chosen in the
order below. We call the norm $\norm{\cdot}_{\cB^{s_*+k_0}}$ of
\eqref{app:global:low-distance} the \textit{base norm}. The closure of
$\mathfrak U_{\rm inv}$ in the base norm lies inside
$\{0<\rho<\rho_0\}$ and inside the image of the preliminary
neighborhood $\mathcal U_0$ of \cref{prop:fixed-tame-map} under the
normalized chart of \cref{subsec:normalized-slice}.
\begin{enumerate}[label=\textup{(\roman*)},leftmargin=*,itemsep=2pt]
\item We fix the analytic widths $\sigma_0>\sigma_->0$ of
\cref{app:minus}, the orientation $\alpha_0$ and the parameter
$\lambda_*$ of the constant reference state of \cref{sec:block}, and a
compact neighborhood in the finite-dimensional parameters. All
constants of the reference inverse of \cref{app:reference} may depend
on the orientation $\alpha_0$.
\item We choose the strength $\gamma$ of the analytic weight and then
the parameter $\delta_A$ of \eqref{app:cross:cell-weights}, so that the
center and affine damping conditions
\eqref{app:minus:eq:phase-choices} and the positive-form condition
\eqref{app:positive:eq:positive-choices} remain strict. These conditions
are uniform up to the circle.
\item We choose $\rho_0<1/4$ so that the bounds for the determinant,
the contours around the center roots, the positive form, the cap, the
quotients, and the traces hold for $0\leq\rho\leq\rho_0$. The Cartesian
inverse of \cref{appref:raw-constant-ellipse} is used only for
$0<\rho<\rho_0$, as explained in \cref{app:lowref}.
\item We choose a preliminary neighborhood in the state and curvature and
take the cell cutoff $J$ large enough for
\eqref{app:minus:eq:state-cutoff-choices}. All subsequent constants may
depend on $J$. We then invert the low block at the circle by
\cref{lem:fixed-J-native}, obtain its constants $K_{0,J}$ and $C_J$,
and put
\begin{equation}\label{sec:inverse:rho-low}
 \rho_{L,J}:=\min\left\{\rho_0,
   \frac{1}{2\max\{1,K_{0,J}C_J\}}\right\}>0.
\end{equation}
\item The constants for $\widehat G_\pm$ and for the off-diagonal
blocks are uniform for $0\leq\rho\leq\rho_0$, and the part at low
frequency of the correction $K_{\alpha,\mathbf b}$ of
\eqref{app:cross:true-diagonal-correction-block} vanishes at the constant
reference state, because the reference is diagonal in the cell
frequency. We may therefore choose
\begin{equation}\label{sec:inverse:rho-star-choice}
 0<\rho_*<\rho_{L,J}
\end{equation}
so that first the two norms
$\norm{\widehat G_{\alpha,\mathbf b_e}K_{\alpha,\mathbf b_e}}$,
$\alpha\in\{-,+\}$, of \eqref{app:cross:diagonal-correction-margin},
and then the norm of the operator $\mathcal R_{-,\mathbf b_e}$ of
\eqref{app:global:return-operator} at the reference state, are strictly
below one.
\item Only after $J$ and $\rho_*$ are fixed do we shrink the final
neighborhood $\mathfrak U_{\rm inv}$ in the state, in the unperturbed
solution, and in the curvature.
\end{enumerate}
The four conditions below control, respectively, the diagonal
correction, the coupling of the two high-frequency blocks, the
perturbation of the low block, and the coupling between low and high
frequencies. On the closure of $\mathfrak U_{\rm inv}$, there are
constants $q_{\rm corr}^*,q_H^*,q_A^*,q_L^*<1$ such that
\begin{align}
 \max_{\alpha\in\{-,+\}}
 \norm{\widehat G_{\alpha,\mathbf b}K_{\alpha,\mathbf b}}
   &\leq q_{\rm corr}^*,                              \label{sec:inverse:q-diagonal-correction}\\
 \norm{D_{-,\mathbf b}^{-1}E_{-+,\mathbf b}
       D_{+,\mathbf b}^{-1}E_{+-,\mathbf b}}
   &\leq q_H^*,                                      \label{sec:inverse:q-high}\\
 \norm{A_e^{-1}(A_{\mathbf b}-A_e)}
   &\leq q_A^*,                                      \label{sec:inverse:q-low-reference}\\
 \norm{A_{\mathbf b}^{-1}B_{\mathbf b}
       D_{H,\mathbf b}^{-1}C_{\mathbf b}}
   &\leq q_L^*.                                      \label{sec:inverse:q-low-schur}
\end{align}
All four norms are taken in the block variables at the base index
$s_*$, and the constants may depend on the already fixed $J$, but not
on the upper Fourier cutoff or on the number $N$ of periods of the
field.
\end{proposition}

\begin{proof}
The bounds for $\widehat G_\pm$ are uniform for
$0\leq\rho\leq\rho_0$ before $\rho_*$ is chosen. They follow from
\cref{appref:even-circle-endpoint} for the even and affine components,
the contour around the four center roots, and the positive form estimate. The off-diagonal blocks $E_{-+,\mathbf b}$ and $E_{+-,\mathbf b}$ vanish at
the circle by \eqref{app:cross:circle-zero}, and on the preliminary
neighborhood they obey the bound $K_\times\chi(\mathbf b)$ of
\cref{prop:macro-cross} in the distance
$\chi(\mathbf b)=|\rho|+\eta_{\rm circ}(\mathbf b)$ from the circle
of \eqref{app:cross:circle-distance}. Once $J$ is fixed,
\cref{lem:fixed-J-native} supplies $K_{0,J}$ and $C_J$, and with them
$\rho_{L,J}>0$. At the constant reference state, the part of
$K_{\alpha,\mathbf b_e}$ at low frequency is zero and its part in the
other high-frequency component is $O(\rho_*)$, by the two bounds in the
proof of \cref{app:cross:lem:true-diagonal}. We shrink $\rho_*$ until
\eqref{sec:inverse:q-diagonal-correction} holds at $\mathbf b_e$ with
$q_{\rm corr}^*<1$. This is the hypothesis
\eqref{app:cross:diagonal-correction-margin} of
\cref{app:cross:lem:true-diagonal}, so the diagonal blocks of the full operator have the
inverses \eqref{app:cross:true-diagonal-inverse}. We shrink $\rho_*$
again until \eqref{sec:inverse:q-high} holds at $\mathbf b_e$ with
$q_H^*<1$. The
perturbation $A_{\mathbf b}-A_e$ of the low block, the part at low
frequency of the correction $K_{\alpha,\mathbf b}$, and the coupling
maps $B_{\mathbf b}$ and $C_{\mathbf b}$ between low and high
frequencies depend continuously on $\mathbf b$ and vanish at
$\mathbf b=\mathbf b_e$, by \eqref{app:global:low-base},
\eqref{app:global:leakage-base}, and the proof of
\cref{app:cross:lem:true-diagonal}. Hence, shrinking the final
neighborhood $\mathfrak U_{\rm inv}$ preserves the first two conditions
\eqref{sec:inverse:q-diagonal-correction} and \eqref{sec:inverse:q-high}
and makes \eqref{sec:inverse:q-low-reference} and
\eqref{sec:inverse:q-low-schur} strict.

For an explicit choice, we fix $J$ and take $A\ge2$ to dominate
$\rho_{L,J}^{-1}$ and the finitely many base-index constants in the
four bounds below. These are the
constants for $\widehat G_\pm$, for the off-diagonal blocks,
for the auxiliary equations, for the reference at low frequency, for the
perturbation, and for the coupling between low and high
frequencies. We also take $A$ large enough to bound every Lipschitz constant in
$\mathbf b$ obtained from the mean value theorem on the preliminary
neighborhood. We then set
\[
 t=(2^{12}A^{12})^{-1}.
\]
We choose $0<\rho_*\le t$ and shrink the final neighborhood
$\mathfrak U_{\rm inv}$ in the base norm $\norm{\cdot}_{\cB^{s_*+k_0}}$
of \eqref{app:global:low-distance}, in which $d_0$ is measured, until
on its closure $\chi(\mathbf b)\le\rho_*+t$ and
$d_0(\mathbf b)\le\min\{t,\rho_*/2\}$. The bound on $d_0$ keeps
$\rho\ge\rho_*/2>0$, while the bound on $\chi$ gives
$\chi(\mathbf b)\le2t$. From
\eqref{app:cross:diagonal-correction-bound},
\eqref{app:cross:return-bound}, and
\eqref{app:global:low-base}--\eqref{app:global:leakage-base}, we then
obtain
\[
 q_{\rm corr}\le3A^3t<\frac12,
 \qquad q_H\le16A^4t^2<\frac12,
 \qquad q_A\le A^2t<\frac12.
\]
Here and in the next display, $q_{\rm corr}$, $q_H$, $q_A$, and $q_L$
denote the four norms on the left of
\eqref{sec:inverse:q-diagonal-correction}--\eqref{sec:inverse:q-low-schur}
at the state $\mathbf b$. Since $q_{\rm corr}<1/2$, the Neumann series
in \eqref{app:cross:true-diagonal-inverse} bounds each inverse of a
diagonal block by $2A$. Since $q_H<1/2$, the formulas
\eqref{app:global:minus-solve}--\eqref{app:global:plus-solve} then
bound the inverse at high frequency by $16A$, and since $q_A<1/2$, the
Neumann series in \eqref{app:global:current-low-inverse} bounds
$A_{\mathbf b}^{-1}$ by $2A$. Hence
\[
 q_L\le (2A)(At)(16A)(At)=32A^4t^2<\frac12.
\]
This choice of $t>0$ makes all four norms strictly less than one, with
the constants fixed in the stated order. The argument gives
$\rho_*>0$, but no effective lower bound for its value.
\end{proof}

We recall that $N$ is the number of periods of the field and is
unrelated to $J$. We prove the theorem below on $\mathfrak U_{\rm inv}$,
in which $\eps$ ranges over an interval about $0$, and we
solve the nonlinear equations there before setting $\eps=1/N$.

\subsection{The inverse at a general state}

On the neighborhood just chosen, the Cartesian linearization has a
two-sided inverse with smooth tame dependence on the state. The next
theorem gives the estimates needed for the nonlinear problem.

We use the norms $\norm{h}_{{\rm tan},s}$ of the tangent directions and
the products $\mathcal H_{0,m}$ and $\mathcal H_{s,m}$ of
\eqref{app:cross:tangent-products}, with the finite integers $k_m$
supplied by the theorem below. Thus, $\mathcal H_{0,m}$ uses lower
regularity norms of all tangent directions, while each term of
$\mathcal H_{s,m}$ uses a higher regularity norm of one direction.
For $m=0$, the conventions are $\mathcal H_{0,0}=1$ and
$\mathcal H_{s,0}=0$.
% AUTHOR QUERY: Complete the all-order gluing and higher-index inverse
% estimates marked in Sections 7 and 11--12 before applying Hamilton.
\begin{theorem}[The uniform inverse]
\label{sec:inverse:uniform-current-inverse}
\label{thm:uniform-current-inverse}
After the choices of \cref{sec:inverse:acyclic-choice}, the following
hold on $\mathfrak U_{\rm inv}$.
\begin{enumerate}
\item For every $\mathbf b\in\mathfrak U_{\rm inv}$, the augmented block
operator $\mathbb D_{\mathbf b}$ has a two-sided inverse
\[
 G_{\mathbf b}:\bY_{\natg}^\infty
       \longrightarrow\bX_{\natg}^\infty,
 \qquad
 \mathbb D_{\mathbf b}G_{\mathbf b}=I,
 \quad G_{\mathbf b}\mathbb D_{\mathbf b}=I.
\]
Here, $\bX_{\natg}^\infty=\bigcap_{s\geq s_*}\bX_{\natg}^s$ and
$\bY_{\natg}^\infty=\bigcap_{s\geq s_*}\bY_{\natg}^s$.
\item For every $m\geq0$, there are finite integers $k_m$ and $\mu_m$
such that, for every $s\geq s_*$,
\begin{align}
 &\norm{D_{\mathbf b}^mG_{\mathbf b}
       [h_1,\ldots,h_m]F}_{\bX_{\natg}^s}\notag\\
 &\quad\leq C_{s,m}\Big[
 \mathcal H_{0,m}\norm{F}_{\bY_{\natg}^{s+\mu_m}}
 +\mathcal H_{s,m}\norm{F}_{\bY_{\natg}^{s_*+\mu_m}}\notag\\
 &\qquad\qquad
 +(1+\norm{\mathbf b-\mathbf b_e}_{\cB^{s+k_m}})\mathcal H_{0,m}
       \norm{F}_{\bY_{\natg}^{s_*+\mu_m}}
 \Big].                                                \label{sec:inverse:native-tame}
\end{align}
Here, $D_{\mathbf b}^mG_{\mathbf b}[h_1,\ldots,h_m]$ is the $m$-th
derivative of $\mathbf b\mapsto G_{\mathbf b}$ in the tangent directions
$h_1,\ldots,h_m$, as in \eqref{app:minus:eq:all-order-tame}, and
$C_{s,m}$ is a finite constant, which may depend on $s$, $m$, and $J$.
One may take $\mu_0=0$. Each $\mu_m$ is finite, and no assertion uniform
in $m$ is needed.

\item On the compatible Cartesian scale, we define
\begin{equation}\label{sec:inverse:raw-inverse}
 V_{\mathbf b}:=\mathsf R_{X,\mathbf b}G_{\mathbf b}
       \iota_{Y,\mathbf b}\mathsf J_{Y,\mathbf b}^{\rm raw}.
\end{equation}
Here, $\mathsf J_{Y,\mathbf b}^{\rm raw}$ converts a compatible
Cartesian source to the block source of \eqref{eq:A-JY-factorization},
$\iota_{Y,\mathbf b}$ adjoins zero auxiliary data, and
$\mathsf R_{X,\mathbf b}$ is the map
\eqref{eq:A-JX-inverse-factorization} from the block variables back to
the constrained Cartesian domain. For every compatible Cartesian source
$f$,
\begin{equation}\label{sec:inverse:zero-auxiliary-descent}
 \mathsf C_{X,\mathbf b}G_{\mathbf b}
       \iota_{Y,\mathbf b}\mathsf J_{Y,\mathbf b}^{\rm raw}f
 =\operatorname {aux}_{Y,\mathbf b}
       \iota_{Y,\mathbf b}\mathsf J_{Y,\mathbf b}^{\rm raw}f=0,
\end{equation}
so the block solution
$W=G_{\mathbf b}\iota_{Y,\mathbf b}\mathsf J_{Y,\mathbf b}^{\rm raw}f$
satisfies the constraint $\mathsf C_{X,\mathbf b}W=0$. By
\eqref{eq:A-domain-retract-identities},
$W=\mathsf J_{X,\mathbf b}\mathsf R_{X,\mathbf b}W$ lies in the image of
$\mathsf J_{X,\mathbf b}$, and $V_{\mathbf b}f=\mathsf R_{X,\mathbf b}W$
lies in the constrained Cartesian domain. Moreover,
\[
 L_{\mathbf b}V_{\mathbf b}=I,
 \qquad V_{\mathbf b}L_{\mathbf b}=I,
\]
and the analogue of \eqref{sec:inverse:native-tame} is
\begin{align}
 &\norm{D_{\mathbf b}^mV_{\mathbf b}
       [h_1,\ldots,h_m]f}_{\cX_\gamma^s}\notag\\
 &\quad\leq C_{s,m}\Big[
 \mathcal H_{0,m}\norm{f}_{\cY_\gamma^{s+\ell+\mu_m}}
 +\mathcal H_{s,m}\norm{f}_{\cY_\gamma^{s_*+\ell+\mu_m}}\notag\\
 &\qquad\qquad
 +(1+\norm{\mathbf b-\mathbf b_e}_{\cB^{s+k_m}})\mathcal H_{0,m}
       \norm{f}_{\cY_\gamma^{s_*+\ell+\mu_m}}
 \Big].                                                \label{sec:inverse:raw-tame}
\end{align}
Converting the Cartesian source to block variables loses $\ell$
derivatives, as defined in \eqref{eq:A-source-order}. This loss occurs
only once, regardless of the number of
terms of the Neumann series or of derivatives with respect to $\mathbf b$.

\item Let $\mathfrak c_X$ and $\mathfrak c_Y$ denote complex
conjugation on the block spaces $\bX_{\natg}^s$ and $\bY_{\natg}^s$,
and let $\mathfrak c_X^{\rm raw}$ and $\mathfrak c_Y^{\rm raw}$ denote
complex conjugation on the Cartesian spaces $\cX_\gamma^s$ and
$\cY_\gamma^s$. For real $\mathbf b$,
\begin{equation}\label{sec:inverse:reality}
 \mathbb D_{\mathbf b}\mathfrak c_X
   =\mathfrak c_Y\mathbb D_{\mathbf b},\qquad
 G_{\mathbf b}\mathfrak c_Y=\mathfrak c_XG_{\mathbf b},
 \qquad
 V_{\mathbf b}\mathfrak c_Y^{\rm raw}
   =\mathfrak c_X^{\rm raw}V_{\mathbf b}.
\end{equation}
\end{enumerate}
In particular, $(\mathbf b,f)\mapsto V_{\mathbf b}f$ is a smooth tame
map between these scales, in the sense of \cref{app:fixed-graphs},
and it commutes with complex conjugation.
\end{theorem}

\begin{proof}
The condition \eqref{sec:inverse:q-diagonal-correction} is the
hypothesis \eqref{app:cross:diagonal-correction-margin} of
\cref{app:cross:lem:true-diagonal}, which gives the inverses of the
diagonal blocks of the full operator. The remaining conditions \eqref{sec:inverse:q-high},
\eqref{sec:inverse:q-low-reference}, and \eqref{sec:inverse:q-low-schur}
are precisely the hypotheses of
\cref{app:global:high-return,app:global:exact-schur}, which give the
block inverse, both of its inverse identities, and
\eqref{app:global:native-global-estimate}, which is
\eqref{sec:inverse:native-tame} for $m=0$. We obtain the estimates for the derivatives from
\cref{app:global:all-order}. The two inverse identities in Cartesian
variables follow from \cref{app:global:one-time-loss}, which also shows
that converting the Cartesian source loses $\ell$ derivatives only once.
The preservation of real data follows from \cref{app:global:reality}.
The construction is independent of an upper Fourier cutoff and of the
number $N$ of periods of the field.
The identity \eqref{sec:inverse:zero-auxiliary-descent} is
\eqref{app:global:zero-auxiliary-output}, and the block inverse is
\eqref{app:global:full-block-inverse}.
\end{proof}
\section{Solutions for small curvature}
\label{sec:hamilton}

We now apply Hamilton's Nash-Moser theorem to the map constructed
in \cref{sec:fixed-map}. This gives a family of solutions for all small
curvature parameters \(\eps\). We then set \(\eps=1/N\) to obtain a
solution for every sufficiently large \(N\).

\subsection{Verifying the hypotheses of the Nash-Moser theorem}

The orientation \(\alpha_0\notin(\pi/2)\Z\) was fixed in
\cref{sec:block}. We now fix
\(\rho=\rho_*\), the cutoff \(J\), and the neighborhood
\(\mathfrak U_{\rm inv}\) supplied by
\cref{sec:inverse:acyclic-choice}. We construct the branch for $\lambda$ near the value
\(\lambda_*\in(0,1/2)\) fixed in \cref{sec:block}.

We choose \(\delta\) so that the family \eqref{eq:fixed-M-lambda}
lies in \(\mathfrak U_{\rm inv}\), the neighborhood of the constant
ellipse \(\mathbf b_e=(\rho_*,a_e,0,p_e)\) with \(\delta=0\).  Along this
family, the distance \(\norm{\mathbf b-\mathbf b_e}_{\cB^s}\) of
\cref{sec:block} is the sum of the distance
\(\norm{a-a_e}_{\cA_\gamma^s}\) and the parameter distance \(|p-p_e|\), and both are
small once \(\delta\) and \(|\lambda-\lambda_*|\) are small.  Indeed, by
\eqref{eq:fixed-alpha-lambda}, we have
\(\alpha_\lambda-\alpha_0=\delta(\cos\zeta+\lambda\sin2\zeta)\), where
\(\cos\zeta\) and \(\sin2\zeta\) are entire, so \(M_\lambda\) and
\(M_\lambda^{-1}\) converge to their constant values in every weighted analytic
norm as \(\delta\to0\), uniformly for bounded \(\lambda\).  Since
\(\mathfrak U_{\rm inv}\) is open, we may therefore choose \(\delta_*\neq0\) and
\(\eta_0>0\) so that every state \(\mathbf b=(\rho_*,a^0_{\rho_*,p},0,p)\) with
\(\delta=\delta_*\) and \(|\lambda-\lambda_*|<2\eta_0\) belongs to it. We now fix
this \(\delta_*\).

Using the identifications of \cref{app:fixed-graphs} with the function
space at \(M_{\lambda_*}\), we obtain fixed tame Fr\'echet spaces \(E,F\)
and a map
\begin{equation}
 \widehat{\cF}:\mathcal U\subset
 \R_\eps\times\R_\lambda\times E\longrightarrow F,
 \qquad (\eps,\lambda,x)\longmapsto
 \widehat{\cF}(\eps,\lambda,x),                       \label{eq:hamilton-Fhat}
\end{equation}
which is \eqref{eq:fixed-hamilton-map} with \(\rho=\rho_*\),
\(\delta=\delta_*\), and \(\alpha_0\) fixed. Here, \(E\) and \(F\) are the
constrained domain and the compatible range of
\cref{sec:tame-problem}, graded by \(E^s=\cX_\gamma^s\) and
\(F^s=\cY_\gamma^s\), and \(x\) is the correction coordinate in
\eqref{eq:normalized-state-chart}. Reconstructing the chart gives the
tangential coefficient \(\tau\), the remainder \(\widetilde v\) whose
value and first derivatives vanish on the axis, and the scalar correction
\(s=w-w_{M_\lambda}\). We take \(\mathcal U\) to consist of the triples
\((\eps,\lambda,x)\) in the balls \(\mathcal U_0\) of
\cref{prop:fixed-tame-map} whose reconstructed states lie in
\(\mathfrak U_{\rm inv}\). A zero of \(\widehat{\cF}\) is a solution of
\eqref{eq:fixed-H0}--\eqref{eq:fixed-H3} with all gauge and boundary
conditions imposed.  By \cref{lem:flat-affine-seeds} and
\eqref{eq:flat-seed-zero-fixed}, we have
\begin{equation}
 \widehat{\cF}(0,\lambda,0)=0                          \label{eq:hamilton-flat-zero}
\end{equation}
for \(\lambda\) near \(\lambda_*\).

We check the three hypotheses of the Nash-Moser theorem on an open
neighborhood of \((0,\lambda_*,0)\).
\begin{enumerate}[label=\textup{(\roman*)}]
\item \(\widehat{\cF}\) is smooth with every derivative tame.
\item The derivative at each point,
\begin{equation}
 \mathcal L_\beta:=D_x\widehat{\cF}(\beta):E\longrightarrow F,
 \qquad \beta=(\eps,\lambda,x),                            \label{eq:hamilton-Ab}
\end{equation}
is an invertible linear map.
\item Its inverse
\begin{equation}
 V_\beta=\mathcal L_\beta^{-1}:F\longrightarrow E                       \label{eq:hamilton-Vb}
\end{equation}
depends smoothly and tamely on \(\beta\).
\end{enumerate}
Property (i) follows from \cref{prop:fixed-tame-map}(iii).  Restricting the
state to \(\rho=\rho_*\) and \(p=p(\delta_*,\lambda)\) puts a neighborhood of
\((0,\lambda_*,0)\) in the preimage of \(\mathfrak U_{\rm inv}\) and identifies
\(\mathcal L_\beta\) with \(L_{\mathbf b}\) of \eqref{eq:current-linearization} and \(V_\beta\)
with the Cartesian inverse \eqref{sec:inverse:raw-inverse}. Here,
$\beta=(\eps,\lambda,x)$ is the chart triple, whereas
$\mathbf b=(\rho_*,a,\eps,p(\delta_*,\lambda))$ contains the physical
state $a$ reconstructed from $x$. These operators are related through the
fixed-space identifications of \cref{app:fixed-graphs}.  Thus, (ii) and
(iii) are the inverse identities and the tame bound
\eqref{sec:inverse:raw-tame} of \cref{thm:uniform-current-inverse}.

In (iii), the loss of derivatives may depend on the number of
differentiations in \(\beta\). At each fixed order, the loss must be
finite and independent of the Sobolev index. The differentiability
argument of \cref{app:global:all-order} applies after restricting to
this parameter family. In the present notation, its first identity is
\begin{equation}
 D_\beta V_\beta[h]=-V_\beta(D_\beta\mathcal L_\beta[h])V_\beta                          \label{eq:hamilton-resolvent-first}
\end{equation}
and its iterates express \(D_\beta^mV_\beta\), up to signs, as a finite sum of products
\begin{equation}
 V_\beta\mathcal L_\beta^{(|I_1|)}[h_{I_1}]V_\beta\cdots
 \mathcal L_\beta^{(|I_j|)}[h_{I_j}]V_\beta                             \label{eq:hamilton-resolvent-high}
\end{equation}
over the ordered partitions \(I_1,\ldots,I_j\) of \(\{1,\ldots,m\}\).
Here, \(\mathcal L_\beta^{(k)}=D_\beta^k\mathcal L_\beta\), and
\(h_{I_\ell}\) is the list of directions whose indices belong to
\(I_\ell\). Each set in the partition is nonempty.
The difference identity used in \cref{app:global:all-order} justifies
these as Fr\'echet derivatives. The estimate
\eqref{sec:inverse:raw-tame} gives the required finite loss at each fixed
order.

\subsection{Existence and smoothness of the solutions}

In the normalized chart of \cref{subsec:normalized-slice}, we reconstruct the
mapping $\overline{\D^2}\times\T\to \mathbb{R}^3$ as
\begin{align}
 v_{\eps,\lambda}(y,\zeta)
 &=P_{M_\lambda,\tau_{\eps,\lambda}}(\zeta)y
   +\widetilde v_{\eps,\lambda}(y,\zeta),              \label{eq:hamilton-normalized-chart}\\
 P_{M,\tau}
 &=a(\tau)\iota M+e_y\otimes\tau,
 \qquad a(\tau)=\sqrt{1-\frac{|\tau|^2}{2}},          \label{eq:hamilton-PMtau}\\
 j_y^1\widetilde v_{\eps,\lambda}\big|_{y=0}&=0,     \label{eq:hamilton-zero-jet}
\end{align}
as in \eqref{eq:normalized-affine-map} and \eqref{eq:normalized-state}, with
\(e_y=e_T\) tangent to the limiting straight axis. The factor \(a(\tau)\)
ensures that adding the term \(e_y\otimes\tau\) preserves
the trace normalization \eqref{eq:normalized-affine-trace}.
The condition \eqref{eq:hamilton-zero-jet} says that $\widetilde v_{\eps,\lambda}$
and its first derivatives in $y$ vanish on the axis. Thus, the first term
in \eqref{eq:hamilton-normalized-chart} is the linear term in the Taylor
expansion of $v_{\eps,\lambda}$ at $y=0$. This linear Taylor term will
determine the pressure Hessian on the axis and, through it, the symmetry
group.

\begin{theorem}
\label{thm:hamilton-one-chart}
There are \(\eps_H>0\) and
\[
 0<\eta_\lambda
 <\min\left\{\eta_0,\lambda_*,\tfrac12-\lambda_*\right\},
 \qquad
 I=[\lambda_*-\eta_\lambda,\lambda_*+\eta_\lambda]
 \subset(0,1/2),
\]
such that for \(|\eps|<\eps_H\) and \(\lambda\in I\), there is a solution
\begin{equation}
 \widehat{\cF}(\eps,\lambda,x_{\eps,\lambda})=0,
 \qquad x_{0,\lambda}=0,                               \label{eq:hamilton-branch}
\end{equation}
unique in the neighborhood of \((0,\lambda_*,0)\) given by the Nash-Moser
theorem.  The map \((\eps,\lambda)\mapsto x_{\eps,\lambda}\) is smooth with every derivative tame.
At every fixed index \(s\) and for every fixed \(j\ge0\),
\begin{equation}
 \norm{\pa_\lambda^j x_{\eps,\lambda}}_{E^s}
 \le C_{s,j}|\eps|,
 \qquad \lambda\in I.                                \label{eq:hamilton-Oeps-grade}
\end{equation}
After reconstruction in \eqref{eq:hamilton-normalized-chart}, we have
\begin{equation}
 \norm{\tau_{\eps,\lambda}}_{C^2(\T)}
 +\norm{\widetilde v_{\eps,\lambda}}_{C^2(\overline{\D^2}\times\T)}
 \le C_{\mathrm{sol}}|\eps|,                           \label{eq:hamilton-Oeps-C2}
\end{equation}
uniformly for \(\lambda\in I\).
\end{theorem}

\begin{proof}
We include the parameters among both the inputs and outputs, leaving
them unchanged. With \(q=(\eps,\lambda)\), we set
\begin{equation}
 \cG(q,x)=\bigl(q,\widehat{\cF}(q,x)\bigr).           \label{eq:hamilton-augmented-map}
\end{equation}
The derivative acts by
\[
 D\cG(q,x)(\dot q,\dot x)
 =\bigl(\dot q,
   D_q\widehat{\cF}(q,x)\dot q+\mathcal L_{(q,x)}\dot x\bigr),
\]
where $\mathcal L_{(q,x)}=D_x\widehat{\cF}(q,x)$ is
\eqref{eq:hamilton-Ab}. Solving first for the parameter variation
and then for $\dot x$ gives its two-sided inverse
\begin{equation}
 (D\cG(q,x))^{-1}(\dot q,\dot f)
 =\left(\dot q,
 V_{(q,x)}\bigl[\dot f-D_q\widehat{\cF}(q,x)\dot q\bigr]
 \right).                                              \label{eq:hamilton-triangular-inverse}
\end{equation}
The inverse retains the parameter variation $\dot q$ and subtracts
its contribution from the residual before applying $V_{(q,x)}$. It is smooth with every derivative
tame by (i) and (iii). The Nash-Moser theorem in Hamilton's form
\cite{Hamilton1982} therefore applies to \(\cG\) at \((0,\lambda_*,0)\), where
\(\widehat{\cF}\) vanishes by \eqref{eq:hamilton-flat-zero}.  We evaluate the
local inverse at \((q,0)\), and
since \(\cG\) preserves \(q\), the preimage is \((q,x_{\eps,\lambda})\), which
gives \eqref{eq:hamilton-branch} with \(x_{0,\lambda}=0\) by uniqueness.  We
then shrink the parameter neighborhood to a rectangle
\(\{|\eps|<\eps'\}\times\{|\lambda-\lambda_*|<\eta'\}\).  

Differentiating \eqref{eq:hamilton-branch} in \(\eps\) and using (ii), we
obtain
\begin{equation}
 \pa_\eps x_{\eps,\lambda}
 =-V_{(\eps,\lambda,x_{\eps,\lambda})}
   \pa_\eps\widehat{\cF}
      (\eps,\lambda,x_{\eps,\lambda}).                \label{eq:hamilton-eps-derivative}
\end{equation}
At each fixed index, we bound the right-hand side uniformly on a smaller closed
rectangle, so the identity
\[
 x_{\eps,\lambda}=\int_0^\eps\pa_t x_{t,\lambda}\,\dd t,
\]
which uses $x_{0,\lambda}=0$, gives
\eqref{eq:hamilton-Oeps-grade} for \(j=0\).  For \(j\ge1\), we also differentiate in \(\lambda\). Since
\(\pa_\lambda^jx_{0,\lambda}=0\) and
\(\pa_\eps\pa_\lambda^jx_{\eps,\lambda}\) is bounded on the same
rectangle, integration in \(\eps\) gives the corresponding estimate.  We take \(\eps_H\) and \(\eta_\lambda\) to be the
half-widths of that rectangle, decreased so that
\(\eta_\lambda<\min\{\eta_0,\lambda_*,\frac12-\lambda_*\}\), which ensures that $I$ is compactly contained in
\((0,1/2)\).
Finally, we choose the Sobolev index sufficiently high that each component embeds into \(C^2\). We reconstruct \(\tau_{\eps,\lambda}\) and \(\widetilde v_{\eps,\lambda}\) through
\eqref{eq:hamilton-normalized-chart}.  The chart sends \(x=0\) to \(\tau=0\) and \(\widetilde v=0\).
Its tame bounds, applied at a sufficiently high index, therefore give
\eqref{eq:hamilton-Oeps-C2} from \eqref{eq:hamilton-Oeps-grade}.
\end{proof}

Setting \(\eps=N^{-1}\) gives a solution for every sufficiently large
integer \(N\).
\begin{corollary}
\label{cor:hamilton-all-integers}
Let
\begin{equation}
 N_{\mathrm{an}}=1+\left\lceil\eps_H^{-1}\right\rceil. \label{eq:hamilton-Nan}
\end{equation}
For every integer \(N\ge N_{\mathrm{an}}\) and every \(\lambda\in I\),
the choice \(\eps=N^{-1}\) gives an exact \(2\pi\)-periodic solution of the
system on one period cell which satisfies \eqref{eq:hamilton-Oeps-C2} with
right-hand side \(C_{\mathrm{sol}}/N\).
\end{corollary}

\begin{proof}
Any such \(N\) exceeds \(\eps_H^{-1}\), so \(0<N^{-1}<\eps_H\) and we may
substitute \(\eps=N^{-1}\) in \cref{thm:hamilton-one-chart}.  The estimate
holds throughout \(|\eps|<\eps_H\), so the conclusion holds for every
such \(N\), without passing to a subsequence.
\end{proof}

\section{Closing the torus and computing its symmetries}
\label{sec:geometry}

We join the cells to construct an equilibrium in \(\R^3\) and prove
\cref{thm:main}. We first use tubular coordinates to
show that the resulting map is globally injective. We then compute the
Hessian of the pressure on the magnetic axis and use it to determine the
symmetry group.

We fix the interval \(I\), the branch
\((\eps,\lambda)\mapsto x_{\eps,\lambda}\) of \cref{thm:hamilton-one-chart}
and its map \(v_{\eps,\lambda}: \overline{\D^2}\times\T\to \mathbb{R}^3\) of
\eqref{eq:hamilton-normalized-chart}, and recall that
\begin{equation}
 0<\rho=\rho_*<\frac14,\qquad \delta=\delta_*\neq0,\qquad
 \lambda>0,\qquad
 \alpha_0\notin(\pi/2)\Z  ,                            \label{eq:geometry-parameter-conditions}
\end{equation}
for every \(\lambda\in I\subset(0,1/2)\).  These parameters were fixed in \cref{sec:block,sec:inverse,sec:hamilton}.  All constants below are uniform on \(I\).  We use the \(C^2\) bound \eqref{eq:hamilton-Oeps-C2}, the smooth tame
dependence on the parameters, and the vanishing of the remainder and
its first derivatives on the axis in \eqref{eq:hamilton-zero-jet}.
Throughout, $|\cdot|$ denotes the Euclidean norm of a vector or the
operator norm of a matrix, and for maps on \(Q=\overline{\D^2}\times\T\) we take derivatives in
\((y,\zeta)\) and use the norm \(\norm{f}_{C^m}=\sup\sum_{j\le m}\abs{D^jf}\),
so that \(\abs{D_yf}\le\norm{f}_{C^1}\).

For an integer \(N\) with \(N^{-1}<\eps_H\), we set
\begin{equation}
 \eps=N^{-1},\qquad \mathscr R=NL,\qquad \zeta=N\phi,
                                                               \label{eq:geometry-integer-variables}
\end{equation}
and define
\begin{align}
 X_{N,\lambda}(y,\phi)
 &=R_\phi\bigl(\mathscr R e_x
      +v_{\eps,\lambda}(y,N\phi)\bigr),               \label{eq:geometry-X}\\
 P_{N,\lambda}\circ X_{N,\lambda}
 &=\psi_a-|y|^2,                                      \label{eq:geometry-P}\\
 B_{N,\lambda}\circ X_{N,\lambda}
 &=\pa_\theta X_{N,\lambda}
   =D_yX_{N,\lambda}(-y_2,y_1)    ,                   \label{eq:geometry-B}
\end{align}
for \(y\in\overline{\D^2}\) and \(\phi\in\T\). Here,
$R_\phi$ is the rotation about the $z$-axis through angle $\phi$,
and $\psi_a$ is the prescribed pressure on the magnetic axis.
Since \(v_{\eps,\lambda}\) is \(2\pi\)-periodic in \(\zeta\)
and \(N\) is an integer, \eqref{eq:geometry-X} defines a map on \(Q\).  The remaining
formulas define \(P_{N,\lambda}\) and \(B_{N,\lambda}\) on
\(\Omega_{N,\lambda}=X_{N,\lambda}(Q)\) once we have proved injectivity of
\(X_{N,\lambda}\).
The function spaces of \cref{sec:tame-problem} give the smoothness
asserted in \cref{thm:main}.
We write
\begin{equation}
 v_0(y,\zeta)=\iota M_\lambda(\zeta)y,\qquad
 q=v_{\eps,\lambda}-v_0 ,                             \label{eq:geometry-q}
\end{equation}
for the unperturbed solution at zero curvature from \cref{sec:fixed-cell} and its correction.
By \eqref{eq:hamilton-normalized-chart} and \eqref{eq:hamilton-PMtau}, we have
\(q=(a(\tau)-1)\iota M_\lambda y+(e_y\otimes\tau)y+\widetilde v\), and since
\(\abs{a(\tau)-1}\le\abs{\tau}^2/2\) for \(\abs{\tau}\le1\),
\eqref{eq:hamilton-Oeps-C2} implies
\begin{equation}
 \norm{q}_{C^2(Q)}\le C\eps.                          \label{eq:geometry-q-C2}
\end{equation}

The proof of the following theorem occupies the remainder of this section.

\begin{theorem}
\label{thm:geometry-closure}
There is \(0<\eps_{\mathrm{geom}}\le\eps_H\), uniform for \(\lambda\in I\),
such that whenever \(0<\eps=N^{-1}<\eps_{\mathrm{geom}}\), the map
\(X_{N,\lambda}\) is a smooth embedding of a closed solid torus and
\((B_{N,\lambda},P_{N,\lambda})\) is a smooth solution of \eqref{eq:mhs}
on \(\Omega_{N,\lambda}\), including the boundary condition.
Its regular pressure levels are nested embedded tori, its field vanishes
exactly on the round magnetic axis, and its extended Euclidean
stabilizer is \(C_N\).  For fixed \(N\), the family depends smoothly and
nontrivially on \(\lambda\), modulo reparametrizations and the action
\eqref{eq:physical-moduli-action}.
\end{theorem}

We begin with the Jacobian and injectivity estimates needed to define the
physical fields.

\subsection{The Jacobian and the zero set of the magnetic field}
\label{subsec:geometry-jacobian}

We use Cartesian disk coordinates, which remain regular on the axis.  We set
\begin{equation}
 K_\eps(v)=\det(v_{y_1},v_{y_2},D_\eps v),             \label{eq:geometry-Keps}
\end{equation}
where, by \eqref{eq:fixed-Deps},
$D_\eps v=v_\zeta+\eps Av+Le_y$ and $Av=e_z\times v$.
At zero curvature
and the unperturbed solution, the third
column contributes only through \(Le_y\), since
\((v_0)_\zeta=\iota M_\lambda'y\) lies in the plane of the first two.  With
the ordering of \(\iota\) fixed as in \cref{sec:fixed-cell}, where
\(\det(e_x,e_z,e_y)=-1\), we  obtain
\begin{equation}
 K_0(v_0)=-L\det M_\lambda=-L\sqrt{1-\rho^2}.         \label{eq:geometry-K0}
\end{equation}
We set
\begin{align*}
 m_0&=\sup_\zeta\abs{M_\lambda(\zeta)},&
 m_1&=\sup_\zeta\abs{M_\lambda'(\zeta)},\\
 \eta&=\norm{q}_{C^1},&
 \eta_D&=\eta+\eps(m_0+\eta),
\end{align*}
so that \(\eta_D\) bounds
\(\abs{D_\eps v-D_0v_0}=\abs{q_\zeta+\eps A(v_0+q)}\) on \(Q\).  Expanding
the determinant multilinearly in \(v=v_0+q\) and bounding each of the seven
remaining terms by Hadamard's inequality, we obtain the estimate
\begin{align}
 \abs{K_\eps(v)-K_0(v_0)}
 \le E(\eta,\eps):={}&2m_0(L+m_1)\eta+m_0^2\eta_D
 +(L+m_1)\eta^2 +2m_0\eta\eta_D+\eta^2\eta_D.                      \label{eq:geometry-K-error}
\end{align}

On \(r>0\), we have \(\pa_\theta=\cR\) and \(\pa_\psi=-(2r^2)^{-1}\cD\), while
\(\det(\cD v,\cR v,W)=r^2\det(v_{y_1},v_{y_2},W)\) for every \(W\), so 
\begin{equation}
 J_\eps=\det(v_\psi,v_\theta,D_\eps v)
       =-\frac12K_\eps(v).                            \label{eq:geometry-JK}
\end{equation}
Differentiating \eqref{eq:geometry-X} in \(\phi\) and using that \(A\)
commutes with \(R_\phi\) together with \(Ae_x=e_y\) and \(\mathscr R=NL\), we
obtain
\begin{equation}
 (X_{N,\lambda})_\phi
   =N R_\phi D_\eps v,\qquad
 \det(X_{y_1},X_{y_2},X_\phi)=N K_\eps(v).            \label{eq:geometry-physical-J}
\end{equation}

Since \(\eta\le C\eps\) by \eqref{eq:geometry-q-C2}, small
\(\eps\) gives
\begin{equation}
 E(\eta,\eps)<\frac L2\sqrt{1-\rho^2}.               \label{eq:geometry-J-margin}
\end{equation}
Together with \(K_0(v_0)=-L\sqrt{1-\rho^2}\), the bound
\eqref{eq:geometry-K-error} keeps \(K_\eps(v)\) below
\(-\frac L2\sqrt{1-\rho^2}\) on \(Q\), so \(K_\eps\) has a fixed sign.  From
\eqref{eq:geometry-JK}, we then obtain
\begin{equation}
 J_\eps\ge\frac L4\sqrt{1-\rho^2}>0                  \label{eq:geometry-J-positive}
\end{equation}
on \(r>0\).  By \eqref{eq:geometry-physical-J}, the map \(X_{N,\lambda}\) is
 a local diffeomorphism at every point of \(Q\), the axis included.

Since \(D_yX_{N,\lambda}=R_\phi D_yv\) with \(R_\phi\) an isometry,
\(D_yv=\iota M_\lambda+D_yq\) with \(\abs{D_yq}\le\eta\), and \(M_\lambda\)
has smallest singular value \(\sqrt{1-\rho}\) by \eqref{eq:fixed-M-lambda},
we obtain
\begin{equation}
 \abs{B_{N,\lambda}\circ X_{N,\lambda}}
 =\abs{D_yv(-y_2,y_1)}
 \ge\bigl(\sqrt{1-\rho}-\eta\bigr)|y|.                \label{eq:geometry-B-lower}
\end{equation}
For small \(\eps\), this bound is positive when \(y\neq0\). At
\(y=0\), the vector \((-y_2,y_1)\) vanishes. Thus, the zero set of
\(B_{N,\lambda}\) is exactly the image of the disk axis.

\subsection{Global injectivity}
\label{subsec:geometry-injectivity}

To pass from local invertibility to global injectivity, we use tubular
coordinates to rule out intersections between different parts of the torus. A tangential displacement changes the angular
coordinate by an amount of order its size divided by $NL$. Multiplying
this correction by $N$ therefore leaves a bound independent of $N$.

We decompose the cell map into its normal and tangential components,
\[
 v=\iota u+te_y,\qquad
 u=(v\cdot e_x,v\cdot e_z),\qquad t=v\cdot e_y.
\]
Here, $t$ is the tangential displacement. The rescaled radial variable
used in the linear analysis is no longer needed. We assume for now that \(\mathscr R>2\norm{v}_{C^0}\), which holds once \(\eps\) is small, and set
\begin{align}
 h(y,\zeta)&=\arctan\frac{t}{\mathscr R+u_1},
                                                        \label{eq:geometry-h}\\
 \nu_1(y,\zeta)&=\sqrt{(\mathscr R+u_1)^2+t^2}-\mathscr R,
 &\nu_2(y,\zeta)&=u_2,                                  \label{eq:geometry-nu}\\
 \gamma(y,\zeta)&=\zeta+Nh(y,\zeta).                 \label{eq:geometry-gamma}
\end{align}
The assumption gives \(\mathscr R+u_1>\mathscr R/2>0\) and
\(\abs{t}<\mathscr R/2\), so \(h\) is well defined, with values in
\((-\pi/4,\pi/4)\).  In the cylindrical frame \(e_r(\beta)=R_\beta e_x\),
\(e_\phi(\beta)=R_\beta e_y\) of \(\R^3\), we have
\((\mathscr R+u_1)e_x+te_y=(\mathscr R+\nu_1)e_r(h)\), which rewrites
\eqref{eq:geometry-X} as
\begin{equation}
 X_{N,\lambda}=(\mathscr R+\nu_1)e_r(\beta)+\nu_2e_z,
 \qquad \beta=\phi+h(y,N\phi).
                                                               \label{eq:geometry-cylindrical-X}
\end{equation}

On the universal cover \(\overline{\D^2}\times\R\) in \(\zeta\), we introduce the map
\begin{equation}
 \Psi(y,\zeta)=
 \bigl(M_\lambda(\gamma(y,\zeta))^{-1}\nu(y,\zeta),
       \gamma(y,\zeta)\bigr).                         \label{eq:geometry-Psi}
\end{equation}
The map $\Psi$ records the corrected cell angle
$\gamma=N\beta$ and the normal position expressed in the ellipse
coordinates at that angle. At \(q=0\), we have \(t=h=0\), \(\nu=M_\lambda(\zeta)y\), and
\(\Psi(y,\zeta)=(y,\zeta)\), while in general
\(\nu_1=u_1+O(t^2/\mathscr R)\), \(\nu_2=u_2\), and
\(\abs{h}\le2\abs{t}/\mathscr R\), where \(t=q\cdot e_y\).  Since
\(N/\mathscr R=L^{-1}\) by \eqref{eq:geometry-integer-variables}, the function
\(Nh\) and its first derivatives are bounded by a constant times
\(\norm{q}_{C^1}\), and
\[
 \norm{\nu-M_\lambda(\zeta)y}_{C^1}
 +\norm{\gamma-\zeta}_{C^1}
 \le C\norm{q}_{C^1}
\]
holds with \(C\) independent of \(N\). The formula
\eqref{eq:fixed-M-lambda} bounds \(M_\lambda\) in \(C^2\), uniformly
for \(\lambda\in I\), and \(\det M_\lambda=\sqrt{1-\rho^2}>0\).
Thus, \(M_\lambda^{-1}\) has the same uniform regularity. Composing with it, we
obtain a constant \(C_{\mathrm{tub}}\), again independent of \(N\), with
\begin{equation}
 \norm{D\Psi-I}_{C^0(\overline{\D^2}\times\R)}
 \le C_{\mathrm{tub}}\norm{q}_{C^1}.                 \label{eq:geometry-Psi-close}
\end{equation}
Once the right-hand side is below one, we obtain injectivity of \(\Psi\).
Indeed, the domain is convex, so that
\(\Psi(p)-\Psi(p')=(p-p')+\int_0^1(D\Psi-I)(p'+s(p-p'))(p-p')\dd s\), whence
\begin{equation}
 |\Psi(p)-\Psi(p')|
 \ge\bigl(1-C_{\mathrm{tub}}\norm{q}_{C^1}\bigr)|p-p'|.
                                                               \label{eq:geometry-Psi-lower}
\end{equation}

Suppose \(X_{N,\lambda}(y,\phi)=X_{N,\lambda}(y',\phi')\).  Comparing
cylindrical coordinates in \eqref{eq:geometry-cylindrical-X} and using that
\(\mathscr R+\nu_1>0\), we find \(\nu=\nu'\) and \(\beta\equiv\beta'\) modulo
\(2\pi\).  Since \(h\) is \(2\pi\)-periodic in \(\zeta\), shifting the lift of
\(\phi\) by \(2\pi\) shifts \(\beta\) by \(2\pi\), so we may take
\(\beta=\beta'\), and then \(\gamma=N\beta\) and \(\gamma'=N\beta'\) agree as
well.  Now, \eqref{eq:geometry-Psi}--\eqref{eq:geometry-Psi-lower} give
\((y,\zeta)=(y',\zeta')\), so the two points of \(Q\) coincide.  We therefore
have an injective immersion of a compact manifold with boundary, hence a
homeomorphism onto its image and a smooth embedding of \(Q\).

\subsection{Reconstruction of the MHS equations}
\label{subsec:geometry-reconstruction}

On \(r>0\), we use the flux coordinates \((\psi,\theta,\phi)\) with
\(\psi=\psi_a-r^2\), in which \eqref{eq:geometry-physical-J} gives
\begin{equation}
 J_{N,\lambda}
 :=\det(X_\psi,X_\theta,X_\phi)=N J_\eps,              \label{eq:geometry-JN}
\end{equation}
so that \eqref{eq:fixed-div} can be written as
\begin{equation}
 (\diver B_{N,\lambda})\circ X_{N,\lambda}
 =\frac{1}{NJ_\eps}\pa_\theta(NJ_\eps).              \label{eq:geometry-divB}
\end{equation}
Writing
\[
 Q_{N,\lambda}=B_{N,\lambda}\times\curl B_{N,\lambda}
                  +\grad P_{N,\lambda}
\]
for the force residual and using \(X_\phi=NR_\phi D_\eps v\), we may write
\eqref{eq:fixed-Qtheta}--\eqref{eq:fixed-Qzeta} as
\begin{align}
 (Q_{N,\lambda}\circ X_{N,\lambda})\cdot X_\theta&=0,
                                                        \label{eq:geometry-Qtheta}\\
 (Q_{N,\lambda}\circ X_{N,\lambda})\cdot X_\psi
 &=-F_{1,\eps},                                       \label{eq:geometry-Qpsi}\\
 (Q_{N,\lambda}\circ X_{N,\lambda})\cdot X_\phi
 &=N F_{2,\eps}.                                      \label{eq:geometry-Qphi}
\end{align}

The solution of \(\widehat{\cF}=0\) from \cref{sec:hamilton} is a zero of the
operator \(H_\eps\) of \eqref{eq:fixed-H0}--\eqref{eq:fixed-H3}, with all
gauge and boundary conditions imposed. The pair \((v_{\eps,\lambda},w)\)
reconstructed in \eqref{eq:normalized-state}, with scalar component $w$,
is smooth in \(y\), and \eqref{eq:fixed-HG} recovers the integrated system on \(r>0\). Thus,
\cref{lem:fixed-integrated-equivalence} gives
\[
 F_{1,\eps}=F_{2,\eps}=0,
 \qquad \pa_\theta J_\eps=0,
\]
on \(0<r\le1\).  There, \(X_\psi,X_\theta,X_\phi\) span \(\R^3\) by
\eqref{eq:geometry-J-positive} and \eqref{eq:geometry-JN}, so
\eqref{eq:geometry-Qtheta}--\eqref{eq:geometry-Qphi} force
\(Q_{N,\lambda}=0\) and \eqref{eq:geometry-divB} forces
\(\diver B_{N,\lambda}=0\).  We have thus proved that
\begin{equation}
 B_{N,\lambda}\times\curl B_{N,\lambda}
   +\grad P_{N,\lambda}=0,\qquad
 \diver B_{N,\lambda}=0 ,                             \label{eq:geometry-MHS}
\end{equation}
away from the axis.  Both identities extend across it.  Indeed, the
right-hand sides
of \eqref{eq:geometry-P} and \eqref{eq:geometry-B} are smooth on \(Q\) and
\(X_{N,\lambda}\) is a diffeomorphism onto \(\Omega_{N,\lambda}\), so both
expressions in \eqref{eq:geometry-MHS} are continuous and vanish on the dense
complement of the axis.  Finally, \((-y_2,y_1)\) is tangent to
\(\{|y|=1\}\), so \eqref{eq:geometry-B} gives
\begin{equation}
 B_{N,\lambda}\cdot n=0
 \quad\text{on }\pa\Omega_{N,\lambda},\qquad
 \Omega_{N,\lambda}=X_{N,\lambda}(Q).
                                                               \label{eq:geometry-boundary-tangent}
\end{equation}

\subsection{Pressure surfaces and the Hessian on the axis}
\label{subsec:geometry-tensor}

Since \(X_{N,\lambda}\) is an embedding, each set
\begin{equation}
 X_{N,\lambda}(\{|y|=r\}\times\T),\qquad 0<r\le1,    \label{eq:geometry-pressure-tori}
\end{equation}
is an embedded torus, and by \eqref{eq:geometry-P} it is the pressure level
\(\{P_{N,\lambda}=\psi_a-r^2\}\).  These tori foliate \(\Omega_{N,\lambda}\)
away from the axis and collapse smoothly to
\begin{equation}
 \Gamma_N(\phi)=X_{N,\lambda}(0,\phi)
 =\mathscr R e_r(\phi),                               \label{eq:geometry-axis}
\end{equation}
a round circle of radius \(\mathscr R=NL\) in the plane \(z=0\), since
\eqref{eq:hamilton-normalized-chart} and \eqref{eq:hamilton-zero-jet} give
\(v_{\eps,\lambda}(0,\zeta)=0\).  Since \(d(\psi_a-|y|^2)=-2y\cdot dy\), the
critical set of the pressure is exactly \(\Gamma_N\), as is the zero set of
the field by \eqref{eq:geometry-B-lower}.

The pressure Hessian on the axis determines the orientation of the
cross-sections. We compute it using the linear Taylor term
$P_{M_\lambda,\tau}y$ in \eqref{eq:hamilton-normalized-chart}, then
divide its inverse by its trace to remove the scalar factor involving
$\tau$. We fix \(p=\Gamma_N(\phi)\) and let \(H_\phi\) be the Euclidean Hessian of
\(P_{N,\lambda}\) restricted to the plane normal to \(\Gamma_N\), written in
the orthonormal basis \((e_r(\phi),e_z)\) of that plane.

By \eqref{eq:hamilton-normalized-chart} and \eqref{eq:hamilton-zero-jet}, the
remainder \(\widetilde v\) contributes neither a value nor a first derivative
at \(y=0\), so \(D_yv(0,\zeta)=P_{M_\lambda,\tau}(\zeta)\), and the same
normalization gives \(D_\eps v(0,\zeta)=Le_y\).  Hence, on the axis, we have
\begin{equation}
 D_yX_{N,\lambda}
 =R_\phi\bigl(a(\tau)\iota M_\lambda+e_y\otimes\tau\bigr),
 \qquad X_\phi=\mathscr R e_\phi(\phi),               \label{eq:geometry-axis-first-jet}
\end{equation}
where each function of \(\zeta\) is evaluated at \(N\phi\).

Let \(w=q_0^1e_r(\phi)+q_0^2e_z=R_\phi\iota q_0\) be a normal vector at
\(p\), with coordinates \(q_0\in\R^2\) in the basis above.  We realize \(w\)
as a derivative of \(X_{N,\lambda}\) by combining a disk variation with an
angle variation, since the first identity in \eqref{eq:geometry-axis-first-jet}
alone also produces a tangential component.  For \(\xi\in\R^2\) and
\(s\in\R\), we have
\begin{equation*}
 D_yX_{N,\lambda}\,\xi+sX_\phi
 =R_\phi\bigl(a(\tau)\iota M_\lambda\xi
   +(\tau\cdot\xi+s\mathscr R)e_y\bigr),
\end{equation*}
so the choice
\[
 \xi=(a(\tau)M_\lambda)^{-1}q_0,
 \qquad s=-\frac{\tau\cdot\xi}{\mathscr R}
\]
gives \(w\).  For two normal vectors \(w,w'\) at \(p\), we let \(Z=(\xi,s)\)
and \(Z'=(\xi',s')\) be the corresponding directions in \((y,\phi)\), and we
set \(F=P_{N,\lambda}\circ X_{N,\lambda}=\psi_a-|y|^2\).  Differentiating
twice, we obtain
\[
 \pa_Z\pa_{Z'}F
 =\operatorname{Hess}P_{N,\lambda}(w,w')
   +dP_{N,\lambda}\bigl(\pa_Z\pa_{Z'}X_{N,\lambda}\bigr),
\]
where the second term vanishes because \(dP_{N,\lambda}=0\) on the axis,
while \(F\) is a quadratic function of \(y\) alone, so that
\(\pa_Z\pa_{Z'}F=-2\xi\cdot\xi'\).  Therefore, we have
\(\operatorname{Hess}P_{N,\lambda}(w,w')
=-2a(\tau)^{-2}(M_\lambda^{-1}q_0)\cdot(M_\lambda^{-1}q_0')\).  In the
basis \((e_r(\phi),e_z)\), this can be written as
\begin{equation}
 H_\phi=-2a(\tau)^{-2}M_\lambda^{-T}M_\lambda^{-1}.
                                                               \label{eq:geometry-axis-Hessian}
\end{equation}

We normalize the inverse Hessian to remove \(a(\tau)^{-2}\).  Since
\(\tr(M_\lambda^TM_\lambda)=2\) by \eqref{eq:fixed-M-invariants}, we have
\((-H_\phi)^{-1}=\frac12a(\tau)^2M_\lambda M_\lambda^T\) and
\(\tr((-H_\phi)^{-1})=a(\tau)^2\), so we obtain the \emph{normalized inverse Hessian}
\begin{equation}
 \cS_\phi
 :=\frac{(-H_\phi)^{-1}}{\tr((-H_\phi)^{-1})}
 =\frac12M_\lambda(N\phi)M_\lambda(N\phi)^T.          \label{eq:geometry-axis-tensor}
\end{equation}
Thus, \(\cS_\phi\) is independent of \(\tau\) and
\(\widetilde v\). It is also unchanged when \(P_{N,\lambda}\) is replaced by \(\kappa^2P_{N,\lambda}+C\),
which multiplies \(H_\phi\) by \(\kappa^2\).  Since \(M_\lambda\) is
symmetric, \eqref{eq:fixed-M-lambda} gives
\[
 \cS_\phi=\frac12\mathsf R_{\alpha_\lambda(N\phi)}
   \begin{pmatrix}1+\rho&0\\0&1-\rho\end{pmatrix}
   \mathsf R_{-\alpha_\lambda(N\phi)},
\]
whose eigenvalues \((1\pm\rho)/2\) are distinct because \(\rho=\rho_*>0\).
The eigenline of the larger eigenvalue is therefore well defined, with angle
\(\alpha_\lambda(N\phi)\) modulo \(\pi\) in the moving frame
\((e_r(\phi),e_z)\).

\subsection{Computing the stabilizer}
\label{subsec:geometry-stabilizer}

We work with the extended stabilizer \(\Stab_\pm\) of
\eqref{eq:stabilizer-pm}, whose elements preserve \(\Omega_{N,\lambda}\) and
\(P_{N,\lambda}\) and send \(B_{N,\lambda}\) to \(\pm B_{N,\lambda}\).  The argument uses the pressure alone and therefore treats both signs
at once.

Let \(g\in\Stab_\pm(\Omega_{N,\lambda},B_{N,\lambda},P_{N,\lambda})\).  Since
\(P_{N,\lambda}\circ g^{-1}=P_{N,\lambda}\), the map \(g\) preserves the
critical set of the pressure, which by \cref{subsec:geometry-tensor} is the
round circle \(\Gamma_N\).  An isometry preserving a round circle fixes its
center, which here is the origin, so \(g\) is linear.  It preserves the affine
plane of the circle, giving \(ge_z=\chi e_z\) for \(\chi\in\{\pm1\}\), and on
that plane it is orthogonal, so that in the angular variable \(g\) acts by
\begin{equation}
 \phi\longmapsto\sigma\phi+c,\qquad
 e_z\longmapsto\chi e_z,\qquad
 \sigma,\chi\in\{\pm1\}.                              \label{eq:geometry-circle-action}
\end{equation}

We now use \eqref{eq:geometry-circle-action} to obtain a condition on
\(\alpha_\lambda\).  Since \(g\Gamma_N(\phi)=\Gamma_N(\sigma\phi+c)\), we
have \(ge_r(\phi)=e_r(\sigma\phi+c)\), so \(g\) carries the normal plane at
\(\Gamma_N(\phi)\) to the normal plane at \(\Gamma_N(\sigma\phi+c)\).  In the
bases \((e_r(\phi),e_z)\) and \((e_r(\sigma\phi+c),e_z)\), it acts by
\(D=\diag(1,\chi)\).  We obtain \(H_{\sigma\phi+c}=DH_\phi D\) from the
invariance of the pressure.  Since \(\tr(DSD)=\tr(S)\), the normalization in
\eqref{eq:geometry-axis-tensor} is preserved by conjugation with
\(D\), giving \(\cS_{\sigma\phi+c}=D\cS_\phi D\).  Finally,
\(D\mathsf R_\alpha D=\mathsf R_{\chi\alpha}\), so conjugation by \(D\)
replaces \(\alpha_\lambda\) by \(\chi\alpha_\lambda\) in
\eqref{eq:geometry-axis-tensor}.  Since the eigenvalues are distinct,
equality of these tensors is equivalent to equality of their principal
eigenlines.  With \(\zeta=N\phi\) and \(c_0=Nc\), this gives
\begin{equation}
 \alpha_\lambda(\sigma\zeta+c_0)
 \equiv\chi\alpha_\lambda(\zeta)\pmod\pi.
                                       \label{eq:geometry-angle-congruence}
\end{equation}

The following lemma determines when this relation can hold. Recall
from \eqref{eq:fixed-alpha-lambda} that
\[
 \alpha_\lambda(\zeta)
 =\alpha_0+\delta_*(\cos\zeta+\lambda\sin2\zeta).
\]
We allow two values of $\lambda$ so that the lemma also applies when
comparing distinct equilibria.

\begin{lemma}
\label{lem:geometry-congruence}
Let \(\lambda,\lambda'\in I\), let \(\sigma,\chi\in\{\pm1\}\), and let
\(c_0\in\R\).  If
\[
 \alpha_{\lambda'}(\sigma\zeta+c_0)
 \equiv\chi\alpha_\lambda(\zeta)\pmod\pi
 \qquad\text{for every }\zeta\in\R,
\]
then \(\chi=1\), \(c_0\in2\pi\Z\), \(\sigma=1\), and \(\lambda'=\lambda\).
\end{lemma}

\begin{proof}
The difference of the two sides is continuous and takes values in
\(\pi\Z\), so it is a constant \(\pi k\), \(k\in\Z\).  Averaging over a period, we are left with
\((1-\chi)\alpha_0=\pi k\), since the two terms in the brackets in
\eqref{eq:fixed-alpha-lambda} have zero mean.  If \(\chi=-1\), this says
that \(2\alpha_0\in\pi\Z\), which \eqref{eq:geometry-parameter-conditions}
excludes.  Hence, \(\chi=1\), and then \(k=0\).

With \(\chi=1\) and \(k=0\), division by \(\delta\neq0\) turns the hypothesis
into
\[
 \bigl[\cos(\sigma\zeta+c_0)-\cos\zeta\bigr]
 +\bigl[\lambda'\sin(2\sigma\zeta+2c_0)-\lambda\sin2\zeta\bigr]=0 .
\]
The terms of frequency one and frequency two must vanish separately.  The first gives \(\cos(\sigma\zeta+c_0)=\cos\zeta\) for every
\(\zeta\), which for either sign of \(\sigma\) forces \(c_0\in2\pi\Z\).  The
second then reads as \(\lambda'\sigma\sin2\zeta=\lambda\sin2\zeta\), that is,
\(\lambda'=\sigma\lambda\).  Since both parameters are positive,
\(\sigma=1\) and \(\lambda'=\lambda\).
\end{proof}

We apply \cref{lem:geometry-congruence} to
\eqref{eq:geometry-angle-congruence}, whose two parameters are equal, and
obtain
\begin{equation}
 \sigma=\chi=1,\qquad c=\frac{2\pi j}{N} ,             \label{eq:geometry-only-CN}
\end{equation}
for some integer \(j\), since \(c_0=Nc\). Thus, \(g=R_{2\pi j/N}\).
Conversely, these rotations preserve both the pressure and the field.
Indeed, \(v_{\eps,\lambda}\) is
\(2\pi\)-periodic in \(\zeta\), so \eqref{eq:geometry-X} gives
\(X_{N,\lambda}(y,\phi+2\pi/N)=R_{2\pi/N}X_{N,\lambda}(y,\phi)\).  Moreover,
the right-hand side of \eqref{eq:geometry-P} does not depend on \(\phi\), and
\(\pa_\theta\) commutes with \(R_{2\pi/N}\), so that \eqref{eq:geometry-B}
gives \((R_{2\pi/N})_*B_{N,\lambda}=B_{N,\lambda}\).  Hence, we obtain
\begin{equation}
 \Stab_\pm(\Omega_{N,\lambda},B_{N,\lambda},P_{N,\lambda})=C_N.
                                                       \label{eq:geometry-stabilizer}
\end{equation}

\subsection{Distinct equilibria and non-isolation}
\label{subsec:geometry-separation}

We fix \(N\) and suppose that the equilibria at \(\lambda,\lambda'\in I\) are
equivalent under \eqref{eq:physical-moduli-action}. Reparametrization leaves the physical fields and $\cS_\phi$ unchanged,
so it suffices to consider the similarity, amplitude rescaling, and
pressure shift.

Such an operation maps the critical set of one pressure to that of the other,
and by \eqref{eq:geometry-axis} both critical circles have radius \(NL\), so
the scale factor is one and we are left with an isometry \(g\).  As in
\cref{subsec:geometry-stabilizer}, \(g\) acts on the axis by
\eqref{eq:geometry-circle-action}. The normalized inverse Hessian
$\cS_\phi$ in \eqref{eq:geometry-axis-tensor} is unchanged by the
amplitude rescaling and the pressure shift. Applying the transformation
rule from \cref{subsec:geometry-stabilizer} gives
\begin{equation}
 \alpha_{\lambda'}(\sigma\zeta+c_0)
 \equiv\chi\alpha_\lambda(\zeta)\pmod\pi,             \label{eq:geometry-two-lambda-congruence}
\end{equation}
with \(c_0=Nc\) as before.  Hence, \cref{lem:geometry-congruence} gives
\(\lambda'=\lambda\), so that distinct parameters give inequivalent
equilibria under \eqref{eq:physical-moduli-action}.

We recall the topology fixed in \cref{ss:moduli}.  For an integer \(k\ge3\), a
configuration is a triple
\begin{equation}
 (X,b,p)\in\Emb^k(Q,\R^3)\times C^{k-1}(Q,\R^3)\times C^k(Q),
 \qquad b=B\circ X,\quad p=P\circ X,                  \label{eq:geometry-configuration-space}
\end{equation}
and we identify two configurations differing by a reparametrization
\begin{equation}
 (X,b,p)\longmapsto(X\circ\Phi,b\circ\Phi,p\circ\Phi),
 \qquad\Phi\in\Diff^k(Q),                             \label{eq:geometry-reparametrization}
\end{equation}
or by \eqref{eq:physical-moduli-action}.  We write
\(\mathscr M^k_{\rm emb}(Q)\) for the quotient, in which the embedded domain
is part of the data, so the domains may move with \(\lambda\).

By \cref{thm:hamilton-one-chart}, for each
\(\eps=N^{-1}\) the solution depends smoothly on \(\lambda\) and so does the reconstruction
\eqref{eq:geometry-X}--\eqref{eq:geometry-B}. For each \(k\), derivatives of order at most \(k\) of \(W_\gamma^{-1}\) are bounded by a
polynomial in \(\langle n\rangle\) times
\(e^{-(\sigma_0-\gamma)\langle n\rangle}\). Thus, removal of the weight
is continuous into every fixed Sobolev space.  Sobolev embedding controls each \(C^k\) norm by
a sufficiently high Sobolev index. Composition with \(\zeta=N\phi\)
is continuous for fixed \(N\). Hence, the curve
\[
 \lambda\longmapsto
 [X_{N,\lambda},B_{N,\lambda}\circ X_{N,\lambda},
                P_{N,\lambda}\circ X_{N,\lambda}]
 \in\mathscr M^k_{\rm emb}(Q)
\]
is continuous, also for the smooth quotient topology of
\cref{ss:moduli}. It is injective by
\eqref{eq:geometry-two-lambda-congruence}.  For any interior $\lambda\in I$, choose distinct
$\lambda_j\in I$ with $\lambda_j\to\lambda$. The corresponding
classes are distinct by injectivity and converge to the class at
$\lambda$ by continuity. Thus, this equilibrium is not isolated in
$\mathscr M^k_{\rm emb}(Q)$.

\subsection{Proof of the main theorem}
\label{subsec:geometry-closure}
We complete the proof of \cref{thm:geometry-closure} and then deduce
\cref{thm:main}.
\begin{proof}[Proof of \cref{thm:geometry-closure}]
The preceding arguments required the four strict inequalities
\[
 E(\eta,\eps)<\frac L2\sqrt{1-\rho^2},\qquad
 \eta<\sqrt{1-\rho},\qquad
 \mathscr R>2\norm{v}_{C^0},\qquad
 C_{\mathrm{tub}}\eta<1 ,
\]
respectively in \eqref{eq:geometry-J-margin},
\eqref{eq:geometry-B-lower}, the construction of
\eqref{eq:geometry-h}--\eqref{eq:geometry-nu}, and
\eqref{eq:geometry-Psi-lower}.  All four hold after one uniform decrease of
\(\eps\), since \(\eta\le C\eps\) by \eqref{eq:geometry-q-C2}, the constants
\(m_0\), \(m_1\), and \(C_{\mathrm{tub}}\) are independent of \(N\), and
\(\norm{v}_{C^0}\le m_0+\eta\) while \(\mathscr R=L/\eps\).  We let
\(\eps_{\mathrm{geom}}\le\eps_H\) be the resulting threshold, uniform on
\(I\) because every constant entering it is.  Then,
\cref{subsec:geometry-jacobian,subsec:geometry-injectivity} give the
embedding of the closed solid torus, and
\cref{subsec:geometry-reconstruction} gives \eqref{eq:geometry-MHS} with
\eqref{eq:geometry-boundary-tangent}.  Moreover, \cref{subsec:geometry-tensor}
gives the nested tori \eqref{eq:geometry-pressure-tori} and the round axis
\eqref{eq:geometry-axis}, on which the field vanishes and nowhere else by
\eqref{eq:geometry-B-lower}.  We computed the stabilizer in
\eqref{eq:geometry-stabilizer}.  Finally, the branch is smooth in
\(\lambda\) by \cref{thm:hamilton-one-chart}, while
\eqref{eq:geometry-two-lambda-congruence} and
\cref{lem:geometry-congruence} make the \(\lambda\)-family nonconstant
modulo reparametrizations, similarities, and the amplitude action.
\end{proof}

We set
\begin{equation}
 N_0=1+\left\lceil
 \max\{2,\eps_H^{-1},\eps_{\mathrm{geom}}^{-1}\}
 \right\rceil.                                        \label{eq:geometry-N0}
\end{equation}
Since \(\eps_{\mathrm{geom}}\le\eps_H\), this is at least the
\(N_{\mathrm{an}}\) of \eqref{eq:hamilton-Nan}, so
\cref{cor:hamilton-all-integers} and \cref{thm:geometry-closure} apply to
every integer \(N\ge N_0\).

\begin{proof}[Proof of \cref{thm:main}]
The orientation \(\alpha_0\) and the ellipse parameter \(\rho_*\) were fixed in
\cref{sec:block,sec:inverse}, the amplitude \(\delta_*\) was fixed in
\cref{sec:hamilton}, and the interval \(I\subset(0,1/2)\) was fixed in
\cref{thm:hamilton-one-chart}. The integer \(N_0\) is defined by
\eqref{eq:geometry-N0}. For
\(N\ge N_0\) and \(\lambda\in I\), \cref{cor:hamilton-all-integers} permits
\(\eps=N^{-1}\), and \cref{thm:geometry-closure} supplies the embedding
\(X_{N,\lambda}\) and the smooth fields on \(\Omega_{N,\lambda}\), which
solve \eqref{eq:mhs} by \eqref{eq:geometry-MHS} and
\eqref{eq:geometry-boundary-tangent}.
Equations \eqref{eq:geometry-pressure-tori} and
\eqref{eq:geometry-axis} give the nested pressure tori and the round axis
\eqref{eq:main-axis}, proving (i). The field vanishes precisely on this
axis by \eqref{eq:geometry-B-lower} and \eqref{eq:geometry-B}, which
proves (ii). Its extended stabilizer is exactly the group of rotations
$R_{2\pi j/N}$ by \eqref{eq:geometry-only-CN} and
\eqref{eq:geometry-stabilizer}. This proves (iii), including the absence
of continuous symmetries and orientation-reversing elements. Finally,
\cref{thm:hamilton-one-chart} gives smooth dependence on $\lambda$, and
\cref{subsec:geometry-separation} proves that the curve
\eqref{eq:main-moduli-curve} is continuous and injective in the quotient.
Thus, the equilibria form a nontrivial family and are not isolated at
interior parameter values, proving (iv).
\end{proof}

\section{Lean formalization}\label{app:lean}

The main result of this paper, Theorem~\ref{thm:main}, has been verified in Lean~4 using its mathematical library Mathlib. The objective of this section is to explain the correspondence between the mathematical statement and its Lean formulation without assuming broad familiarity with Lean. We repeat the relevant mathematical definitions alongside their formal counterparts. The Lean code is available at the accompanying repository:
\begin{center}   \url{https://github.com/lukasliehr/Grad-Conjecture}.
\end{center}

In this repository, the file \lean{Showcase.lean} contains the main theorem in Lean. It is self-contained, imports only Mathlib, and defines every additional notion needed to state the result. The companion file \lean{Showcase_WithProofs.lean} proves the same statements using the verified library. The distinction between the presentation and the proofs is explained in Subsection~\ref{app:lean-verification}.

\subsection{Main theorem}\label{app:lean-result}

In this section, we recall the main theorem of the present paper and state its corresponding formalization in Lean. To do so, we recall the reference domain and the representatives used to compare equilibria on different embedded domains:
\[
 \mathbb T=\mathbb R/(2\pi\mathbb Z),\qquad
 Q=\overline{\mathbb D^2}\times\mathbb T,\qquad
 (X,b,p)=(X,B\circ X,P\circ X).
\]
Here, $X:Q\to\mathbb R^3$ is an embedding and $\Omega=X(Q)$ is the
closed solid torus. The notation $\mathscr M^k_{\mathrm{emb}}(Q)$ denotes
the space of such configurations modulo reference reparametrizations,
spatial similarities, magnetic rescaling, and pressure shifts, with the
quotient topology. Their corresponding Lean definitions are spelled out and explained in
Section~\ref{app:lean-moduli}. The theorem we verified in Lean is the following (see also Theorem \ref{thm:main}).

\begin{leanmainresult}
Fix $L>0$. There exist $0<a<b<1/2$, an integer $N_0\geq1$, and
representatives $(X_{N,\lambda},b_{N,\lambda},p_{N,\lambda})$ for
$N\geq N_0$ and $\lambda\in I=[a,b]$ with the following properties.
Each $X_{N,\lambda}$ is a smooth embedding and there are fields
$B_{N,\lambda},P_{N,\lambda}$, smooth on a neighborhood of
$\Omega_{N,\lambda}=X_{N,\lambda}(Q)$, whose pullbacks are
$b_{N,\lambda},p_{N,\lambda}$, and which satisfy
\[
 B\times\operatorname{curl}B+\nabla P=0,\qquad
 \operatorname{div}B=0\quad\hbox{in }\operatorname{int}\Omega,
 \qquad B\cdot n=0\quad\hbox{on }\partial\Omega.
\]
Suppressing $N,\lambda$ on the physical fields and domain, one has:
\begin{enumerate}
\renewcommand{\labelenumi}{(\roman{enumi})}
\item The critical set of $P$ in $\Omega$ is the circle
$
 \Gamma_N=\{NL(\cos\phi,\sin\phi,0):\phi\in\mathbb T\}
 \subset\operatorname{int}\Omega.
$
Every nonempty regular pressure level is an embedded torus, and these
levels smoothly foliate $\Omega\setminus\Gamma_N$.
\item For $x\in\Omega$, one has $B(x)=0$ if and only if $x\in\Gamma_N$.
\item The extended Euclidean stabilizer, allowing a global reversal of
$B$, consists exactly of the rotations
$R_{2\pi j/N}$, $j=0,\ldots,N-1$, about the third coordinate axis.
\item For each fixed $N\geq N_0$, the representatives depend smoothly on
$\lambda$. The curve of their classes in
$\mathscr M^k_{\mathrm{emb}}(Q)$ is continuous and injective for every
finite $k\geq3$ and for $k=\infty$. Every member with
$\lambda\in(a,b)$ is non-isolated in that moduli space.
\end{enumerate}
\end{leanmainresult}

The interval $I$ and threshold $N_0$ are common to all four conclusions.
The Lean statement also makes explicit the boundary regularity and the
boundary leaf of the pressure foliation. We use the closed body $X(Q)$
throughout; the differential equations are imposed on its interior.

The formalization in Lean of the latter theorem takes the following form:

\par
\noindent\begin{minipage}{\linewidth}
\begin{leancode}
theorem equilibria_with_exact_cyclic_symmetry (L : ℝ) (hL : 0 < L) :
    ∃ a b : ℝ, 0 < a ∧ a < b ∧ b < 1 / 2 ∧
      ∃ N₀ : ℕ, 1 ≤ N₀ ∧ ∃ family : ℕ → Set.Icc a b → Representative,
      ∀ N : ℕ, N₀ ≤ N →
          SmoothFamily (Set.Icc a b) (family N) ∧
          (∀ lambda : Set.Icc a b, Equilibrium L N (family N lambda)) ∧
          (∀ k : ℕ, 3 ≤ k → ModuliCurve (.finite k) (Set.Icc a b) (family N)) ∧
          ModuliCurve .smooth (Set.Icc a b) (family N)
\end{leancode}
\end{minipage}
\par

In the following subsections we provide a precise description of all clauses in the previous Lean theorem.

\subsection{The reference domain and regularity}\label{app:lean-domain}

We recall that Lean records the type of each object. For example, writing \lean{L : ℝ} means mathematically that $L\in\mathbb R$, and \lean{X : A → B} denotes a
map $X$ from $A$ to $B$. A definition is
introduced by \lean{def} or \lean{abbrev}. For reading the following Lean statements, an abbreviation can be treated as a definition. We begin with the types corresponding to $\mathbb R^3$, $\mathbb R^2$, $\mathbb T$,
$\overline{\mathbb D^2}$, $Q$, and $\mathbb T^2$. In Lean, we define them as follows:

\par
\noindent\begin{minipage}{\linewidth}
\begin{leancode}
abbrev Vec := EuclideanSpace ℝ (Fin 3)
abbrev Plane := EuclideanSpace ℝ (Fin 2)
abbrev PeriodicCircle := AddCircle (2 * Real.pi)
abbrev ClosedDisk := {x : Plane // ‖x‖ ≤ 1}
abbrev ReferenceDomain := ClosedDisk × PeriodicCircle
abbrev Torus := PeriodicCircle × PeriodicCircle
\end{leancode}
\end{minipage}
\par

The type \lean{EuclideanSpace ℝ (Fin 3)} represents $\mathbb R^3$
endowed with its Euclidean norm, where \lean{Fin 3} is the index set $\{0,1,2\}$. The type \lean{AddCircle (2 * Real.pi)} is the quotient $\mathbb R/(2\pi\mathbb Z)$. The notation \lean{x : Plane // ‖x‖ ≤ 1} restricts the points of the plane to those satisfying $\|x\|\leq1$. Thus, \lean{ClosedDisk} is precisely the closed unit disk.

% For derivatives on $Q$, we use the covering cylinder
% \[
%  C=\overline{\mathbb D^2}\times\mathbb R,\qquad
%  \pi_C(y,z)=(y,[z])\in Q.
% \]
% The definitions \lean{planarPart}, \lean{cylinder}, and
% \lean{quotientPoint} implement $(y,z)\mapsto y$, the set $C$, and
% the map $\pi_C$, respectively. For $f:Q\to E$, the function
% \lean{periodicLift f} agrees with $f\circ\pi_C$ on $C$ and is
% assigned the value zero outside $C$. This outside value is only a way to
% represent the lift as a function on all of $\mathbb R^3$. Regularity is
% required through local extensions from $C$, and derivatives on $C$ are
% taken within $C$.

For derivatives on $Q$, we use the covering cylinder
$
 C=\overline{\mathbb D^2}\times\mathbb R
$
and the projection
$
 \pi_C:C\longrightarrow Q
$
defined by $\pi_C(y,z)=(y,[z])$, where $[z]$ denotes the class of $z$ in
$\mathbb T=\mathbb R/(2\pi\mathbb Z)$.
The projection onto the first two coordinates and the cylinder $C$
are defined by

\begin{leancode}
def planarPart (x : Vec) : Plane := WithLp.toLp 2 ![x 0, x 1]
def cylinder : Set Vec := {x | ‖planarPart x‖ ≤ 1}
\end{leancode}

In the above code, \lean{x 0}, \lean{x 1}, and \lean{x 2} are the three
coordinates of $x\in\mathbb R^3$. The expression
\lean{![x 0, x 1]} collects its first two coordinates, and
\lean{WithLp.toLp 2} regards this pair as an element of the
Euclidean plane. Thus, \lean{planarPart} is the map
$(y,z)\mapsto y$, and \lean{cylinder} is precisely the set
$\{(y,z):\|y\|\leq1\}$. The covering projection $\pi_C$ is represented in Lean by

\begin{leancode}
def quotientPoint (x : Vec) (hx : x ∈ cylinder) :
    ReferenceDomain := (⟨planarPart x, hx⟩, (x 2 : PeriodicCircle))
\end{leancode}

The argument \lean{hx} is a proof that $x\in C$, or equivalently
that $\|\operatorname{planarPart}(x)\|\leq1$.
The expression \lean{⟨planarPart x, hx⟩} therefore gives a point
of \lean{ClosedDisk}. The second component,
\lean{(x 2 : PeriodicCircle)}, takes the third coordinate modulo
$2\pi$. Together, these components give a point of
$Q=\overline{\mathbb D^2}\times\mathbb T$.

For $f:Q\to E$, where $E$ has a specified zero element, define
\[
 \widetilde f(x)=
 \begin{cases}
   f(\pi_C(x)),&x\in C,\\
   0,&x\notin C.
 \end{cases}
\]
Thus, $\widetilde f$ is a function on $\mathbb R^3$ whose
restriction to $C$ is the lift of $f$. In Lean, we define the periodic lift as follows:

\begin{leancode}
def periodicLift {E : Type*} [Zero E]
    (f : ReferenceDomain → E) (x : Vec) : E := by
  classical
  exact if hx : x ∈ cylinder then
    f (quotientPoint x hx)
  else
    0
\end{leancode}

The notation \lean{E : Type*} allows the target $E$ to be
inferred from $f$, and \lean{[Zero E]} specifies its zero element.
The keywords \lean{by} and \lean{exact} introduce a construction
and give its resulting value.
The command \lean{classical} permits the use of classical logic.
Here it allows the case distinction $x\in C$ or $x\notin C$,
without requiring a procedure/algorithm that decides membership in $C$.
In \lean{if hx : x ∈ cylinder then}, the name \lean{hx}
records the proof of membership in the first case.
This proof is passed to \lean{quotientPoint}, since $\pi_C(x)$
is defined only for $x\in C$. The two branches therefore give
exactly the two values in the displayed formula.
The values outside $C$ serve only to represent the lift as a
function on all of $\mathbb R^3$. The zero extension need not be
continuous across the boundary of $C$. Regularity on $Q$ is
expressed through local extensions of $\widetilde f|_C$ to open
neighborhoods in $\mathbb R^3$, with the required differentiability.
Derivatives are taken within $C$ and therefore do not depend on
the assigned exterior values.

We encode the regularity of a configuration by a type
\lean{Regularity}, whose values are \lean{finite k}, for
$k\in\mathbb N$, and \lean{smooth}. These represent $C^k$ and
$C^\infty$ regularity, respectively. The keyword \lean{inductive}
declares the type by listing these alternatives and each line beginning with \lean{|} introduces one of them.

\par
\noindent\begin{minipage}{\linewidth}
\begin{leancode}
inductive Regularity where
  | finite (k : ℕ) 
  | smooth        
\end{leancode}
\end{minipage}
\par

We associate two operations with these regularity labels.
The first assigns the differentiability order used by Mathlib.
The second lowers finite regularity by one derivative, as needed
for the magnetic component of a configuration, while preserving
smooth regularity.

\par
\noindent\begin{minipage}{\linewidth}
\begin{leancode}
def Regularity.order : Regularity → ℕ∞ω
  | .finite k => k
  | .smooth => ∞
def Regularity.predecessor : Regularity → Regularity
  | .finite k => .finite (k - 1)
  | .smooth => .smooth
\end{leancode}
\end{minipage}
\par

Both functions are defined by cases: each \lean{|} introduces
a possible input, and \lean{=>} specifies the corresponding
output. The notation \lean{.finite k} abbreviates
\lean{Regularity.finite k}, and \lean{.smooth} abbreviates
\lean{Regularity.smooth}.
The type \lean{ℕ∞ω} contains the natural differentiability
orders, the smooth order $\infty$, and an additional order
$\omega$ for analyticity. Only the natural orders and $\infty$
are used here. Thus, \lean{Regularity.order} sends
\lean{.finite k} to $k$ and \lean{.smooth} to $\infty$.
The operation \lean{Regularity.predecessor} replaces $C^k$
by $C^{k-1}$ for $k\geq1$ and leaves $C^\infty$ unchanged.
At $k=0$ it returns $C^0$, since subtraction in $\mathbb N$
is truncated at zero.

Let $E,F$ be real normed vector spaces, let $S\subset E$, and let
$f:E\to F$. We express regularity on $S$ through local extensions:
for every $x\in S$, there must be an open neighborhood $U$ of $x$
and a $C^r$ map $g$ on $U$ such that
$
 g(y)=f(y)
$
for all $y\in U\cap S$. In our application, $S$ is the closed covering cylinder $C$.
This formulation includes regularity at its boundary. The
corresponding Lean definition is

\par
\noindent\begin{minipage}{\linewidth}
\begin{leancode}
def HasLocalExtensions {E F : Type*}
    [NormedAddCommGroup E] [NormedSpace ℝ E]
    [NormedAddCommGroup F] [NormedSpace ℝ F]
    (r : Regularity) (f : E → F) (S : Set E) : Prop :=
      ∀ x ∈ S, ∃ U : Set E, IsOpen U ∧ x ∈ U ∧ ∃ g : E → F,
      ContDiffOn ℝ r.order g U ∧ EqOn g f (U ∩ S)
\end{leancode}
\end{minipage}
\par

The first lines specify the spaces and their structures.
The declaration \lean{NormedAddCommGroup E} provides addition,
a zero element, additive inverses, and a norm on $E$; addition
is commutative, and the distance is given by $\|x-y\|$.
The declaration \lean{NormedSpace ℝ E} provides real scalar
multiplication compatible with these operations and with the
norm.
Together, these declarations express that $E$ is a real normed
vector space. The corresponding declarations for $F$ have the
same meaning.
% Square brackets indicate structures that Lean
% can infer from the spaces being used. For example, the Euclidean spaces introduced above already carry them.
The remaining arguments specify the regularity label $r$,
the function $f$, and the set $S$. The notation \lean{S : Set E}
means that $S$ is a subset of $E$, and \lean{: Prop}
indicates that \lean{HasLocalExtensions r f S} is a mathematical
assertion about these data.
The body of the definition follows the stated extension condition.
The symbols \lean{∀}, \lean{∃}, and \lean{∧} have the usual meaning. Thus, \lean{IsOpen U ∧ x ∈ U} says that $U$ is an open
neighborhood of $x$. The expression
\lean{ContDiffOn ℝ r.order g U} states that $g$ is $C^r$
on $U$, with derivatives taken over $\mathbb R$, and
\lean{EqOn g f (U ∩ S)} states that $g=f$ on $U\cap S$.
Although $g$ is represented as a function on all of $E$,
its regularity is required only on $U$.

For a function $f:Q\to E$, we apply this condition to its
periodic lift $\widetilde f$ on the cylinder $C$:

\par
\noindent\begin{minipage}{\linewidth}
\begin{leancode}
def HasRegularity {E : Type*}
    [NormedAddCommGroup E] [NormedSpace ℝ E]
    (r : Regularity) (f : ReferenceDomain → E) : Prop :=
      HasLocalExtensions r (periodicLift f) cylinder
\end{leancode}
\end{minipage}
\par

Thus, \lean{HasRegularity r f} means that the lift of $f$
admits local $C^r$ extensions near every point of $C$,
including its boundary. Since $\pi_C:C\to Q$ is a covering
map, this gives $C^r$ regularity on the reference solid torus
in local coordinates.

We also require the position map $X:Q\to\mathbb R^3$ to be
an embedding. At differentiability orders at least one, this
means that $X$ is $C^r$, is a homeomorphism onto its image,
and has injective derivative at every point. In covering
coordinates, these requirements are expressed by

\par
\noindent\begin{minipage}{\linewidth}
\begin{leancode}
def IsEmbeddingOfRegularity
    (r : Regularity) (X : ReferenceDomain → Vec) : Prop :=
      HasRegularity r X ∧ Topology.IsEmbedding X ∧ ∀ x ∈ cylinder,
      Function.Injective (fderivWithin ℝ (periodicLift X) cylinder x)
\end{leancode}
\end{minipage}
\par

The three conditions have distinct roles.
The first supplies the $C^r$ regularity of $X$.
The second, \lean{Topology.IsEmbedding X}, says that $X$
is a homeomorphism from $Q$ onto its image, equipped with
the topology inherited from $\mathbb R^3$.
The third is the immersion condition:
\lean{Function.Injective} requires the derivative to be
injective.
%, or equivalently to have rank three.
More precisely,
\lean{fderivWithin ℝ (periodicLift X) cylinder x}
is the Fr\'echet derivative of $\widetilde X$ at $x$,
computed using points of $C$. It is a continuous real-linear map
$
 D_C\widetilde X(x):\mathbb R^3\longrightarrow\mathbb R^3.
$
At a boundary point, this is the derivative of any local
$C^1$ extension. It is independent of the chosen extension:
two such extensions agree on the interior of $C$, hence have
the same derivatives there, and continuity gives equality of
their derivatives at the boundary. Thus, the definition imposes
the full immersion condition at boundary points as well as
at interior points.
\subsection{Configurations and their moduli space}\label{app:lean-moduli}
The representative $(X,b,p)=(X,B\circ X,P\circ X)$ is stored as three
functions on $Q$. The keyword \lean{structure} collects these functions
into one object.

\par
\noindent\begin{minipage}{\linewidth}
\begin{leancode}
structure Representative where
  position : ReferenceDomain → Vec 
  magnetic : ReferenceDomain → Vec 
  pressure : ReferenceDomain → ℝ   
def IsConfiguration (r : Regularity) (c : Representative) : Prop :=
  IsEmbeddingOfRegularity r c.position ∧
    HasRegularity r.predecessor c.magnetic ∧ HasRegularity r c.pressure
abbrev Configuration (r : Regularity) := {c : Representative // IsConfiguration r c}
\end{leancode}
\end{minipage}
\par

For a representative \lean{c}, the fields \lean{c.position},
\lean{c.magnetic}, and \lean{c.pressure} are $X,b,p$, respectively.
In particular, the last two are the \emph{pulled-back} fields, not the
ambient functions $B,P$. For finite $k\geq3$, the configuration space is
\[
 \operatorname{Emb}^k(Q,\mathbb R^3)
 \times C^{k-1}(Q,\mathbb R^3)\times C^k(Q,\mathbb R).
\]
The subtype \lean{Configuration r} attaches the required regularity
conditions to the triple. It does not itself impose the MHS equations.

We equip this space with its usual product topology. The implementation
uses the compact fundamental cylinder
$C_0=\overline{\mathbb D^2}\times[0,2\pi]$, called
\lean{fundamentalCylinder}. The map \lean{jets r f} records the
derivatives $D^j\widetilde f|_{C_0}$ for every permitted order $j$;
\lean{Regularity.admits} expresses $j\leq k$ for finite regularity
and allows every $j\in\mathbb N$ for smooth regularity. Each derivative
space carries Mathlib's uniform convergence topology, \lean{UniformFun}.
Equivalently, for two maps $f,g$ the topology is generated by
\[
 d_j(f,g)=\sup_{x\in C_0}
 \bigl\|D^j\widetilde f(x)-D^j\widetilde g(x)\bigr\|,
\]
using $0\leq j\leq k$ or all finite $j$, as appropriate.
The instance \lean{configurationTopology} induces the topology from
the three jet maps, with orders $(k,k-1,k)$. Thus, the smooth case uses
the full $C^\infty$ topology.

Two configurations are equivalent if they differ by the action
\begin{equation}\label{app:lean-action}
 (X,b,p)\longmapsto
 \bigl(sOX\circ\Phi+c,\ \kappa Ob\circ\Phi,
       \ \kappa^2p\circ\Phi+c_0\bigr),
\end{equation}
where $s>0$, $\kappa\ne0$, $O\in O(3)$, $c\in\mathbb R^3$,
$c_0\in\mathbb R$, and $\Phi$ is a $C^r$ diffeomorphism of $Q$.
To express the last condition, \lean{HasRegularLocalLifts r Φ}
requires, near each $x\in C$, a $C^r$ map $F$ in Euclidean coordinates
such that $F(y)\in C$ and
$\pi_C(F(y))=\Phi(\pi_C(y))$ for $y$ in that neighborhood within $C$.
The predicate \lean{IsReparametrization} requires such lifts for
both $\Phi$ and its inverse. With this convention, the full relation is

\par
\noindent\begin{minipage}{\linewidth}
\begin{leancode}
def ModuliEquivalent (r : Regularity) (u v : Configuration r) : Prop :=
  ∃ (s κ : ℝ) (O : Vec ≃ₗᵢ[ℝ] Vec) (c : Vec) (c₀ : ℝ)
      (Φ : ReferenceDomain ≃ ReferenceDomain),
    0 < s ∧ κ ≠ 0 ∧ IsReparametrization r Φ ∧ ∀ q : ReferenceDomain,
      v.val.position q = s • O (u.val.position (Φ q)) + c ∧
      v.val.magnetic q = κ • O (u.val.magnetic (Φ q)) ∧
      v.val.pressure q = κ ^ 2 * u.val.pressure (Φ q) + c₀
\end{leancode}
\end{minipage}
\par

Here, \lean{Vec ≃ₗᵢ[ℝ] Vec} is the type of real linear isometric
isomorphisms, hence orthogonal transformations. The notation
\lean{ReferenceDomain ≃ ReferenceDomain} specifies a bijection; its
regularity and that of its inverse are imposed separately.
The symbol \lean{•} denotes scalar multiplication, and
\lean{u.val} is the representative underlying a configuration.
Notice that the scale $s$, the possibly negative magnetic factor
$\kappa$, and the pressure shift $c_0$ are all included.

The moduli space is the quotient by this relation, with its quotient
topology:

\par
\noindent\begin{minipage}{\linewidth}
\begin{leancode}
def ModuliSpace (r : Regularity) := Quot (ModuliEquivalent r)
def moduliClass (r : Regularity) : Configuration r → ModuliSpace r :=
  Quot.mk (ModuliEquivalent r)
instance moduliTopology (r : Regularity) : TopologicalSpace (ModuliSpace r) :=
  TopologicalSpace.coinduced (moduliClass r) inferInstance
\end{leancode}
\end{minipage}
\par

Mathlib's \lean{TopologicalSpace.coinduced} constructs the quotient
topology from the projection \lean{moduliClass}. In mathematical terms,
a subset of the quotient is open precisely when its inverse image under
this projection is open in the configuration space. The following
theorem ensures that equality in the quotient means exactly the action
in~\eqref{app:lean-action}:

\par
\noindent\begin{minipage}{\linewidth}
\begin{leancode}
theorem same_moduli_class_iff (r : Regularity) (u v : Configuration r) :
    moduliClass r u = moduliClass r v ↔ ModuliEquivalent r u v
\end{leancode}
\end{minipage}
\par

The symbol \lean{↔} means ``if and only if''. Lean's general quotient
\lean{Quot} may be formed from a relation without first attaching a
proof that the relation is an equivalence relation. The displayed theorem states that two configurations have the
same moduli class if and only if they differ by the reference
reparametrization and physical transformations specified in
\lean{ModuliEquivalent}.

\subsection{The equations and the boundary condition}\label{app:lean-mhs}

The physical fields are represented by functions
$B:\mathbb R^3\to\mathbb R^3$ and $P:\mathbb R^3\to\mathbb R$.
Their values away from a neighborhood of the closed body play no role.
The required smoothness is expressed by the following definition

\par
\noindent\begin{minipage}{\linewidth}
\begin{leancode}
def SmoothNear {E : Type*} [NormedAddCommGroup E] [NormedSpace ℝ E]
    (Ω : Set Vec) (f : Vec → E) : Prop :=
  ∃ U : Set Vec, IsOpen U ∧ Ω ⊆ U ∧ ContDiffOn ℝ ∞ f U
\end{leancode}
\end{minipage}
\par

This says that there is an open $U\supset\Omega$ on which the field is
$C^\infty$. The showcase defines \lean{basisVector i} as the $i$th
coordinate vector and defines \lean{gradient}, \lean{curl}, and
\lean{divergence} by their usual Cartesian formulas, using Mathlib's
Fr\'echet derivative \lean{fderiv}. Explicitly,
\[
 (\nabla P)_i=\partial_iP,\qquad
 \operatorname{curl}B=
 (\partial_2B_3-\partial_3B_2,\partial_3B_1-\partial_1B_3,
  \partial_1B_2-\partial_2B_1),\qquad
 \operatorname{div}B=\sum_{i=1}^3\partial_iB_i.
\]
The helper \lean{vector} forms a vector from its three coordinates,
and \lean{cross} is the usual cross product. Lean numbers these
coordinates $0,1,2$ rather than $1,2,3$.

Boundary tangency is expressed through curves. For a smooth surface $S$,
the condition $v\in T_xS$ means that $v=\gamma'(0)$ for a smooth curve
in $S$ with $\gamma(0)=x$. The definition used here is

\par
\noindent\begin{minipage}{\linewidth}
\begin{leancode}
def TangentTo (S : Set Vec) (x v : Vec) : Prop :=
  ∃ γ : ℝ → Vec, ContDiff ℝ ∞ γ ∧ γ 0 = x ∧ (∀ t, γ t ∈ S) ∧ fderiv ℝ γ 0 1 = v
\end{leancode}
\end{minipage}
\par

The last expression evaluates the derivative of $\gamma$ at $0$ on the
real number $1$, and hence gives $\gamma'(0)$. For a smooth boundary,
this is the usual tangency condition $B(x)\cdot n(x)=0$. Consequently,
the equations and their boundary condition are combined as follows:

\par
\noindent\begin{minipage}{\linewidth}
\begin{leancode}
def IsMHS (Ω : Set Vec) (B : Vec → Vec) (P : Vec → ℝ) : Prop :=
  (∀ x ∈ interior Ω, cross (B x) (curl B x) + gradient P x = 0 ∧ divergence B x = 0) ∧
  ∀ x ∈ frontier Ω, TangentTo (frontier Ω) x (B x)
\end{leancode}
\end{minipage}
\par

Here, \lean{interior Ω} and \lean{frontier Ω} mean
$\operatorname{int}\Omega$ and $\partial\Omega$. The first line of the
condition imposes the unforced MHS equations at every interior point;
the second imposes tangency at every boundary point.

\subsection{The axis and nested pressure surfaces}\label{app:lean-foliation}

The round magnetic axis and pressure levels are defined directly by
\[
 \Gamma_R=\{(R\cos\theta,R\sin\theta,0):\theta\in\mathbb R\},
 \qquad \Sigma_p=\{x\in\Omega:P(x)=p\}.
\]
In Lean these are \lean{roundAxis R} and \lean{pressureLevel Ω P p}:

\par
\noindent\begin{minipage}{\linewidth}
\begin{leancode}
def roundAxis (R : ℝ) : Set Vec :=
  Set.range fun θ : ℝ => vector (R * Real.cos θ) (R * Real.sin θ) 0
def pressureLevel (Ω : Set Vec) (P : Vec → ℝ) (p : ℝ) : Set Vec :=
  {x | x ∈ Ω ∧ P x = p}
def IsRegularLevel (Ω : Set Vec) (P : Vec → ℝ) (p : ℝ) : Prop :=
  ∀ x ∈ pressureLevel Ω P p, fderiv ℝ P x ≠ 0
\end{leancode}
\end{minipage}
\par

The notation \lean{Set.range} denotes the image of a map. A pressure
level is regular when $dP_x\ne0$ at each of its points. In Euclidean
space this is equivalent to $\nabla P(x)\ne0$. Nonemptiness is imposed
separately when asserting that a regular level is a torus.

The predicate \lean{IsEmbeddedTorus S} means that $S$ is the image
of a smooth embedding $Y:\mathbb T^2\to\mathbb R^3$. More precisely,
\lean{torusLift Y} is the map
$(\theta_1,\theta_2)\mapsto Y([\theta_1],[\theta_2])$ on
$\mathbb R^2$. The definition requires this lift to be smooth with
injective derivative, $Y$ to be a topological embedding, and $Y(\mathbb
T^2)=S$. These are the same embedding conditions as above, now for a
two-dimensional torus.

To specify a foliation rather than merely a collection of toroidal
levels, we require a smooth embedding
\begin{equation}\label{app:lean-foliation-map}
 F:(0,1]\times\mathbb T^2\longrightarrow\mathbb R^3,
 \qquad F((0,1]\times\mathbb T^2)=\Omega\setminus\Gamma.
\end{equation}
For every $r\in(0,1]$, the image $F(\{r\}\times\mathbb T^2)$ must
be an entire regular pressure level, and the slice $r=1$ must be exactly
$\partial\Omega$.
% The definition \lean{FoliationDomain} is
% $(0,1]\times\mathbb T^2$. Its covering domain, \lean{foliationCylinder},
% is $\{x\in\mathbb R^3:0<x_0\leq1\}$; \lean{foliationLift}
% reduces the other two coordinates modulo $2\pi$ and is assigned zero
% outside this covering domain.
The parameter domain of the foliation is
$(0,1]\times\mathbb T^2$. We represent it and its covering
domain as follows:

\par
\noindent\begin{minipage}{\linewidth}
\begin{leancode}
abbrev FoliationDomain := Set.Ioc (0 : ℝ) 1 × Torus
def foliationCylinder : Set Vec :=
  {x | x 0 ∈ Set.Ioc (0 : ℝ) 1}
\end{leancode}
\end{minipage}
\par

The notation \lean{Set.Ioc (0 : ℝ) 1} denotes the interval
$(0,1]$: the left endpoint is excluded and the right endpoint
is included. Since \lean{Torus} is
$\mathbb T\times\mathbb T$, the first line defines precisely
$(0,1]\times\mathbb T^2$. Its first coordinate labels the
leaves, and the two circle coordinates parametrize each leaf.

The set \lean{foliationCylinder} is
\[
 C_{\mathrm{fol}}
 =\{x=(x_0,x_1,x_2)\in\mathbb R^3:0<x_0\leq1\}
 =(0,1]\times\mathbb R^2.
\]
The corresponding covering projection is
\[
 \pi_{\mathrm{fol}}:C_{\mathrm{fol}}
 \longrightarrow(0,1]\times\mathbb T^2,\qquad
 \pi_{\mathrm{fol}}(x_0,x_1,x_2)
 =(x_0,[x_1],[x_2]).
\]
Thus, the first coordinate is unchanged, while the other two
are taken modulo $2\pi$.

For a parametrization
$F:(0,1]\times\mathbb T^2\to\mathbb R^3$, we use the lift
\[
 \widetilde F(x)=
 \begin{cases}
 F(x_0,[x_1],[x_2]),&x\in C_{\mathrm{fol}},\\
 0,&x\notin C_{\mathrm{fol}}.
 \end{cases}
\]
Its Lean definition is the following, with \lean{X} denoting
the map $F$:

\par
\noindent\begin{minipage}{\linewidth}
\begin{leancode}
def foliationLift (X : FoliationDomain → Vec)
    (x : Vec) : Vec := by
  classical
  exact if hx : x ∈ foliationCylinder then
    X (⟨x 0, hx⟩,
       (x 1 : PeriodicCircle),
       (x 2 : PeriodicCircle))
  else
    0
\end{leancode}
\end{minipage}
\par

In the first branch, \lean{hx} is a proof that
$0<x_0\leq1$. The expression \lean{⟨x 0, hx⟩} therefore
specifies an element of $(0,1]$. The other two entries
give the classes of $x_1,x_2$ in
$\mathbb R/(2\pi\mathbb Z)$.
As in \lean{periodicLift}, the command \lean{classical}
permits the case distinction according to membership in the
covering domain.
The exterior value zero serves only to define the lift on all
of $\mathbb R^3$. Smoothness is imposed through local extensions
from $C_{\mathrm{fol}}$, and derivatives are taken within that
set. In particular, local extensions at $x_0=1$ express
smoothness up to the outermost leaf. The endpoint $x_0=0$
does not belong to the parameter domain.
The precise foliation condition is

\par
\noindent\begin{minipage}{\linewidth}
\begin{leancode}
def IsPressureFoliation (Ω Γ : Set Vec) (P : Vec → ℝ) : Prop :=
  ∃ X : FoliationDomain → Vec,
    HasLocalExtensions .smooth (foliationLift X) foliationCylinder ∧
    Topology.IsEmbedding X ∧ Set.range X = Ω \ Γ ∧
    (∀ x ∈ foliationCylinder,
      Function.Injective
        (fderivWithin ℝ (foliationLift X) foliationCylinder x)) ∧
    (∀ r : Set.Ioc (0 : ℝ) 1, ∃ p : ℝ, IsRegularLevel Ω P p ∧
      Set.range (fun θ : Torus => X (r, θ)) = pressureLevel Ω P p) ∧
    Set.range (fun θ : Torus => X (⟨1, zero_lt_one, le_rfl⟩, θ)) = frontier Ω
\end{leancode}
\end{minipage}
\par

The local extension condition includes smoothness at the outer leaf
$r=1$. The equality with \lean{Ω \ Γ} covers every point off the
axis. The variable $r$ here labels the leaves; it need not equal the
pressure value $p$.

We combine the critical-set assertion, regularity off the axis, toroidal
level sets, and the foliation into the predicate corresponding to
clause~(i):

\par
\noindent\begin{minipage}{\linewidth}
\begin{leancode}
def NestedPressureTori (Ω Γ : Set Vec) (P : Vec → ℝ) : Prop :=
  Γ ⊆ interior Ω ∧ {x | x ∈ Ω ∧ fderiv ℝ P x = 0} = Γ ∧
    (∀ p, (pressureLevel Ω P p).Nonempty → IsRegularLevel Ω P p →
      IsEmbeddedTorus (pressureLevel Ω P p)) ∧
    (∀ x ∈ Ω \ Γ, IsRegularLevel Ω P (P x)) ∧ IsPressureFoliation Ω Γ P
\end{leancode}
\end{minipage}
\par

In this definition, the set equation says precisely
$\{x\in\Omega:dP_x=0\}=\Gamma$. In the main theorem we take
$\Gamma=\Gamma_{NL}$. Clause~(ii), namely
$B(x)=0\Longleftrightarrow x\in\Gamma_{NL}$ for $x\in\Omega$,
will appear explicitly in the combined equilibrium condition below.

\subsection{The extended stabilizer}\label{app:lean-symmetry}

For $g(x)=Ox+c$, with $O\in O(3)$ and $c\in\mathbb R^3$, membership
in the extended stabilizer means
\[
 g\Omega=\Omega,\qquad P(gx)=P(x),\qquad
 B(gx)=\sigma\,OB(x)\quad(x\in\Omega),
 \qquad \sigma\in\{1,-1\}.
\]
There is one sign $\sigma$ for the entire field. This is equivalent to
the definition using $g_*B=\pm B$ and $P\circ g^{-1}=P$, since $g$
maps $\Omega$ onto itself. In Lean it is written as follows:

\par
\noindent\begin{minipage}{\linewidth}
\begin{leancode}
def SignedStabilizes (Ω : Set Vec) (B : Vec → Vec) (P : Vec → ℝ)
    (O : Vec ≃ₗᵢ[ℝ] Vec) (c : Vec) : Prop :=
  let g := fun x => O x + c
  g '' Ω = Ω ∧ (∀ x ∈ Ω, P (g x) = P x) ∧
    ∃ σ : ℝ, (σ = 1 ∨ σ = -1) ∧ ∀ x ∈ Ω, B (g x) = σ • O (B x)
\end{leancode}
\end{minipage}
\par

Here \lean{g '' Ω} denotes the image $g(\Omega)$. The orthogonal map $O$ is unrestricted, so
orientation-reversing isometries are included in the test.

The helper \lean{rotation θ} is the map
\[
 R_\theta(x_1,x_2,x_3)=
 (\cos\theta\,x_1-\sin\theta\,x_2,
  \sin\theta\,x_1+\cos\theta\,x_2,x_3).
\]
The exact stabilizer assertion in clause~(iii) is
\[
 \operatorname{Stab}_{\pm}(\Omega,B,P)
 =\{R_{2\pi j/N}:0\leq j<N\}\cong C_N.
\]
It is formalized in Lean as follows:

\par
\noindent\begin{minipage}{\linewidth}
\begin{leancode}
def ExactCyclicSymmetry (N : ℕ) (Ω : Set Vec) (B : Vec → Vec) (P : Vec → ℝ) : Prop :=
  ∀ (O : Vec ≃ₗᵢ[ℝ] Vec) (c : Vec), SignedStabilizes Ω B P O c ↔
    ∃ j : ℕ, j < N ∧ ∀ x : Vec, O x + c = rotation (2 * Real.pi * j / N) x
\end{leancode}
\end{minipage}
\par

Thus every signed symmetry is one of the stated rotations, and every
stated rotation is a signed symmetry. The equality in the last line is
required on all of $\mathbb R^3$, so it identifies the Euclidean motion
itself. This finite group contains neither a nontrivial continuous group
of Euclidean symmetries nor an orientation-reversing symmetry.

\subsection{The equilibrium and its family}\label{app:lean-families}

For a representative $c=(X,b,p)$, put $\Omega=X(Q)$ and
$\Gamma=\Gamma_{NL}$. The predicate \lean{Equilibrium L N c}
collects smooth embeddedness, the physical fields and their pullback
identities, the MHS equations, and all of clauses~(i)--(iii):

\par
\noindent\begin{minipage}{\linewidth}
\begin{leancode}
def Equilibrium (L : ℝ) (N : ℕ) (c : Representative) : Prop :=
  let Ω := Set.range c.position
  let Γ := roundAxis (N * L)
  IsConfiguration .smooth c ∧ ∃ (B : Vec → Vec) (P : Vec → ℝ),
    SmoothNear Ω B ∧ SmoothNear Ω P ∧
    (∀ q : ReferenceDomain,
      B (c.position q) = c.magnetic q ∧ P (c.position q) = c.pressure q) ∧
    IsMHS Ω B P ∧ NestedPressureTori Ω Γ P ∧
    (∀ x ∈ Ω, B x = 0 ↔ x ∈ Γ) ∧ 
    ExactCyclicSymmetry N Ω B P
\end{leancode}
\end{minipage}
\par

The notation \lean{let} introduces the local abbreviations $\Omega$
and $\Gamma$. The two pullback identities are
\[
 B(X(q))=b(q),\qquad P(X(q))=p(q)\qquad(q\in Q).
\]
They ensure that the physical fields in the equations and stabilizer
are those represented by $c$. In particular, the definition does not
allow an unrelated choice of fields when forming the moduli curve.
The equivalence \lean{B x = 0 ↔ x ∈ Γ} is the magnetic-axis
assertion in clause~(ii). Thus, \lean{Equilibrium} has a stronger
meaning here than merely being a solution of the MHS equations.

It remains to explain the smooth family and non-isolation in
clause~(iv). A family $c_\lambda=(X_\lambda,b_\lambda,p_\lambda)$
on $I=[a,b]$ is smooth when it extends to an open neighborhood
$U\supset I$ and its three components are jointly smooth in the parameter
and covering coordinates, including the boundary of $Q$. This is
expressed by

\par
\noindent\begin{minipage}{\linewidth}
\begin{leancode}
def SmoothFamily (I : Set ℝ) (c : I → Representative) : Prop :=
  ∃ U : Set ℝ, IsOpen U ∧ I ⊆ U ∧ ∃ e : ℝ → Representative,
    (∀ lambda : I, e lambda.val = c lambda) ∧
    HasLocalExtensions .smooth
      (fun z : ℝ × Vec => periodicLift (e z.1).position z.2) (U ×ˢ cylinder) ∧
    HasLocalExtensions .smooth
      (fun z : ℝ × Vec => periodicLift (e z.1).magnetic z.2) (U ×ˢ cylinder) ∧
    HasLocalExtensions .smooth
      (fun z : ℝ × Vec => periodicLift (e z.1).pressure z.2) (U ×ˢ cylinder)
\end{leancode}
\end{minipage}
\par

The map $e$ extends the representatives from $I$; only its restriction
to $U$ is used in the smoothness conditions. The variable
\lean{z : ℝ × Vec} contains the parameter and spatial coordinates,
so \lean{z.1} and \lean{z.2} select these two entries. The set
\lean{U ×ˢ cylinder} is their Cartesian product. The local extensions
are smooth in all these variables together. This condition therefore
includes mixed derivatives and extension across both endpoints of $I$.

For each regularity $r$, the associated curve is
\[
 f:I\longrightarrow\mathscr M^r_{\mathrm{emb}}(Q),
 \qquad f(\lambda)=[X_\lambda,b_\lambda,p_\lambda].
\]
The following predicate asserts that the configurations have regularity
$r$, that $f$ is continuous and injective, and that every member with
$\lambda\in\operatorname{int}I$ is non-isolated:

\par
\noindent\begin{minipage}{\linewidth}
\begin{leancode}
def ModuliCurve (r : Regularity) (I : Set ℝ) (c : I → Representative) : Prop :=
  ∃ h : ∀ lambda, IsConfiguration r (c lambda),
    let f := fun lambda => moduliClass r ⟨c lambda, h lambda⟩
    Continuous f ∧ Function.Injective f ∧
      ∀ lambda : I, lambda.val ∈ interior I →
        ¬ IsOpen ({f lambda} : Set (ModuliSpace r))
\end{leancode}
\end{minipage}
\par

An element \lean{lambda : I} is a real parameter together with its
membership in $I$ and \lean{lambda.val} is the underlying real number.
The witness $h$ gives the regularity needed to form the moduli class.
The last line says that the singleton $\{f(\lambda)\}$ is not open
in the full moduli space, which is exactly non-isolation.
Mathematically, this also follows from continuity and injectivity:
if that singleton were open, its inverse image would be the open
singleton $\{\lambda\}$ in $I$, impossible for an interior parameter.
The Lean predicate includes this conclusion explicitly.

\subsection{The formal theorem}\label{app:lean-main}

With these definitions, the theorem recalled in
Subsection~\ref{app:lean-result} has the following Lean statement.
The name \lean{lambda} is the paper's parameter $\lambda$.

\par
\noindent\begin{minipage}{\linewidth}
\begin{leancode}
theorem equilibria_with_exact_cyclic_symmetry (L : ℝ) (hL : 0 < L) :
    ∃ a b : ℝ, 0 < a ∧ a < b ∧ b < 1 / 2 ∧
      ∃ N₀ : ℕ, 1 ≤ N₀ ∧ ∃ family : ℕ → Set.Icc a b → Representative,
        ∀ N : ℕ, N₀ ≤ N →
          SmoothFamily (Set.Icc a b) (family N) ∧
          (∀ lambda : Set.Icc a b, Equilibrium L N (family N lambda)) ∧
          (∀ k : ℕ, 3 ≤ k → ModuliCurve (.finite k) (Set.Icc a b) (family N)) ∧
          ModuliCurve .smooth (Set.Icc a b) (family N)
\end{leancode}
\end{minipage}
\par

The arguments \lean{(L : ℝ)} and \lean{(hL : 0 < L)} are the
cell length and its positivity hypothesis. All other data are conclusions:
there are real endpoints $0<a<b<1/2$, a threshold $N_0\geq1$, and one
family of representatives. Mathlib's \lean{Set.Icc a b} is the closed
interval $[a,b]$. Thus the first line gives precisely a compact interval
inside $(0,1/2)$ with nonempty interior.

The family is indexed by natural numbers, and its stated properties are
required whenever $N\geq N_0$. Since $N_0\geq1$, this expresses the
paper's assertion for all sufficiently large integers. The line
\lean{Equilibrium L N (family N lambda)} supplies the equations
and properties~(i)--(iii) for every parameter in that same interval.
Finally, the two \lean{ModuliCurve} clauses give property~(iv) for
every finite $k\geq3$ and for smooth regularity. Neither the interval nor
the family is allowed to depend on $k$.

\subsection{Source files and verification}\label{app:lean-verification}

All displayed declarations come from \lean{Showcase.lean}. The two occurrences of \lean{sorry} omit the proofs for exposition.
The file \lean{Showcase_WithProofs.lean} proves both statements. The proof files were checked with Lean~4.33.1. The proofs depend only on Lean's three standard axioms,
\lean{propext}, \lean{Classical.choice}, and \lean{Quot.sound}.

\begingroup
\renewcommand{\baselinestretch}{1}
\selectfont
\bibliographystyle{plain}
\bibliography{references}
\endgroup
\end{document}